%% file: main.tex
\pdfoutput=1
\documentclass{amsbook}
\usepackage{macros}
 
\usepackage{setspace}



\newcommand{\fronttitle}[1]{%
  \clearpage
  \phantomsection
  \addcontentsline{toc}{chapter}{#1}%
  {\normalfont\huge\bfseries\noindent #1\par}%
  \vspace{3em}
}

\begin{document}

\frontmatter

\input{front_matter_ist}

\cleardoublepage
\newgeometry{margin=2.5cm}
\fronttitle{Abstract}
\begin{spacing}{1.5}
\input{support-files/abstract_en.tex}
\end{spacing}
\medskip

\fronttitle{Resumo}
\begin{spacing}{1.5}
\input{support-files/abstract_pt.tex}
\end{spacing}
\fronttitle{Resumo Alargado}
\begin{spacing}{1.5}
\input{support-files/abstract_extended_pt.tex}
\end{spacing}

\restoregeometry
\fronttitle{Extended Abstract}
\begin{spacing}{1.5}

\nid The notion of defects in quantum field theories (QFTs) has amassed an extraordinary amount of attention and importance in physics and, perhaps surprisingly, also in mathematics. Let $M$ be a smooth $n$-dimensional manifold on which a given QFT $\tau$ lives. A defect of dimension $d$ can be heuristically articulated as the data of a submanifold $D\subset M$ of codimension $(n-d)$, together with some discontinuity of $\tau$ along $D$ which produces another quantum field theory. Examples of these phenomena are given by Wilson and ’t Hooft line defects in Yang--Mills theory. Giving a definition and a classification of these general defects is the subject of extensive ongoing work (see \cite{CEG23,CZR25}).

QFTs which lack dependence on any geometric parameter are dubbed \emph{topological QFTs} (TQFTs). In particular, defects whose interaction with the bulk theory is independent of the isotopy class of their embedding are called topological defects. One can organize topological QFTs with topological defects into a single package, called a defect (or stratified) topological quantum field theory. These are symmetric monoidal functors out of the category of (oriented) stratified cobordisms with defect data into the category of finite-dimensional complex vector spaces
\[\cZ \maps \mathsf{Bord}_{d}^{\sf strat} \longra \mathsf{Vect}_{\CC}^{\sf f.d}.\]
Extending these categories upwards — by adding stratified diffeomorphisms and stratified isotopies between them, and higher stratified isotopies between these, etc. — we can contemplate derived versions of these defect TQFTs; whereas by extending them downwards — letting stratified bordisms be cut into stratified cornered pieces of arbitrary depth — we can contemplate more computable versions of these derived defect TQFTs. This thesis is concerned with building such an extended object.

In his sketch of a proof of the \emph{Cobordism Hypothesis} \cite{Lur09}, Lurie introduces the first incarnation of the fully extended $(\oo,d)$-category of framed cobordisms as a (complete) $d$-fold Segal space. In the last section of that paper, he gives a very rough \emph{sketch} of a definition for an $(\oo,d)$-category of cobordisms with singularities through a convoluted notion of $\vec{X}$-manifold with singularities. In principle, this is the object we are after; but in practice, this class of spaces lacks differential-topological foundations that would enable its widespread use. We revise his construction in three ways: 1) we use a different model for higher categories; 2) we make use of the theory of conically smooth stratified spaces and their tangential structures; 3) we limit our construction to stratified bordisms whose underlying space is a topological manifold.

We now provide a general outline of the thesis. In Chapter 2 we introduce a model for $(\oo,1)$-categories that complements the existing literature, which we call \emph{interval Segal spaces} — these are simplicial co-presheaves on $\Int$, a subcategory of the poset of opens of $[0,1]$. In particular, we construct a model structure on $\mathsf{coPsh}(\Int)$ which is Quillen equivalent to that of Rezk’s Segal spaces. The motivation here is geometric: most models of $(\oo,1)$-categories are not tailored for the construction of geometric moduli spaces. A rare exception is the recent work of Ayala--Francis--Rozenblyum \cite{AFR1}, where the authors present \emph{striation sheaves} — a highly bespoke class of constructible sheaves on stratified spaces, designed precisely for constructing $(\oo,1)$-categories of a moduli flavour using point-set data. Our construction attempts a similar feat, while being less ambitious and remaining closer to more traditional Segal spaces, thereby facilitating higher-dimensional versions thereof.

In Chapter 3 we present a new construction of the $d$-dimensional bordism $(\oo,1)$-category. This construction is a spiritual cousin of other works, namely those of Galatius--Madsen--Tillman--Weiss \cite{GMTW09}, which later inspired the work of Lurie \cite{Lur09}, and more recently that of Scheimbauer--Calaque \cite{CS19}.

The pioneering work of Galatius--Madsen--Tillman--Weiss offered the first tools to construct an $(\oo,1)$-category of cobordisms. In their paper on the homotopy type of the cobordism category, they construct a topologically enriched category of $d$-dimensional bordisms using classifying spaces. Building on their work, Lurie, in his sketch of a proof of the cobordism hypothesis \cite{Lur09}, put forth a first construction of a model for an $(\oo,n)$-category of cobordisms using $n$-fold complete Segal spaces. This construction contained a small mistake preventing it from satisfying the \emph{essential constancy} condition. In order to address this, in her PhD thesis and in a subsequent paper, Scheimbauer--Calaque propose an alternative definition, which fixes Lurie’s while retaining the use of $n$-fold complete Segal spaces. Despite this almost successful attempt at a \emph{corrigendum}, the work of Scheimbauer--Calaque contains a separate flaw in its claim that the simplicial sets modelling the spatial levels of its $n$-fold Segal spaces are Kan complexes — a claim which, being false, translates into a loss of fibrancy necessary for the proof of the Segal condition.

To remedy this, we make use of the theory of smooth $\oo$-sheaves, i.e. $\oo$-sheaves on the site of smooth manifolds, and their respective classifying spaces — as articulated in the work of Berwick-Evans--Boavida de Brito--Pavlov \cite{BBEP24} — as a conduit for constructing appropriate moduli spaces of smooth manifolds, smooth bordisms, and smooth families of composed bordisms. This gives our approach a distinct flavour and clarifies the minor pathologies identified in previous constructions.

Let us sketch the path we will take in this section. We start by defining a general construction $\widetilde{\bord}_{d,V}(-,-)$ which is a bifunctor on $\man^{\op}\times \Int$ that, when fixing $T\in \Int$, produces a sheaf of sets on $\man^{\op}$ and captures a type of moduli space construction for $d$-dimensional smooth manifolds embedded in $V$ and fibred over smooth manifolds of varying dimension — designed so that, upon taking its concordance resolution and further restricting along the cosimplicial manifold $\bDelta \ra \man$, we obtain a simplicial set $\bord_{d,V}$, or space, which we show is the classifying space of $(d-1)$-dimensional manifolds, $d$-dimensional bordisms, and composites of $d$-dimensional bordisms. Lastly, the most dense and technical part of this chapter is spent showing that this construction, a pre-interval space, is an interval Segal space which recovers the $d$-dimensional bordism category $\mathsf{Cob}_d$ as its homotopy 1-category.

In Chapter 4 we finally turn to the theory of \emph{conically smooth stratified spaces}. These were first introduced in the work of Ayala--Francis--Tanaka and further elaborated in the subsequent works of Ayala--Francis--Rozenblyum \cite{AFT1, AFR1, AFR2}. They form the geometric foundations for our stratified higher bordism categories, and their adoption for that role is undoubtedly one of the more novel features of this thesis.

This chapter contains a heavier dose of exposition. We introduce the necessary differential-topological toolkit for our construction of the stratified bordism category, namely: we introduce an $(\oo,1)$-category of conically smooth stratified spaces; the exit-path $(\oo,1)$-category; and the constructible tangent bundle — all of which are necessary for a discussion of stratified tangential structures thereon. We prove new theorems about the contractibility of spaces of embeddings of stratified spaces into Euclidean space, and we prove a stratified analogue of the famed Ehresmann fibration theorem — which is the key to unlocking comparisons with the infrastructure developed for the unstratified case.

In Chapter 5 we construct an interval Segal space, i.e. an $(\oo,1)$-category of conically smooth stratified $d$-dimensional bordisms. We begin by introducing the notion of a conically smooth stratified topological $d$-manifold by stipulating a suitable $(\oo,1)$-category of basics $\Bsc_{d}^{\cS_{d-1}}$, whose links are stratified spheres of appropriate codimensions. Secondly, we produce a pre-interval space $\cS \bord_d$ of stratified topological $d$-dimensional bordisms, again using the concordance resolution for smooth $\oo$-sheaves and a restriction along the cosimplicial smooth manifold $\bDelta \longra \man$, and show that it is an interval Segal space.

Lastly, we endow this interval Segal space with a symmetric monoidal structure agreeing with the expected one coming from disjoint union. After defining this symmetric monoidal interval Segal space, we initiate a discussion on suitable notions of tangential structures for stratified bordisms — in particular, we introduce the notion of \emph{oriented discrete defect datum} by way of discrete right fibrations over $\Bsc_{d}^{\cS_{d-1}}$, as an attempt to make contact with the defect data for stratified bordisms introduced in the work of Carqueville et al. (see \cite{CRS19}) and we prove a comparison result for the stratified oriented case $d=2$.

We conclude Chapter~5 by endowing the resulting interval Segal space of conically smooth stratified $d$-bordisms with the symmetric monoidal structure induced by disjoint union. We then initiate a discussion of tangential structures in the stratified setting and, in particular, explain how discrete defect data can be encoded as tangential data via discrete right fibrations over the chosen $(\oo,1)$-category of basics. As a first test of this perspective, we prove a comparison result in the oriented stratified case $d=2$. Finally, we define a (derived) \emph{stratified topological quantum field theory} to be a symmetric monoidal functor from the symmetric monoidal $(\oo,1)$-category of stratified bordisms constructed here to a symmetric monoidal $(\oo,1)$-category $\cC$. 

Chapter 6 is a future outlook. In the first section, following ideas in unpublished work of Manuel Ara\'ujo and Pedro Boavida de Brito, we introduce a higher-dimensional generalization of $\Int$, denoted $\cube_d$ for $d>1$, aimed at extending interval Segal techniques to the study of $(\oo,d)$-categories. We prove that simplicial co-presheaves on $\cube_d$ satisfying suitable higher-dimensional analogues of the Segal and completeness conditions present the $(\oo,1)$-category $\Cat_{(\oo,d)}$, and we outline how this formalism should be used to build a $\cube_d$-space of stratified bordisms and verify the corresponding Segal condition, thus providing a candidate model for fully extended stratified TQFTs \`a la Atiyah--Baez--Dolan--Lurie. In the second section we return to the non-extended setting and sketch how constructible factorization algebras and \emph{stratified factorization homology} can be leveraged as a mechanism for producing stratified (and derived) defect TQFTs from appropriate stratified disk-algebras, in the spirit of the factorization-homological approach of Scheimbauer and Calaque.

\end{spacing}

\fronttitle{Acknowledgments}
\begin{spacing}{1.5}

The effort that has culminated in this thesis is mostly invisible to the casual reader: mere words that blend into sections, which in turn blend into chapters. The hours spent talking to people, reading, writing, absorbing and being absorbed; the time and space for reflection asked of those close to me; the doubts, detours, and conversations that shaped the work — all remain out of sight, except perhaps in this brief section.

I will begin with my advisors, John Huerta and Pedro Boavida de Brito. To John, without whom I would never have considered the possibility of attempting to be a mathematician: you have become a friend, and I hope your courage and attitude towards mathematical and academic freedom are threaded throughout this text — your patience, your mathematical culture, and your razor-blade rigour, too, have saved me from terminal exhaustion. To Pedro Boavida de Brito, your willingness to endure my many unannounced door knocks has helped shape this thesis as much as any theorem or proposition: our endless conversations and your kindness and generosity in sharing your ideas are the lifeblood of this thesis — doing manifold topology near you has been an honour reserved for a select few. I would also like to acknowledge my former co-supervisor, Roger Picken — to you I am indebted for all the hours spent listening to me talk in seminars, all the time you spent engaging me and my ideas, and for your unmatched attentiveness to the well-being of your students.

I am also grateful to the Fundação para a Ciência e a Tecnologia (FCT) for funding my PhD, and for making possible the years of attention and freedom on which this work depended. I thank Pedro Resende, president of CAMGSD, for his constant support, for the kindness with which he accompanied this process, and for the enthusiasm for category theory that he has always made feel not only legitimate, but contagious. I am equally grateful to Filipe Mena, president of the Doctoral Program in Mathematics and a member of the jury, for his steady support of the program and for the care with which he accompanied this final stage. I would also like to leave a note of appreciation to João Pimentel Nunes, from whom I learned so much about geometry and teaching, and who always heeded my calls to attend conferences, workshops, and seminars.

I am deeply indebted to Manuel Araújo for our many hours of conversation about TQFTs, string diagrams, and higher categories, for your kindness in sharing your unpublished work, and for your generosity in calmly explaining it to me. To David Ayala I owe all the courage I had to muster to face the homotopical higher-category theory laid bare in this thesis — your kindness in sharing your work, your encouragement, and your mentorship left an indelible mark on me. Your work, together with that of John Francis, Aaron Mazel-Gee, Hiro Lee Tanaka, and Nick Rozenblyum, on stratified spaces and factorization homology is among the strongest foundations on which this thesis stands. To Nils Carqueville I in some sense owe this thesis: our long conversations during your stay in Lisbon, and your stern and steady words, helped me channel and focus my attention towards the problem of strengthening the foundations of (higher) homotopical defect TQFTs. Your seminal work on defect TQFTs, together with that of Jürgen Fuchs, Catherine Meusburger, and Ingo Runkel, is the raison d'être for this work. Thanks for your generosity and your stewardship. 

Across the many years I spent thinking about higher algebra and topological field theories, I met some of the most passionate, brilliant, and kindest people I have ever met: Chun-Yu — to whom I owe some joint work —, Jennifer Brown, Rodrigo Duarte, Ben Haioun, David Jordan — I will never again witness someone draw skeins with such fiery passion on their kitchen table —, Patrick Kinnear, Theodoros Lagiotis, L\''ukas Muller, Jack Romo — an academic brother of sorts —, Chris Schommer-Pries — it was one of my greatest pleasures to hear about 3-dualizability and the status of the cobordism hypothesis directly from you —, Arber Selimi, Luuk Stehouwer, Daniel Teixeira, Matthew Yu, and Markus Zetto. All of you have made my four-year path through research mathematics among the most enjoyable and intellectually stimulating years I have had. A final note of sincere appreciation to Nino Scalbi — one my academic brothers with whom I shared hours and hours talking about mathematics and higher categories, without your precious inputs and your example to lead the way, I fear for what could have been. To Rui Peixoto, my other academic sibling, to you I owe you the moments of brilliance you shared with me in the form of questions, seminars, lectures, notes, insights, travels and so much more — I wish you well.

To my office mates: Filipe Berletta, for being our office's academic stepfather — without you we would all have been lost; Rafael Ferreira, for all the forró and all the relaxed energy you brought to our office; Lucas Ferrero; Diogo Franquinho, for all we have shared and lived so far, and for all we will still live together — you have become a friend for life; long live Gradient; Jasper Kager; Mena; Christoph Obenaus, thank you for being a friend, for sharing your passion for sports and for learning, and for inspiring me to become a better man; Francisco Pais, you have become a friend with whom I am sure I will still share great memories — thanks for everything; João Rodrigues; Miguel Santiago, thank you for your camaraderie and for our discussions about cylinders and life; Serena; and Leander Stecker, for all the beers and camaraderie.

To those who made my life outside Técnico not only bearable, but joyful, I begin with the people with whom I had the pleasure of teaching at Nova SBE. Frederico Toulson, more than an academic cousin and my shadow through the halls of Nova, you have become a trusted friend; it was a pleasure to do everything with you: mathematics, lunches, teaching, meetings, everything — please stay around. To Maria João Braga, I owe my job at Nova SBE to you; more importantly, I owe you all the compassion, respect, and kindness with which you always treated me — you are, without a shadow of a doubt, my academic mother. To Tânia I owe our many long conversations about mathematics and teaching, and a great deal of gratitude for your patience in hearing me ramble on about what needed to change in the Linear Algebra course. So many more people should be named here: Maria Helena Almeida, Rodrigo Barrela, Ernesto Freitas, Patrícia Ramos, and, I am sure, many others whose names I cannot now recall. To Feliu and Lu, thank you for the warmth with which you made me feel at home in Barcelona when I felt as though the world was ending.

To my family and friends: no words will ever make good on the debts I owe you. To my parents, everything — I speak for the four of us: everything. To my siblings, Madalena, Miguel, and Carolina, I owe the joy you have brought into my life. Discovering the world beside you gave me the confidence, and the reason, to pursue a PhD in the first place.

Finally, to Elisa: I owe you the most of what can be measured, and possibly all of what cannot. Without you, I would have given up on writing this long ago; without you, this journey would have been a fruitless chase for reason. I am happy to have had your hand in mine through it all, and happier still to feel it holding mine as I write these words.

\end{spacing}

\cleardoublepage
\tableofcontents

\cleardoublepage
\chapter*{Index of notation}\

\begin{itemize}
    \item The symbol $\cong$ means isomorphism in whichever category the morphism inhabits.
    \item The symbol $\simeq$ means equivalence. Further context should specify what type of equivalence we are referring to.
    \item $\Cat$ means the category of small categories and functors.
    \item $\bDelta$ means the skeleton of the category of finite totally ordered sets and order-preserving maps.
    \item $\sO_{\pm}$ refers to the category of finite linearly ordered sets with a distinguished minimum and maximum and order- and extrema-preserving maps.
    \item $\Gamma$ refers to the Segal category.
    \item $\Delta^{n}$ is the simplicial set corresponding to the representable of $[n]\in\bDelta$---it is also referred to by the standard $n$-simplex.
    \item $\Delta[n]$ is the simplicial space corresponding to $\Delta^n$ under the image of $\sset\hookra \fun(\bDelta^{\op,}\sset)$.
    \item $\Int$ refers to the \emph{interval} category.
    \item $\cInt$ refers to the \emph{closed interval} category.
    \item $\gamma\colon \Int \longra \bDelta^{\op}$ is the canonical comparison functor.
    \item $\sU(A_1,\ldots,A_n)\in \Int$ denotes the interval object determined by a finite string of closed intervals.
    \item $\core(I)$ denotes the core of an interval object $I\in \Int$.
    \item $\cW_{\pi_0}$ denotes the class of morphisms in $\Int$ inducing isomorphisms on path components.
    \item $\mathsf{ISSp}$ denotes the model category of interval Segal spaces.
    \item $\mathsf{CISSp}$ denotes the model category of complete interval Segal spaces.
    \item $\ho(\cC)$ denotes the homotopy category of an $(\oo,1)$-category or Segal-type model $\cC$.
    \item $\man$ denotes the category, or site, of smooth manifolds.
    \item $\AA^\bullet$ denotes the standard cosimplicial smooth affine simplex.
    \item $\widetilde{\bord}_{d,V}$, $\bord_{d,V}$, and $\bord_d$ denote the smooth bordism pre-interval object, the fixed-ambient interval Segal space, and the stabilized interval Segal space of $d$-dimensional smooth bordisms, respectively.
    \item $\bTheta_n$ is the iterated wreath product $\bDelta \wr \bTheta^{n-1}$.
    \item $\cube_n$ denotes the $n$-cubes category introduced in the outlook chapter.
    \item $\overline{\sC}(X)$ refers to the \emph{closed} cone on a topological stratified space $\big(X \ra P\big)$.
    \item $\strat$ denotes the ordinary category of conically smooth stratified spaces and conically smooth stratified maps.
    \item $\Strat$ denotes the associated $(\oo,1)$-category of conically smooth stratified spaces.
    \item $\Bsc$ denotes the $(\oo,1)$-category of basics.
    \item $\Snglr$ denotes the $(\oo,1)$-category of singularity types.
    \item $\Bun$ denotes the $(\oo,1)$-category of constructible bundles.
    \item $\Exit$ denotes the absolute exit-path $(\oo,1)$-category.
    \item $\Vect^{\inj}$ denotes the $(\oo,1)$-category of constructible vector spaces with injective transition maps.
    \item $\tau$ denotes the tangent classifier.
    \item $\exit(M)$ and $\entr(M)$ denote the exit-path and enter-path $(\oo,1)$-categories of a stratified space $M$.
    \item $\sT_M$ denotes the constructible tangent bundle of a conically smooth stratified space $M$.
    \item $\sT^{\mathrm{fib}}_{\pi}$ denotes the fiberwise constructible tangent bundle of a constructible bundle $\pi$.
    \item $\Mfd(\cB)$ denotes the $(\oo,1)$-category of $\cB$-manifolds, and $\Disk(\cB)\subset \Mfd(\cB)$ the full symmetric monoidal subcategory generated by the basics.
    \item $\int_X A$ denotes factorization homology of $X$ with coefficients in $A$.
    \item $\SMfd_d$ denotes the $(\oo,1)$-category of conically smooth stratified topological $d$-manifolds.
    \item $\SMfd_d^\partial$ denotes the $(\oo,1)$-category of conically smooth stratified topological $d$-manifolds with boundary.
    \item $\cS_d$ and $\cS_{\le d}$ denote the $(\oo,1)$-categories of stratified spheres of dimension exactly $d$, and of dimension at most $d$, respectively.
    \item $\Bsc_d^{\cS_{d-1}}$ denotes the chosen $(\oo,1)$-category of basics used to model stratified topological $d$-manifolds.
    \item $\cS\mathsf{Cob}_d$ denotes the ordinary category of conically smooth stratified bordisms.
    \item $\tilbordstr_{d,V}$, $\bordstr_{d,V}$, and $\bordstr_d$ denote the stratified bordism pre-interval object, the fixed-ambient interval Segal space, and the stabilized interval Segal space of stratified bordisms, respectively.
    \item $\Bsc_d^{\cS_{d-1},\DD}$ denotes the category of basics equipped with discrete defect data.
    \item $\Bsc_d^{\cS_{d-1},\mathsf{or}}$ denotes the oriented version of the chosen category of stratified basics.
    \item $\BSO$ denotes the stable tangential structure of orientations.
    
\end{itemize}

\mainmatter
\raggedbottom

\input{./Chapters/Intro.tex}

\input{./Chapters/Ch2.tex}

\input{./Chapters/Ch3.tex}

\input{./Chapters/Ch4.tex}

\input{./Chapters/Ch5.tex}

\input{./Chapters/FO.tex}

\appendix

\input{./Chapters/Appendix.tex}

\printbibliography

\end{document}

%% file: front_matter_ist.tex
\newcommand{\istthesislogo}{%
  \includegraphics[trim=36.72bp 36.72bp 36.72bp 36.72bp,clip,width=3.8cm]{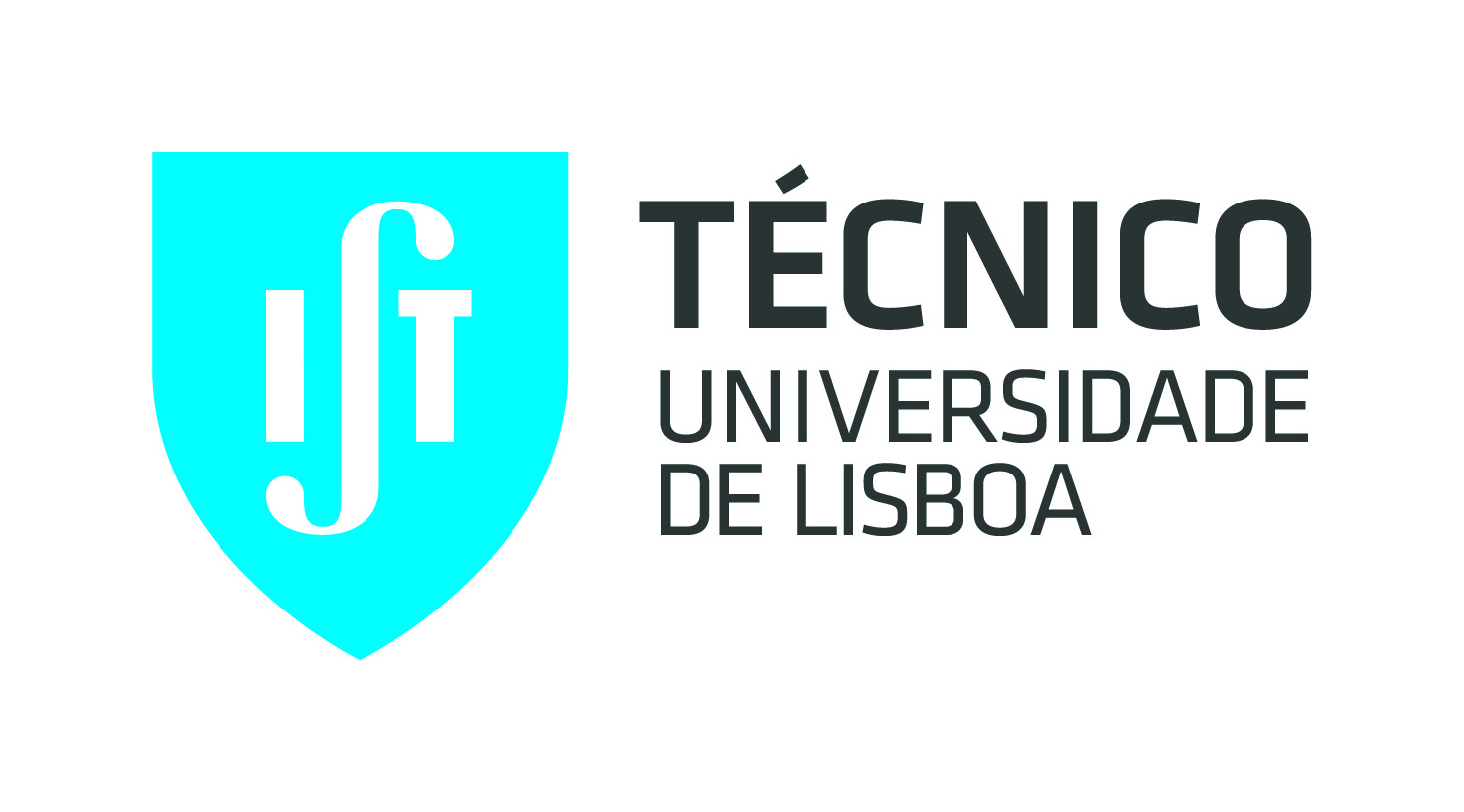}%
}

\newcommand{\istthesistitle}{%
  {\Large\bfseries An $(\infty,1)$-category of conically smooth stratified bordisms\par}
}
\newcommand{\istthesisauthor}{%
  {\large\bfseries Diogo Maria Marques Coelho Freire de Andrade\par}
}
\newcommand{\istthesissupervisors}{%
  \noindent\textbf{Supervisor:} Doctor John Gmerek Huerta\\[0.15cm]
  \textbf{Co-supervisor:} Doctor Pedro Boavida de Brito\par
}
\newcommand{\istthesisapproval}{%
  \begin{center}
    \emph{Thesis approved in public session to obtain the PhD Degree in Mathematics}\\[0.35cm]
    \textbf{Jury final classification: Pass with Distinction}
  \end{center}
}

\begin{titlepage}
  \thispagestyle{empty}

  \noindent\parbox[c][4.2cm][c]{\textwidth}{\istthesislogo}
  \vspace{0.3cm}

  \begin{center}
    \textbf{UNIVERSIDADE DE LISBOA\\[0.25cm]
    INSTITUTO SUPERIOR TÉCNICO}
    \vspace{2.2cm}

    \istthesistitle
    \vspace{1.6cm}
    \istthesisauthor
  \end{center}

  \vspace{1.7cm}
  \hspace*{2.7cm}%
  \begin{minipage}{0.6\textwidth}
    \textbf{Supervisor:} Doctor John Gmerek Huerta\\[0.25cm]
    \textbf{Co-supervisor:} Doctor Pedro Boavida de Brito
  \end{minipage}

  \vspace{1.5cm}
  \begin{center}
    \emph{Thesis approved in public session to obtain the PhD Degree in\\
    Mathematics}
    \vspace{0.9cm}

    \textbf{Jury final classification: Pass with Distinction}
    \vfill

    \textbf{2026}
  \end{center}
\end{titlepage}

\begin{titlepage}
  \thispagestyle{empty}

  \noindent\parbox[c][2.0cm][c]{\textwidth}{\istthesislogo}
  \begin{center}
    \textbf{UNIVERSIDADE DE LISBOA\\[0.15cm]
    INSTITUTO SUPERIOR TÉCNICO}
    \vspace{0.4cm}

    \istthesistitle
    \vspace{0.4cm}
    \istthesisauthor
  \end{center}
  \vspace{0.45cm}
  \istthesissupervisors
  \vspace{0.35cm}
  \istthesisapproval
  \vspace{0.35cm}
  \begin{center}{\large\bfseries Jury}\end{center}
  \begin{flushleft}
    \footnotesize
    \textbf{Chairperson:}\\
    Doctor Filipe Artur Pacheco Neves Carteado Mena, Instituto Superior Técnico, Universidade de Lisboa.\\[0.25cm]
    \textbf{Members of the Committee:}\\
    Doctor David Ayala, Department of Mathematical Sciences, College of Letters \& Science, Montana State University, USA.\\[0.12cm]
    Doctor Filipe Artur Pacheco Neves Carteado Mena, Instituto Superior Técnico, Universidade de Lisboa.\\[0.12cm]
    Doctor Nils Carqueville, Faculty of Physics, University of Vienna, Austria.\\[0.12cm]
    Doctor Roger Francis Picken, Instituto Superior Técnico, Universidade de Lisboa.\\[0.12cm]
    Doctor Marko Stosic, Faculdade de Ciências, Universidade de Lisboa.\\[0.12cm]
    Doctor John Gmerek Huerta, Faculty of Physics, Ludwig Maximilian University of Munich, Germany.
  \end{flushleft}
  \vfill
  \begin{center}
    \textbf{Funding Institution}\\[0.15cm]
    FCT: Fundação para a Ciência e a Tecnologia\\[0.35cm]
    \textbf{2026}
  \end{center}
\end{titlepage}

%% file: support-files/abstract_en.tex
With the goal of initiating a rigorous study of derived defect topological quantum field theories, this thesis constructs an $(\oo,1)$-category of stratified bordisms and sketches its $(\oo,d)$-dimensional analogue.

The first part develops a method for constructing symmetric monoidal $(\oo,1)$-categories, inspired by Segal spaces but better suited to geometric input: interval Segal spaces. Using this framework, we construct the $(\oo,1)$-category of $d$-dimensional smooth manifolds via smooth $\oo$-sheaves of manifolds, correcting a mistake in earlier work of Calaque--Scheimbauer. We then equip this category with a symmetric monoidal structure, compute its homotopy 1-category, and extend the construction to manifolds endowed with tangential structures, obtaining the interval Segal space of $d$-dimensional bordisms with tangential structure.

To model the stratifications supporting topological defects, we work with conically smooth stratified spaces, a suitably smooth framework for conical stratifications. This choice is motivated by two intertwined requirements: the ability to handle families of stratifications and the flexibility to impose tangential structures on stratified spaces.

We define stratified topological manifolds as stratified $\cB$-manifolds for a suitable $(\oo,1)$-category of basics $\cB$, that is, a sieve of singularity types. Using these, we construct an interval Segal space $\tilbordstr_d$ whose objects are compact closed stratified topological $(d-1)$-manifolds and whose 1-morphisms are stratified $d$-dimensional bordisms, and we endow it with a symmetric monoidal structure. We also use tangential structures on conically smooth stratified spaces to formulate a notion of defect labels, aiming toward bordisms with labeled defects in the style of Carqueville et al.

Finally, we sketch a higher-dimensional generalization of $\Int$ to $\cube_d$ for $d>1$, conjecturally presenting $\Cat_{(\oo,d)}$, and outline how stratified factorization homology can be used to produce derived stratified defect TQFTs.

\bigskip

\begin{samepage}
\noindent\textbf{Keywords}
\begin{itemize}[nosep]
    \item conically smooth stratified spaces
    \item stratified bordisms
    \item interval Segal spaces
    \item higher categories
    \item topological quantum field theories
\end{itemize}
\end{samepage}

%% file: support-files/abstract_pt.tex
Com o objectivo de iniciar um estudo rigoroso de versões derivadas de teorias de campo quântico topológicas com defeitos, construímos uma $(\oo,1)$-categoria de bordismos estratificados e esboçamos também uma construção do seu análogo $(\oo,d)$. Na primeira parte da tese introduzimos um método de construção de $(\oo,1)$-categorias simétricas monoidais, inspirado nos espaços de Segal, mas mais adequado a construções geométricas, a que chamamos \emph{espaços de Segal intervalares}.

Usamos estes espaços para construir a $(\oo,1)$-categoria de variedades suaves de dimensão $d$, recorrendo à teoria de $\oo$-feixes suaves e corrigindo um erro no trabalho de Calaque--Scheimbauer. Dotamos depois esta estrutura de uma estrutura simétrica monoidal, calculamos a sua categoria de homotopia 1-truncada e construímos o espaço de Segal intervalar dos bordismos de dimensão $d$ com estruturas tangenciais.

Para captar as estratificações sobre as quais os defeitos topológicos são suportados, recorremos a \emph{espaços estratificados conicamente suaves}, uma versão adequadamente suave das estratificações cónicas. Esta escolha permite trabalhar simultaneamente com famílias de estratificações e com noções flexíveis de estrutura tangencial. Introduzimos então uma noção de \emph{variedade topológica estratificada} como $\cB$-variedade estratificada, onde $\cB$ é uma adequada $(\oo,1)$-categoria de básicos. Com esta classe de espaços construímos um espaço de Segal intervalar $\tilbordstr_d$ cujos objectos são variedades topológicas estratificadas compactas e fechadas de dimensão $(d-1)$ e cujos 1-morfismos são bordismos estratificados de dimensão $d$, dotando-o ainda de uma estrutura simétrica monoidal. Usamos também estruturas tangenciais para introduzir uma noção de rótulos de defeito.

Por fim, generalizamos $\Int$ para $\cube_d$ com $d>1$, delineando um modelo para $\Cat_{(\oo,d)}$, um $\cube_d$-espaço de bordismos estratificados e uma aplicação da homologia de factorização estratificada à construção de TQFTs estratificadas com defeitos.
\bigskip

\begin{samepage}
\noindent\textbf{Palavras-chave:}
\begin{itemize}[nosep]
    \item espaços estratificados conicamente suaves
    \item bordismos estratificados
    \item espaços de Segal intervalares
    \item categorias superiores
    \item teorias de campo quântico topológicas
\end{itemize}
\end{samepage}

%% file: support-files/abstract_extended_pt.tex
A noção de defeitos em teorias de campo quântico (QFTs) tem ganho uma atenção e importância extraordinárias na física e, talvez de modo menos esperado, também na matemática. Seja $M$ uma variedade diferenciável $n$-dimensional na qual vive uma dada QFT $\tau$. Heuristicamente, um defeito de dimensão $d$ pode ser descrito como o dado de uma subvariedade $D\subset M$ de codimensão $(n-d)$, juntamente com alguma descontinuidade de $\tau$ ao longo de $D$ que produz outra teoria de campo quântico. Exemplos paradigmáticos destes fenómenos são os defeitos de linha de Wilson e de ’t Hooft na teoria de Yang--Mills. Dar uma definição geral e uma classificação conceptual destes defeitos é hoje tema de um vasto corpo de trabalho em desenvolvimento.

QFTs que não dependem de quaisquer parâmetros geométricos são chamadas \emph{QFTs topológicas} (TQFTs). Em particular, defeitos cuja interação com a teoria de bulk é independente da classe de isotopia do seu mergulho são chamados defeitos topológicos. Pode então organizar-se uma TQFT topológica com defeitos topológicos num único pacote, chamado teoria topológica de campo quântico com defeitos, ou teoria topológica de campo quântico estratificada. Em linguagem categórica, estas teorias são modeladas por funtores simetricamente monoidais da categoria de cobordismos estratificados com dados de defeito para uma categoria-alvo linear, como a categoria dos espaços vetoriais complexos de dimensão finita. Se enriquecermos estas categorias “para cima”, acrescentando difeomorfismos estratificados, isotopias estratificadas e isotopias de ordem superior, obtemos versões derivadas destas TQFTs; se, por outro lado, as estendermos “para baixo”, permitindo decomposições em peças estratificadas com cantos de profundidade arbitrária, obtemos versões mais computáveis dessas teorias. Esta tese ocupa-se da construção de um tal objeto, numa primeira etapa ao nível de uma $(\oo,1)$-categoria e, de forma prospetiva, ao nível totalmente estendido.

No seu esboço de prova da \emph{Hipótese do Cobordismo}, Lurie introduziu uma primeira encarnação da $(\oo,d)$-categoria completamente estendida de cobordismos framizados como um espaço de Segal completo $d$-vezes iterado. Na secção final desse trabalho, sugere ainda uma definição de uma $(\oo,d)$-categoria de cobordismos com singularidades, baseada numa noção de variedade com singularidades suficientemente geral para capturar defeitos. Em princípio, esse é precisamente o tipo de objeto que aqui pretendemos construir. Contudo, na prática, essa classe de espaços não vinha acompanhada dos fundamentos diferencial-topológicos necessários para uma teoria robusta. A presente tese revê esse programa em três direções complementares: usa um modelo alternativo para categorias superiores, recorre à teoria de espaços estratificados conicamente suaves e às suas estruturas tangenciais, e restringe-se a cobordismos estratificados cujo espaço subjacente é uma variedade topológica.

A tese começa por introduzir, no Capítulo 2, um novo modelo para $(\oo,1)$-categorias, a que chamamos \emph{espaços de Segal intervalares}. Estes são copré-feixes simpliciais na categoria $\Int$, uma subcategoria do poset dos abertos de $[0,1]$, concebida de modo a tornar mais natural a codificação geométrica de composição. Mostramos que a categoria correspondente admite uma estrutura de modelo equivalente de Quillen à estrutura clássica dos espaços de Segal de Rezk. A motivação desta escolha é geométrica: muitos modelos de $(\oo,1)$-categorias são tecnicamente poderosos, mas pouco adequados à construção direta de espaços de módulos geométricos. Uma exceção notável é dada pelas \emph{striation sheaves} de Ayala, Francis e Rozenblyum, desenhadas precisamente para construir $(\oo,1)$-categorias a partir de dados ao nível de conjuntos. Os espaços de Segal intervalares introduzidos nesta tese podem ser vistos como uma alternativa mais próxima da tecnologia tradicional dos espaços de Segal, mas ainda assim suficientemente adaptada a construções geométricas e, por isso, promissora para generalizações a dimensões superiores.

No Capítulo 3 aplicamos este formalismo para construir uma nova versão da $(\oo,1)$-categoria de cobordismos $d$-dimensionais suaves. A abordagem inspira-se em três antecedentes principais: o trabalho de Galatius, Madsen, Tillmann e Weiss sobre o tipo de homotopia da categoria de cobordismos; o programa de Lurie para a hipótese do cobordismo; e a reformulação proposta por Scheimbauer e Calaque. O ponto de partida é que as construções existentes apresentam pequenas falhas técnicas quando reinterpretadas no contexto necessário para uma prova cuidada das condições de Segal e de fibrância. Em particular, a construção de Lurie não satisfaz literalmente uma condição de constância essencial, e a correção posterior de Scheimbauer e Calaque apoia-se numa afirmação problemática sobre a fibrância de certos níveis simpliciais.

Para ultrapassar estas dificuldades, recorremos à teoria de $\oo$-feixes suaves no sítio das variedades diferenciáveis e aos respetivos espaços classificantes, tal como desenvolvidos por Berwick-Evans, Boavida de Brito e Pavlov. A ideia central é construir, a partir de dados geométricos honestos, espaços classificantes de variedades suaves, cobordismos suaves e famílias suaves de composições de cobordismos. Mais concretamente, definimos um bifuntor $\widetilde{\bord}_{d,V}(-,-)$ em $\man^{\op}\times \Int$ que, para cada intervalo $T$, produz um feixe de conjuntos sobre variedades suaves e codifica famílias de variedades $d$-dimensionais mergulhadas num espaço euclidiano auxiliar $V$. Tomando a sua resolução por concordâncias e restringindo depois ao longo da variedade co-simplicial padrão, obtemos um conjunto simplicial $\bord_{d,V}$ que funciona como espaço classificante de variedades $(d-1)$-dimensionais, de cobordismos $d$-dimensionais e das suas composições. A parte técnica central do capítulo consiste em mostrar que este pré-espaço intervalar satisfaz as condições que o tornam num espaço de Segal intervalar e que a sua categoria de homotopia 1-truncada recupera a categoria de cobordismos esperada.

No Capítulo 4 passamos à geometria estratificada propriamente dita. Introduzimos a teoria dos \emph{espaços estratificados conicamente suaves}, desenvolvida por Ayala, Francis, Tanaka e Rozenblyum, como o contexto geométrico adequado para modelar defeitos. Esta escolha é motivada por duas exigências fundamentais: trabalhar com famílias de estratificações e formular noções suficientemente flexíveis de estrutura tangencial em contexto singular. O capítulo desenvolve o ferramental diferencial-topológico necessário para os capítulos seguintes: a $(\oo,1)$-categoria de espaços estratificados conicamente suaves, a categoria de caminhos de saída e o fibrado tangente construtível. Para além disso, demonstramos resultados novos sobre a contratilidade de certos espaços de mergulhos de espaços estratificados em espaços euclidianos, bem como um análogo estratificado do teorema de Ehresmann. Estes resultados são essenciais para transportar para o mundo estratificado parte da infraestrutura já conhecida no caso liso.

O Capítulo 5 contém o núcleo geométrico da tese: a construção de uma $(\oo,1)$-categoria de cobordismos $d$-dimensionais estratificados conicamente suaves. Começamos por introduzir a noção de variedade topológica estratificada conicamente suave, modelada como uma $\cB$-variedade estratificada para uma adequada $(\oo,1)$-categoria de básicos $\Bsc_{d}^{\cS_{d-1}}$, cujas ligações são esferas estratificadas das codimensões relevantes. Esta restrição a espaços subjacentes topológicos permite manter um equilíbrio entre generalidade geométrica e tratabilidade técnica. Construímos então um pré-espaço intervalar de cobordismos estratificados, novamente através da resolução por concordâncias de $\oo$-feixes suaves e da restrição ao longo da variedade co-simplicial padrão, e provamos que ele satisfaz a condição de Segal intervalar.

Em seguida, dotamos este espaço de uma estrutura simetricamente monoidal compatível com a união disjunta, obtendo assim o candidato natural à categoria de cobordismos estratificados necessária para formular TQFTs estratificadas derivadas. Um aspeto importante desta construção é a incorporação de estruturas tangenciais no contexto estratificado. Em particular, introduzimos a noção de \emph{dado de defeito discreto orientado} através de fibrados à direita discretos sobre a categoria de básicos escolhida, com o objetivo de aproximar o formalismo geométrico desenvolvido nesta tese dos dados de defeito considerados por Carqueville e colaboradores. Como primeiro teste desta perspetiva, provamos um resultado de comparação no caso orientado estratificado de dimensão $2$. O capítulo culmina com a definição de uma teoria topológica de campo quântico estratificada derivada como um funtor simétrico monoidal da $(\oo,1)$-categoria simétrica monoidal de cobordismos estratificados construída na tese para uma $(\oo,1)$-categoria simétrica monoidal alvo $\cC$.

O Capítulo 6 tem natureza prospetiva e delineia duas direções de desenvolvimento futuro. Na primeira, seguindo ideias não publicadas de Manuel Araújo e Pedro Boavida de Brito, introduzimos uma generalização de dimensão superior da categoria $\Int$, denotada $\cube_d$, com vista à formulação de um análogo intervalar para $(\oo,d)$-categorias. Argumentamos que copré-feixes simpliciais em $\cube_d$ satisfazendo análogos superiores das condições de Segal e de completude deverão apresentar a $(\oo,1)$-categoria $\Cat_{(\oo,d)}$, e esboçamos como este formalismo poderá ser usado para construir um $\cube_d$-espaço de cobordismos estratificados. Na segunda direção regressamos ao contexto não estendido e indicamos como álgebras de fatorização construtíveis e homologia de fatorização estratificada podem servir de mecanismo de produção de TQFTs estratificadas com defeitos, generalizando para o contexto estratificado e derivado ideias presentes na abordagem de Calaque e Scheimbauer.

Em síntese, a tese propõe um novo quadro para a construção de categorias superiores geométricas adaptadas ao estudo de cobordismos estratificados e de teorias topológicas de campo quântico com defeitos. A introdução dos espaços de Segal intervalares, a nova construção da categoria de cobordismos lisos via $\oo$-feixes suaves, e a passagem para o contexto dos espaços estratificados conicamente suaves constituem os três pilares principais deste programa. O resultado é uma construção rigorosa de uma $(\oo,1)$-categoria de cobordismos estratificados suficientemente rica para servir de domínio a TQFTs estratificadas derivadas, bem como um esboço das extensões futuras necessárias para tratar a versão totalmente estendida do problema.

\bigskip

\begin{samepage}
\noindent\textbf{Palavras-chave:}
\begin{itemize}[nosep]
    \item espaços estratificados conicamente suaves
    \item bordismos estratificados
    \item espaços de Segal intervalares
    \item categorias superiores
    \item teorias de campo quântico topológicas
\end{itemize}
\end{samepage}

%% file: Chapters/Intro.tex
\chapter{Introduction}

The notion of defects in quantum field theories (QFTs) has amassed an extraordinary amount of attention and importance in physics and, perhaps surprisingly, also in mathematics. Let $M$ be a smooth $n$-dimensional manifold on which a given QFT $\tau$ lives. A defect of dimension $d$ can be heuristically articulated as the data of a submanifold $D\subset M$ of codimension $(n-d)$, together with some discontinuity of $\tau$ along $D$ which produces another quantum field theory. Examples of these phenomena are given by Wilson and ’t Hooft line defects in Yang--Mills theory. Giving a definition and a classification of these general defects is the subject of extensive ongoing work (see \cite{CEG23,CZR25}).

QFTs which lack dependence on any geometric parameter are dubbed \emph{topological QFTs} (TQFTs). In particular, defects whose interaction with the bulk theory is independent of the isotopy class of their embedding are called \emph{topological defects}. One can organize topological QFTs with topological defects into a single package, called a defect topological quantum field theory, as they are defined in the work of Carqueville et al. — for an overview of recent developments see \cite{CZR25} and for a thorough introduction see \cite{CRS19}. These are symmetric monoidal functors out of the 1-category of (oriented) stratified cobordisms with defect data into the category of complex vector spaces
\begin{equation}
\cZ \maps \mathsf{Bord}_{d}(\DD) \longra \mathsf{Vect}_{\CC},
\end{equation}
where $\DD=(D_d,\ldots,D_0)$ is the collection of defect labels — when omitting the defect labels $\DD$, we will call these \emph{stratified} TQFTs.

Extending the categories of defect bordisms upwards — by adding stratified diffeomorphisms and stratified isotopies between them, and higher stratified isotopies between these, etc. — we can contemplate derived versions (or $\oo$ versions) of these defect TQFTs; whereas by extending them downwards — letting stratified bordisms be cut into stratified cornered pieces of arbitrary depth — we can contemplate extended and more computable versions of these derived defect TQFTs. In this thesis we will construct the $(\oo,1)$-category of stratified bordisms whose objects are closed $(d-1)$-dimensional topological manifolds endowed with a conically smooth stratified atlas; the 1-morphisms are $d$-dimensional bordisms endowed also with a conically smooth atlas; the higher morphisms encode conically smooth isomorphisms and conically smooth isotopies among them. We consider also versions with defect labels. In particular, we sketch a way of interpreting the defect labels $\DD$ defined in \cite[Definition~2.4]{CRS19} as part of the tangential data associated to a conically smooth stratified space. Lastly, we provide a sketch of a definition of the $(\oo,d)$-category of fully-extended stratified bordisms: objects are points; $k$-morphisms are stratified $k$-bordisms for $k\leq d$ and stratified isotopies for $k > d$.

Despite the growing interest in defect TQFTs, especially derived and extended ones, the literature on such $(\oo,d)$-categories of stratified bordisms (or defect bordisms) is scant still, and mostly limited to specific ranges of codimension for the stratification and to low dimensions. Our construction is based on the following three choices:
\begin{enumerate}
\item We use a different model for higher categories, which we call complete ``interval'' Segal spaces (\Cref{intro:: thm: interval models}).
\item We make use of the theory of conically smooth stratified spaces and their tangential structures — based on the work of Ayala--Francis--Tanaka--Rozenblyum (\Cref{intro:: thm: stratified bordisms}).
\item We limit our construction to stratified bordisms whose underlying space is a topological manifold.
\end{enumerate}
These choices mean that our construction, while close in spirit to Lurie's in \cite[Definition~4.11]{Lur09}, explores a more geometrically-biased variant of Segal spaces, enjoys a more solid differential-topological foundation coming from the comprehensive works of Ayala--Francis--Tanaka \cite{AFT1} and Ayala--Francis--Rozenblyum \cite{AFR1}\cite{AFR2}, and is less general insofar as it concerns only topological spaces whose underlying space is locally Euclidean and of fixed dimension. This last point, we believe, is not a true limitation of our construction and we intend for future work to address it — concretely, the cornerstone of our approach is a version of the stratified Ehresmann fibration theorem, which does not apply only to conically smooth stratified spaces whose underlying topological space is locally Euclidean.

\section*{Main results}

We now state the three principal results of the thesis. The key idea is to
replace the simplicial indexing category $\bDelta^{\op}$ by a geometric
category $\Int$. Informally, an object of $\Int$ is a finite union of open
subintervals of $[0,1]$ containing both endpoints, but not the whole interval;
morphisms are inclusions. Its complement is a non-empty ordered collection
of closed intervals, which specify regions over which a bordism can be cut
and glued. The functor $\gamma\colon\Int\to\bDelta^{\op}$ remembers only
the ordered set of these complementary components.

A \emph{pre-interval space} is a functor $\Int\to\sset$, where $\sset$
denotes simplicial sets. An \emph{interval Segal space} sends refinements
preserving the cutting components to weak equivalences and reconstructs
chains of morphisms from their individual pieces by homotopy pullbacks.
Completeness identifies the space of objects with the space of equivalences,
as for Rezk's complete Segal spaces.

\begin{samepage}
\begin{maintheorem}[Interval models for $(\oo,1)$-categories]
\label{intro:: thm: interval models}
Let $\Int$ be the inclusion poset of finite unions of open subintervals of
$[0,1]$ containing $\{0,1\}$ and different from $[0,1]$. The category
$\fun(\Int,\sset)$ of pre-interval spaces admits a left proper, combinatorial, and
simplicial model structure whose fibrant objects are the complete interval
Segal spaces. Moreover, restriction and right Kan extension along
\[
    \gamma\colon \Int\longrightarrow\bDelta^{\op}
\]
where $\gamma$ records the ordered complementary components, induce a Quillen equivalence
\[
    \mathsf{CSSp}
    \ \underset{\rkan_{\gamma}}{\overset{\gamma^{*}}{\rightleftarrows}}\
    \mathsf{CISSp}.
\]
Here $\mathsf{CSSp}$ is Rezk's model category of complete Segal spaces and
$\mathsf{CISSp}$ is the model category just described.
Consequently, complete interval Segal spaces present the
$(\oo,1)$-category of $(\oo,1)$-categories.
\end{maintheorem}
\end{samepage}

The significance of \Cref{intro:: thm: interval models} is that replacing
$\bDelta^{\op}$ by $\Int$ retains geometric cutting data without changing
the homotopy theory of $(\oo,1)$-categories. The precise statement and proof
are given in \Cref{ch2:: thm: complete interval model}.

For the geometric constructions we use the smooth-sheaf framework of
Berwick-Evans, Boavida de Brito, and Pavlov (BBEP) \cite{BBEP24}.
The \emph{smooth moduli sheaf of embedded bordisms} assigns to a smooth
parameter manifold the smooth families of embedded bordisms over it, with
specified cylindrical cutting regions. Its \emph{concordance resolution}
organizes these families into a simplicial space: vertices are individual
bordisms, paths are smooth one-parameter families, and higher simplices
record families with more parameters. This retains deformation data rather
than taking only a set of diffeomorphism classes. Stabilizing the ambient
Euclidean space means allowing arbitrarily many embedding coordinates.

\begin{samepage}
\begin{maintheorem}[The smooth bordism $(\oo,1)$-category]
\label{intro:: thm: smooth bordisms}
For every $d\geq 0$, form $\bord_d\colon\Int\to\sset$ from smooth families
of embedded $d$-dimensional bordisms with cylindrical cutting regions by
the concordance resolution of Berwick-Evans, Boavida de Brito, and Pavlov
\cite{BBEP24} and stabilization of the ambient Euclidean
space. Then $\bord_d$ is a symmetric monoidal interval Segal space, with
monoidal operation given by disjoint union. Its objects are closed smooth
$(d-1)$-manifolds, its morphisms are smooth $d$-bordisms, and its higher
morphisms encode diffeomorphisms and their isotopies. Its homotopy category
$\ho(\bord_d)$, obtained by identifying homotopic morphisms, is equivalent
as a symmetric monoidal category to the ordinary unoriented bordism category
$\mathsf{Cob}_d$, whose morphisms are diffeomorphism classes of bordisms:
\[
    \ho(\bord_d)\simeq\mathsf{Cob}_d.
\]
\end{maintheorem}
\end{samepage}

Thus the interval model supports a geometric construction of the bordism
$(\oo,1)$-category in which smooth families are retained rather than
prematurely divided out by diffeomorphisms. The interval Segal condition,
symmetric monoidal structure, and homotopy-category comparison are proved
in \Cref{ch3:: thm: interval segal bordism cat,ch3:: thm: symmetric monoidal,ch3:: thm: homotopy category bordism}, respectively.

\begin{samepage}
\begin{maintheorem}[The stratified bordism $(\oo,1)$-category]
\label{intro:: thm: stratified bordisms}
For every $d\geq 0$, form $\bordstr_d\colon\Int\to\sset$ from smooth
families of embedded conically smooth stratified topological $d$-bordisms
with cylindrical cutting regions by the concordance resolution of
Berwick-Evans, Boavida de Brito, and Pavlov \cite{BBEP24} and
stabilization of the ambient Euclidean space. Here the underlying spaces
are topological manifolds and the stratifications have conically smooth
local models. Then $\bordstr_d$ is a symmetric monoidal interval Segal space,
with monoidal operation given by disjoint union. Its objects are closed
stratified topological $(d-1)$-manifolds, its morphisms are stratified
$d$-bordisms, and its higher morphisms encode conically smooth isomorphisms
and their isotopies. Its homotopy category is equivalent as a symmetric
monoidal category to $\cS\mathsf{Cob}_d$, the ordinary unoriented stratified
bordism category whose morphisms are conically smooth isomorphism classes
of bordisms:
\[
    \ho(\bordstr_d)\simeq\cS\mathsf{Cob}_d.
\]
\end{maintheorem}
\end{samepage}

The construction and comparison in \Cref{intro:: thm: stratified bordisms}
are proved in \Cref{ch5:: thm: stratbord interval segal,ch5:: thm: symmetric monoidal stratified,ch5:: thm: homotopy category stratified bord}.
They rest on several geometric results established in
Chapter~4. We prove that every proper conically smooth submersion is a
conically smooth fibre bundle, together with the corresponding
boundary-compatible statement needed for bordisms. We also prove that every
finite conically smooth stratified space admits a proper conically smooth
embedding into some Euclidean space and that the stabilized embedding space
\[
    \Strat^{\sf emb}(M,\RR^{\oo})
\]
is weakly contractible. For compact depth-$1$ spaces of dimension $n$, we
obtain the sharper embedding dimension $\RR^{2n+1}$. These results supply the
local triviality and embedding-independence required in the construction of
$\bordstr_d$.

The remainder of the introduction explains the three choices underlying
these results and places them in relation to earlier constructions.

\emph{Why do we use a different model for higher categories?} The theory of $(\oo,d)$-categories is already filled with an abundance of papers and books covering a variety of different models for $(\oo,d)$-categories. Ours is extremely close to that of complete Segal spaces — first articulated in \cite{Rezk01}.\footnote{A thorough introduction to the zoo of models for the homotopy theory of $(\oo,1)$-categories, we highly recommend the comprehensive volume of Julia Bergner \cite{Ber18}} Recall that a complete Segal space is a pre-sheaf 
\begin{equation}
    \cC \maps \bDelta^{\op} \longra \sset
\end{equation}
which satisfies the Segal condition 
\begin{equation}
    \cC([n]) \,\, \simeq \,\, \cC([1]) \times_{\cC([0])} \cC([1]) \times_{\cC([0])} \cdots \times_{\cC([0])} \cC([1]).
\end{equation}
Since the collection of all objects $\cC([0])$ forms not a set but a space, the \emph{completeness condition} merely ensures that this space agrees with the space of equivalences of $\cC$, so category theory and homotopy theory work in tandem. We introduce a variant of the simplex category $\bDelta$ by replacing finite linearly ordered sets with something slightly more geometric: open intervals in the closed interval $[0,1]$. We introduce a category $\Int$, whose objects are finite disjoint unions of open subintervals containing both endpoints, but different from $[0,1]$; morphisms are inclusions. This admits a functor to $\bDelta^{\op}$ which can be sketched as follows:
\begin{equation}\label{intro:: eq: interval_object}
    I \,\, = \,\, \tikzfig{figures/Ch2/INT1}{\input{figures/Ch2/INT1.tikz}}
\end{equation}
We prove that functors from $\Int$ to $\sset$, mapping inclusions which preserve connected components to weak equivalences, form a model category that is Quillen equivalent to a projective/injective model category of functors from $\bDelta^{\op}$ to $\sset$ — where in both cases the weak equivalences are defined levelwise. This is \emph{almost} equivalent to stating that $\bDelta^{\op}$ is an $\oo$-localization of $\Int$ relative to component-preserving inclusions. From this Quillen equivalence we articulate an adapted Segal condition which we call interval Segal condition, and prove that Segal spaces are Quillen equivalent to interval Segal spaces.

The reason we develop this model for $(\oo,1)$-categories, instead of using complete Segal spaces, is twofold:
\begin{enumerate}
    \item It is intrinsically geometric in nature, which better aligns with our intent of constructing an $(\oo,1)$-category of cobordisms. In past constructions, this interval structure has always played an adjacent role rather than an intrinsic one. It is the author's belief that interval Segal spaces offer a more suitable canvas for building $(\oo,1)$-categories of a geometric nature.
    \item It generalizes well for higher categories, i.e.,\ $(\oo,d)$-categories. As we expand on in the Further Outlook chapter, at the end of the thesis, $\Int$ admits a natural extension to cubes which can just as well parametrize $\btheta_d$ or $\bDelta^{\times d}$ — serving as appropriate shapes to model $(\oo,d)$-categories or cubical versions thereof.
\end{enumerate}

A last --- yet paramount — technical reason for our shift in model, has to do with how well it blends with the spatial structure, i.e.,\ the $\oo$ part of the definition of bordism categories. It has been diagnosed, first in the work of Violeta Marques and Gustavo Granja, that the most state-of-the-art definition of $(\oo,d)$-category of fully-extended cobordisms, given by Claudia Scheimbauer and Damien Calaque \cite{CS19}, suffers from two subtle technical flaws: 

\begin{enumerate}
    \item It assumes that the smooth singular complex of a smooth manifold with corners is a Kan complex — which is proved not to be case in the work of Dan Christensen and Enxin Wu \cite{CW}. See \cite[Section~4.1]{CS19}.
    \item It states that the smooth singular simplices of the topological space of submanifolds satisfy a local triviality condition which would ensure that they form a Kan complex, but they do not satisfy this condition. See \cite[Lemma~5.7]{CS19}.
\end{enumerate}

To provide some additional context for the importance of this issue, we will first discuss how modern constructions of higher categories of bordisms are usually set up. We are used to formulating bordism categories from the \emph{bottom up}: in other words, we define the bordism category as a $k$-category (or $(\oo,k)$-category) whose $k$-morphisms are $(d-k)$-dimensional manifolds with corners and whose composition is given by gluing along the ``right'' boundary components. While this interpretation is retained in the works of Lurie and Calaque--Scheimbauer, their approach can be described as \emph{top down}. As is well documented in the literature, one of the most contrived aspects of constructing higher categories of bordisms is defining composition, since it comes with a compounding pile of choices whose compatibility is easy to check, but very hard to write down. We will first go over a sketch of the definition of Lurie, then we will explain how Calaque--Scheimbauer (CS) address the problems in Lurie's construction and expand on how ours fixes the issues with CS indicated above.

Sticking to the $(\oo,1)$-case, for simplicity, one might describe a point in Lurie's space of $k$-composable chains of 1-morphisms as a choice of a finite-dimensional real vector space $V$ and an embedding $M \hookra V \times \RR$, where $M$ is a $d$-dimensional smooth manifold together with a set $\{t_0\leq \cdots \leq t_k\} \subset \RR$ of $k$ points whose restricted projection  $\pi_{|M}\maps M \fib \RR$ is proper and such that
\begin{equation}
    \mathsf{Crit}(\pi_{|M}) \cap \{t_0\leq \cdots \leq t_k\} = \emptyset.
\end{equation}
Here $\mathsf{Crit}(\pi_{|M})\subset\RR$ denotes the set of critical
\emph{values} of the restricted projection. We depict one such example for $k=2$
\begin{equation}
    \scalebox{0.8}{\tikzfig{figures/Ch3/BOR6}{\input{figures/Ch3/BOR6.tikz}}}
\end{equation} 
\medskip

This set of embeddings can be topologized appropriately and a valid model for the correct $(\oo,1)$-category of bordisms might be extracted. Despite this, the Segal condition proves hard to check, and it does so for reasons which are unsurprising at this stage: cutting manifolds across a single submanifold, or gluing them if you see it as the reverse operation, is easier to do if you first specify collars. Failure to do so begets issues of coherence and pseudo-isotopy — which, although not intractable, increase the burden of working with such models. This issue can be ameliorated by placing the legwork elsewhere, and this is at the core of the variant introduced in the work of Calaque--Scheimbauer. 

To build collars into the definition, Calaque--Scheimbauer introduce an auxiliary $(\oo,1)$-category of closed intervals $\cI$. Calaque--Scheimbauer denote this category by $\mathrm{Int}$; we use $\cI$ to distinguish it from our interval category $\Int$. Assuming, for the sake of simplicity, that there were no issues with this definition, then as a complete Segal space its $k$-levels are defined as $(k+1)$-tuples 
\begin{equation}
    \cI_k = \{([a_0,b_0]\leq \cdots \leq [a_k,b_k])\,:\, [a_i,b_i]\subset \RR \}
\end{equation}
topologized via the standard Euclidean topology on $\RR^{2k+2}$,  and this is where the technical issue appears — in order for these spaces to organize into a Segal space, their smooth singular complex must be a Kan complex, and for $k > 1$ these are not; as $\cI_{k}$ is a smooth manifold with corners. Working with these base intervals, it improves on Lurie's definition by pairing the embedding $M \hookra V\times \RR$ with a choice of a $k$-simplex $(I_0\leq \cdots \leq I_k)\in \cI_k$ such that $\mathsf{Crit}_{|\pi}$ is instead disjoint from $I_0 \cup \cdots \cup I_k$. Visually, we can depict a $k$-simplex as follows 
\begin{equation}
    \scalebox{0.8}{\tikzfig{figures/Ch3/BOR7}{\input{figures/Ch3/BOR7.tikz}}}
\end{equation}
\medskip

\nid where the red shaded regions in $M$ are manifestly the pre-images of $\pi_{|M}$ over the intervals $I_i$, and correspond to the collars along which the bordism decomposes. This definition has a slightly larger amount of overhead than Lurie's, but it offers a clearer way of checking the Segal condition and it fixes a failure of essential constancy in Lurie's model. In particular, it is claimed that the spatial $l$-simplices are obtained as the smooth singular complex of a moduli space of embedded submanifolds crossed with $\cI_{k}$. There are two issues here: first, it is not clear to us that what smooth structure is being assumed for the moduli space of embedded submanifolds, calling into question whether the resulting simplicial set is a Kan complex; second, it is stated that these $l$-simplices can be described as the those satisfying conditions (1)--(3) in \cite[Definition~5.6]{CS19}, but it is not clear how condition 2b. would result from taking smooth singular simplices, and we argue that if condition 2b. is assumed then the Kan condition must fail, as one could find examples of horn conditions which do not admit fillings. This issue lies at the core of our choice to change both the model for $(\oo,1)$-category as well as our approach to the spatial structure.

Our definition tries to reduce the extra machinery in Calaque--Scheimbauer's interval $(\oo,1)$-category $\cI$ by, once again, placing it elsewhere. We use interval Segal spaces and we no longer talk about $k$-simplices as direct encoders of chains of composable morphisms, instead we use a disjoint union of open intervals $I \subset [0,1]$ containing the endpoints. An $I$-level element will be, similarly to the above constructions, an embedding $M \hookra V\times [0,1]$, which we depict as
\begin{equation}
    \scalebox{0.8}{\tikzfig{figures/Ch3/BOR8}{\input{figures/Ch3/BOR8.tikz}}}
\end{equation}
\medskip

\nid where the restricted projection $\pi_{|M}\maps M \fib [0,1]$ must be proper and such that $\mathsf{Crit}(\pi_{|M})$ is disjoint from $[0,1]\setminus I$. Note that the red shaded regions of $M$ correspond to the pre-image of $I$ under the restricted projection $\pi_{|M}$, and they contain the critical points of the map --- it is over these regions that the cobordism will undergo meaningful changes in topology as it may contain any finite number of critical points with respect to the projection --- on the other hand, in its complementary regions the property outlined above implies that the embedding must be diffeomorphic to a cylinder of a $(d-1)$-dimensional manifold. Lastly, we note that Calaque--Scheimbauer need their simplicial sets to be Kan complexes so as to make their proof of the Segal condition work, while we do not. We instead rely on the work of Berwick-Evans--Boavida de Brito--Pavlov (see \cite{BBEP24}) to work around this issue, while still maintaining a similar approach as CS to the proof of the Segal condition. 

\emph{Why do we use conically smooth stratified spaces?}

A central choice in this thesis is to model stratified bordisms using conically smooth stratified spaces, in the sense of Ayala--Francis--Tanaka--Rozenblyum, rather than using more classical theories of stratified spaces, such as Whitney stratified spaces \cite{M12,Tro20}. From our perspective, the point is that conically smooth stratified spaces are singular spaces for which a substantial amount of differential topology survives. They admit a fully developed differential-topological package: one has local models by basics, good notions of smooth map and smooth family, constructible bundles, a robust theory of links, a method called unzipping which allows for functorial resolutions of singularities, tubular neighborhoods, as well as suitably stratified notions of tangent data and tangential structures. In particular, they are especially well tailored to the study of objects in families, which is indispensable for moduli-type constructions like the ones necessitated by modern derived versions of bordism categories. Just as importantly for the purposes of this thesis, this package can be pushed further in directions directly relevant to bordism theory: in Chapter 4 we prove a stratified analogue of the Ehresmann fibration theorem, prove that spaces of conically smooth embeddings into $\mathbb{R}^{\infty}$ are weakly contractible, and establish both weak and depth-$1$ strong Whitney-type embedding results for conically smooth stratified spaces. These results are not ancillary to our later constructions, but rather part of the geometric infrastructure that makes them possible. This is one of the main reasons the conically smooth framework is preferable, for our purposes, to more classical approaches such as Whitney stratifications, where the geometry is often organized through an ambient embedding and is correspondingly less well adapted to intrinsic family-based constructions; see, for instance, the introduction of \cite{AFT1,AFR1}.

Another important point is that conically smooth stratified spaces live within a particularly good categorical context. Rather than just treating stratifications as auxiliary structure on an underlying topological space, they provide a framework in which exit-paths, constructible tangent bundles, and tangential structures fit naturally together. This glue is crucial for the development of stratified bordism categories in the last chapters, where one wants not only to remember how strata meet, but also how singularities vary in families, and how geometric structures behave under pullback. In this respect, the conically smooth formalism is also better suited to our aims than approaches of a more stratified-homotopy-theoretic flavour, such as work of Hughes and others (see \cite{Hu99}), where the emphasis falls instead on homotopy links, path-space constructions, and neighborhood theory. Those approaches are highly natural for questions in stratified homotopy theory, but the present thesis requires a more rigid differential-topological setting, one in which singularities still carry enough geometry to support bordisms, tangential structures, and their families. For this reason, conically smooth stratified spaces provide the most natural ambient category for the constructions undertaken here.

\emph{Why do we restrict to topological manifolds?}
Our construction of stratified bordisms is carried out only for conically smooth stratified spaces whose underlying topological space is a topological manifold of fixed dimension. This is, in part, a conceptual choice. The bordisms relevant for defect topological field theory are meant to model manifold-like bulk geometries together with lower-dimensional strata on which defects live. Requiring the underlying space to remain locally Euclidean ensures that the stratification is encoding controlled singular behaviour inside an ambient topological manifold, rather than replacing the ambient geometry by a more general singular space altogether.

It is also a technical choice. The machinery on which our construction relies — most notably the stratified Ehresmann theorem and the resulting control over families of stratified bordisms — is currently available in this setting, and it is precisely this control that allows the interval Segal machinery to run. We do not view this as a genuine limitation of the theory, only as a first natural range in which the differential-topological and higher-categorical aspects of the construction meet cleanly. Indeed, the formalism developed here appears flexible enough to eventually encompass conically smooth stratified spaces whose underlying space is not a topological manifold.

\section*{Related work}

\nid This thesis lies at the intersection of three bodies of work: the construction of higher bordism categories, the differential topology of conically smooth stratified spaces, and the study of topological field theories with defects. The aim of this section is not to survey these subjects in full, but rather to locate the present construction among the works that are closest either in spirit or in kind.

The foundational backdrop is the standard lineage of higher bordism categories: the embedded topological bordism category of Galatius--Madsen--Tillmann--Weiss \cite{GMTW09}, Lurie's complete $n$-fold Segal-space sketch of the fully extended bordism category \cite{Lur09}, and the detailed Calaque--Scheimbauer construction with tangential structures \cite{CS19}. As explained earlier in this introduction, it is in the Calaque--Scheimbauer framework that the relevant loss-of-fibrancy issue enters: the Segal argument relies on spatial simplicial sets with Kan-complex behaviour, while the smooth singular complexes and moduli spaces involved do not supply this automatically. The construction in this thesis is closest in spirit to this lineage, but uses a new model for $(\oo,1)$-categories which we name interval Segal spaces in contrast to the usual Segal spaces. It makes use of concordance resolutions together with the classifying-space formalism of Berwick-Evans--Boavida de Brito--Pavlov \cite{BBEP24}.

To define interval Segal spaces we introduce a category of intervals $\Int$. This belongs to a broader pattern in which geometric indexing categories replace more rigid combinatorial ones, such as $\Delta$, and one proves that the extra geometric data is homotopically inessential. In stratified factorization homology, Ayala--Francis--Tanaka already relate disk categories to localizations of un-topologized disk categories by inverting isotopy equivalences \cite{AF3}; a close analogue also appears in the work of Karlsson--Scheimbauer--Walde, who prove that the $\oo$-category of disks is a localization of sufficiently fine posets of disks, and relate constructible factorization algebras to the $\oo$-operad of embedded disks \cite{KSW25}. The localization statements for $\Int$ in this thesis should be read in this same spirit: interval data is kept because it is useful for cutting geometric families of bordisms, and then shown not to change the intended homotopy type.

The second source of comparison is the theory of conically smooth stratified spaces developed by Ayala--Francis--Tanaka and Ayala--Francis--Rozenblyum \cite{AFT1,AFR1,AFR2}. Their work provides, not only inspiration, but also the geometric language used throughout: basics, conically smooth maps, constructible bundles, exit-path categories,\footnote{The exit-path category goes back to MacPherson's point of view on constructible sheaves \cite{Mac91}; for published higher-categorical treatments see Lurie \cite{LurieHA} and Treumann \cite{Tre09}.} links, unzipping, and stratified tangential structures. From the perspective of this thesis, conical smoothness serves us not just with a convenient regularity condition but the structure that affords working in families — a primary consideration when extending upwards the ordinary category of stratified bordisms. The theorem of Nocera--Volpe that Whitney stratifications are conically smooth \cite{NV23} further explains why this framework remains close to the classical examples of stratified geometry, while being better adapted to intrinsic moduli constructions. This thesis does not propose a new foundation for stratified spaces. Rather, it uses the conically smooth foundation to construct a symmetric monoidal $(\oo,1)$-category of stratified bordisms and sketch a higher, extended version thereof.

The third source of comparison is the literature on defect topological field theories. Carqueville--Meusburger--Schaumann and Carqueville--Runkel--Schaumann define defect TQFTs as symmetric monoidal functors out of categories of decorated stratified bordisms, and analyze the algebraic structures arising from such theories \cite{CMS20,CRS19}. More recent work, such as the survey of Carqueville--Del Zotto--Runkel \cite{CZR25} and the higher-categorical framework of Müller \cite{Mu25}, places topological defects in a broad field-theoretic context. In comparison, this thesis focuses on the geometric source side. Its goal is not to classify defects algebraically, but to build an $(\oo,1)$-categorical bordism object from which derived defect TQFTs can be defined as symmetric monoidal functors.

Two works are especially close in spirit to the present construction. The first is Will Stewart's thesis on topological domain walls and relative field theories \cite{Ste24}. In that model, the smooth bordism machinery of Calaque--Scheimbauer is adapted to domain walls by replacing ordinary proper submersion conditions with flagged versions thereof. The similarity with our work is that both approaches aim to encode singular or defect-like geometric behavior into bordism constructions. There are, however, two related differences in scope and method. First, Stewart focuses on a restricted and carefully controlled class of domain-wall singularities, whereas the present work builds bordisms directly from conically smooth stratified topological manifolds and allows more general singularity types through the chosen category of basics, with defect labels later encoded through tangential data. Second, Stewart remains close to the flagged version of the Calaque--Scheimbauer model, while this thesis changes the categorical model itself: it uses interval Segal spaces and constructs the spatial levels through concordance resolutions, smooth $\oo$-sheaves, and conically smooth families. Thus the two constructions should be read as closely related but technically distinct approaches to defect-like bordism geometry, rather than as competing treatments of the same model.

The second is the work of Tetik on linked topology \cite{Tet23,Tet22}. Tetik studies depth-$1$ stratified phenomena by replacing a stratified space with a span of ordinary smooth manifolds
\[
M \xleftarrow{\pi} L \xrightarrow{\iota} N,
\]
where the map $\pi$ is a fibre bundle and $\iota$ is a closed embedding. This allows part of the topology of conically smooth stratified spaces, including exit-path quasi-categories and tangential structures, to be expressed using ordinary smooth data. The comparison with conically smooth stratified spaces is already part of Tetik's work: the linked exit-path construction induces a fully faithful functor whose essential image contains the Lurie--MacPherson exit-path quasi-categories of conically smooth spaces over depth-$1$ posets, and the tangential theory of linked smooth manifolds recovers and generalizes the corresponding depth-$1$ conically smooth tangential theory of Ayala--Francis--Rozenblyum. I am not aware of an independent published comparison beyond these results, but this already suggests that linked topology is, at depth $1$, a broader span-based framework in which the conically smooth examples sit naturally. The overlap with this thesis is substantial at the level of motivation: both works are concerned with links, exit paths, tangential structures, and defects. The difference is again one of geometric strategy. Tetik's linked spaces offer a depth-$1$ reduction of stratified geometry to spans of smooth manifolds; this thesis instead keeps conically smooth stratified spaces as the primary geometric objects and constructs an interval Segal space of stratified bordisms from them.

Finally, factorization homology provides a complementary mechanism for producing field theories from local algebraic data. The works of Ayala--Francis, Ayala--Francis--Tanaka, Ayala--Francis--Rozenblyum, and Scheimbauer show how factorization-homological methods organize field theories and observables on manifolds, including stratified manifolds \cite{AF1,AF3,AFR2,Sch14}. In this thesis, factorization homology enters mainly as a future direction: once a suitable stratified bordism category has been constructed, one expects stratified disk-algebras and constructible factorization homology to produce examples of derived defect TQFTs.

I also record several contributions from the unpublished manuscript of Manuel Araújo and Pedro Boavida de Brito, without which the fulfilment of this thesis would not have been possible. The idea of the interval category and its localization as a model for the simplex category underlies Chapter 2. The use of the interval category as a foundation for a variant of complete Segal spaces suited to the $(\oo,1)$-category of bordisms, together with the observation that weak Kan fibrations can be used to transport the proof of \cite[Theorem 5.1]{BBEP24} to the interval Segal condition, enters in Chapter 3. Finally, the notion of the $\cube_n$ category, its localization as a model for the $\bTheta_n$ category, and the use of Segal $\cube_n$-spaces as a foundation for constructing $(\oo,n)$-categories of bordisms enter in Chapter 6.

\section*{Outline}

In Chapter 2 we introduce a model for $(\oo,1)$-categories that complements the existing literature, which we call \emph{interval Segal spaces} — these are simplicial co-presheaves on $\Int$, a subcategory of the poset of opens of $[0,1]$. In particular, we construct a model structure on $\mathsf{coPsh}(\Int)$ which is Quillen equivalent to that of Rezk’s Segal spaces. The motivation here is geometric: most models of $(\oo,1)$-categories are not tailored for the construction of geometric moduli spaces. A rare exception is the recent work of Ayala--Francis--Rozenblyum \cite{AFR1}, where the authors present \emph{striation sheaves} — a highly bespoke class of constructible sheaves on stratified spaces, designed precisely for constructing $(\oo,1)$-categories of a moduli flavour using point-set data. Our construction attempts a similar feat, while being less ambitious and remaining closer to more traditional Segal spaces, thereby facilitating higher-dimensional versions thereof.

In Chapter 3 we present a new construction of the $d$-dimensional bordism $(\oo,1)$-category. This construction is a spiritual cousin of other works, namely those of Galatius--Madsen--Tillman--Weiss \cite{GMTW09}, which later inspired the work of Lurie \cite{Lur09}, and more recently that of Scheimbauer--Calaque \cite{CS19}.

The pioneering work of Galatius--Madsen--Tillman--Weiss offered the first tools to construct an $(\oo,1)$-category of cobordisms. In their paper on the homotopy type of the cobordism category, they construct a topologically enriched category of $d$-dimensional bordisms using classifying spaces. Building on their work, Lurie, in his sketch of a proof of the cobordism hypothesis \cite{Lur09}, put forth a first construction of a model for an $(\oo,n)$-category of cobordisms using $n$-fold complete Segal spaces. This construction contained a small mistake preventing it from satisfying the \emph{essential constancy} condition. In order to address this, in her PhD thesis and in a subsequent paper, Scheimbauer--Calaque propose an alternative definition, which fixes Lurie’s while retaining the use of $n$-fold complete Segal spaces. Despite this almost successful attempt at a \emph{corrigendum}, the work of Scheimbauer--Calaque contains a separate flaw in its claim that the simplicial sets modelling the spatial levels of its $n$-fold Segal spaces are Kan complexes — a claim which, if false, translates into a loss of fibrancy necessary for the proof of the Segal condition.

To remedy this, we make use of the theory of smooth $\oo$-sheaves, i.e. $\oo$-sheaves on the site of smooth manifolds, and their respective classifying spaces — as articulated in the work of Berwick-Evans--Boavida de Brito--Pavlov \cite{BBEP24} — as a conduit for constructing appropriate moduli spaces of smooth manifolds, smooth bordisms, and smooth families of composed bordisms. This gives our approach a distinct flavour and clarifies the minor pathologies identified in previous constructions.

Let us sketch the path we will take in this section. We start by defining a general construction $\widetilde{\bord}_{d,V}(-,-)$ which is a bifunctor on $\man^{\op}\times \Int$ that, when fixing $T\in \Int$, produces a sheaf of sets on $\man^{\op}$ and captures a type of moduli space construction for $d$-dimensional smooth manifolds embedded in $V$ and fibred over smooth manifolds of varying dimension — designed so that, upon taking its concordance resolution and further restricting along the cosimplicial manifold $\bDelta \ra \man$, we obtain a simplicial set $\bord{d,V}$, or space, which we show is the classifying space of $(d-1)$-dimensional manifolds, $d$-dimensional bordisms, and composites of $d$-dimensional bordisms. Lastly, the most dense and technical part of this chapter is spent showing that this construction, a pre-interval space, is an interval Segal space which recovers the $d$-dimensional bordism category $\mathsf{Cob}_d$ as its homotopy 1-category.

In Chapter 4 we finally turn to the theory of \emph{conically smooth stratified spaces}. These were first introduced in the work of Ayala--Francis--Tanaka and further elaborated in the subsequent works of Ayala--Francis--Rozenblyum \cite{AFT1, AFR1, AFR2}. They form the geometric foundations for our stratified higher bordism categories, and their adoption for that role is undoubtedly one of the more novel features of this thesis.

This chapter contains a heavier dose of exposition. We introduce the necessary differential-topological toolkit for our construction of the stratified bordism category, namely: we introduce an $(\oo,1)$-category of conically smooth stratified spaces; the exit-path $(\oo,1)$-category; and the constructible tangent bundle — all of which are necessary for a discussion of stratified tangential structures thereon. We prove new theorems about the contractibility of spaces of embeddings of stratified spaces into Euclidean space, and we prove a stratified analogue of the famed Ehresmann fibration theorem — which is the key to unlocking comparisons with the infrastructure developed for the unstratified case.

In Chapter 5 we construct an interval Segal space, i.e. an $(\oo,1)$-category of conically smooth stratified $d$-dimensional bordisms. We begin by introducing the notion of a conically smooth stratified topological $d$-manifold by stipulating a suitable $(\oo,1)$-category of basics $\Bsc_{d}^{\cS_{d-1}}$, whose links are stratified spheres of appropriate codimensions. Secondly, we produce a pre-interval space $\cS \bord_d$ of stratified topological $d$-dimensional bordisms, again using the concordance resolution for smooth $\oo$-sheaves and a restriction along the cosimplicial smooth manifold $\bDelta \longra \man$, and show that it is an interval Segal space.

Lastly, we endow this interval Segal space with a symmetric monoidal structure agreeing with the expected one coming from disjoint union. After defining this symmetric monoidal interval Segal space, we initiate a discussion on suitable notions of tangential structures for stratified bordisms — in particular, we introduce the notion of \emph{oriented discrete defect datum} by way of discrete right fibrations over $\Bsc_{d}^{\cS_{d-1}}$, as an attempt to make contact with the defect data for stratified bordisms introduced in the work of Carqueville et al. (see \cite{CRS19}) and we prove a comparison result for the stratified oriented case $d=2$.

We conclude Chapter 5 by defining a (derived) \emph{stratified topological quantum field theory} as a symmetric monoidal functor from the $(\oo,1)$-category of stratified bordisms to a symmetric monoidal $(\oo,1)$-category $\cC$.

Chapter 6 is instead prospective in nature. In its first part, following ideas in unpublished work of Manuel Araújo and Pedro Boavida de Brito, we sketch a higher-dimensional generalization of $\Int$, denoted $\cube_d$, with the aim of extending interval Segal techniques to the study of $(\oo,d)$-categories. We conjecture that simplicial co-presheaves on $\cube_d$ satisfying a suitable higher-dimensional Segal condition should present the $(\oo,1)$-category $\Cat_{(\oo,d)}$, and we outline how this formalism could be used to construct a $\cube_d$-space of stratified bordisms, hence providing a candidate model for fully extended stratified topological quantum field theories \`a la Atiyah--Baez--Dolan--Lurie. In the second part of the chapter, we return to the non-extended setting and sketch how stratified factorization homology could provide a mechanism for constructing stratified topological quantum field theories from appropriate stratified disk-algebras, in the spirit of Scheimbauer's factorization-homological approach to topological field theory.

%% file: Chapters/Ch2.tex
\chapter{Models for symmetric monoidal \texorpdfstring{$(\oo,1)$}-categories}
\label{ch2:: title}

The goal of this section is to introduce a model for $(\oo,1)$-categories that complements the existing literature. The motivation here is geometric: most models of $(\oo,1)$-categories are not necessarily tailored for the construction of geometric moduli spaces. The rare exception is the recent work of Ayala--Francis--Rozenblyum \cite{AFR1}, where the authors present \emph{striation sheaves} — a very bespoke class of constructible sheaves on stratified spaces, designed exactly with the intent of constructing $(\oo,1)$-categories of a moduli-flavour using point-set data. Our stated intent is to construct an $(\oo,1)$-category of stratified cobordisms, and to sketch its $(\oo,n)$-version; doing so with complete Segal spaces proves to be fraught with complications.

In \Cref{ch2:: sec: complete_segal_spaces} we introduce the construction of $(\oo,1)$-categories using simplicial spaces i.e., presheaves of simplicial sets, as a way of motivating our construction using open intervals — namely we introduce \emph{complete Segal spaces}.  In \Cref{ch2:: subsec: interval_spaces} we introduce \emph{interval spaces} and provide some intuition for how to think about them. In \Cref{ch2:: subsec: interval_spaces_model_infinity_cats} we prove that they offer a model of $(\oo,1)$-categories, by proving the existence of an explicit Quillen equivalence of model categories, and we discuss important homotopy-categorical invariants associated to interval Segal spaces, like the homotopy category. Lastly, in \Cref{ch2:: sec: symmetric monoidal structure} we discuss symmetric monoidal structures for $(\oo,1)$-categories and construct a model structure whose fibrant objects are symmetric monoidal Segal interval spaces. The ideas in this chapter, namely the reformulation of complete Segal spaces using the interval category, are based on unpublished work by Manuel Araújo and Pedro Boavida de Brito.

\section{Preliminaries on simplicial homotopy theory}
\label{ch2:: sec: preliminaries on simplicial homotopy theory}

We will introduce fundamental objects, definitions, and results and a few important notation conventions which we will stick to throughout the remainder of the thesis.

\subsection{The simplicial indexing category}

Throughout the thesis we will use simplicial sets as avatars or models for spaces, so we begin by introducing them.

\begin{definition}\label{ch2:: def: simplex category}
    We define the \textbf{simplex category} $\bDelta$ as the category whose:
    \begin{itemize}
        \item Objects are finite linearly ordered sets of the form 
        \begin{equation}
        [n] = \{0 < 1 < \cdots < n\}
        \end{equation} for $n\geq 0$.
        \item A morphism from $[n]$ to $[m]$ is a weakly order-preserving function
        \begin{equation}
            \alpha \maps [n] \longra [m].
        \end{equation} That is, a function such that 
        \begin{equation}
            \alpha(i) \leq \alpha(j)
        \end{equation}
        for every $0\leq i < j \leq n$.
    \end{itemize}
\end{definition}

\begin{observation}
    Less frequently, the simplex category is defined as the category whose objects include every finite non-empty linearly ordered set. This category is bigger and manifestly equivalent, as it is its skeleton. More can be said, since it is canonically equivalent: given a finite nonempty linearly ordered set $I$, there is a \emph{unique} weakly order preserving function $I \xra{\cong} [n]$, for some $n\geq 0$.  
\end{observation}

\begin{definition}\label{ch2:: def: join}
    The simplex category carries a (strict) non-unital monoidal structure, the \textbf{join}
    \begin{equation}
        \star \maps \bDelta \times \bDelta \longra \bDelta, 
        \qquad ([m],[n]) \longmapsto [m+n+1].
    \end{equation}
    For morphisms $\alpha \maps [m] \ra [n]$ and $\beta \maps [p] \ra [q]$, their join is the morphism
    \begin{equation}
        \alpha \star \beta \maps [m]\star [p] \longra [n] \star [q], \quad
        k \longmapsto
        \begin{cases}
            \alpha(k), & \text{if } k = 0,\ldots, m, \\
            n + \beta(k-m-1), & \text{if } k = m + 1, \ldots, m + p + 1.
        \end{cases}
    \end{equation}
    If we adjoin an initial object $[-1] = \emptyset$ to obtain the augmented simplex category $\bDelta_{+}$, then this non-unital monoidal structure extends to a unital monoidal structure with unit $[-1]$.
\end{definition}

\begin{proposition}\cite[Section~5, Proposition~1]{Mac}
    The augmented simplex category $(\bDelta_{+},\star,[-1])$ is the free monoidal category generated by a single object.
\end{proposition}

\begin{observation}
    Each $[n]$ is a finite linearly ordered set, hence a partially ordered set, hence a small category, and weakly order-preserving maps are exactly functors between these categories. Thus $\bDelta$ may be viewed as a full subcategory
    \begin{equation}
        \bDelta \longhookra \Cat.
    \end{equation}
\end{observation}

\begin{definition}\label{ch2:: def: simplicial objects}
    Let $\sC$ be a category, a \textbf{simplicial object} in $\sC$ is a functor 
    \begin{equation}
        K_\bullet \maps \bDelta^{\op} \longra \sC.
    \end{equation}
    The analog functors whose source is $\bDelta$ are called \textbf{cosimplicial objects.}
\end{definition}
\begin{notation}
    The $\bullet$ in the subscript of a simplicial object is meant as a placeholder. The image of $K$ on $[n]\in \bDelta^{\op}$ will be denoted $K_n$. We will not always be consistent with the notation $K_\bullet$ for simplicial objects, at some point we will drop it and use it only under specific circumstances.
\end{notation}
\begin{notation}
    Whenever we want to talk about the category of simplicial objects in $\cF$, we will write 
    \begin{equation}
        \mathsf{s}\sC := \fun(\bDelta^{\op},\sC).
    \end{equation}
    for the functor category from $\bDelta^{\op}$ to $\sC$. In particular, the case of simplicial sets $\sset$ will show up very frequently. When we want to emphasize that the target is already $\sset$, we will write
    \begin{equation}
        \mathsf{s}\sset := \fun(\bDelta^{\op},\sset)
    \end{equation}
    and refer to its objects as \textbf{simplicial spaces}.
\end{notation}

\begin{definition}\label{ch2:: def: face and deg maps}
    The simplex category $\bDelta$ has two important sets of distinguished morphisms. Let $n>0$, the first class $\{d^{i}\maps [n-1] \longra [n]\}_{0\leq i \leq n}$ is defined as 
    \begin{equation}
        d^{i}(j) =  
        \begin{cases}
        j & \text{ if } j < i \\
        j + 1 & \text{ if } j \geq i,    
        \end{cases}
    \end{equation} 
    and we call these the \textbf{face maps}.

    The second distinguished class of maps $\{s^{i}\maps [n+1] \longra [n]\}_{0\leq i \leq n}$ are defined as 
    \begin{equation}
        s^{i}(j) = 
        \begin{cases}
            j & \text{ if } j \leq i \\
            j - 1 & \text{ if } j > i,
        \end{cases}
    \end{equation}
    and we call these \textbf{degeneracy maps.}
\end{definition}

Note that we abuse notation when we write $d^{i}$ and $s^{i}$ without explicit reference to the source and target objects. We do this for mere notational parsimony and with the understanding that it is unlikely to hamper readability. The face and degeneracy maps satisfy some long but easy-to-check relations, which we will now enumerate for completion, together with a characterization result.

\begin{proposition} 
    The simplex category can be presented as the category whose objects are the nonempty finite totally ordered sets $[n]$, and the morphisms are generated by the face and degeneracy maps $d^{i},s^{i}$, together with the following set of relations:
    \begin{align*}
        d^{i}\circ d^{j} &= d^{j+1} \circ d^{i} \quad \text{for } i\leq j, \\
        s^{j}\circ s^{i} &= s^{i} \circ s^{j+1} \quad \text{for } i\leq j, \\
    \end{align*}
    and
    \begin{equation}
        s^{j}\circ d^{i} = 
        \begin{cases}
           d^{i} \circ s^{j-1} & \text{ if }  i < j \\
           \id_n, & \text{ if } i - j \leq 1 \\
            d^{i-1} \circ s^{j} & \text{ if } j + 1 < i,
        \end{cases}
    \end{equation}
    which we will call the \textbf{simplicial relations or identities}.
\end{proposition}

\begin{notation}
We can describe maps in $\bDelta$ using a more verbose and explicit notation. Let $f\maps [n] \longra [m]$ be an order-preserving map which sends $k \mapsto i_k$, for $0\leq k \leq n$ and $0\leq i_k \leq n$. We must have that $i_0 \leq i_1 \leq \cdots \leq i_n$, so we define 
\begin{equation}
    f \, \overset{\sf def}{=} \, \langle i_0, i_1, \cdots, i_n \rangle.
\end{equation}  
As an example, the face map $d_1 \maps [1] \longra [2]$ is written as $\langle 0,2 \rangle$ and the degeneracy map $s \maps [2] \longra [1]$ is written as $\langle 001 \rangle$, which in our opinion often improve readability. 
\end{notation}

\begin{remark}\label{ch2:: rmk: inert active}
    There is a useful alternative factorization system on $\bDelta$ — common-knowledge for users of $\oo$-operads and alike $\oo$-structures, but less so for more pragmatic users of $\oo$-categories. The two classes of morphisms are called \textbf{inert} and \textbf{active}. \emph{Inert} morphisms $\alpha \maps [m] \hookra [n]$ are those that satisfy $\alpha(i+1) = \alpha(i) + 1$ for all $i = 0,\ldots, m-1$, i.e., they are like embeddings of linear subposets. \emph{Active} morphisms $\beta \maps [m] \ra [n]$ are those that preserve extrema $\beta(0) = 0$ and $\beta(m) = n$.
\end{remark}

\subsection{Simplicial sets}

We now recall the basic comparison between simplicial sets and topological spaces. We introduce the Kan--Quillen model structure on $\sset$ and the adjunction relating it to the usual homotopy theory of spaces.

\begin{definition}\label{ch2:: def: simplicial set}
    A \textbf{simplicial set} is a simplicial object in $\set$, i.e.\ a functor
    \begin{equation}
        K_\bullet\maps \bDelta^{\op} \longra \set.
    \end{equation} 
\end{definition}

Simplicial sets are ubiquitous in modern homotopy theory. A set $X$ is a simplicial set via the \emph{constant} construction: we set $X_n = X$ for all $n$ and take all face and degeneracy maps to be the identity on $X$. Graphs can also be realised as simplicial sets, with $0$-simplices given by vertices, $1$-simplices encoding edges, and appropriate face maps recording the endpoints.

\begin{definition}\label{ch2:: def: standard simplex}
    Let $n\geq 0$. The \textbf{standard $n$-simplex}, denoted $\Delta^{n}$, is the representable functor
    \begin{equation}
        \Hom_{\bDelta}(-, [n])\maps \bDelta^{\op} \longra \set.
    \end{equation} 
    Its $k$-simplices are given by 
    \begin{equation}
        (\Delta^{n})_k = \Hom_{\bDelta}([k],[n]).
    \end{equation}
\end{definition}

\begin{observation}
    The Yoneda embedding gives a fully faithful functor
    \begin{equation}
        \bDelta \hookra \fun(\bDelta^{\op},\set)
    \end{equation} 
    sending $[n]$ to $\Delta^n$ for every $n\geq 0$, and an order-preserving $f\maps [n] \longra [m]$ in $\bDelta$ to the natural transformation $f^\ast\colon \Delta^{n}\to\Delta^{m}$, which acts on $k$-simplices by 
    \begin{align*}
        \Delta^n &\longra \Delta^m \\
        (\alpha\maps [k] \longra [n]) &\longmapsto ( f\circ \alpha \maps [k] \longra [m]).
    \end{align*}
    In particular, there is a natural isomorphism 
    \begin{equation}
        \Hom_{\sset}(\Delta^n, K) \xra{\,\,\cong\,\,} K_n
    \end{equation} 
    for every simplicial set $K$ and every $n\geq 0$.
\end{observation}

\begin{observation}
    There is a fully faithful inclusion 
    \begin{equation}\label{ch2:: def: set in ssets}
        \set \hookra \sset 
    \end{equation}
    sending a set $S$ to the constant simplicial set $([n] \mapsto S)$.
\end{observation}

\begin{definition}\label{ch2:: def: sset-mapping-space}
    The category $\sset$ is cartesian closed. Given simplicial sets $K,L\in\sset$, we denote by
    \begin{equation}
        \Map_{\sset}(K,L)
    \end{equation}
    the \textbf{simplicial mapping space} whose $k$-simplices are
    \begin{equation}
        \Map_{\sset}(K,L)_k \;\overset{\sf def}{=}\;
        \Hom_{\sset}\big(K \times \Delta^k,\, L\big),
    \end{equation}
    where $\Delta^k$ is the standard $k$-simplex. This construction makes $\sset$ enriched, tensored, and cotensored over itself.
\end{definition}

There are a few additional simplicial sets derived from the standard simplices that will be important for us.

\begin{definition}\label{ch2:: def: boundary n-simplex}
    Let $n\geq 0$. The \textbf{boundary of} $\Delta^{n}$ is the simplicial set $\partial \Delta^{n}$ defined by
    \begin{equation}
        (\partial \Delta^n)([m]) = \{\alpha \in \Hom_{\bDelta}([m],[n])\,\maps\, \alpha \text{ is not surjective}\}.
    \end{equation}
\end{definition}

\begin{definition}\label{ch2:: def: k-horn}
    Let $n > 0$ and let $0\leq k \leq n$. The \textbf{$k$-th horn of the standard $n$-simplex} is the simplicial set $\Lambda_{k}^{n}$ defined by
    \begin{equation}
        (\Lambda_{k}^{n})([m]) = \{\alpha\in \Hom_{\bDelta}([m],[n])\, : \, ([n]\setminus \{k\}) \not\subseteq \alpha([m])\}.
    \end{equation} 
    If $0<k<n$ we call $\Lambda_k^n$ an \textbf{inner horn}, and if $k\in \{0,n\}$ we call it an \textbf{outer horn}. 
\end{definition}

Both the boundary and the $k$-th horn are sub-simplicial sets of the standard $n$-simplex. In particular, the corresponding inclusions of sets give rise to monomorphisms
\begin{equation}
    \partial \Delta^n \hookra \Delta^n,
    \qquad
    \Lambda^{n}_k \hookra \Delta^n.
\end{equation}

As mentioned above, simplicial sets arise as a combinatorial model for topological spaces. We now make this comparison explicit by constructing a topological analogue of the standard $n$-simplex.

\begin{definition}\label{ch2:: def: topological n-simplex}
    Let $n\geq 0$. The \textbf{standard topological $n$-simplex} is the subset
    \begin{equation}
        \Delta^n_{t} = \Big\{(t_0,\ldots,t_n)\in \RR^{n+1}\,\maps \, t_{i} \geq 0,\ \sum_{i=0}^n t_i = 1\Big\}
    \end{equation}
    endowed with the subspace topology from $\RR^{n+1}$.
\end{definition}

\begin{proposition}
    There is a cosimplicial object in topological spaces
    \begin{equation}
        \Delta^{\bullet}_{t} \maps \bDelta \longra \Top
    \end{equation}
    which sends $[n]$ to the standard topological $n$-simplex $\Delta^n_{t}$ and sends an order-preserving map $\alpha\maps [m] \longra [n]$ to the linear map
    \begin{equation}
        \alpha(t_0,\ldots,t_m)_i = 
        \begin{cases}
            0, & \text{if } \alpha^{-1}(i) = \emptyset, \\
            \displaystyle\sum_{j\in \alpha^{-1}(i)} t_j, & \text{if } \alpha^{-1}(i) \neq \emptyset,    
        \end{cases}
    \end{equation}
    for $i=0,\ldots,n$.
\end{proposition}

\begin{definition}
    Given a simplicial set $K\maps \bDelta^{\op} \longra \set$, its \textbf{category of simplices} $\bDelta / K$ is its category of elements.
\end{definition}

Recall from \Cref{ap:: def: cat of els} that the category of simplices associated to $K$ has objects pairs $([n],x)$ with $[n]\in \bDelta^{\op}$ and $x\in K([n])$, and morphisms $([n],x)\to([m],y)$ given by order-preserving maps $\alpha \maps [n] \longra [m]$ in $\bDelta^{\op}$ such that $\alpha^\ast(y) = x$. Crucially, knowing the category of simplices of $K$ is enough to reconstruct $K$.

\begin{lemma}
    Given a simplicial set $K$, there is a natural isomorphism 
    \begin{equation}
        K \,\cong\, \mathsf{colim}\big( \bDelta / K \longra \bDelta \xra{\,\,\Delta^\bullet\,\,} \sset \big).
    \end{equation}
\end{lemma}

The previous lemma states that every simplicial set can be constructed by gluing together copies of standard simplices. This suggests how to construct a topological space from a simplicial set. The next definition makes this precise and provides the bridge between the homotopy theory of spaces and that of simplicial sets.

\begin{definition}
    Let $K$ be a simplicial set. We define its \textbf{geometric realization} $|K|$ as the following colimit of topological spaces:
    \begin{equation}
        |K| = \mathsf{colim}\big(\bDelta / K \longra \bDelta \xra{\,\,\Delta^\bullet_{t}\,\,} \Top\big).
    \end{equation}
    The construction is manifestly functorial, and we write
    \begin{equation}\label{ch2:: eq: geometric real}
        |-| \maps \sset \longra \Top
    \end{equation}
    for the associated functor.
\end{definition}

\begin{observation}
    The category of simplices $\bDelta/\Delta^{n}$ is canonically equivalent to $\bDelta/[n]$ via the Yoneda embedding. Moreover, $\bDelta/[n]$ has a terminal object, namely $[n] \xra{\id} [n]$. Together, these observations imply that 
    \begin{equation}
        |\Delta^n| \simeq \Delta^n_{t}.
    \end{equation} 
\end{observation}

Geometric realization constructs a topological space from a simplicial set. Conversely, we often wish to associate a simplicial set to a given space $X$ that captures its homotopy type. For this we use the singular complex.

\begin{definition}\label{ch2:: def: singular complex}
    Let $X$ be a topological space. We define a simplicial set $\sing(X)$ by
    \begin{equation}
        \sing(X)_n = \Hom_{\Top}(\Delta^{n}_t, X),
    \end{equation}
    and for a morphism $f\maps [m]\longra [n]$ in $\bDelta$ we set
    \begin{equation}
        \sing(X)(f) \colon \sing(X)_n = \Hom_{\Top}(\Delta^{n}_{t},X)
        \longra \Hom_{\Top}(\Delta^{m}_{t},X) = \sing(X)_m
    \end{equation}
    to be precomposition with the induced map $\Delta^m_t \to \Delta^n_t$.
    We call $\sing(X)$ the \textbf{singular complex of} $X$, and this defines a functor
    \begin{equation}
        \sing \maps \Top \longra \sset.
    \end{equation}
\end{definition}

\begin{proposition}
    The geometric realization functor is left adjoint to the singular complex functor:
    \begin{equation}
        \begin{tikzcd}[ampersand replacement=\&]
            \sset \&\& {\Top}
            \arrow[""{name=0, anchor=center, inner sep=0}, "{|-|}", curve={height=-6pt}, from=1-1, to=1-3]
            \arrow[""{name=1, anchor=center, inner sep=0}, "{\sing}", curve={height=-6pt}, from=1-3, to=1-1]
            \arrow["\dashv"{anchor=center, rotate=-90}, draw=none, from=0, to=1]
        \end{tikzcd}
    \end{equation} 
\end{proposition}

\begin{proof}
    The proof is an instance of abstract nonsense. If $K$ is a simplicial set and $X$ is a topological space, then
    \begin{align}
        \Hom_{\Top}(|K|,X)\,\, 
        &\cong \,\, \Hom_{\Top}\Big(\mathsf{colim}_{\bDelta/ K}|\Delta^{n}|, X\Big) \\
        &\cong\,\, \Limit[\bDelta / K]\Hom_{\Top}\big(|\Delta^n|,X\big) \\
        &\cong\,\, \Limit[\bDelta / K]\Hom_{\sset}\big(\Delta^n,\sing(X)\big) \\
        &\cong \,\, \Hom_{\sset}\big(K, \sing(X)\big),
    \end{align}
    where every isomorphism is natural. We have used the fact that $\Hom$ out of a colimit is a limit of $\Hom$-sets, the identification $|\Delta^n|\simeq\Delta^n_t$, and the Yoneda lemma to pass from the third line to the fourth.
\end{proof}


\begin{observation}\label{ch2:: obs: pi_0 simp sets}
    Just as spaces admit a functor $\pi_0 \maps \Top \longra \set$ which records its path components, so do simplicial sets:
    \begin{equation}
        \pi_0 \maps \sset \xra{\,\,\, |-|\,\,\,} \Top \xra{\,\,\, \pi_0 \,\,\,} \set.
    \end{equation}
    For $X\in \sset$, we have that $[x] = [y]$ in $\pi_0(X)$ if and only if there exists a 1-simplex $\gamma \maps \Delta^1 \ra X$ such that its two 0-faces 
    \begin{equation}
    \partial \gamma = \gamma_{|\partial \Delta^1} = \{x,y\}
    \end{equation}
    connecting the two.
\end{observation}

\begin{observation}
    If $S$ and $T$ are discrete simplicial sets, then they are weakly equivalent if and only if $\pi_k(|S|) \cong \pi_k(|T|)$ for every $k \geq 0$. All of these are equivalent to $\pi_0(|S|) \cong \pi_0(|T|)$ being bijective i.e., $S \cong T$.
\end{observation}

\subsection{A primer on model categories}
\label{subsec:model-cats-primer}

The language of model categories provides an axiomatic framework for doing
homotopy theory in settings more general than topological spaces.  Our main
references are Hirschhorn's book \cite{Hir14} and Hovey's monograph on model
categories.  In this subsection we collect the basic definitions and results
we will use later, together with a few guiding examples.

\subsubsection{Model categories and homotopy categories}

\begin{definition}
  A \emph{model category} is a complete and cocomplete category
  $\cM$ equipped with three distinguished classes of morphisms:
  \begin{itemize}
    \item \define{weak equivalences} $\cW$,
    \item \define{cofibrations} $\cC\mathsf{of}$,
    \item \define{fibrations} $\cF\mathsf{ib}$,
  \end{itemize}
  such that the following axioms hold:
  \begin{enumerate}
    \item[(MC1)] (\emph{2-out-of-3}) If $f$ and $g$ are composable morphisms
      in $\cM$ and two of the three maps $f$, $g$ and $g \circ f$ are weak
      equivalences, then so is the third.
    \item[(MC2)] (\emph{Retracts}) The classes $\cW$, $\cC\mathsf{of}$ and
      $\cF\mathsf{ib}$ are closed under retracts.
    \item[(MC3)] (\emph{Lifting}) Cofibrations have the left lifting property
      with respect to acyclic fibrations, and fibrations have the right
      lifting property with respect to acyclic cofibrations, where
      \begin{equation}
        \cC\mathsf{of} \cap \cW = \{\text{acyclic cofibrations}\}, \qquad
        \cF\mathsf{ib} \cap \cW = \{\text{acyclic fibrations}\}.
      \end{equation}
      Concretely, $i\colon A \to B$ is a cofibration if and only if for every
      trivial fibration $p\colon X \to Y$, any solid commutative square
      \begin{equation}
        \begin{tikzcd}[ampersand replacement=\&]
            A \&\& X \\
            \\
            B \&\& Y
            \arrow[from=1-1, to=1-3]
            \arrow["i"', from=1-1, to=3-1]
            \arrow["p", from=1-3, to=3-3]
            \arrow[from=3-1, to=3-3]
            \arrow[dashed, from=3-1, to=1-3]
        \end{tikzcd}
      \end{equation}
      admits a dashed diagonal lift.  Dually for trivial cofibrations and
      fibrations.
    \item[(MC4)] (\emph{Factorization}) Every morphism $f \colon X \to Y$ in
      $\cM$ admits two functorial factorizations
      \begin{equation}
        X \xrightarrow{i} Z \xrightarrow{p} Y,
      \end{equation}
      with $i$ a cofibration and $p$ an acyclic fibration, and
      \begin{equation}
        X \xrightarrow{j} Z' \xrightarrow{q} Y,
      \end{equation}
      with $j$ an acyclic cofibration and $q$ a fibration.
  \end{enumerate}
  Since $\cM$ is complete and cocomplete it has an initial and a terminal
  object.  An object $X \in \cM$ is called \define{cofibrant} if the unique map
  $\varnothing \to X$ is a cofibration, and \define{fibrant} if the unique map
  $X \to *$ is a fibration.
\end{definition}

\begin{lemma}\label{ch2:: lem: fibrant-cofibrant}
    Let $\cM$ be a model category. For every object $X$ there is
    \begin{itemize}
      \item a cofibrant object $\bQ X$ together with a weak equivalence
        $\bQ X \xra{\sim} X$;
      \item a fibrant object $\bR X$ together with a weak equivalence
        $X \xra{\sim} \bR X$.
    \end{itemize}
    We call these \define{cofibrant} and \define{fibrant} replacements of $X$.
\end{lemma}

\begin{proof}
    Factor the unique map $\varnothing \to X$ as
    \begin{equation}
        \varnothing \xra{i} \bQ X \xra{p} X
    \end{equation}
    as in~\textup{(MC4)}, with $i$ a cofibration and $p$ an acyclic fibration.
    Then $\bQ X$ is cofibrant and $p$ is a weak equivalence.  The construction
    of $\bR X$ is dual, factoring $X \to *$ as an acyclic cofibration followed
    by a fibration.
\end{proof}

\begin{definition}[Cofibrantly generated model category]
  A model category $\cM$ is \define{cofibrantly generated} if there exist sets
  of maps
  \begin{equation}
    I, J \subseteq \Mor(\cM)
  \end{equation}
  such that:
  \begin{enumerate}
    \item A map is an acyclic fibration if and only if it has the right
          lifting property with respect to every map in \(I\).
    \item A map is a fibration if and only if it has the right lifting
          property with respect to every map in \(J\).
    \item The domains of maps in \(I\) are small relative to \(I\)-cell
          complexes, and the domains of maps in \(J\) are small relative to
          \(J\)-cell complexes (see \cite[Def.~10.1.2]{Hir14} for the precise
          smallness condition).
  \end{enumerate}
  We call \(I\) a \define{set of generating cofibrations} and \(J\) a
  \define{set of generating acyclic cofibrations}.
\end{definition}

\begin{definition}[Combinatorial model category]
  A model category $\cM$ is \define{combinatorial} if it is cofibrantly
  generated and its underlying category is locally presentable (equivalently,
  is cocomplete and generated under filtered colimits by a set of objects).
\end{definition}

\begin{definition}[Left and right properness]\label{ch2:: def: proper}
  A model category $\cM$ is
  \begin{itemize}
    \item \define{left proper} if every pushout of a weak equivalence along a
          cofibration is again a weak equivalence;
    \item \define{right proper} if every pullback of a weak equivalence along a
          fibration is again a weak equivalence;
    \item \define{proper} if it is both left and right proper.
  \end{itemize}
\end{definition}

\begin{remark}
  In practice we will almost always work with \emph{cofibrantly generated} and
  \emph{combinatorial} model categories, because these admit good existence
  results for localizations and diagram model structures; see
  \cite[Ch.~11, Ch.~13]{Hir14}.  Properness will be important when we want
  homotopy (co)limits to behave well with respect to weak equivalences.
\end{remark}

\begin{theorem}[Quillen model structure on topological spaces]
  \label{ch2:: thm: Quillen model structure}
  There is a proper, cofibrantly generated model structure on the category
  $\Top$ of topological spaces such that:
  \begin{enumerate}
    \item weak equivalences are the \emph{weak homotopy equivalences}, i.e.\
      maps $f\colon X \to Y$ inducing isomorphisms
      \begin{equation}
        \pi_0(X) \xrightarrow{\;\cong\;} \pi_0(Y),
        \qquad
        \pi_n(X,x_0) \xrightarrow{\;\cong\;} \pi_n(Y,f(x_0))
      \end{equation}
      for all $n\geq 1$ and all basepoints $x_0\in X$;
    \item fibrations are the Serre fibrations;
    \item cofibrations are the retracts of relative cell complexes.
  \end{enumerate}
  A set of generating cofibrations is given by
  \begin{equation}
    I = \{S^{n-1} \hookra D^{n}\}_{n\geq 0},
  \end{equation}
  where $S^{n-1}$ is the $(n-1)$-sphere and $D^n$ the $n$-disk, and a set of
  generating acyclic cofibrations is
  \begin{equation}
    J = \{D^{n}\times \{0\} \hookra D^{n} \times [0,1]\}_{n \geq 0}.
  \end{equation}
  We denote this model structure by $\Top_{\mathrm Q}$ and refer to it as the
  \textbf{Quillen model structure}.
\end{theorem}

\begin{theorem}[Kan--Quillen model structure on simplicial sets]
  \label{ch2:: thm: Kan-quillen model}
  There is a proper, cofibrantly generated model structure on the category
  $\sset$ of simplicial sets such that:
  \begin{enumerate}
    \item weak equivalences are those maps $f\colon K \to L$ whose geometric
      realization $|f|\colon |K| \to |L|$ is a weak homotopy equivalence of
      topological spaces;
    \item fibrations are the Kan fibrations;
    \item cofibrations are the monomorphisms.
  \end{enumerate}
  In particular, every simplicial set is cofibrant.  A set of generating
  cofibrations is
  \begin{equation}
    I = \{\partial \Delta^n \hookra \Delta^{n}\}_{n\geq 0},
  \end{equation}
  and a set of generating acyclic cofibrations is
  \begin{equation}
    J = \{\Lambda^{n}_{k} \hookra \Delta^{n}\,:\, 0 \leq k \leq n\}_{n \geq 0}.
  \end{equation}
  We denote this model structure by $\sset_{\mathrm{QK}}$ and call it the
  \textbf{Kan--Quillen model structure}.
\end{theorem}

\begin{remark}
  Theorems~\ref{ch2:: thm: Quillen model structure} and
  \ref{ch2:: thm: Kan-quillen model} provide our basic examples of model
  categories.  We will see in \Cref{ch2:: thm: homotopy hyp} that
  $\Top_{\mathrm Q}$ and $\sset_{\mathrm{QK}}$ model the same homotopy theory.
\end{remark}

The homotopy theory encoded by a model category $\cM$ is captured by its
\emph{homotopy category}.

\begin{definition}
  Let $\cM$ be a model category.  Its \emph{homotopy category}
  $\ho(\cM)$ is obtained from $\cM$ by formally inverting the weak
  equivalences:
  \begin{equation}
    \ho(\cM) \,\overset{\sf def}{=} \, \cM[\cW^{-1}].
  \end{equation}
\end{definition}

This definition is entirely independent of the fact that $\cM$ is a model
category.  In particular, it makes no reference to the axioms.  In general,
localization (the process of formally inverting a class of maps) is
ill-behaved.  One way to think about model categories is precisely as a
framework where localization behaves well.

We can therefore understand $\ho(\cM)$ through a more concrete description
using homotopies.  Let $X$ be cofibrant and $Y$ fibrant.  A
\define{cylinder object} for $X$ is a factorization of the codiagonal
\begin{equation}
  X \sqcup X \longrightarrow X
\end{equation}
as an acyclic cofibration $X \sqcup X \to \mathsf{Cyl}(X)$ followed by a
fibration $\mathsf{Cyl}(X) \to X$.  A \emph{left homotopy} between two
morphisms $f,g \colon X \to Y$ is then a map
\begin{equation}
  H \colon \mathsf{Cyl}(X) \longrightarrow Y
\end{equation}
restricting to $f$ and $g$ on the two copies of $X$.  If $X$ is cofibrant and
$Y$ is fibrant, left homotopy induces an equivalence relation $\sim_l$ on
$\Hom_{\cM}(X,Y)$; see \cite[Proposition~4.7]{Hir14}.

\begin{theorem}\label{ch2:: thm: homotopy cat def}
  Let $\cM$ be a model category and let $X$ be cofibrant and $Y$ fibrant.
  Then the canonical map
  \begin{equation}
    \Hom_{\cM}(X,Y)/\sim_{l}
    \longrightarrow
    \ho(\cM)(X,Y)
  \end{equation}
  is a bijection.  In particular, to compute morphisms in $\ho(\cM)$ it
  suffices to replace objects by weakly equivalent cofibrant–fibrant objects.
\end{theorem}

\begin{cor}[Fundamental theorem of model categories]
  \label{ch2:: cor: fundamental theorem model categories}
  Let $\cM_{\mathrm{cf}} \subseteq \cM$ denote the full subcategory of
  cofibrant–fibrant objects, and let
  \begin{equation}
    \gamma \colon \cM \longra \ho(\cM)
  \end{equation}
  be the localization functor.  Then $\gamma$ induces an equivalence of
  categories
  \begin{equation}
    \cM_{\mathrm{cf}}/\sim_l \xra{\simeq} \ho(\cM),
  \end{equation}
  where $\cM_{\mathrm{cf}}/\sim_l$ has the same objects as
  $\cM_{\mathrm{cf}}$ and morphisms given by left-homotopy classes of maps.
\end{cor}

\begin{remark}
  \Cref{ch2:: cor: fundamental theorem model categories} gives the concrete
  description of the localization promised above.  In particular, for every
  pair of objects $X,Y\in \cM$ there is a natural bijection
  \begin{equation}
    \Hom_{\ho(\cM)}(X,Y)
    \;\cong\;
    \Hom_{\cM}(\bQ X,\bR Y)/\sim_l.
  \end{equation}
  Moreover, if $X$ and $Y$ are cofibrant–fibrant, then a morphism
  $f\colon X \to Y$ is a weak equivalence if and only if there exists a map
  $g\colon Y \to X$ such that
  \begin{equation}
    g f \sim_l \id_X,
    \qquad
    f g \sim_l \id_Y.
  \end{equation}
  In other words, weak equivalences between cofibrant–fibrant objects are
  precisely homotopy equivalences.
\end{remark}

\begin{example}
  The model category $\Top_{\mathrm Q}$ of
  \Cref{ch2:: thm: Quillen model structure} has a homotopy category equivalent
  to the classical homotopy category of CW-complexes and homotopy classes of
  maps.
\end{example}

\begin{example}
  The model category $\sset_{\mathrm{QK}}$ of
  \Cref{ch2:: thm: Kan-quillen model} has a homotopy category equivalent to
  $\ho(\Top_{\mathrm Q})$; this follows from the Quillen equivalence of
  \Cref{ch2:: thm: homotopy hyp} below.
\end{example}

\subsubsection{Quillen adjunctions and Quillen equivalences}

Model categories are meant to be compared by functors that preserve the
homotopy theory they encode.

\begin{definition}
  Let $\cM$ and $\cN$ be model categories.  An adjunction
  \[
    F \colon \cM \rightleftarrows \cN \colon G
  \]
  is called a \emph{Quillen adjunction} (or \emph{Quillen pair}) if the left
  adjoint $F$ preserves cofibrations and acyclic cofibrations.  Equivalently,
  the right adjoint $G$ preserves fibrations and acyclic fibrations.
\end{definition}

Given a Quillen adjunction, one can form its \emph{derived} adjunction on the
homotopy categories.

\begin{definition}
  Let $F \dashv G$ be a Quillen adjunction between model categories.  Choose a
  cofibrant replacement functor $Q \colon \cM \to \cM$ and a fibrant
  replacement functor $R \colon \cN \to \cN$.  The \emph{total left derived
  functor} and \emph{total right derived functor} are the functors
  \begin{equation}
    \LL F \colon \ho(\cM) \longrightarrow \ho(\cN),
    \qquad
    \RR G \colon \ho(\cN) \longrightarrow \ho(\cM),
  \end{equation}
  defined on objects by $\LL F(X) = F(\bQ X)$ and $\RR G(Y) = G(\bR Y)$, and on
  morphisms by applying $F$ or $G$ to representatives between cofibrant or
  fibrant replacements.  The pair $(\LL F,\RR G)$ is again an adjoint pair.
\end{definition}

\begin{definition}
  A Quillen adjunction $F \dashv G$ is called a \define{Quillen equivalence} if
  any (hence all) of the following equivalent conditions hold:
  \begin{enumerate}
    \item For every cofibrant object $X \in \cM$ and every fibrant object
      $Y \in \cN$, a morphism
      \begin{equation}
        F(X) \longrightarrow Y
      \end{equation}
      is a weak equivalence in $\cN$ if and only if its adjoint
      $X \to G(Y)$ is a weak equivalence in $\cM$.
    \item The derived adjunction
      \begin{equation}
        \LL F \colon \ho(\cM) \rightleftarrows \ho(\cN) \colon \RR G
      \end{equation}
      is an equivalence of categories.
    \item The left adjoint $F$ reflects weak equivalences between cofibrant
      objects and, for every fibrant $Y \in \cN$, the composite
      \begin{equation}
        \bQ F(G(Y)) \longrightarrow F(G(Y)) \xrightarrow{\;\varepsilon_Y\;} Y
      \end{equation}
      is a weak equivalence in $\cN$, where $\bQ$ is a cofibrant replacement
      and $\varepsilon$ is the counit.
  \end{enumerate}
\end{definition}

\begin{example}
  The adjoint pair
  \[
    |{-}| \colon \sset_{\mathrm{QK}}
    \rightleftarrows
    \Top_{\mathrm Q} \colon \Sing
\]
  of geometric realization and the singular complex is a Quillen adjunction,
  and indeed a Quillen equivalence (Theorem~\textup{2.1.30}).  Its derived
  adjunction identifies $\ho(\sset_{\mathrm{QK}})$ with $\ho(\Top_{\mathrm Q})$.
\end{example}

\begin{theorem}[Projective and injective model structures on functor categories, cf.\ {\cite[Ch.~11]{Hir14}}]
  \label{ch2:: thm: functor model structures}
  Let $\cF$ be a small category and let $\cM$ be a combinatorial model
  category. Then $\Fun(\cF,\cM)$ admits both the projective and the injective
  model structures, each with weak equivalences defined levelwise. Moreover, if
  $\cM$ is left proper or simplicial, then the corresponding functor model
  structures inherit the same property.
\end{theorem}

\begin{example}
  Let $u \colon \cF \to \cG$ be a functor between small categories, and let
  $\cM$ be a combinatorial model category.  Then precomposition
  \[
    u^{\ast} \colon \Fun(\cF,\cM) \longrightarrow \Fun(\cG,\cM)
  \]
  admits both a left and a right adjoint, given by left and right Kan
  extension along $u$, which we denote by $\lkan_u$ and $\rkan_u$:
    \begin{equation}
        \begin{tikzcd}[ampersand replacement=\&]
        {\fun(\cF,\cM)} \\
        \\
        {\fun(\cG,\cM)}
        \arrow[""{name=0, anchor=center, inner sep=0}, "{u^*}"{description}, from=1-1, to=3-1]
        \arrow[""{name=1, anchor=center, inner sep=0}, "{\lkan_{u}}", bend left=50, from=3-1, to=1-1]
        \arrow[""{name=2, anchor=center, inner sep=0}, "{\rkan_{u}}"', bend right=50, from=3-1, to=1-1]
        \arrow["\dashv"{anchor=center}, draw=none, from=0, to=2]
        \arrow["\dashv"{anchor=center}, draw=none, from=1, to=0]
        \end{tikzcd}
    \end{equation}
  For the \emph{projective} model structures on functor categories,
  $(\lkan_u,u^{\ast})$ is a Quillen pair, and for the \emph{injective}
  model structures $(u^{\ast},\rkan_u)$ is a Quillen pair; to see why, note that if $F\maps \cF'\longra \cM$ is a level-wise weak equivalence and $\alpha \maps c \longra d$ is a morphism in $\cF$, then $u^{*}F(\alpha) = F\big(u(\alpha)\big)$ will be a level-wise weak equivalence i.e., $u^*$ preserves weak equivalences. Moreover, if $F$ is a projective fibration, then $u^* F$ is a projective fibration by the same argument. Dually, $u^*$ preserves injective cofibrations. We will use this repeatedly when passing
  between model structures on diagram categories.
\end{example}

\begin{theorem}[Homotopy hypothesis for spaces vs.\ simplicial sets]
  \label{ch2:: thm: homotopy hyp}
  The adjoint pair
  \begin{equation}
    |{-}| \colon \sset_{\mathrm{QK}}
    \rightleftarrows
    \Top_{\mathrm Q} \colon \Sing
  \end{equation}
  of geometric realization and singular complex is a Quillen adjunction, and
  indeed a Quillen equivalence of model categories.
\end{theorem}

\begin{remark}
  As a consequence, the homotopy category of $\sset_{\mathrm{QK}}$ is
  equivalent to that of $\Top_{\mathrm Q}$.  In particular, the familiar
  homotopy theory of topological spaces can be studied equally well using
  simplicial sets with the notion of weak equivalence from
  \Cref{ch2:: thm: Kan-quillen model}.  From now on we will often refer to
  simplicial sets simply as \emph{spaces}, and we will only distinguish
  between topological spaces and simplicial sets when the context requires it.
\end{remark}

\subsubsection{Homotopy limits, colimits and left Bousfield localization}
\label{ch2:: subsubsec: homotopy-limits-bousfield}

We will need homotopy limits and colimits of diagrams valued in a simplicial
model category.

\begin{definition}[Homotopy (co)limits]\label{ch2:: def: homotopy-(co)limits}
  Let $\cM$ be a model category and let $\cF$ be a small category.
  We endow the diagram category $\Fun(\cF,\cM)$ with one of the standard
  projective or injective model structures from
  \Cref{ch2:: thm: functor model structures}, so that weak equivalences are
  levelwise in~$\cM$.

  The constant diagram functor
  \begin{equation}
    \const \colon \cM \longra \Fun(\cF,\cM)
  \end{equation}
  admits both a left adjoint and a right adjoint, namely the colimit and
  limit functors
  \begin{equation}
    \Colimit[\cF] \colon \Fun(\cF,\cM) \longra \cM,
    \qquad
    \Limit[\cF] \colon \Fun(\cF,\cM) \longra \cM,
  \end{equation}
  and for the chosen diagram model structure:
  \begin{enumerate}
    \item $\Colimit[\cF]$ is a left Quillen functor;
    \item $\Limit[\cF]$ is a right Quillen functor.
  \end{enumerate}

  Let $X \in \Fun(\cF,\cM)$.
  \begin{itemize}
    \item The \define{homotopy limit} of $X$ is the value of the total
      right derived functor of $\lim_{\cF}$:
      \begin{equation}
        \hoLimit[\cF] X
        \;\overset{\mathsf{def}}{=}\;
        \RR\!\lim_{\cF} X.
      \end{equation}
      Equivalently, if $\bR X$ is a fibrant replacement of $X$
      in $\Fun(\cF,\cM)$ for the chosen diagram model structure, then
      \begin{equation}
        \hoLimit[\cF] X
        \;\simeq\;
        \lim_{\cF} (\bR X).
      \end{equation}
    \item The \define{homotopy colimit} of $X$ is the value of the total
      left derived functor of $\colim_{\cF}$:
      \begin{equation}
        \hoColimit[\cF] X
        \;\overset{\mathsf{def}}{=}\;
        \LL\!\Colimit[\cF] X.
      \end{equation}
      Equivalently, if $\bQ X$ is a cofibrant replacement of $X$
      in $\Fun(\cF,\cM)$, then
      \begin{equation}
        \hoColimit[\cF] X
        \;\simeq\;
        \Colimit[\cF] (\bQ X).
      \end{equation}
  \end{itemize}
\end{definition}

When $\cM$ is simplicial, one can also describe homotopy (co)limits via
Bousfield--Kan formulas \cite[Ch.~18]{Hir14}. In many examples of interest
(e.g.\ $\cM=\sset$ and $\cF$ direct), the Reedy model structure on $\Fun(\cF,\cM)$ exists and coincides with the injective model structure,
but we will not use this explicitly.

\begin{definition}[Homotopy Kan extensions]\label{ch2:: def: homotopy-kan}
  Let $\cM$ be a model category, let $\cF,\cG$ be small categories, and let
  \[
    u \colon \cF \longra \cG
  \]
  be a functor.

  \begin{enumerate}
    \item \textbf{Homotopy left Kan extension.}
    Endow $\Fun(\cF,\cM)$ and $\Fun(\cG,\cM)$ with the \emph{projective}
    diagram model structures. Precomposition with $u$,
    \begin{equation}
      u^{\ast} \colon \Fun(\cG,\cM) \longra \Fun(\cF,\cM),
    \end{equation}
    has a left adjoint
    \begin{equation}
      \lkan_u \colon \Fun(\cF,\cM) \longra \Fun(\cG,\cM),
    \end{equation}
    the (ordinary) left Kan extension along $u$.
    In the projective model structure the functor $\lkan_u$ is left
    Quillen, so we can form its total left derived functor
    \begin{equation}
      \Llkan \colon
      \ho\big(\Fun(\cF,\cM)\big) \longra \ho\big(\Fun(\cG,\cM)\big).
    \end{equation}
    For $X\in \Fun(\cF,\cM)$ we call
    \begin{equation}
      \Llkan_u X
    \end{equation}
    the \define{homotopy left Kan extension} of $X$ along $u$.

    Equivalently, if $\bQ X$ is a cofibrant replacement of $X$ in
    $\Fun(\cF,\cM)$ (projective), then
    \begin{equation}
      \Llkan_u X \;\simeq\; \lkan_u(\bQ X)
    \end{equation}
    in $\ho\big(\Fun(\cG,\cM)\big)$.

    \item \textbf{Homotopy right Kan extension.}
    Endow $\Fun(\cF,\cM)$ and $\Fun(\cG,\cM)$ with the \emph{injective}
    diagram model structures. Then $u^{\ast}$ has a right adjoint
    \begin{equation}
      \rkan_{u} \colon \Fun(\cF,\cM) \longra \Fun(\cG,\cM),
    \end{equation}
    the (ordinary) right Kan extension along $u$.
    In the injective model structure the functor $\rkan_{u}$ is right
    Quillen, so we can form its total right derived functor
    \begin{equation}
      \Rrkan_{u} \colon
      \ho\big(\Fun(\cF,\cM)\big) \longra \ho\big(\Fun(\cG,\cM)\big).
    \end{equation}
    For $X\in \Fun(\cF,\cM)$ we call
    \begin{equation}
      \Rrkan_{u} X
    \end{equation}
    the \emph{homotopy right Kan extension} of $X$ along $u$.

    Equivalently, if $\bR X$ is a fibrant replacement of $X$ in
    $\Fun(\cF,\cM)$ (injective), then
    \begin{equation}
      \Rrkan_{u} X \;\simeq\; \rkan_{u}(\bR X)
    \end{equation}
    in $\ho\big(\Fun(\cG,\cM)\big)$.
  \end{enumerate}
\end{definition}

\begin{remark}\label{ch2:: rem: pointwise-kan}
 For $X\colon \cF\to\cM$ and
  $g\in \cG$ one has the usual formulas
  \begin{equation}
    (\lkan X)(g) \;\cong\;
    \mathsf{colim}\big(\cF/g \xra{\,\, \sf forget \,\,} \cG \xra{u} \cM \big),
    \qquad
    (\rkan_{u}X)(g) \;\cong\;
    \mathsf{lim}\big(g/\cF \xra{\,\, \sf forget \,\,} \cG \xra{u} \cM \big),
  \end{equation}
  where $(u/g)$ and $(g/u)$ are the usual comma categories.
  Passing to derived functors and using
  Definition~\textup{\ref{ch2:: def: homotopy-(co)limits}} we obtain
  natural weak equivalences
  \begin{equation}
    (\Llkan_u X)(g)
    \;\simeq\;
    \mathsf{hocolim}\big(\cF/g \xra{\,\, \sf forget \,\,} \cG \xra{u} \cM \big),
    \qquad
    (\Rrkan_{u} X)(g)
    \;\simeq\;
    \mathsf{holim}\big(g/\cF \xra{\,\, \sf forget \,\,} \cG \xra{u} \cM \big).
  \end{equation}
  The existence of this pointwise forumalae for homotopy Kan extensions is asserted and proved in \cite[Proposition~1.14]{Cis09}
\end{remark}

\begin{remark}\label{ch2:: rem: holim-hocolim-as-kan}
  Let $\cF$ be a small category and let
  \begin{equation}
    p \colon \cF \longra *
  \end{equation}
  be the unique functor to the terminal category.
  Then $\lkan_p$ is the colimit functor and $\rkan_{p}$ is the limit functor.
  In particular, for $X\in\Fun(\cF,\cM)$ there are canonical identifications
  \begin{equation}
    \hoColimit[\cF] X
    \;\simeq\;
    \big(\LL\mathsf{Kan}_p X\big)(*),
    \qquad
    \hoLimit[\cF] X
    \;\simeq\;
    \big(\Rrkan_{p} X\big)(*).
  \end{equation}
  Thus, homotopy (co)limits are precisely the (left) right homotopy Kan extensions along the terminal functor.
\end{remark}

\begin{proposition}\label{ch2:: prop: homotopy-limit-homotopically-constant}
  Let $\cM$ be a simplicial model category, let $\cF$ be a small category
  whose nerve is contractible, and let
  $X \colon \cF \to \cM$ be a homotopically constant diagram (i.e.~every map
  $X(f)$ is a weak equivalence) which is fibrant as an $\cF$-diagram.  Then
  for every $i \in \cF$ the canonical maps
  \[
    \holim_{\cF} X \longrightarrow X(i),
    \qquad
    X(i) \longrightarrow \hocolim_{\cF} X
  \]
  are weak equivalences in $\cM$.
\end{proposition}

\begin{example}
  Let $X \colon \Delta^{\op} \to \sset$ be an injectively fibrant simplicial object in $\sset$, which is levelwise weakly equivalent to a constant diagram on a
  Kan complex $K$.  Then
  $\holim_{\Delta^{\op}} X \simeq K \simeq \hocolim_{\Delta^{\op}} X$.  This
  is a typical situation when simplicial resolutions are used to compute
  derived mapping spaces.
\end{example}

\begin{definition}[Homotopy initial and homotopy final functors]
  Let $u \colon \cF \to \cG$ be a functor between small categories and let
  $\cM$ be a model category.

  \begin{enumerate}
    \item We say that $u$ is \define{homotopy final} if, for every diagram
      $X \colon \cG \to \cM$, the canonical map
      \[
        \hoColimit[\cF] (X \circ u)
        \longrightarrow
        \hoColimit[\cG] X
      \]
      is a weak equivalence in $\cM$, whenever these homotopy colimits are
      defined.

    \item We say that $u$ is \define{homotopy initial} if, for every diagram
      $X \colon \cG \to \cM$, the canonical map
      \[
        \hoLimit[\cG] X
        \longrightarrow
        \hoLimit[\cF] (X \circ u)
      \]
      is a weak equivalence in $\cM$, whenever these homotopy limits are
      defined.
  \end{enumerate}
\end{definition}

The preceding definition is conceptually clean but not very practical.  The
usual criterion is phrased in terms of comma categories and their nerves.

\begin{proposition}\label{ch2::prop:homotopy-final-initial}
  Let $u \colon \cF \to \cG$ be a functor between small categories.
  \begin{enumerate}
    \item The functor $u$ is homotopy final if and only if, for every object
      $g \in \cG$, the nerve of the over-category $(g / \cF)$ is
      weakly contractible:
      \begin{equation}
        \sN\bigl(g / \cF\bigr) \simeq *.
      \end{equation}
    \item The functor $u$ is homotopy initial if and only if, for every object
      $g \in \cG$, the nerve of the under-category $(\cF / g)$ is
      weakly contractible:
      \begin{equation}
        \sN\bigl(\cF / g\bigr) \simeq *.
      \end{equation}
  \end{enumerate}
  In particular, these conditions do not depend on the choice of model
  category $\cM$.
\end{proposition}

\begin{example}
  If $\cF$ has a terminal object $t$ and $u \colon \cF \to \cG$ sends $t$ to a
  terminal object of $\cG$, then $u$ is homotopy final.  Indeed, for any
  $g \in \cG$ the over-category $(g / \cF)$ has a terminal object (if
  it is nonempty), so its nerve is contractible.  Dually, if $\cF$ has an
  initial object that $u$ preserves, then $u$ is homotopy initial.
\end{example}

Finally, we recall left Bousfield localization, which will be used to
construct the model structures for (complete) Segal spaces and (complete) interval Segal spaces.

\begin{definition}
  Let $\cM$ be a model category and let $K$ be a set of maps in $\cM$.
  \begin{itemize}
    \item A fibrant object $Z \in \cM$ is called \define{$K$-local} if for every
      $f \colon A \to B$ in $K$, the induced map on derived mapping spaces
      \begin{equation}
        \Map^{h}_{\cM}(B,Z) \longrightarrow \Map^{h}_{\cM}(A,Z)
      \end{equation}
      is a weak equivalence of simplicial sets.
    \item A map $g \colon X \to Y$ in $\cM$ is a \define{$K$-local equivalence}
      if for every $K$-local $Z$ the induced map
      \begin{equation}
        \Map^{h}_{\cM}(Y,Z) \longrightarrow \Map^{h}_{\cM}(X,Z)
      \end{equation}
      is a weak equivalence.
  \end{itemize}
  A \define{left Bousfield localization} of $\cM$ with respect to $K$ is a model
  structure $L_K\cM$ on the same underlying category, with the same
  cofibrations, whose weak equivalences are the $K$-local equivalences and
  whose fibrant objects are the $K$-local fibrant objects.
\end{definition}

\begin{theorem}[Left Bousfield localization, {\cite[Thm.~3.3.1, Thm.~4.1.1]{Hir14}}]
  \label{ch2:: thm: left Bousfield localization}
  Let $\cM$ be a left proper, cellular or combinatorial model category, and
  let $K$ be a set of maps in $\cM$.  Then the left Bousfield localization
  $L_K\cM$ exists.  Moreover:
  \begin{itemize}
    \item $L_K\cM$ is left proper and cellular or combinatorial (accordingly);
    \item if $\cM$ is a simplicial model category, then so is $L_K\cM$.
  \end{itemize}
\end{theorem}

\begin{proposition}
  Let $\cM$ be a model category and let $K$ be a set of maps in $\cM$.  Then
  a map between $K$-local fibrant objects is a $K$-local weak equivalence if
  and only if it is a weak equivalence in $\cM$.
\end{proposition}

\begin{theorem}\label{ap:: thm: transfer of localizations}
    Suppose that $\cM$ and $\cN$ are model categories and $F\maps \cM \leftrightarrows \cN\maps G$ is a Quillen pair. Suppose that $T$ is a class of maps in $\cM$ such that the localized model structure $L_{T}\cM$ exists, and furthermore that the localized model structure $L_{\mathbb{L}F(T)}\cN$ also exists. Then:
    \begin{enumerate}
        \item The adjoint pair $F\dashv G$ still defines a Quillen pair
        \begin{equation}
            F\maps L_{T}\cM \leftrightarrows L_{\mathbb{L}F(T)}\cN \maps G.
        \end{equation}
        \item If the pair $F\dashv G$ defines a Quillen equivalence between $\cM$ and $\cN$, then it also defines a Quillen equivalence on the localized model structures.
    \end{enumerate}
\end{theorem}


\section{A primer on complete Segal spaces}
\label{ch2:: sec: complete_segal_spaces}

The main goal of the present chapter is to propose a variant model for $(\oo,1)$-categories, yet, so far, we have not made any progress into introducing what these are. This section not only addresses that missing exposition, but it introduces Segal spaces after which we will model our own variant of $(\oo,1)$-categories. The goal is not to give a precise account of the subject, which contains many intricate bells and whistles — details about them can be found in Rezk's original paper \cite{Rezk01} or Julie Bergner's comprehensive book on the homotopy theory of $(\oo,1)$-categories \cite{Ber18} — these are our main references for this section.

An $(\oo,1)$-category is a generalization of the usual notion of category, where the set of morphisms is endowed with more structure, namely that of a homotopy type modeled by a topological space or a simplicial set. Of course, this conceptual change comes with modifications to the way one compares and classifies $(\oo,1)$-categories, especially since there is an obvious expectation that the collection of $(\oo,1)$-categories forms itself an $(\oo,1)$-category. With these objects in mind, we finally turn to the work of Rezk.

\subsection{Simplicial spaces}

Segal spaces, despite their suggestive name, are not actually spaces, they are simplicial spaces. Before we delve into why these are the right objects to use, let us quickly go over what is a category again. category $\sC$ is a collection of objects $\sC_{0}$, together with a collection of morphisms $\sC_{1}$, an additional pair of maps 
\begin{equation}
    \sC_0 \xleftarrow{\,\,\,\mathsf{source}\,\,\,} \sC_{1} \xra{\,\,\,\mathsf{target}\,\,\,} \sC_{0},
\end{equation} 
together with a composition rule 
\begin{equation}
    \sC_{1} \times_{\sC_{0}} \sC_{1} \xra{\,\,\ \circ \,\,} \sC_{1} 
\end{equation}
satisfying an additional associativity condition, encoded by the commutativity of the square (\ref{ch2:: eq: assoc})
\begin{equation}\label{ch2:: eq: assoc}
    \begin{tikzcd}
        {\sC_1\times_{\sC_{0}} \sC_{1}\times_{\sC_{0}} \sC_{1}} && {\sC_1\times_{\sC_{0}} \sC_{1}} \\
        \\
        {\sC_1\times_{\sC_{0}} \sC_{1}} && {\sC_1}\, .
        \arrow["{ \id \times \circ}", from=1-1, to=1-3]
        \arrow["\circ", from=1-3, to=3-3]
        \arrow["{\circ \times \id}", from=3-1, to=1-1]
        \arrow["{ \circ}"', from=3-1, to=3-3]
    \end{tikzcd}
\end{equation}
Additionally, there is a map $\sC_{0} \xra{e} \sC_{1}$ sending $x$ to its corresponding identity morphism. 

Importantly, all of this data and structure could just as well be packaged as a simplicial set: 

\begin{equation}
    \sN(\sC) \maps \bDelta^{\op} \longra \mathsf{Set}.
\end{equation} 
Let $\sN(\sC)_0 = \sC_0$, $\sN(\sC)_{1} = \sC_1$, and 
\begin{equation}
\sN(\sC)_{n} = \underbrace{\sC_{1} \times_{\sC_{0}} \sC_{1} \times_{\sC_{0}} \cdots \times_{\sC_{0}} \sC_{1}}_{n \text{ times}}.
\end{equation} 

Moreover, let $\langle 0\rangle^{*},\langle 1\rangle^{*}$ get sent to the source and target maps, respectively. Let the unique map $\langle 00 \rangle^{*}$ get sent to $\sC_{0} \xra{e} \sC_{1}$. Immediately we note that the simplicial relations enforce a condition — the fact that the source and target of the identity morphism are the same. Composition is defined as the map $\langle 02 \rangle^{*}\maps \sN(\sC)_{2} \longra \sN(\sC)_{1}$, which means that associativity, i.e., a condition on the data, also gets enforced or encoded by the simplicial relations. 

This construction highlights the fact that categories naturally give rise to simplicial sets. More is true, since they do so in a functorial way.

\begin{proposition} \label{ch2:: prop: nerve}
    There exists a fully faithful functor, which we will call the \textbf{nerve}
    \begin{equation}
        \sN \maps \Cat \longra \sset.
    \end{equation}
\end{proposition}

An interesting question is whether we can characterize the essential image of the nerve construction. In other words, can we find a condition on a simplicial set $X_\bullet$ which ensures it is isomorphic to the nerve $\sN(\sC)$ of some category $\sC$. The answer is in the positive, and the condition is recorded in the theorem below.
\begin{theorem}\label{ch2:: thm: Kan--quillen}
    Let $\cC\maps \bDelta^{\op} \longra \set$ be a simplicial set satisfying the following condition: for every $n,m\geq 1$, the diagram
    \begin{equation}
        \begin{tikzcd}
            {\cC_{n+m}} && {\cC_{n}} \\
            \\
            {\cC_{m}} && {\cC_{0}}
            \arrow["{\langle 01\cdots n \rangle^{*} }", from=1-1, to=1-3]
            \arrow["{\langle n\cdots n+m\rangle^{*}}"', from=1-1, to=3-1]
            \arrow["\lrcorner"{anchor=center, pos=0.125}, draw=none, from=1-1, to=3-3]
            \arrow["{\langle n\rangle^{*}}", from=1-3, to=3-3]
            \arrow["{\langle m \rangle^{*}}"', from=3-1, to=3-3]
        \end{tikzcd}
    \end{equation}
    is a pullback square in $\mathsf{Set}$. We call this condition the \textbf{Segal condition}. A simplicial set $\cC_\bullet$ is isomorphic to the nerve of some category $\sC$ if and only if it satisfies the Segal condition.
\end{theorem}

\begin{remark}
 Note that we can consider an equivalent variant of the Segal condition where instead of decomposing $[n+m]$ into linear subposets $[n]$ and $[m]$, we decompose them into more linear subposets of length 1. We can consider the inert maps $[1] \xra{\langle i, i+1\rangle} [n]$ and re-define the Segal condition by demanding that for every $n\geq 2$ the map
 \begin{equation}
    \cC_{n} \longra \cC_{1} \times_{\cC_{0}} \cdots \times_{\cC_{0}} \cC_{1}
 \end{equation}
 is an isomorphism of sets.
\end{remark}

\begin{observation}
The Segal condition says that the $n$-simplices of $\cC$ are, up to isomorphism, strings of $n$-composable 1-simplices, which in turn are manifestly interpreted as the 1-morphisms of some category whose nerve is isomorphic to $\cC$.
\end{observation}

This theorem informs us that simplicial sets house the two ingredients which are important to build $(\oo,1)$-categories: 1) categories and 2) homotopy types. Even more explicitly, we want our classes of morphisms to become homotopy types and we want that associativity and unitality be modified accordingly — associativity need only hold up to a weak homotopy equivalence. Simplicial spaces, or bisimplicial sets, offer an appropriate receptacle for both structures at once as we are interested in going from sets to spaces, it makes sense that we study functors
\begin{equation}
    \bDelta^{\op} \longra \sset,
\end{equation}
which can equivalently be defined and studied as bisimplicial sets. The former perspective is the one we will take when we refer to simplicial spaces — with exceptions, all of which are duly flagged.

\begin{remark}
    The category $\fun(\bDelta^{\op},\sset)$ of simplicial spaces contains two fully faithful copies of $\sset$, and their respective images have very different interpretations. Given a simplicial set $K$, we define the \textbf{constant} simplicial space on $K$ by $([n] \mapsto K\in \sset)$ and every $[n] \ra [m]\in \bDelta^{\op}$ gets sent to the identity. Alternatively, we have the case that $K$ factors through the inclusion $\mathsf{Set} \hookra \sset$ so that $K_n$ is a discrete simplicial set i.e., a set for every $n\geq 0$. What makes this distinction so relevant is the fact that it is only one of the versions of $\Delta^{n}$ in $\fun(\bDelta^{\op},\sset)$ which is the stand-in for the standard $n$-simplex, namely the discrete one. We denote this simplicial space by $\Delta[n]$ in order to contrast with the simplicial set $\Delta^n$. In particular we have that for every simplicial space $X\maps \bDelta^{\op} \ra \sset$, there is a weak equivalence of simplicial sets
    \begin{equation}
        \Map(\Delta[n], X) \simeq X_n 
    \end{equation}
    for every $n\geq 0$.
\end{remark}

We will shift our focus to simplicial spaces now. 

\begin{observation}\label{ch2:: obs: reedy is injective}
In particular, the category $\fun(\bDelta^{\op},\sset)$ can be endowed with
either the injective model structure or the projective model structure. In fact
it can even be endowed with a third one, called the \emph{Reedy model
structure}. All of them are Quillen equivalent, and in this case the injective
and Reedy model structures agree. For simplicity, we will choose the injective
model structure.
\end{observation}

\begin{theorem}
    There exists a model structure on the category of simplicial spaces $\fun(\bDelta^{\op},\sset)$ which we will denote $\fun(\bDelta^{\op},\sset)_{\sf inj}$, with the following properties:
    \begin{enumerate}
        \item The weak equivalences are levelwise weak equivalences in $\sset_{\sf QK}$.
        \item The cofibrations are levelwise cofibrations in $\sset_{\sf QK}$.
        \item The model structure is left proper combinatorial and simplicial.
    \end{enumerate}
\end{theorem}

\begin{observation}
    \label{ch2:: obs: discreteness}
A pair $X$ and $Y$ of discrete simplicial spaces are weakly equivalent as injective simplicial spaces if and only if $X_{k} \simeq Y_{k}$ for every $k$. In other words, if the sets $X_k$ and $Y_k$ are bijective for all $k\geq 0$. 
\end{observation}

\subsection{Imposing the Segal condition}

We want a definition of $(\oo,1)$-category which resembles the nerve of an
ordinary category but relaxes the requirements of strict composition,
associativity, and unitality.

\begin{construction}
Consider the class of inert morphisms
\[
  \big\{[1] \xra{\langle i, i+1 \rangle} [n]\,:\, 0\leq i < n\big\}_{n \geq 1}
\]
in $\bDelta$. Let $\bDelta^{\sf int}\subset\bDelta$ denote the subcategory of
inert morphisms and let $\bDelta^{\sf int}_{\leq 1}\subset\bDelta^{\sf int}$ be
the full subcategory on $[0]$ and $[1]$. For each $[n]\in\bDelta$ we can form
the slice category $\bDelta^{\sf int}_{\leq 1}/[n]$. We define the
\emph{$n$-spine} to be the homotopy colimit
\begin{equation}
  G(n) \,\overset{\sf def}{=}\,
  \hocolim\Big(
    \bDelta^{\sf int}_{\leq 1}/[n]
      \longra \bDelta_{\leq 1}^{\sf int}
      \xra{\,\bDelta[\bullet]\,}
      \Fun(\bDelta^{\op},\sset)
  \Big),
\end{equation}
where $\bDelta[\bullet]$ is the Yoneda embedding.
\end{construction}

\begin{remark}
    The objects of $\bDelta^{\sf int}_{\leq 1}/[n]$ are precisely the inert
    maps $[1] \xra{\langle i,i+1\rangle} [n]$ for $0\leq i<n$.
\end{remark}

\begin{observation}
    The $2$-spine $G(2)$ happens to coincide with the inner horn
    $\Lambda_{1}[2]$, but this is no longer true for $n\geq 3$. For higher $n$,
    the $n$-spine is a string of $1$-simplices and does not contain higher
    non-degenerate simplices, unlike the higher inner horns.
\end{observation}

Let $\cF\colon \bDelta^{\op} \to \sset$ be a simplicial space, which we will
assume to be injectively fibrant. For $n\geq 2$ there is an inclusion of
simplicial spaces $G(n) \hookra \bDelta[n]$, which induces a map
\begin{equation}
    \Map(\bDelta[n], \cF) \longra \Map(G(n), \cF)
      \,\simeq\,
      \hoLimit{(\bDelta^{\sf int}_{\leq 1}/[n])}
      \Map(\bDelta[1], \cF).
\end{equation}
Using the standard identifications $\Map(\bDelta[n],\cF)\cong\cF_n$ and
$\Map(\bDelta[1],\cF)\cong\cF_1$, this can equivalently be recast as a map
\begin{equation}
    \varphi_{n} \colon
    \cF_{n} \longra
    \underbrace{\cF_1 \times_{\cF_{0}}^{h} \cdots
      \times_{\cF_{0}}^{h} \cF_{1}}_{n \text{ times}}.
\end{equation}
These maps $\varphi_n$ for all $n\geq 2$ are called the \textbf{Segal maps}.\footnote{Segal first defined
these maps in the reduced case $\cF_{0} = *$, as we will re-encounter when we
discuss symmetric monoidal $(\oo,1)$-categories.}

Recall that a simplicial set is (isomorphic to) the nerve of an ordinary
category if and only if its Segal maps are isomorphisms of (discrete) simplicial
sets. For a simplicial space we want $\cF_0$ to play the role of the space of
objects, $\cF_1$ the space of morphisms, and $\cF_n$ the space of $n$ composable
morphisms; asking the Segal maps to be weak equivalences is exactly the
homotopy-theoretic analogue of strict composition.

\begin{definition}\label{ch2:: def: segal space}
    Let $\cF \colon \bDelta^{\op} \to \sset$ be an injectively fibrant simplicial
    space. We say that $\cF$ is a \emph{Segal space} if for all $n\geq 2$ the
    Segal maps
    \begin{equation}
        \varphi_{n} \colon
        \cF_{n} \longra
        \underbrace{\cF_1 \times_{\cF_{0}}^{h} \cdots
          \times_{\cF_{0}}^{h} \cF_{1}}_{n \text{ times}}
    \end{equation}
    are weak equivalences of simplicial sets in the Quillen--Kan model
    structure, i.e.\ weak homotopy equivalences. Since $\cF$ is fibrant, the
    ordinary fibre products compute these homotopy fibre products, so we may
    omit the superscript $h$ when convenient.
\end{definition}

\begin{observation}
    Observe that we can also define Segal spaces to be the injectively fibrant simplicial spaces $\cF\maps \bDelta^{\op} \longra \sset$ such that their restriction $\cF_{|\bDelta^{\mathsf{int}, \op}}$ corresponds to the right Kan extension 
    \begin{equation}
    \begin{tikzcd}
        {\bDelta^{\mathsf{int},\op}_{\leq 1}} &&& \sset \\
        \\
        {\bDelta^{\mathsf{int},\op}}
        \arrow["{\cF_{|\bDelta^{\mathsf{int}, \op}_{\leq 1}}}", from=1-1, to=1-4]
        \arrow[hook, from=1-1, to=3-1]
        \arrow["{\cF_{|\bDelta^{\mathsf{int}, \op}}}"', curve={height=12pt}, dashed, from=3-1, to=1-4]
    \end{tikzcd}
    \end{equation}
    We find this a useful reframing of the Segal conditions, since it admits a more immediate generalization when passing from $(\oo,1)$- to $(\oo,n)$-categories via presheaves of simplicial sets.
\end{observation}

\begin{proposition}\label{ch2:: prop: Segal-model-structure}
    There exists a model structure on the category of simplicial spaces
    $\fun(\bDelta^{\op},\sset)$, which we denote by $\mathsf{SSp}$, obtained by
    left Bousfield localizing the injective model structure
    $\fun(\bDelta^{\op},\sset)_{\mathsf{inj}}$ at the set of spine inclusions
    \begin{equation}
        \{G(n)\hookrightarrow \Delta[n]\}_{n\geq 2}.
    \end{equation}
    It has the following properties:
    \begin{enumerate}
        \item The cofibrations are the injective cofibrations, i.e., the monomorphisms.
        \item The fibrant objects are the Segal spaces.
        \item The weak equivalences between Segal spaces are the levelwise weak equivalences.
        \item The model structure is left proper, combinatorial, and simplicial.
    \end{enumerate}
\end{proposition}

\begin{proof}
    Existence follows from \Cref{ch2:: thm: left Bousfield localization},
    applied to the injective model structure on simplicial spaces from the
    theorem above and the set of spine inclusions. Since cofibrations are
    unchanged under left Bousfield localization, (1) holds. By construction,
    the fibrant objects are precisely the injectively fibrant simplicial spaces
    that are local with respect to the spine inclusions, which is exactly the
    Segal condition from \Cref{ch2:: def: segal space}; this gives (2). Finally,
    (3) is the usual description of weak equivalences between fibrant objects in
    a left Bousfield localization, and (4) is inherited from the injective
    model structure.
\end{proof}

\subsection{Homotopy theory of Segal spaces}

Segal spaces, just as any good homotopical gadget, admit a few important
invariants which facilitate their study using algebraic tools. Namely:

\begin{itemize}
    \item Its set of objects.
    \item Its mapping spaces between pairs of objects.
    \item Homotopy equivalence of objects (and morphisms).
    \item Homotopy category.
\end{itemize}

\subsubsection{Set of objects and mapping space}

The best encapsulating property of what it means to be an $(\oo,1)$-category is that the collection of morphisms between any fixed pair of objects organizes as a \emph{space}, with composition being a map of spaces — thereby defining all higher invertible $k$-morphisms as homotopies, in some suitable sense.

\begin{construction}
    \label{ch2:: cons: objects mapping space}
    Let $\cF$ be a Segal space. Let us define its set of \define{objects} to be the set of $0$-simplices of $\cF_0$; and given a pair of objects $x,y\in \cF_{0,0}$, let us define the mapping space $\map_{\cF}(x,y)$ to be the pullback
    \begin{equation}
        \begin{tikzcd}
            {\Map_{\cF}(x,y)} && {\cF_{1}} \\
            \\
            {*} && {\cF_{0}\times \cF_{0}}
            \arrow[from=1-1, to=1-3]
            \arrow[from=1-1, to=3-1]
            \arrow["\lrcorner"{anchor=center, pos=0.125}, draw=none, from=1-1, to=3-3]
            \arrow["\lag 0\rag^* \times \lag 1 \rag^*", from=1-3, to=3-3]
            \arrow["{(x,y)}"', from=3-1, to=3-3]
        \end{tikzcd}
    \end{equation}
    in $\sset$. Recall our assumption that $\cF$ is injectively fibrant — it implies that the right vertical map of simplicial sets is a fibration; as observed in \Cref{ch2:: obs: reedy is injective}, the injective and Reedy model structure agree and since $\partial \Delta[1] \longhookra \Delta[1]$ is cofibration therein, we have that the following induced map is a Kan fibration
    \[\Hom_{\sset}(\Delta[1], \cF) \longra \Hom_{\sset}(\partial \Delta[1], \cF).\]
    In particular, this means that $\Map_{\cF}(x,y)$ is homotopy equivalent to the homotopy fiber of \[\cF_1 \longra \cF_0 \times \cF_0,\] which is desirable, since we want the mapping space $\Map_{\cF}(x,y)$ to be well-defined only up to homotopy, i.e., only depend on the path-component of $x,y\in \cF_0$. We will use the notation
    \begin{equation}
        f \maps x \longra y
    \end{equation}
    for a 0-simplex of $\Map_{\cF}(x,y)$.
\end{construction}

\begin{definition}
    Let $\cF$ be a Segal space. Let $f\maps x \ra y$ and $g\maps x \ra y$ be 0-simplices of the mapping space $\Map_{\cF}(x,y)$. We declare $f \sim g$ as \define{homotopic} if they are in the same component, i.e., if $([f] = [g]) \in \pi_0 \Map_{\cF}(x,y)$.
\end{definition}

\subsubsection{Homotopy equivalences and homotopy category}

Given a Segal space $\cF\maps \bDelta^{\op} \longra \sset$, in \Cref{ch2:: cons: objects mapping space} we associated to it: 1) a set of objects; 2) a mapping space for each pair of objects. Thus, a sensible question is the following: how do we compose morphisms in a Segal space? Let $(x,y,z) \in \cF_{0,0}^{\times 3}$ be a triple of objects of $\cF$. Let $f \maps x \longra y$ and $g \maps y \longra z$ be morphisms in $\Map_{\cF}(x,y)$ and $\Map_{\cF}(y,z)$ respectively. We define the pullback
\begin{equation}
    \begin{tikzcd}[ampersand replacement=\&]
            {\Map_{\cF}(x,y,z)} \&\& {\cF_{2}} \\
            \\
            \ast \&\& {\cF_0\times \cF_0\times \cF_0}
            \arrow[from=1-1, to=1-3]
            \arrow[from=1-1, to=3-1]
            \arrow["\lrcorner"{anchor=center, pos=0.125}, draw=none, from=1-1, to=3-3]
            \arrow["{\lag 0 \rag^* \times \lag 1 \rag^* \times \lag 2 \rag^*}", from=1-3, to=3-3]
            \arrow["{(x,y,z)}"', from=3-1, to=3-3]
        \end{tikzcd}
\end{equation}
to be the space of composable pairs of morphisms $(x \xra{f} y \xra{g} z)$ with fixed $x,y,z$. Using the fact that the Segal map $\cF_2 \longra \cF_1 \times_{\cF_0} \cF_1$ is a weak equivalence, we assert that the dashed map in the following diagram is a weak equivalence:
\begin{equation}
\begin{tikzcd}[ampersand replacement=\&]
	\&\&\& {\Map_{\cF}(x,y,z)} \& {\Map_{\cF}(x,y) \times \Map_{\cF}(y,z)} \\
	\\
	\&\& {\cF_{2}} \& {\cF_1 \times_{\cF_0} \cF_1} \& {*} \\
	\\
	{} \&\&\& {\cF_0\times \cF_0\times \cF_0} \& {.}
	\arrow["\simeq"', dashed, from=1-4, to=1-5]
	\arrow[from=1-4, to=3-3]
	\arrow[""{name=0, anchor=center, inner sep=0}, between={0.5}{1}, from=1-4, to=3-5]
	\arrow["\lrcorner"{anchor=center, pos=0.125, rotate=-45}, draw=none, from=1-4, to=5-4]
	\arrow[from=1-5, to=3-4]
	\arrow[from=1-5, to=3-5]
	\arrow["\lrcorner"{anchor=center, pos=0.125, rotate=-90}, draw=none, from=1-5, to=5-4]
	\arrow["\simeq"', from=3-3, to=3-4]
	\arrow[from=3-3, to=5-4]
	\arrow[from=3-4, to=5-4]
	\arrow[from=3-5, to=5-4]
	\arrow[between={0}{0.75}, no head, from=1-4, to=0]
\end{tikzcd}
\end{equation}
\smallskip

\nid Given a pair $\big(x \xra{f} y,y \xra{g} z\big)\in \Map_{\cF}(x,y)\times \Map_{\cF}(y,z)$, we can form the following diagram:
\smallskip

\begin{equation}\label{comp}
    \begin{tikzcd}[ampersand replacement=\&]
        \&\&\& {\cF_2} \&\& {\cF_1} \\
        \&\& {\Map_{\cF}(x,y,z)} \&\& {\Map_{\cF}(x,z)} \\
        \&\&\& {\cF_1\times_{\cF_0}\cF_1} \&\& {\cF_1} \\
        {*} \&\& {\Map_{\cF}(x,y)\times \Map_{\cF}(y,z)}
        \arrow["{\lag 0,2 \rag^*}", from=1-4, to=1-6]
        \arrow[between={0}{0.4}, no head, from=1-4, to=3-4]
        \arrow["{\rotatebox[origin=c]{90}{$\sim$}}"'{pos=0.75}, between={0.6}{1}, from=1-4, to=3-4]
        \arrow[equals, from=1-6, to=3-6]
        \arrow[from=2-3, to=1-4]
        \arrow[from=2-3, to=2-5]
        \arrow["{\rotatebox[origin=c]{90}{$\sim$}}",from=2-3, to=4-3]
        \arrow[from=2-5, to=1-6]
        \arrow[from=3-4, to=3-6]
        \arrow["h", dashed, from=4-1, to=2-3]
        \arrow["{(f,g)}"', from=4-1, to=4-3]
        \arrow[curve={height=24pt}, dotted, from=4-3, to=2-5]
        \arrow[""{name=0, anchor=center, inner sep=0}, from=4-3, to=3-4]
        \arrow["\lrcorner"{anchor=center, pos=0.125}, draw=none, from=2-3, to=0]
    \end{tikzcd}
\end{equation}
\smallskip

\nid where $h$ is lift of $(f,g)$. We define a composite for $(f,g)$ to be the image of $h$ under $\Map_{\cF}(x,y,z) \longra \Map_{\cF}(x,y)\times \Map_{\cF}(y,z)$. We indicate this ``composition map'' with the curved dotted arrow. We note that this arrow is only well-defined up to homotopy since it depends on the choice $h$ of a lift of $(f,g)$.

\begin{definition}\label{ch2:: def: homotopy category segal space}
    Let $\cF$ be a Segal space. We define its \textbf{homotopy category} $\ho(\cF)$ to be the category with the following structure:
    \begin{itemize}
        \item \textbf{Objects:} The set of 0-simplices of the simplicial set $\cF_0$.
        \item \textbf{Morphisms:} Given a pair of objects $x,y \in \cF_{0,0}$ we define the set of morphisms of $x$ to $y$ to be 
        \begin{equation}
            \Hom_{\ho(\cF)}(x,y) \, \overset{\sf def}{=} \, \pi_0 \Map_{\cF}(x,y).
        \end{equation}
        \item \textbf{Composition:} Let \[([f],[g])\in \Hom_{\ho(\cF)}(x,y) \times \Hom_{\ho(\cF)}(y,z)\]
        be a pair of composable morphisms. Passing diagram (\ref{comp}) to $\pi_0$, we get a well-defined composition map given by the dashed arrow in the following commutative diagram: 
        \begin{equation}
            \begin{tikzcd}[ampersand replacement=\&]
            {*} \& {\pi_{0}\Map_{\cF}(x,y,z)} \&\& {\pi_0\Map_{\cF}(x,z)} \\
            \\
            {\pi_0\Map_{\cF}(x,y)\times \pi_0\Map_{\cF}(y,z)} \& {\pi_0\Map_{\cF}(x,y)\times \pi_0\Map_{\cF}(y,z)}
            \arrow["{([f],[g])}"{description}, from=1-1, to=3-1]
            \arrow[from=1-2, to=1-4]
            \arrow["\cong", from=1-2, to=3-2]
            \arrow["{\cong} "', from=3-1, to=3-2]
            \arrow["\circ"'{pos=0.6}, curve={height=24pt}, dashed, from=3-2, to=1-4]
        \end{tikzcd}
        \end{equation} 
        The key here is the fact that the weak equivalence became an isomorphism of sets.\footnote{The well-definedness of the composition map comes from the fact that any two choices $\lag 0,2\rag^* h$ and $\lag 0,2\rag^* h'$ of composition for a fixed pair $(f,g)$, are homotopic.} 
    \end{itemize}
\end{definition}

\begin{observation}
    For every object $x\in \ho(\cF)$, there is a 1-morphism $\lag 00\rag ^* x \in \cF_1$ which is the identity. We note that for $f\maps x \ra y$, we can contemplate the following 2-simplex $\lag 0,0,1 \rag^* f \in \cF_2$, whose image under $\lag 0,2\rag^*$ is
    \begin{equation}
        \lag 0,2 \rag ^* \circ \lag 0,0,1\rag ^* f = \left(\lag 0,0,1 \rag \circ \lag 0,2 \rag \right)^*f = \lag 0,1\rag ^{*}f = f,
    \end{equation}
    which confirms that a composite of $f$ with $\id_x$ is again $f$.
\end{observation}

\medskip

Given the fact that composition is not defined on the nose, we cannot (or should not) ask that a given morphism $f\maps x \ra y$ be an isomorphism. For that would force us to make a choice of composite and then ask that choice be the identity —  all the while ignoring all other possible composites. We avoid this by using the homotopy category.

\begin{definition}\label{ch2:: def: hom eq}
Let $\cF$ be a Segal space. We define a morphism $f\maps x \ra y$ to be a \define{homotopy equivalence}\footnote{It is also common to call homotopy equivalences in a Segal space just equivalences.} if $[f]\in \ho(\cF)$ is an isomorphism.
\end{definition}

The next step is to understand to what extent an equivalence of Segal spaces, encoded as a levelwise equivalence of simplicial spaces, captures the correct notion of equivalence between $(\oo,1)$-categories. The short answer is that it does not, but we must first articulate what we expect an equivalence of $(\oo,1)$-categories to be. To do that, we suitably emulate the usual definitions from ordinary category theory. Namely those of fully-faithful and essentially surjective functor.

\begin{definition}
    A morphism $F\maps \cF \longra \cG$ of Segal spaces (i.e., a natural transformation $\bDelta^{\op}\times [1] \longra \sset$) is \define{fully faithful} if the commutative square 
    \[\begin{tikzcd}[ampersand replacement=\&]
        {\cF_1} \&\& {\cG_1} \\
        \\
        {\cF_0\times \cF_0} \&\& {\cG_0 \times \cG_0}
        \arrow[from=1-1, to=1-3]
        \arrow[from=1-1, to=3-1]
        \arrow[from=1-3, to=3-3]
        \arrow[from=3-1, to=3-3]
    \end{tikzcd}\]
    is a homotopy pullback square.
\end{definition}

\begin{observation}
Unpacking the definition above, we see that for every $(x,y) \in \cF_0^{\times 2}$ by the universal property of the pullback, there exists the following dashed map  
\begin{equation}
    \begin{tikzcd}[ampersand replacement=\&]
        \& {\Map_{\cG}(Fx,Fy)} \&\& {\cG_1} \\
        {\Map_{\cF}(x,y)} \&\& {\cF_1} \\
        \& {*} \&\& {\cG_0 \times \cG_0} \\
        {*} \&\& {\cF_0\times \cF_0}
        \arrow[from=1-2, to=1-4]
        \arrow[""{name=0, anchor=center, inner sep=0}, between={0.6}{1}, from=1-2, to=3-2]
        \arrow["\lrcorner"{anchor=center, pos=0.125}, draw=none, from=1-2, to=3-4]
        \arrow[from=1-4, to=3-4]
        \arrow[dashed, from=2-1, to=1-2]
        \arrow[from=2-1, to=2-3]
        \arrow[from=2-1, to=4-1]
        \arrow[""{name=1, anchor=center, inner sep=0}, "\lrcorner"{anchor=center, pos=0.125}, draw=none, from=2-1, to=4-3]
        \arrow[from=2-3, to=1-4]
        \arrow["\lrcorner"{anchor=center, pos=0.125}, draw=none, from=2-3, to=3-4]
        \arrow[from=2-3, to=4-3]
        \arrow[between={0.7}{1}, from=3-2, to=3-4]
        \arrow[from=4-1, to=3-2]
        \arrow["\lrcorner"{anchor=center, pos=0.125, rotate=45}, draw=none, from=4-1, to=3-4]
        \arrow["{(x,y)}"{description}, from=4-1, to=4-3]
        \arrow[from=4-3, to=3-4]
        \arrow[between={0}{0.9}, no head, from=1-2, to=0]
        \arrow["{(Fx,Fy)}"{description, pos=0.3}, between={0.1}{0.6}, no head, from=1, to=3-4]
    \end{tikzcd}
\end{equation}
\medskip

\nid and this map is an equivalence — matching the usual notion of fully-faithfulness in ordinary category theory. In particular, the induced map between homotopy categories
\[\ho(F)\maps \ho(\cF) \longra \ho(\cG) \]
is fully-faithful. 
\end{observation}

\begin{observation}
    Let $\cF$ be a Segal space. We define an equivalence relation on the set of path components of the space of objects $\pi_0\cF_0$, by determining 
    \begin{equation}
        x \sim y, \text{ if and only if there is } f\maps x \ra y  
    \end{equation} 
    with $f$ a homotopy equivalence (\Cref{ch2:: def: hom eq}). Note the subtle distinction between being path connected in $\cF_{0}$ and being connected by an equivalence. They do not always agree, and that is at the heart of our next topic of discussion: \emph{completeness}.
\end{observation}

Before we can discuss completeness, what it means, and why we need it, let us first define what it means for a map of Segal spaces to be essentially surjective.

\begin{definition}
    Let $F \maps \cF \longra \cG$ be a map of Segal spaces. We declare $F$ to be \define{essentially surjective} every object in $\cG$ is equivalent to one in the image of $F$ or if the induced map 
    \begin{equation}
    \begin{tikzcd}[ampersand replacement=\&]
        {\cF_0} \&\& {\cG_0} \\
        {\pi_0\cF_0} \&\& {\pi_{0}\cG_{0}} \\
        {\pi_{0}\cF_{0}{/\sim}} \&\& {\pi_{0}\cG_{0}/\sim}
        \arrow[from=1-1, to=1-3]
        \arrow[from=1-1, to=2-1]
        \arrow[from=1-3, to=2-3]
        \arrow[from=2-1, to=2-3]
        \arrow[from=2-1, to=3-1]
        \arrow[from=2-3, to=3-3]
        \arrow[two heads, from=3-1, to=3-3]
    \end{tikzcd}
    \end{equation}
    in the bottom of this diagram is surjective.
\end{definition}

\begin{definition}
    We define a map $F\maps \cF \longra \cG$ to be a \define{Dwyer--Kan equivalence} if it is fully-faithful and essentially surjective.
\end{definition}

Since we want $(\oo,1)$-categories to be equivalent if and only if they are Dwyer--Kan equivalent, an immediate question is: 

\begin{enumerate}
    \item[\textbf{Q:}] To which extent do levelwise equivalences between Segal spaces agree with Dwyer--Kan equivalences? 
\end{enumerate}

If the answer to this question were affirmative we could confidently assert that
$\mathsf{SSp}$ was a correct model for the homotopy theory of
$(\oo,1)$-categories. In other words, it would serve as a model for the
$(\oo,1)$-category of $(\oo,1)$-categories. Unfortunately, the answer is no.
We construct in \Cref{ch2:: ex: need complete} a pair of Segal spaces which are
Dwyer--Kan equivalent but fail to be so levelwise. To add some context for how
these examples are constructed, let us introduce a generalization of the
simplicial nerve $\sN$, called the classifying diagram, due to Rezk
(\cite{Rezk01}).

\begin{definition}
    \label{ch2:: def: rezk nerve}
    Let $\sC$ be an ordinary category, and let $\sW \subseteq \sC$ be a wide subcategory of $\sC$. We define the following simplicial space which we call \define{Rezk's classifying diagram} of $\sC$ relative to $\sW$: 
    \begin{equation}
        \sN^{\cR}(\sC,\sW) \maps \bDelta^{\op} \longra \sset,\quad [m] \longmapsto \sN\big(\fun^{\sW}([m],\sC) \big).
    \end{equation} 
    Here $\fun^{\sW}([m],\sC)$ denotes the category whose objects are all
    functors $[m]\to \sC$, and whose morphisms are the natural transformations
    whose components lie in $\sW$.
\end{definition}

\begin{remark}
    We can describe the $n$-simplices of the $m$th space of $\sN^{\cR}(\sC,\sW)$ as an $m\times n$ commutative grid:
    \begin{equation}
       \begin{tikzcd}[ampersand replacement=\&]
        {c_{0,0}} \&\& {c_{0,1}} \&\& \cdots \&\& {c_{0,n}} \\
        {c_{1,0}} \&\& {c_{1,1}} \&\& \cdots \&\& {c_{1,n}} \\
        \vdots \&\& \vdots \&\& \vdots \&\& \vdots \\
        {c_{m,0}} \&\& {c_{m,1}} \&\& \cdots \&\& {c_{m,n}}
        \arrow[from=1-1, to=1-3]
        \arrow["\shortmid"{marking}, from=1-1, to=2-1]
        \arrow[from=1-3, to=1-5]
        \arrow["\shortmid"{marking}, from=1-3, to=2-3]
        \arrow[from=1-5, to=1-7]
        \arrow["\shortmid"{marking}, from=1-5, to=2-5]
        \arrow["\shortmid"{marking}, from=1-7, to=2-7]
        \arrow[from=2-1, to=2-3]
        \arrow["\shortmid"{marking}, from=2-1, to=3-1]
        \arrow[from=2-3, to=2-5]
        \arrow["\shortmid"{marking}, from=2-3, to=3-3]
        \arrow[from=2-5, to=2-7]
        \arrow["\shortmid"{marking}, from=2-5, to=3-5]
        \arrow["\shortmid"{marking}, from=2-7, to=3-7]
        \arrow["\shortmid"{marking}, from=3-1, to=4-1]
        \arrow["\shortmid"{marking}, from=3-3, to=4-3]
        \arrow["\shortmid"{marking}, from=3-5, to=4-5]
        \arrow["\shortmid"{marking}, from=3-7, to=4-7]
        \arrow[from=4-1, to=4-3]
        \arrow[from=4-3, to=4-5]
        \arrow[from=4-5, to=4-7]
    \end{tikzcd} 
    \end{equation}
    where each $c_{i,j} \in \sC$ for $(i,j)\in [m] \times [n]$, and the vertical marked arrows are in $\sW$, while the horizontal unmarked ones are arbitrary morphisms in $\sC$.
\end{remark}

\begin{observation}\label{ch2:: obs: nerves agree}
If $\sC_{\mathrm{id}}\subseteq \sC$ denotes the wide subcategory containing
only identity morphisms, then Rezk's classifying diagram
$\sN^{\cR}(\sC,\sC_{\mathrm{id}})$ recovers the simplicial nerve $\sN(\sC)$.
\end{observation}

\begin{definition}
    Let $\sC$ be a category, as a way of simplifying notation and naming, let us define the Rezk nerve to be \[\sN^{\cR}(\sC) \,\, \overset{\sf def}{=} \,\, \sN^{\cR}(\sC,\sC^{\simeq})\] the Rezk classifying diagram of $\sC$ relative to its maximal underlying groupoid.
\end{definition}

\begin{proposition}\label{ch2:: prop: Rezk nerve equiv disc nerve}\cite[Proposition 2.2.16]{Zha13}
Given a category $\sC$, there is a canonical inclusion of Segal spaces 
\begin{equation}
    \sN(\sC) \xlonghookra{\,\,\,\,} \sN^{\cR}(\sC)
\end{equation}    
which is a Dwyer--Kan equivalence.
\end{proposition}

\begin{example}\label{ch2:: ex: need complete}
    Let us consider $G$ a finite group. Let $\sB G$ denote the category with a single object $*$ and $\Aut_{\sB G}(*) = G$. There is a pair of Segal spaces we can construct out of this category:
\begin{enumerate}
    \item \textbf{Discrete nerve}: We can define the simplicial nerve of $\sB G$, or as noted in \Cref{ch2:: obs: nerves agree} the Rezk classifying diagram of $\sB G$ relative to $\{\id_G\}$. This is a constant simplicial space, i.e., every $m$th level is given by the discrete set \[\sN(\fun^{\{\id_* \}}([m],\sB G))\,\simeq\, G^{\times m}.\] In particular, the $0$th space is the one-point set.
    \item \textbf{Rezk nerve}: Consider the Rezk nerve $\sN^{\cR}(\sB G)$. Its $0$th space is given by $\sN(\fun^{\sB G}([0],\sB G))$, in particular 
    \begin{equation}
        \sN(\fun^{\sB G}([0],\sB G))_{\bullet} \,\, \simeq \,\, G^{\times \bullet}.
    \end{equation} 
\end{enumerate}

\nid By \Cref{ch2:: prop: Rezk nerve equiv disc nerve} we have that 
\begin{equation}
    \sN(\sB G) \hookra \sN^{\cR}(\sB G, \sB G)
\end{equation}
is a Dwyer--Kan equivalence. On the other hand, it is manifest that they are not levelwise equivalent.
\end{example}

\Cref{ch2:: ex: need complete} suggests that we should apply a left Bousfield
localization of $\mathsf{SSp}$ along the Dwyer--Kan equivalences to obtain the
correct homotopy theory of $(\oo,1)$-categories. What Rezk shows in
\cite{Rezk01}, is that this is a \emph{reflective localization}, meaning that
the localization functor admits a fully-faithful right adjoint. The image of
this right adjoint contains the Segal spaces whose space of objects captures
the correct homotopy type i.e., the homotopy type of homotopy equivalences.
These Segal spaces will be called complete, and the next sequence of
definitions gives a more precise characterization of these objects.

\begin{definition}\label{ch2:: def: complete Segal}
    Let $\cF$ be a Segal space. We say that $\cF$ is a \define{complete Segal space} if the functor
    \begin{equation}
        \Map^{h}(-, \cF) \maps \mathsf{SSp}^{\op} \longra \sset
    \end{equation}
    sends the following commutative diagram 
    \begin{equation}
        \label{eq complete}
        \begin{tikzcd}[ampersand replacement=\&]
        {\bDelta[\{0<2\}] \coprod \bDelta[\{1<3\}]} \&\& {\bDelta[\{0<1<2<3\}]} \\
        \\
        {\bDelta[0] \coprod\bDelta[0]} \&\& {\bDelta[0]}
        \arrow[from=1-1, to=1-3]
        \arrow[from=1-1, to=3-1]
        \arrow[from=1-3, to=3-3]
        \arrow[from=3-1, to=3-3]
    \end{tikzcd}
    \end{equation}
    to a homotopy pullback square in $\sset$.
\end{definition}

To make sense of \Cref{ch2:: def: complete Segal}, we must define what we mean by the space of homotopy equivalences. 

\begin{definition}
    Let $\cF$ be a Segal space. We define its \define{space of homotopy equivalences} to be
    \begin{equation}
        \cF_{1}^{\sf eq} \,\, \overset{\sf def}{=} \,\, \Map_{\mathsf{SSp}}^h(J,\cF)
    \end{equation}
    where $J$ is the following Segal space\footnote{Note that every object is cofibrant and that $\cF$ is fibrant, so there is no need for replacements to achieve the correct mapping space.} 
    \begin{equation}
        J \,\, \overset{\sf def}{=} \,\, \Delta[\{0<1<2<3\}] \coprod_{\bDelta[0]\coprod \bDelta[0]} \left( \bDelta[\{0<2\}] \sqcup \bDelta[\{1<3\}] \right),
    \end{equation}
    defined by contracting the $(0\ra 2)$ and the $(1\ra 3)$ edges of the 3-simplex. It is the pushout of diagram (\ref{eq complete}). The resulting simplicial space can be visualized as follows:
    \begin{equation}
        \begin{tikzcd}[ampersand replacement=\&]
            \&\& 2 \\
            \\
            0 \&\&\&\& 3 \&\& 0 \& 1 \& 0 \& 1 \\
            \&\& 1
            \arrow[from=1-3, to=3-5]
            \arrow[from=3-1, to=1-3]
            \arrow[between={0.6}{1}, dotted, from=3-1, to=3-5]
            \arrow[between={0}{0.4}, dotted, no head, from=3-1, to=3-5]
            \arrow[from=3-1, to=4-3]
            \arrow[between={0.2}{0.8}, squiggly, from=3-5, to=3-7]
            \arrow[from=3-7, to=3-8]
            \arrow[curve={height=-20pt}, equals, from=3-7, to=3-9]
            \arrow[from=3-8, to=3-9]
            \arrow[curve={height=20pt}, equals, from=3-8, to=3-10]
            \arrow[from=3-9, to=3-10]
            \arrow[from=4-3, to=1-3]
            \arrow[from=4-3, to=3-5]
        \end{tikzcd}
    \end{equation}
\end{definition}

\begin{observation}
    Note how giving a map $J \longra \cF$ is equivalent to specifying a morphism $f\maps x \ra y$ together with a left inverse $l$ and a right inverse $r$
    \begin{equation}
        \begin{tikzcd}[ampersand replacement=\&]
	y \& x \& y \& x.
        \arrow["l", from=1-1, to=1-2]
        \arrow[curve={height=-20pt}, equals, from=1-1, to=1-3]
        \arrow["f", from=1-2, to=1-3]
        \arrow[curve={height=20pt}, equals, from=1-2, to=1-4]
        \arrow["r", from=1-3, to=1-4]
    \end{tikzcd}
    \end{equation}
    thereby asserting that $f$ must be a homotopy equivalence. Note also that fixing $f\maps x \ra y$, the homotopy fiber of the map $ \cF^{\sf eq}_{1} \longra \cF_1$ induced by the inclusion of $\bDelta[\{1<2\}] \hookra J$:
    \begin{equation}
        \cF^{\sf eq}_{1} \overset{h}{\underset{\cF_1}{\times}} \{f \maps x \ra y\},
    \end{equation}   
    which contains all left and right homotopy inverses to $f$, might well be empty. This simplicial set can only be either contractible or empty, as one would hope. A proof of this fact can be found, for example, in \cite[Corollary 2.44]{NR18}.
\end{observation}

\begin{observation}
    Revisiting \Cref{ch2:: def: complete Segal}, we see that a complete Segal space $\cF$ is one for which the dashed arrow 
    \begin{equation}
        \begin{tikzcd}[ampersand replacement=\&]
            {\cF([0])} \\
            \& {\cF_1^{\sf eq}} \&\& {\cF([\{0<1<2<3\}])} \\
            \\
            \& {\cF([0])\times \cF([0])} \&\& {\cF([\{0<2\}])\times \cF([\{1<3\}])}
            \arrow["{\rotatebox[origin=c]{-35}{$\simeq$}}", dashed, from=1-1, to=2-2]
            \arrow[curve={height=-18pt}, from=1-1, to=2-4]
            \arrow[curve={height=18pt}, from=1-1, to=4-2]
            \arrow[from=2-2, to=2-4]
            \arrow[from=2-2, to=4-2]
            \arrow["\lrcorner"{anchor=center, pos=0.125}, draw=none, from=2-2, to=4-4]
            \arrow[from=2-4, to=4-4]
            \arrow[from=4-2, to=4-4]
        \end{tikzcd}
    \end{equation}
    is a weak equivalence of simplicial sets. Complete Segal spaces $\cF$ are thus Segal spaces for which the inclusion
    \begin{equation}
        \cF_0 \longra \cF_{1}^{\sf eq}
    \end{equation} 
    is an equivalence.
\end{observation}

Finally, we have the following result of Rezk.

\begin{proposition}\label{ch2:: prop: Rezk-model-structure}
    There is a model structure $\mathsf{CSSp}$, obtained by left Bousfield
    localizing $\mathsf{SSp}$ along the terminal map
    \[\label{local map complete} J \longra \bDelta[0],\]
    and it has the following properties:
    \begin{enumerate}
        \item The cofibrations are the monomorphisms.
        \item The fibrant objects are the complete Segal spaces.
        \item The weak equivalences among complete Segal spaces are the levelwise equivalences.  
    \end{enumerate}
    Additionally, the model structure is left proper, combinatorial and Cartesian closed.
\end{proposition}

\begin{theorem}\cite[Theorem 7.7]{Rezk01}
    Let $\cF$ and $\cG$ be Segal spaces. A morphism $F \maps \cF \longra \cG$
    is a weak equivalence in $\mathsf{CSSp}$ if and only if it is a
    Dwyer--Kan equivalence. 
\end{theorem}

In particular, this means that despite $\sN(\sB G)$ not being a complete Segal
space, therefore not fibrant in $\mathsf{CSSp}$, the Dwyer--Kan equivalence
\begin{equation}
    \sN(\sB G) \longra \sN^{\cR}(\sB G, \sB G)
\end{equation}
into its Rezk nerve is a weak equivalence, as desired.

\begin{remark}
    If necessary, we write $\cC \mathsf{SSp}\subset \mathsf{SSp}$ for the full
    subcategory of complete Segal spaces, as a way of distinguishing it from
    the localized model category $\mathsf{CSSp}$.
\end{remark}

\begin{remark}
    As a consequence of the fact that $\mathsf{CSSp}$ is Cartesian closed we
    have that given a Segal space $\cF$ and a complete Segal space $\cG$, the
    simplicial space $\fun(\cF,\cG)$ is also a complete Segal space. We will
    refer to it as the $(\oo,1)$-category of functors between $\cF$ and $\cG$.
\end{remark}

\section{Interval spaces}\label{ch2:: subsec: interval_spaces}

When working with complete Segal spaces as a model for $(\oo,1)$-categories, the combinatorial engine, the object doing the work of encoding all the higher homotopical data, is provided by the simplex category $\bDelta$. Yet, as the literature on $(\oo,1)$-categories has evolved, the community has found it useful to adopt variants, which facilitate their direct application. As we will unravel throughout the remainder of the thesis, there is a variant of $\bDelta$, called the \emph{interval category}, which provides a more flexible framework for working with $(\oo,1)$-categories of cobordisms and higher versions thereof. As such, this subsection is dedicated to introducing its mechanics.

\begin{definition}\label{ch2:: def: interval_category}
    The interval category $\mathsf{Int}$ is the category whose objects are finite disjoint unions of open intervals $I \subseteq [0,1]$ containing the endpoints $\{0,1\}$ but are \textbf{not} the whole interval $[0,1]$, and whose morphisms are the inclusions $(I \longra J):= I\subseteq J$ of open intervals. 
\end{definition}

\begin{observation}\label{ch2:: obs: poset_observation}
    The category $\mathsf{Int}$ is a poset. Manifestly, $\Int$ is a full subcategory spanned by a subset of objects of $\cU([0,1])$, the poset of opens of the unit interval. 
\end{observation}

For convenience of the reader and as a means to stipulate notation for the rest of the chapter, we introduce the following diagrammatic means of depicting the objects of $\mathsf{Int}$:
\medskip
\begin{equation}\label{ch2:: eq: interval_object}
    I \,\, = \,\, \tikzfig{figures/Ch2/INT1}{\input{figures/Ch2/INT1.tikz}}
\end{equation}
\medskip

\noindent The diagram above depicts an object $I$ in $\mathsf{Int}$. It
consists of a finite disjoint union of four open intervals which contain both
endpoints but do not exhaust the whole interval. Unless we need to specify the
start- and end-points of the intervals, we will refrain from doing so, which,
given the topological nature of our applications, will prove to be a reasonable
default.

Additionally, we will often depict morphisms in this category by writing diagrams like the following one: 
\medskip
\begin{equation}\label{ch2:: eq: interval_morphism}
    \tikzfig{figures/Ch2/INT2}{\input{figures/Ch2/INT2.tikz}}
\end{equation}
\medskip

\noindent which codifies a unique morphism from the object $I$ to the object $J$. Every dashed edge is merely meant to indicate how these objects relate to one another and are not data associated to the objects or the morphisms themselves. In fact, if the morphism is clear we will sometimes avoid the dashed edges.

\begin{definition}\label{ch2:: def: face_maps_intervals}
    Let $I\subseteq [0,1]$ be an object of $\mathsf{Int}$. We can identify its complement with a disjoint union of finitely many closed intervals of $[0,1]$:
    \begin{equation}
        [0,1]\setminus I = \coprod_{k = 0}^{n} A_k
    \end{equation}
    where $A_k$ are closed subsets of $[0,1]$ which do not contain $0$ or $1$. Let $\sU(\{A_i\}_{i=0}^n)$ denote the open subset which is the complement of the union of closed subsets $A_i$ in $[0,1]$. We define 
    \begin{equation}
        \alpha_{i}^{I}\maps \sU(\{A_1,\cdots,A_n\}) \longrightarrow \sU(\{A_1,\cdots,\widehat{A_{i}},\cdots,A_n\})
    \end{equation} for every $i=0,\cdots,n$ as the map which \emph{forgets} the $i$-th closed interval $A_i$. We will call these maps \textbf{gluings}.

    \medskip

    \nid Let us define 
    \begin{equation}\label{ch2:: eq: degeneracy_map_interval}
        \beta_{i}^{I}\maps \sU(\{A_1,\ldots, A_n\}) \longleftarrow \sU(\{A_1,\ldots, \mathsf{Hull}(A_i,A_{i+1}),\ldots, A_n\}) 
    \end{equation}
    for every $i=1,\ldots, n-1$, as the inclusion which \emph{merges} the intervals $A_i$ and $A_{i+1}$, replacing it with its \emph{hull}. We will call these maps \textbf{refinements}.
\end{definition}

\medskip

\begin{remark}
    By $\mathsf{Hull}(A_{1},\cdots, A_n)$ we mean the smallest closed interval containing all the closed intervals $A_i$ for $i=1,\cdots,n$. In particular, if $A_i = [a_i,b_i]$ and $A_{i+1} = [a_{i+1},b_{i+1}]$, then $\mathsf{Hull}(A_i,A_{i+1}) = [\min(a_i,a_{i+1}),\max(b_i,b_{i+1})]$.
\end{remark}

\begin{example}
    Let us pictorially exhibit examples of face maps in $\mathsf{Int}$, so as to get the reader a better feel for how they relate to their analogues in $\bDelta^{\op}$. Let $I\in\mathsf{Obj}(\mathsf{Int})$ be depicted as
    \begin{equation}\label{ch2:: eq: face_map_1}
        \medskip
        \sU(\{A_0,A_1,A_2\}) \,\, = \,\, \tikzfig{figures/Ch2/INT3}{\input{figures/Ch2/INT3.tikz}}.
        \medskip
    \end{equation}
    Then we have the following triple of face maps: 
    \begin{equation}
        \medskip
        \tikzfig{figures/Ch2/INT5}{\input{figures/Ch2/INT5.tikz}}
        \medskip
    \end{equation}
    which forgets $A_0$, and
    \begin{equation}
        \medskip
        \tikzfig{figures/Ch2/INT4}{\input{figures/Ch2/INT4.tikz}}
        \medskip
    \end{equation}
    which forgets $A_1$, and finally
    \begin{equation}
        \medskip
        \tikzfig{figures/Ch2/INT6}{\input{figures/Ch2/INT6.tikz}}
        \medskip
    \end{equation}
    which forgets $A_2$.
\end{example}

\noindent We introduce the notation above also as a means of aiding our description of the morphisms in $\mathsf{Int}$, given how inherently geometric they are — we like to think of them as a discrete version of the  configuration category of points in the open interval —  but also to initiate the discussion of how it relates to the simplex category $\bDelta$.

\begin{notation}
Let $[i,j] \subseteq \{1,\ldots,n\}$ for $0\leq i < j\leq n$ denote the subset which contains all consecutive integers from $i$ through $j$. For $I = \sU(A_1,\cdots,A_n)$, we denote the maps 
\begin{equation}
    \alpha_{[i,j]}^{I}\maps \sU(A_{1},\ldots, A_{n}) \longra \sU(A_1, \ldots,  A_{i-1}, A_{j+1}, \ldots, A_{n}) 
\end{equation}
which forget the consecutive closed intervals $A_i,\ldots,A_j$.
\end{notation}

\begin{observation}\label{ch2:: obs: interval_category_small}
The interval category $\Int$ is a poset category, so it is small. 
\end{observation}
\begin{remark}
    The opposite $\Int^{\op}$ of the interval category is equivalent to a category with an equally geometric interpretation. Let $\cInt$ be the category whose objects are finite non-empty disjoint unions of connected closed intervals $K\subset [0,1]$ which do not intersect $\{0,1\}$. A morphism $K\longra L$ is given by inclusion — so $\cInt$ is obviously also a poset. The equivalence is defined as 
    \begin{align}
        (-)^{c}\maps \Int^{\op} &\xra{\,\,\,\cong \,\,\,} \cInt \\
        I &\longmapsto \big( I^{c} = [0,1]\setminus I\big).
    \end{align}
    where an inclusion $I\ra J$ naturally induces an inclusion $J^{c} \ra I^{c}$. This equivalence is manifest, since it is bijective poset map.
\end{remark}

\medskip

\nid When convenient, we will use this isomorphism to shift perspectives which suit a given narrative — minding, of course, the fundamental change in variance.

\begin{definition}
    Let 
    \begin{equation}
        \gamma \maps \Int \longra \bDelta^{\op}, \quad (I \subset [0,1])  \longmapsto \pi_{0}(I^{c}),
    \end{equation}
    be the functor which sends an interval $I$ to the set of connected components of its complement in $[0,1]$. Given an inclusion a morphism $I \subseteq J$, $\gamma$ maps 
    \begin{equation}
        \gamma (I \subseteq J) = \big(\,\pi_{0}(J^{c}) \longra \pi_0(I^{c})\,\big).
    \end{equation} 
\end{definition}

\begin{observation}
    We have the following canonical commutative diagram of categories:
    \begin{equation}
        \begin{tikzcd}[ampersand replacement=\&]
            \Int \&\& {\bDelta^{\op}} \\
            \\
            {\cInt^{\op}} \&\& {.}
            \arrow["\gamma", from=1-1, to=1-3]
            \arrow["{\rotatebox[origin=c]{90}{$\sim$}}"', from=1-1, to=3-1]
            \arrow["{\pi_0}"', curve={height=6pt}, dashed, from=3-1, to=1-3]
        \end{tikzcd}
    \end{equation}
\end{observation}

\begin{definition}\label{ch2:: def: updelta cat}
    Let $\sO_{\pm}$ be the following category:
    \begin{itemize}
        \item \textbf{Objects}: Linearly ordered finite sets $S$ with a distinct minimum and maximum.
        \item \textbf{Morphisms}: Order-preserving maps $f \maps S \longra T$ that preserve extrema. 
    \end{itemize}
\end{definition}

\begin{remark}
    Let $[p]\in \bDelta$, note that the set $\Hom_{\bDelta}([p],[1])$ is a linearly ordered finite set with a distinct minimum and maximum. For simplicity, we will denote an object in this set as 
    \begin{equation}
        \lag i_0, \ldots, i_p \rag,
    \end{equation}
    where $0 \leq i_{k} \leq i_{k+1} \leq 1$ for every $0\leq k < p$. We declare
    \begin{equation}
    \lag i_0, \ldots, i_p \rag \leq \lag j_0, \ldots, j_p \rag, \quad \text{ if } \quad i_k \leq j_k        
    \end{equation}
    for every $0\leq k \leq p$. This is a linear order precisely because $[p]$ is linearly ordered with the minimum and maximum elements identified respectively as
    \begin{equation}
        \lag 0, \ldots, 0\rag, \quad \lag 1,\ldots, 1 \rag. 
    \end{equation}
    Equivalently, these maps can be identified with splittings
    \begin{equation}
        [p] \,\, \simeq \,\, [q]\star [p-1-q],
    \end{equation}
    where $q = -1, 0,\ldots, p$ — by identifying $\alpha \maps [p] \ra [1]$ with $\alpha^{-1}([0])\star \alpha^{-1}([1])$.\footnote{We establish the convention that $[-1] = \emptyset$.}
\end{remark}

\begin{notation}
    Let $P$ be a poset. We define $P^{\tr}$ to be the poset obtained by adjoining $P \star [0]$ a maximal element. Analogously, we define $P^{\,\tl}$ to be the poset obtained by adjoining $[0] \star P$ a minimal element. 
\end{notation}

\begin{observation}
There is an identification in $\sO_{\pm}$
\begin{equation}
\Hom_{\bDelta}([p],[1]) \,\, \simeq \,\, [p]^{\,\tl},    
\end{equation}
defined by sending the splitting $[q] \star [p-1-q]$ to $q \in [p]^{\,\tl}$.
\end{observation}

\begin{lemma}
    The functor
    \begin{equation}
        \bDelta^{\op} \longra \sO_{\pm}, \quad \big(\,[p] \mapsto \Hom_{\bDelta}([p],[1]) \, \big), \quad \big(\, \bDelta \ni ([p] \xra{\alpha} [q]) \mapsto \Hom_{\bDelta}([q],[1]) \xra{\alpha^*} \Hom_{\bDelta}([p],[1]) \,\big)
    \end{equation}
    is an equivalence of categories.
\end{lemma}
\begin{proof}
    Surjectivity on objects follows from the fact that given a finite linearly ordered set with distinct minimum and maximum $S$, there is canonical isomorphism $[p]^{\,\tl}$ for some $p\geq 0$. Fully-faithfullness is manifest from the definition.
\end{proof}

\begin{observation}
    We have the following commutative diagram of categories:
    \smallskip
    \begin{equation}\label{diagram combinatorics}
        \begin{tikzcd}[ampersand replacement=\&]
            \&\& {\sO_{\pm}} \\
            \\
            \Int \&\& {\bDelta^{\op}} \\
            \\
            {\cInt^{\op}}
            \arrow["{\pi_{0}}", curve={height=-12pt}, from=3-1, to=1-3]
            \arrow["\gamma", from=3-1, to=3-3]
            \arrow["{\rotatebox[origin=c]{90}{$\sim$}}"', from=3-1, to=5-1]
            \arrow["{(-)^{\,\tl}}"', from=3-3, to=1-3]
            \arrow["{\rotatebox[origin=c]{90}{$\sim$}}", draw=none, from=3-3, to=1-3]
            \arrow["{\pi_0}"', curve={height=6pt}, from=5-1, to=3-3]
        \end{tikzcd}
    \end{equation}
    In particular, given a $J \in \cInt$ we have that $\pi_0(J^c) = \gamma_{0}(J^{c})^{\, \tl}$. As innocuous and simple as this fact is, it will serve us well in the upcoming sections.
\end{observation}

Understanding the yoga of~(\ref{diagram combinatorics}) in more depth and in a graphical way will pay dividends in the next sections, so let us add some elucidating remarks.

\begin{remark}\label{ch2:: rmk:: graph calc}
    The simplex category $\bDelta$ admits a particularly simple graphical calculus. We depict an object $[n]$ as a collection of $[n]$-marked points in a horizontal array:
    \begin{equation}
        \tikzfig{figures/Ch2/INT10}{\input{figures/Ch2/INT10.tikz}}
    \end{equation}
    \smallskip

    \nid A morphism $\alpha \maps [n] \ra [m]$ is depicted by drawing connecting 1-dimensional open intervals between the marked points. This is best understood through examples:
    
     \begin{equation}
       \tikzfig{figures/Ch2/INT9}{\input{figures/Ch2/INT9.tikz}}, \qquad \quad  \tikzfig{figures/Ch2/INT11}{\input{figures/Ch2/INT11.tikz}}
    \end{equation}
    \smallskip
    
    \nid We might call this a graphical calculus since if we allow the 1-dimensional arcs to be isotoped, while keeping them anchored at the marked points, we get equivalent representations of the same morphism. To facilitate a better visual interpretation for how the functor $(-)^{\sf \tl}\maps \bDelta^{\op} \longra \sO_{\pm}$ acts on objects and morphisms, we extend the graphical calculus of $\bDelta$ by embedding the $[n]$-marked points into the open interval $(0,1)$.\footnote{The choice of the open interval $(0,1)$ is for practical reasons — any other space homeomorphic to $\RR$ would have done the job.}
    \begin{equation}
        \tikzfig{figures/Ch2/INT10}{\input{figures/Ch2/INT10.tikz}} \,\,\, \xlonghookra{\,\,\,\,\,\,\,} \,\,\,\tikzfig{figures/Ch2/INT12}{\input{figures/Ch2/INT12.tikz}}.
    \end{equation}
    \smallskip

    \nid Additionally, we label the gaps, the 1-dimensional open intervals, between the $[n]$-marked points with the totally ordered set $[n]^{\, \tl}$ as follows
    \begin{equation}
        \tikzfig{figures/Ch2/INT13}{\input{figures/Ch2/INT13.tikz}}
    \end{equation}
    \smallskip  
    
    \nid This makes clear how the functor $(-)^{\, \tl}$ acts on morphisms. The simple idea is that it tracks the induced mapping on the 1-dimensional open intervals, as follows:
    \begin{equation}
         \tikzfig{figures/Ch2/INT14}{\input{figures/Ch2/INT14.tikz}} \qquad \qquad \tikzfig{figures/Ch2/INT15}{\input{figures/Ch2/INT15.tikz}}
    \end{equation} 
\end{remark}

\subsection{Interval spaces}

When working with complete Segal spaces as a model for $(\oo,1)$-categories, all the combinatorics is encoded in the simplex category $\bDelta$. In this chapter we instead use the interval category $\Int$ as our indexing category. One should think of $\Int$ as keeping track of \emph{where} composition of morphisms is allowed to happen (via unions of open subintervals of $[0,1]$), while the functor
\begin{equation}
  \gamma \colon \Int \longra \bDelta^{\op}
\end{equation}
forgets this extra geometric information and remembers only the number and order of “slots” for composition (the number of inner connected components).

\begin{definition}
A \define{pre-interval space} is a functor
\begin{equation}
  \cF \colon \Int \longra \sset.
\end{equation}
\end{definition}

The fibres of $\gamma$ consist of intervals which have the same underlying combinatorics in $\bDelta^{\op}$ but differ by refining or sliding the closed pieces in their complements. We single out the class of maps
\begin{equation}
  \cW_{\gamma} \subset \Fun(\Int,\sset)
\end{equation}
consisting of morphisms of co-representables
\begin{equation}
  \Hom_{\Int}(J,-) \xra{\,f_*\,} \Hom_{\Int}(I,-)
\end{equation}
induced by morphisms $f\colon I \to J$ in $\Int$ such that $\gamma(f)=\id_{[n]}$ for some $n$. Equivalently, these are the maps coming from morphisms which induce an isomorphism on path-components. A pre-interval space which sends all such $f$ to weak equivalences is, heuristically, “constant along the
fibres of $\gamma$”.

\begin{definition}
A pre-interval space $\cF\colon \Int \to \sset$ is called an \define{interval space} if for every morphism $f\colon I \to J$ in $\Int$ with $\gamma(f)=\id$ the induced map
\begin{equation}
  \cF(I) \xra{\ \cF(f)\ } \cF(J)
\end{equation}
is a weak equivalence of simplicial sets.
\end{definition}

Thus, an interval space is a pre-interval space whose homotopy type only depends on the image of $I$ in $\bDelta^{\op}$. In §\ref{ch2:: subsec: interval_spaces_model_infinity_cats} we will make this precise by
left Bousfield localizing $\Fun(\Int,\sset)$ along $\cW_{\gamma}$ and showing that the resulting model category is Quillen equivalent to Rezk's model structure on simplicial spaces.

The next definition is the analogue of the Segal condition for simplicial spaces, but phrased geometrically in terms of cutting an interval into smaller pieces.

\begin{definition}\label{ch2:: def: interval Segal condition}
Let $\cF\colon \Int \to \sset$ be a pre-interval space. We say that $\cF$ satisfies the \define{interval Segal condition} if for every $n \geq 2$ and every $\sU(A_1,\ldots,A_n)\in\Int$ the induced
commuting square
\begin{equation}
    \begin{tikzcd}[ampersand replacement=\&,column sep=2.25em]
        {\cF(\sU(A_1,\ldots, A_n))} \&\& {\cF(\sU(A_1,\ldots, A_k))} \\
        \\
        {\cF(\sU(A_k,\ldots, A_n))} \&\& {\cF(\sU(A_k))}
        \arrow["{(\alpha_{[1,k]})_*}", from=1-1, to=1-3]
        \arrow["{(\alpha_{[k,n]})_*}"', from=1-1, to=3-1]
        \arrow["{(\alpha_{[k]})_*}", from=1-3, to=3-3]
        \arrow["{(\alpha_{[k]})_*}"', from=3-1, to=3-3]
    \end{tikzcd}
\end{equation}
is a homotopy pullback square in $\sset$, for every $1 \le k \le n-1$.
\end{definition}

Informally, one should think of $\sU(A_1,\ldots,A_n)$ as encoding a string of $n$ composable 1-morphisms, with the closed intervals $A_i$ marking where the topology is allowed to change. The interval Segal condition says that gluing along any intermediate cut $A_k$ recovers $\cF(\sU(A_1,\ldots,A_n))$
as the appropriate homotopy pullback of the two shorter strings.

\begin{definition}[Complete interval Segal space]
An \define{interval Segal space} is a pre-interval space which is both an interval space and satisfies the interval Segal condition. We will say that an interval Segal space $\cF$ is \define{complete} if its associated
Segal space
\begin{equation}
  \Rrkan_{\gamma}(\cF)\colon \bDelta^{\op} \longra \sset
\end{equation}
is a complete Segal space in the sense of Rezk.   
\end{definition}

We will justify this definition and construct a model structure whose fibrant objects are complete interval Segal spaces in §\ref{ch2:: subsec: interval_spaces_model_infinity_cats}, where we prove that this model structure is Quillen equivalent to the Rezk model structure on complete Segal spaces.

\section{Interval Segal spaces as a model for \texorpdfstring{$(\oo,1)$}{(infty,1)}-categories}
\label{ch2:: subsec: interval_spaces_model_infinity_cats}

In this subsection we construct a model structure on the category
$\fun(\Int,\sset)$ of pre-interval spaces and prove that it is Quillen
equivalent to the Rezk model structure \cite{Rezk01} on complete Segal spaces
described in \Cref{ch2:: prop: Rezk-model-structure}. In particular, we obtain
a left proper, combinatorial model structure on the category of pre-interval
spaces whose fibrant objects are complete interval Segal spaces and whose weak
equivalences are the Dwyer--Kan equivalences.

\medskip

\nid The structure of this subsection can be understood through progressively more refined model structures on the category of pre-interval spaces. We begin by localizing pre-interval spaces along $\pi_0$-isomorphisms, which yields a model category Quillen equivalent to simplicial spaces; we then further localize along the derived image of the spine inclusions to recover a Segal condition for interval spaces.

\subsection{Equivalence with simplicial spaces}

\begin{definition}\label{ch2:: def: pi_0 equiv}
    Let $\cW_{\pi_0}$ be the set of morphisms in $\fun(\Int,\sset)$ which are the image of $\pi_0$-isomorphisms under the co-Yoneda embedding 
    \begin{equation}
        \Int^{\op} \hookrightarrow \fun(\Int,\sset).
    \end{equation}
\end{definition}

\medskip

Recall the discussion in
\Cref{ch2:: subsubsec: homotopy-limits-bousfield} about left Bousfield
localizations. We are interested in studying pre-interval spaces
$\cF\colon \Int \to \sset$ which send every $\pi_0$-isomorphism $f\colon I \to
J$ to a weak equivalence $\cF(I) \xra{\simeq} \cF(J)$ in $\sset_{\sf QK}$, and
we would like to design a model structure for which these are precisely the
fibrant objects.

\begin{observation}
    A morphism in $\cW_{\pi_0}$ is a map of co-representables
    \begin{equation}
        \Hom_{\Int}(J,-) \xra{\,f_*\,} \Hom_{\Int}(I,-)
    \end{equation}
    such that 
    \[
      \gamma(I) \xra{\gamma(f)} \gamma(J)
    \]
    is an isomorphism, i.e., the identity since $\bDelta^\op$ is gaunt.
\end{observation}

\begin{proposition}\label{ch2:: prop: model structure interval}
    There exists a model structure on the category of pre-interval spaces $\fun(\Int,\sset)$, which we will denote $\fun(\Int,\sset)_{\mathsf{inj},h}$, with the following properties:
    \begin{enumerate}
        \item The weak equivalences are the $\cW_{\pi_0}$-local equivalences.
        \item The cofibrations are the injective cofibrations.
        \item The fibrant objects are the $\cW_{\pi_0}$-local objects, i.e.\ injectively fibrant functors $\cF$ such that $\cF(I) \xra{\cF(f)} \cF(J)$ is a level-wise weak equivalence of simplicial sets for every $f\colon I \to J$ in $\cW_{\pi_0}$.
        \item The model structure is left proper, combinatorial, and simplicial.
    \end{enumerate} 
\end{proposition}

\begin{proof}
We construct $\fun(\Int,\sset)_{\mathsf{inj},h}$ as the left Bousfield
localization of $\fun(\Int,\sset)_{\mathsf{inj}}$ with respect to the set
$\cW_{\pi_0}$ of inclusions inducing an isomorphism on $\pi_0$.

To see that this localization exists, note first that
$\sset_{\mathsf{QK}}$ is a cofibrantly generated, left proper, simplicial
model category. Since $\Int$ is small, the general existence theorem for
injective model structures on functor categories (see, for instance,
\Cref{ch2:: thm: functor model structures}) implies that
$\fun(\Int,\sset)_{\mathsf{inj}}$ exists and is again cofibrantly generated,
left proper, combinatorial, and simplicial. In particular, $\fun(\Int,\sset)$
is locally presentable. The set $\cW_{\pi_0}$ consists of maps between
co-representables, so it is indeed a \emph{set} of morphisms in
$\fun(\Int,\sset)$. We may therefore apply the general existence theorem for
left Bousfield localizations (see \Cref{ch2:: thm: left Bousfield localization})
to obtain the localized model structure
$\fun(\Int,\sset)_{\mathsf{inj},h}$ with the claimed weak equivalences. This
establishes (1) and (4).

Assertion (2) is straightforward from the properties of left Bousfield
localization: the class of cofibrations is unchanged.

For (3), recall that $\cW_{\pi_0}$-local objects are precisely those
injectively fibrant pre-interval spaces $\cF$ such that, for every morphism
$\eta\colon \cG\to\cH$ in $\cW_{\pi_0}$, the induced map
\[
  \Hom_{\fun(\Int,\sset)}(\cH,\cF)
  \xlongrightarrow{\ \eta^*\ }
  \Hom_{\fun(\Int,\sset)}(\cG,\cF)
\]
is a weak equivalence of simplicial sets. By definition of $\cW_{\pi_0}$, any
$\eta$ is of the form
\[
  \Hom_{\Int}(J,-)\xrightarrow{\,f_*\,}\Hom_{\Int}(I,-)
\]
for some $f\colon I\to J$ in $\Int$ inducing an isomorphism on $\pi_0$. In
this case the Yoneda lemma identifies the induced map above with
\[
  \cF(J)\xrightarrow{\;\cF(f)\;}\cF(I).
\]
Thus $\cF$ is $\cW_{\pi_0}$-local if and only if it is injectively fibrant and
$\cF(f)$ is a weak equivalence for every $f\in\cW_{\pi_0}$, which is exactly
the description in (3).
\end{proof}

Having defined a model structure on $\fun(\Int, \sset)$ with the fibrant objects having the homotopical properties we originally desired, we now construct a Quillen pair between the injective model structure on simplicial spaces $\fun(\bDelta^{\op},\sset)_{\sf inj}$ and $\fun(\Int, \sset)_{\mathsf{inj},h}$. We will use the injective model structure to obtain a naive Quillen pair from which we will build the desired one between $\fun(\bDelta^{\op},\sset)_{\mathsf{inj}}$ and $\fun(\Int,\sset)_{\mathsf{inj},h}$.

\begin{proposition}\label{ch2:: prop: Quillen_pair_interval_simplicial_naive}
There is a Quillen pair
\begin{equation}
    \begin{tikzcd}[ampersand replacement=\&]
	{\fun(\bDelta^{\op},\sset)_{\mathsf{inj}}} \&\& {\fun(\Int,\sset)_{\mathsf{inj}}}
	\arrow[""{name=0, anchor=center, inner sep=0}, "{\gamma^*}", bend left=10, from=1-1, to=1-3]
	\arrow[""{name=1, anchor=center, inner sep=0}, "{\rkan_{\gamma}}", bend left=10, from=1-3, to=1-1]
	\arrow["\dashv"{anchor=center, rotate=-90}, draw=none, from=1, to=0]
    \end{tikzcd}
\end{equation}
induced by the functor $\gamma\colon \Int \to \bDelta^{\op}$. 
\end{proposition}

\begin{proof}
Precomposition with $\gamma$ preserves levelwise cofibrations and levelwise
weak equivalences. Since these are precisely the cofibrations and weak
equivalences in the injective model structures, $\gamma^*$ is a left Quillen
functor. Its right adjoint is the right Kan extension $\rkan_\gamma$.
\end{proof}

We want to show that this Quillen pair extends to a Quillen pair between
$\fun(\bDelta^{\op},\sset)_{\mathsf{inj}}$ and
$\fun(\Int,\sset)_{\mathsf{inj},h}$. Using
\Cref{ch2:: obs: interval_category_small} together with the general existence
results for injective model structures on functor categories (again
\Cref{ch2:: thm: functor model structures}), we know that
$\fun(\Int,\sset)_{\mathsf{inj}}$ is left proper and combinatorial, so its left
Bousfield localization at $\cW_{\pi_0}$ exists.

\begin{proposition}
    There is a Quillen pair 
    \begin{equation}
    \begin{tikzcd}[ampersand replacement=\&]
	{\fun(\bDelta^{\op},\sset)_{\mathsf{inj}}} \&\& {\fun(\Int,\sset)_{\mathsf{inj}}[\cW_{\pi_0}^{-1}]}
	\arrow[""{name=0, anchor=center, inner sep=0}, "{\gamma^*}", bend left=10, from=1-1, to=1-3]
	\arrow[""{name=1, anchor=center, inner sep=0}, "{\rkan_{\gamma}}", bend left=10, from=1-3, to=1-1]
	\arrow["\dashv"{anchor=center, rotate=-90}, draw=none, from=1, to=0]
    \end{tikzcd}
    \end{equation}
\end{proposition}

\begin{proof}
    Building on the Quillen pair from \Cref{ch2:: prop: Quillen_pair_interval_simplicial_naive}, this is a formal consequence of left Bousfield localization: the localization leaves the class of cofibrations unchanged and enlarges the class of weak equivalences. In particular, $\gamma^*$ continues to send level-wise cofibrations to injective cofibrations, and level-wise weak equivalences to level-wise weak equivalences, hence in particular to weak equivalences in the localized model structure. Therefore the adjunction descends to the localized model structure.
\end{proof}

Having shown that the adjoint pair $\gamma^{*} \dashv \rkan_{\gamma}$ extends to the localized model structure, we now proceed to show that it is in fact a Quillen equivalence. Before we delve into the details of the proof, it is helpful to offer a brief blueprint of the argument:
\begin{enumerate}
    \item[\textbf{\color{brown}{Step \ding{192}}}] We first check that if $\cF$ is fibrant, then the right Kan extension agrees with the homotopy right Kan extension:
    \begin{equation}
        \rkan_{\gamma}(\cF) \xra{\,\simeq\,} \Rrkan_{\gamma}(\cF).
    \end{equation}
    \item[\textbf{\color{brown}{Step \ding{193}}}] The homotopy right Kan extension is computed pointwise by 
    \begin{equation}
        \Rrkan_{\gamma}(\cF)([n]) \,\simeq\, \underset{I \in ([n]/ \Int)}{\mathsf{holim}} \, \cF(I).
    \end{equation}
    \item[\textbf{\color{brown}{Step \ding{194}}}] The full subcategory $\cJ_n$ of $([n]/ \Int)$ on the objects $\big([n] \xra{\id} \gamma(I)\big)$ is shown to be homotopy initial (see \Cref{ch2:: prop: hom initial} below).
    \item[\textbf{\color{brown}{Step \ding{195}}}] The previous observations imply that we can compute the homotopy right Kan extension as a homotopy limit over the category $\cJ_n$.
    \item[\textbf{\color{brown}{Step \ding{196}}}] Since $\cF$ is $\cW_{\pi_0}$-local, every morphism in $\cJ_n$ is mapped to a weak equivalence. Together with the fact that the nerve of $\cJ_n$ is weakly contractible, this implies that $\Rrkan_{\gamma}(\cF)([n]) \to \cF(I)$ is a weak equivalence for every $I\in \Int$ (cf.\ \Cref{ch2:: prop: homotopy-limit-homotopically-constant}).
\end{enumerate}

\subsubsection{Digression on configuration spaces}

Before we move on to proving that $\gamma^* \dashv \rkan_\gamma$ is a Quillen equivalence, let us first provide a brief recollection on configuration spaces of $\RR^{n}$. Let $r \geq 0$ be an integer, and let $\underline{r} = \{1,\ldots,r\}$. We define the space
\begin{equation}
    \mathsf{Conf}_{r}(\RR^{n}) \,\overset{\sf def}{=}\, \big\{\underline{r} \xlonghookra{\phi} \RR^{n}\, : \, \phi \text{ is an injection} \big\}
\end{equation}
of injections from the finite set $\underline{r}$ into $\RR^{n}$, topologized by removing all the closed diagonal subspaces 
\begin{equation}
    (\RR^{n})^{\underline{r}} \setminus \Big(\coprod_{0\leq i\neq j\leq r} \Delta_{i,j}\Big).
\end{equation}
We call this the \emph{ordered configuration space} of $r$ particles in $\RR^{n}$. The ordered configuration space of $r$ points in $\RR^{n}$ admits an action by $\Sigma_r$, the symmetric group on $r$ letters, defined by
\begin{align}
 \Sigma_r \times \mathsf{Conf}_r(\RR^{n}) &\longra \mathsf{Conf}_r(\RR^{n}) \\
 \big(\sigma,(x_1,\ldots,x_r) \big) &\longmapsto (x_{\sigma(1)},\ldots, x_{\sigma(r)}).
\end{align}

\begin{lemma}
    The action of the symmetric group $\Sigma_r \curvearrowright \mathsf{Conf}_r(\RR^{n})$ is free and properly discontinuous. As a result, the quotient map 
    \begin{equation}
        q \colon \mathsf{Conf}_{r}(\RR^{n}) \longra \mathsf{UConf}_{r}(\RR^{n})
    \end{equation} 
    is a regular covering map (of degree $r!$).
\end{lemma}

\begin{definition}
    We call the space $\mathsf{UConf}_{r}(\RR^{n})$ the \emph{unordered configuration space of $r$ points in $\RR^{n}$}.
\end{definition}

In general, configuration spaces of $\RR^{n}$ for arbitrary $n$ can have a complicated and interesting homotopy type, but not when $n=1$.

\begin{lemma}
    The space $\mathsf{UConf}_{r}(\RR) \cong \RR^{r}$. In particular, the unordered configuration space of $r$ points in $\RR$ is contractible.
\end{lemma}

The point of the preceding digression is that in dimension one the unordered
configuration spaces we will encounter are \emph{contractible}: they are open,
convex subsets of a Euclidean space. We now explain how these configuration
spaces appear naturally in the proof \textbf{\color{brown}{Step \ding{194}}}, stating that 
\begin{equation}
  \iota_n \maps \cJ_n \xlonghookra{\,\,\,\,\,\,\,} ([n]/\Int)
\end{equation}
is homotopy initial, which is the key combinatorial input in our computation of
the homotopy right Kan extension along $\gamma$.

Roughly speaking, fixing an object
\begin{equation}
  \alpha \maps [n] \longra \gamma(J)
\end{equation}

in the slice category $([n]/\Int)$, we will associate to each object of the
comma category $(\alpha / \cJ_n)$ an open region in a suitable configuration
space of points in $J^c \subset (0,1)$. The Lurie--Seifert--van Kampen theorem
then allows us to identify the nerve of $(\alpha / \cJ_n)$ with this
configuration space. Since the latter is built out of unordered configuration
spaces of open intervals in~$\RR$, the contractibility results above imply that
each $(\alpha / \cJ_n)$ has contractible nerve, hence that $\iota_n$ is
homotopy initial.

\begin{definition}
Let $\cJ_n$ denote the full subcategory of $([n]/\Int)$ on the identity arrows $\big([n] \xra{\id} \gamma(I)\big)\in \bDelta^{\op}$, for $I\in \Int$.    
\end{definition}

\begin{observation}
    The objects of $\cJ_n$ are parameterized by the intervals $I\in \Int$ with $n$ interior connected components. Moreover, given a fixed pair of objects 
    \[
      [n] \xra{\id} \gamma(I)
      \quad\text{and}\quad
      [n] \xra{\id} \gamma(I')
    \]
    in $\cJ_n$, the set of morphisms 
    \[
        \Hom_{\cJ_n}\big([n] \xra{\id} \gamma(I), [n] \xra{\id} \gamma(I')\big)
    \]
    is naturally identified with
    \[
        \{\,f\in\Hom_{\Int}(I,I') \mid \gamma(f) = \id_{[n]} \,\},
    \]
    i.e.\ with the inclusions $I \hookrightarrow I'$ which induce isomorphisms on $\pi_0$.
\end{observation}

The next two results are key propositions. As mentioned in the outline of the proof above, they ensure that the structure maps of the homotopy limit are weak equivalences — the proof presented here is a generalization of a similar kind in the work of Karlsson--Scheimbauer--Walde \cite{KSW25}. 

\begin{proposition}\label{ch2:: prop: hom initial}
    The inclusion functor 
    \begin{equation} 
        \iota_n\maps \cJ_n \xlonghookra{\,\,\,\,\,} [n]/\Int 
    \end{equation} 
    is homotopy initial. 
\end{proposition} 
\begin{proof}
     It is useful to start by outlining the strategy of the proof, that will greatly simplify exposition of its technical part. The first disclaimer is that we are going to show its opposite: we will prove the theorem using the $\cInt$ perspective. What we mean by that, is that we will instead prove the functor 
     \begin{equation} 
        \iota_{n}^{\op} \maps \cJ_{n}^{\op} \xlonghookra{\,\,\,\,} ([n]/\Int)^{\op} = \cInt/[n] 
    \end{equation} is homotopy final. For simplification, we will just call $\cJ_n^{\op} = \cJ_n$ for the remainder of the proof — assuring the reader that this is only local notation.

    Secondly, the way we will prove that $\left|\alpha / \cJ_n\right|$ is contractible is through the usage of Lurie--Seifert-van-Kampen theorem \Cref{ap:: thm: lurie Seifert}. Roughly, proofs of contractibility using this theorem work as follows: one shows that the category $\cD$ is weakly equivalent to a space $C$, which is known to be contractible, by constructing a functor 
    \begin{equation} 
        f\maps \cD \longra \cU(C) 
    \end{equation} 
    from $\cD$ into the poset of open subsets of $C$ and satisfying two additional properties: 
    \begin{enumerate} 
        \item For every $d\in \cD$ the open subset $f(d) \subset C$ is contractible. 
        \item For every point $x\in C$, the subcategory $\cD_{x} \hookra \cD$ containing the objects $d\in \cD$ for which $x \in f(d)$ is weakly contractible — in particular it is non-empty. 
    \end{enumerate} 
    The contractibility of $\cD$ then follows from the following sequence of equivalences: 
    \begin{align*} 
        \left| \cD \right| \,\,\, &\overset{(1)}{\simeq} \,\,\, \hoColimit[d \in \cD] * \\ 
        &\overset{(2)}{\simeq} \,\,\, \hoColimit[d \in \cD] f(d) \\
        &\overset{(3)}{\simeq} \,\,\, C \\ 
        &\overset{(4)}{\simeq} \,\,\, *. 
    \end{align*} 
    Equivalence (1) follows from ..., (2) results from the contractibility of $f(d)$ and the homotopy invariance of homotopy colimits, (3) follows from the Lurie--Seifert-van-Kampen theorem, if assuming contractibility of $\cD_x$ for every $x$, and finally (4) follows from the contractibility of $C$. \smallskip 
    
    We begin by constructing the space $C$, in the outline above. Fix $\alpha \maps \pi_{0}(J) \ra [n]$ a morphism in $\bDelta$, i.e., an object of $\cInt/[n]$. We define: 
    \begin{equation} 
        \mathsf{Conf}(J;\alpha)\,\, \overset{\sf def}{=} \,\, \left\{ 
            \begin{tikzcd}[ampersand replacement=\&] 
                {[n]^{\, \tl} } \& {(0,1)} \\ \& {J^{c}} \arrow[hook, from=1-1, to=1-2] \arrow["p"', curve={height=6pt}, dashed, from=1-1, to=2-2] \arrow[hook, from=2-2, to=1-2] \end{tikzcd} \text{ s.t. } p \text{ is injective and } \begin{tikzcd}[ampersand replacement=\&] {\pi_{0}([n]^{\, \tl})} \& {\pi_0(J^{c})} \\ { [n]^{\, \tl}} \& {\gamma(J^{c})^{\, \tl} } 
                    \arrow["{\pi_{0}(p)}", from=1-1, to=1-2] 
                    \arrow[equals, from=1-1, to=2-1] 
                    \arrow[equals, from=1-2, to=2-2] 
                    \arrow["{(\alpha^{\op})^{\, \tl}}"', from=2-1, to=2-2] 
                \end{tikzcd} \right\} \end{equation} This is the ordered space of $(n+2) = \left| [n]^{\, \tl} \right|$ points in $(0,1)$ satisfying two additional conditions: 
                
                \begin{enumerate} 
                    \item The injection factors through the disjoint union of open intervals $J^c$. 
                    \item Upon passage to $\pi_0$, agrees with $(\alpha^\op)^{\, \tl}$ in $\sO_{\pm}$. 
                \end{enumerate} 
                
                \nid We let $\mathsf{UConf}(J;\alpha)$ be its space of $\Sigma_{n+2}$-coinvariants, i.e., the unordered analogue. This is an open and convex subspace of $\mathsf{UConf}_{n+2}((0,1))$; therefore contractible. Before continuing with the proof, it is useful to indicate an example of a point in this space. Let $J$ be the following object of $\cInt$ 
                \smallskip 
                \begin{equation} 
                    \tikzfig{figures/Ch2/INT16}{\input{figures/Ch2/INT16.tikz}}, 
                \end{equation} 
                \smallskip 
                
                \nid and let $\alpha \maps \pi_{0}(J) \simeq [1] \xra{\lag 0,2 \rag} [2]$. A point $p \maps [2]^{\, \tl} \ra J^{c}$ in $\mathsf{UConf}(J;\alpha)$ is thus an unordered configuration of $4$ points in $J^{c} = \sqcup_{i} U_i$ — which is the disjoint union of three ordered open intervals — satisfying an additional condition: 
                \begin{enumerate} 
                    \item[$(*)$] The number of points in each connected open interval is determined by the cardinality of $\pi_{0}(p)^{-1}([U_i])$. 
                \end{enumerate} 
                Which in this case is 1 for $U_1$, 2 for $U_2$, and $1$ for $U_3$. Therefore the configuration on the left is a point in $\mathsf{UConf}(J;\alpha)$ and the one on the right is not 
                \begin{equation} 
                    \tikzfig{figures/Ch2/INT18}{\input{figures/Ch2/INT18.tikz}}, \quad \quad \tikzfig{figures/Ch2/INT19}{\input{figures/Ch2/INT19.tikz}} \, . 
                \end{equation} 
                \smallskip 
                
                \nid We indicate the combinatorics at play with the following diagram 
                \begin{equation} 
                    \tikzfig{figures/Ch2/INT17}{\input{figures/Ch2/INT17.tikz}} 
                \end{equation} 
                \smallskip 
                
                \nid This example also clarifies how $\mathsf{UConf}(J;\alpha)$ is homeomorphic to a product of unordered configuration spaces of open intervals with fixed numbers of particles \begin{equation} \mathsf{UConf}_{k_1}(U_1) \times \mathsf{UConf}_{k_2}(I_2) \times \mathsf{UConf}_{k_3}(U_3), \end{equation} where $k_i = \pi_{0}(p)^{-1}([U_i])$, for $i \in \pi_{0}(J^{c}) = [- < 0 < 1]$. 
                \medskip 
                
                \nid Next, we want to define the subspaces $f(d)$ in the outline. In particular, for every object 

                \[\big(\pi_0(J) \ra \pi_0(I) \xra{\id} [n]\big) = \alpha,\]
                in $\alpha /\cJ_{n}$, we seek to construct an open subspace $\mathsf{UConf}(J;\alpha; I) \subset \mathsf{UConf}(J;\alpha)$. We do so as follows: 
                
                \begin{equation} 
                    \mathsf{UConf}(J;\alpha;I)\,\, \overset{\sf def}{=} \,\, \left\{ 
                        \begin{tikzcd}[ampersand replacement=\&]
                         {[n]^{\, \tl} } \& {(0,1)} \\ {I^{c}} \& {J^{c}} \arrow["\phi", hook, from=1-1, to=1-2] \arrow["p"', dashed, from=1-1, to=2-1] \arrow[hook, from=2-1, to=2-2] \arrow[hook, from=2-2, to=1-2] 
                        \end{tikzcd} 
                    \text{ s.t. } p\text{ is injective and } 
                    \begin{tikzcd}[ampersand replacement=\&] {\pi_{0}([n]^{\, \tl})} \& {\pi_0(I^{c})} \& {\pi_0(J^{c})} \\ { [n]^{\, \tl}} \& {\gamma(I^{c})^{\, \tl} } \& {\gamma(J^{c})^{\, \tl} } 
                    \arrow["{{\pi_{0}(p)}}", from=1-1, to=1-2] 
                    \arrow[equals, from=1-1, to=2-1] 
                    \arrow[from=1-2, to=1-3] 
                    \arrow[equals, from=1-2, to=2-2] 
                    \arrow[equals, from=1-3, to=2-3] 
                    \arrow[from=2-1, to=2-2] 
                    \arrow["{\alpha^\op}"', curve={height=18pt}, from=2-1, to=2-3] 
                    \arrow[from=2-2, to=2-3] 
                \end{tikzcd} \right\} 
                \end{equation} 
                \smallskip 
                
                \nid Note that for $(p\maps [2]^{\, \tl} \ra I^{c})$ to be an object in $\mathsf{UConf}(J;\alpha;I)$ we need only to ensure that $\pi_{0}(p)=\id$. In particular, this means that we can interpret points in this open subspace as configurations of $n+2$ points in $I^{c}$ which have a unique particle per connected component. We again depict an example of a point in this space (implicitly making a choice of $I\in \cInt$): 
                \begin{equation} 
                    \tikzfig{figures/Ch2/INT20}{\input{figures/Ch2/INT20.tikz}} 
                \end{equation} 
                \smallskip 
                
                \nid This subspace is, again, convex for similar reasons and therefore contractible. All that is left showing to finish the proof is that, given an arbitrary point in $\mathsf{UConf}(J;\alpha)$: 
                \[([n]^{\, \tl} \xra{p} J^{c} \hookra (0,1)) \in \mathsf{UConf}(J;\alpha).\] 
                the subposet 
                \begin{equation}
                     \left\{(I\subseteq J) \in \alpha/\cJ_n \, : (p \maps [n]^{\, \tl} \ra J^{c}) \in \mathsf{UConf}(J;\alpha; I) \right\} \,\, \subseteq \,\, \cJ_{n} 
                \end{equation}
                is contractible. To show that is non-empty we offer a simple procedure which produces an object $I\in \cInt$ such that $p\in \mathsf{UConf}(J;\alpha;I)$. For each connected open interval \begin{equation} 
                    J^{c} = U_1 \sqcup \cdots \sqcup U_{k} 
                \end{equation} 
                you do one of two things: if it contains no points, an impossibility for the extremal components, you close it by taking the hull of the two adjacent closed intervals; if it contains a single point, you do nothing; if it contains $m$ points $(p_i < \cdots < p_{i+m})$ with $m\geq 1$, you add $m-1$ closed connected intervals in $(p_{j},p_{j+1})$ for $j=i,\ldots,m+i-1$. We essentially used the quantity of particles on a given open edge of $J^c$, as a way of determining how many connected open intervals an object $I$ should have in that region, in order for it to receive an inclusion from $J$ in $\cInt$. 
                
                Finally, contractibility of this subposet comes from the fact that it is isomorphic to a product of posets, each given by the poset of connected closed intervals in $(p_i, p_{i+1})$ for $i=0,\ldots, n+2$ — which is filtered and therefore contractible. \end{proof}

\begin{proposition}\label{ch2:: prop: weak contractibility of pi_0}
    The nerve of the full subcategory $\iota_n \colon \cJ_n \hookrightarrow ([n]/\Int)$ is weakly contractible for every $n\geq 0$.
\end{proposition}

\begin{proof}
   This is the special case of \Cref{ch2:: prop: hom initial} obtained by taking $\alpha = \id$.
\end{proof}

\begin{proposition}\label{ch2:: prop: structure_maps_weak_equivalences}
    Let $\cF\colon \Int \to \sset$ be a $\cW_{\pi_{0}}$-local object. Then for every $n\geq 0$ and every object $I\in \cJ_{n}$, the structure map
    \begin{equation}
        \underset{J\in ([n]/\Int)}{\mathsf{holim}} \,\, \cF(J) \xlongrightarrow{\ \sim\ } \cF(I)
    \end{equation}
    is a weak equivalence.
\end{proposition}

\begin{proof}
    Since $\cF$ is $\cW_{\pi_0}$-local, in particular it is injectively fibrant and sends all morphisms in $\cJ_n$ to weak equivalences, so the restriction $\cF|_{\cJ_n}$ is a homotopically constant diagram of Kan complexes. By \Cref{ch2:: prop: hom initial} the inclusion $\iota_n\colon \cJ_n \hookrightarrow ([n]/\Int)$ is homotopy initial, hence the canonical map
    \[
      \holim_{([n]/\Int)} \cF \;\xrightarrow{\ \sim\;}\; \holim_{\cJ_n} \cF
    \]
    is a weak equivalence. On the other hand, by \Cref{ch2:: prop: weak contractibility of pi_0} the nerve of $\cJ_n$ is weakly contractible, so we can apply the general result on homotopy limits of homotopically constant diagrams over contractible indexing categories (see \Cref{ch2:: prop: homotopy-limit-homotopically-constant}). This yields, for every object $I\in \cJ_n$, a weak equivalence
    \[
      \holim_{\cJ_n} \cF \;\xrightarrow{\ \sim\;}\; \cF(I).
    \]
    Composing these two weak equivalences gives the desired statement.
\end{proof}

\begin{theorem}
    The Quillen pair 
    \begin{equation}
    \begin{tikzcd}[ampersand replacement=\&]
	{\fun(\bDelta^{\op},\sset)_{\mathsf{inj}}} \&\& {\fun(\Int,\sset)_{\mathsf{inj},h}}
	\arrow[""{name=0, anchor=center, inner sep=0}, "{\gamma^*}", bend left=10, from=1-1, to=1-3]
	\arrow[""{name=1, anchor=center, inner sep=0}, "{\rkan_{\gamma}}", bend left=10, from=1-3, to=1-1]
	\arrow["\dashv"{anchor=center, rotate=-90}, draw=none, from=1, to=0]
    \end{tikzcd}
    \end{equation}
    is a Quillen equivalence.
\end{theorem}

\begin{proof}
    We use the standard criterion for Quillen equivalences which says that it is enough to show:
    \begin{enumerate}
        \item \label{ch2:: item: Quillen_equivalence_conditions_1} The right adjoint $\rkan_\gamma$ reflects weak equivalences between fibrant objects.
        \item \label{ch2:: item: Quillen_equivalence_conditions_2} For every cofibrant object $X$ in $\fun(\bDelta^{\op},\sset)_{\mathsf{inj}}$, the derived unit 
        \begin{equation}
            X \longra \rkan_{\gamma}(\bR \gamma^{*}X)
        \end{equation}
        is a weak equivalence.
    \end{enumerate} 

    To check \eqref{ch2:: item: Quillen_equivalence_conditions_1}, let $\eta\colon \cF \to \cG$ be a natural transformation between fibrant objects in $\fun(\Int,\sset)_{\mathsf{inj}, h}$ and assume that
    \begin{equation}
        \rkan_\gamma(\cF) \xra{\rkan_\gamma (\eta)} \rkan_\gamma(\cG)
    \end{equation}
    is a weak equivalence in $\fun(\bDelta^{\op},\sset)_{\mathsf{inj}}$. Since weak equivalences in $\fun(\bDelta^{\op}, \sset)_{\mathsf{inj}}$ are level-wise, for every $[n]\in \bDelta^{\op}$ the map
    \begin{equation}
        \rkan_\gamma(\cF)([n]) \longra \rkan_\gamma(\cG)([n])
    \end{equation}
    is a weak equivalence in $\sset$.

    Fix $I\in \Int$ and write $[n]=\gamma(I)$. By \Cref{ch2:: prop: structure_maps_weak_equivalences}, the canonical maps
    \[
      \rkan_\gamma(\cF)([n]) \longra \cF(I),
      \qquad
      \rkan_\gamma(\cG)([n]) \longra \cG(I)
    \]
    are weak equivalences. Since the square
    \[
        \begin{tikzcd}[ampersand replacement=\&]
            {\rkan_\gamma(\cF)([n])} \&\& {\rkan_\gamma(\cG)([n])} \\
            \\
            {\cF(I)} \&\& {\cG(I)}
            \arrow[from=1-1, to=1-3]
            \arrow["\sim"', from=1-1, to=3-1]
            \arrow["\sim", from=1-3, to=3-3]
            \arrow[from=3-1, to=3-3]
        \end{tikzcd}
    \]
    commutes, the map $\eta_I\colon \cF(I)\to\cG(I)$ is a weak equivalence by the $2$-out-of-$3$ property. This holds for every $I\in\Int$, so $\eta$ is a level-wise weak equivalence, hence a weak equivalence in $\fun(\Int,\sset)_{\mathsf{inj},h}$.

    To verify \eqref{ch2:: item: Quillen_equivalence_conditions_2}, let $X$ be cofibrant in $\fun(\bDelta^{\op},\sset)_{\mathsf{inj}}$ and choose an injective fibrant replacement
    \[
      \gamma^*X \longra \cF
    \]
    in $\fun(\Int,\sset)_{\mathsf{inj}}$. Since $\gamma^*X$ sends every morphism in $\cW_{\pi_0}$ to an isomorphism, the fibrant object $\cF$ is automatically $\cW_{\pi_0}$-local, hence fibrant in $\fun(\Int,\sset)_{\mathsf{inj},h}$. The derived unit is represented by the adjoint map
    \[
      u_X\colon X \longra \rkan_\gamma(\cF).
    \]
    Fix $[n]\in \bDelta^{\op}$ and choose $I\in\Int$ with $\gamma(I)=[n]$. Then the composite
    \[
      X([n]) \xra{\,u_X([n])\,} \rkan_\gamma(\cF)([n]) \longra \cF(I)
    \]
    is exactly the component at $I$ of the fibrant replacement map $\gamma^*X\to\cF$, hence is a weak equivalence in $\sset$. By \Cref{ch2:: prop: structure_maps_weak_equivalences}, the second map is also a weak equivalence. Therefore $u_X([n])$ is a weak equivalence by the $2$-out-of-$3$ property. Since this holds for every $[n]$, the derived unit is a level-wise weak equivalence, as required.
\end{proof}

\subsection{Imposing the Segal condition}
\smallskip

The next step in constructing a model structure on $\Fun(\Int,\sset)$ which
models the homotopy theory of $(\oo,1)$-categories is to encode the Segal
condition. Concretely, we want to further localize the model structure
$\Fun(\Int,\sset)_{\mathsf{inj},h}$ so that the fibrant objects not only send
$\pi_0$-isomorphisms to weak equivalences, but also satisfy the interval Segal
condition introduced in \Cref{ch2:: subsec: interval_spaces}. The heavy lifting
here is provided by \Cref{ap:: thm: transfer of localizations}.

First, let us recall what we want from an \emph{interval Segal condition}. If an
interval space $\cF\colon \Int \to \sset$ is to model an $(\oo,1)$-category,
then its space of morphisms should be encoded by $\cF\big(\sU(A,B)\big)$, with
the evident maps to $\cF\big(\sU(B)\big)$ and $\cF\big(\sU(A)\big)$ playing the
role of target and source, respectively. Under this interpretation, we expect
the Segal-type condition to express that $\cF\big(\sU(A,B,C)\big)$ captures the
data of composable pairs of 1-morphisms. More explicitly, the Segal-type
condition requires that the canonical map
\begin{equation}
    \cF\big(\sU(A,B,C)\big)
      \longra
      \cF\big(\sU(A,B)\big)
        \underset{\cF\big(\sU(B)\big)}{\overset{h}{\times}}
      \cF\big(\sU(B,C)\big)
\end{equation}
is a levelwise weak equivalence of simplicial sets.

\begin{definition}
    For each $n\geq 2$, let
    \begin{equation}
      \varphi_n\colon G(n) \hookrightarrow \Delta[n]
    \end{equation}
    denote the $n$-dimensional spine inclusion of simplicial spaces, and let
    \begin{equation}
      \mathbb{L}\gamma^{*}(\varphi_n)
    \end{equation}
    be its image under the left derived functor of the left Quillen functor
    $\gamma^{*}$. We write
    \begin{equation}
      \varphi := \{\varphi_n\}_{n\geq 2}
    \end{equation}
    for the resulting family of maps, and similarly
    $\mathbb{L}\gamma^{*}(\varphi)$ for their images.
\end{definition}

\begin{observation}
    Since we are working with the injective model structure on the category of
    simplicial spaces, cofibrations are precisely monomorphisms. In particular,
    each spine inclusion $\varphi_n\colon G(n)\hookrightarrow\Delta[n]$ is a
    cofibration, so its left derived image satisfies
    \begin{equation}
        \mathbb{L}\gamma^*(\varphi_n)
        \simeq
        \gamma^*\big(G(n)\big) \longra \gamma^*\big(\Delta[n]\big).
    \end{equation}
    In other words, the cofibrant replacement is trivial. The same remark
    applies to finite coproducts of spine inclusions, since these are again
    monomorphisms.
\end{observation}

The next result is the key observation of this subsection: it shows that we can
recover the interval Segal condition for interval spaces by localizing
simultaneously on the $\Int$-side and on the simplicial side.

\begin{lemma}\label{ch2:: lem: fibrant Segal condition}
Let $\cF\in \fun(\Int,\sset)_{\mathsf{inj},h}$. Then $\cF$ is
$\mathbb{L}\gamma^{*}(\varphi)$-local if and only if $\cF$ is fibrant and
satisfies the interval Segal condition.
\end{lemma}

\begin{proof}
    Suppose first that $\cF\colon \Int \to \sset$ is
    $\mathbb{L}\gamma^{*}(\varphi)$-local. By definition of locality, $\cF$ is
    fibrant in $\Fun(\Int,\sset)_{\mathsf{inj},h}$, hence injectively fibrant
    and homotopically constant along $\cW_{\pi_0}$. It remains to check that it
    satisfies the interval Segal condition, i.e.\ that for every
    $\sU(A_1,\ldots,A_n)$ with $n\geq 2$, the canonical map
    \begin{equation}\label{eq segal}
        \cF\big(\sU(A_1,\ldots, A_n)\big)
          \longra
        \cF\big(\sU(A_1,A_2)\big)
          \times^{h}_{\cF\big(\sU(A_2)\big)}
        \cdots
          \times_{\cF\big(\sU(A_{n-1})\big)}^{h}
        \cF\big(\sU(A_{n-1},A_n)\big)
    \end{equation}
    is a weak equivalence.

    By $\mathbb{L}\gamma^*(\varphi)$-locality, for every $n\geq 2$ the map
    \begin{equation}
      \map^{h}\big(\gamma^* \Delta[n], \cF\big)
        \xrightarrow{\ \simeq\ }
      \map^{h}\big(\gamma^* G(n), \cF\big)
    \end{equation}
    is a weak equivalence of simplicial sets. Since every object is cofibrant
    and $\cF$ is fibrant, the derived mapping spaces agree with the ordinary
    ones, and we can identify
    \begin{equation}
      \map^{h}(-,-)\;\simeq\;
      \Hom_{\fun(\Int,\sset)}(-,-).
    \end{equation}
    Using the adjunction $\gamma^\ast \dashv \rkan_{\gamma}$, this yields a
    commutative diagram
    \begin{equation}
        \begin{tikzcd}[ampersand replacement=\&]
            {\Hom_{\fun(\Int,\sset)}\big(\gamma^* \Delta[n], \cF\big)} \&\&
            {\Hom_{\fun(\Int,\sset)}\big(\gamma^* G(n), \cF\big)} \\
            \\
            {\Hom_{\fun(\bDelta^{\op},\sset)}\big(\Delta[n], \rkan_{\gamma}(\cF)\big)} \&\&
            {\Hom_{\fun(\bDelta^{\op},\sset)}\big(G(n), \rkan_{\gamma}(\cF)\big)} \\
            \\
            {\rkan_{\gamma}(\cF)([n])} \&\&
            {\rkan_{\gamma}(\cF)([1])\underset{\rkan_{\gamma}(\cF)([0])}{\times}\cdots
              \underset{\rkan_{\gamma}(\cF)([0])}{\times}\rkan_{\gamma}(\cF)([1])}
            \arrow["\sim", from=1-1, to=1-3]
            \arrow["\cong"', from=1-1, to=3-1]
            \arrow["\cong", from=1-3, to=3-3]
            \arrow["\sim",  from=3-1, to=3-3]
            \arrow["\cong"', from=3-1, to=5-1]
            \arrow["\cong", from=3-3, to=5-3]
            \arrow["\sim",from=5-1, to=5-3]
        \end{tikzcd}
    \end{equation}
    where:
    \begin{itemize}
      \item the top horizontal arrow is a weak equivalence by
            $\mathbb{L}\gamma^*(\varphi)$-locality;
      \item the first pair of vertical arrows are the natural isomorphisms of
            the adjunction;
      \item the bottom vertical arrows are the identifications
            $\Hom(\bDelta[k],X)\cong X_k$ coming from the (enriched) Yoneda
            lemma;
      \item the middle horizontal arrow is therefore the Segal map for the
            simplicial space $\rkan_\gamma(\cF)$ at $[n]$, and is a weak
            equivalence by the $2$-out-of-$3$ property.
    \end{itemize}

    Recall that the $n$-dimensional spine $G(n)$ is obtained as the colimit of
    the diagram
    \begin{equation}
        \begin{tikzcd}[ampersand replacement=\&,column sep=small]
            {\Delta[\{0<1\}]} \&\& {\Delta[\{1<2\}]} \&\& \cdots \&\&
            {\Delta[\{n-1<n\}]} \\
            \\
            \& {\Delta[\{1\}]} \&\& {\Delta[\{2\}]} \&\& {\,\Delta[\{n-1\}]}
            \arrow[from=1-1, to=3-2]
            \arrow[from=1-3, to=3-2]
            \arrow[from=1-3, to=3-4]
            \arrow[from=1-5, to=3-4]
            \arrow[from=1-5, to=3-6]
            \arrow[from=1-7, to=3-6]
        \end{tikzcd}
    \end{equation}
    in $\Fun(\bDelta^{\op},\sset)$, so that
    $\Hom_{\fun(\bDelta^{\op},\sset)}\big(G(n), \rkan_{\gamma}(\cF)\big)$ can be
    identified with the iterated fibre product in the bottom right corner of
    the previous diagram.

    Now fix an object $\sU(A_0,\ldots,A_n)\in \Int$ with $n$ interior connected
    components. By \Cref{ch2:: prop: structure_maps_weak_equivalences}, the
    structure maps
    \begin{equation}
      \rkan_\gamma(\cF)([k])
        \longrightarrow
      \cF\big(\text{corresponding interval with $k$ ‘‘slots’’}\big)
    \end{equation}
    are weak equivalences for $k=0,1,n$ and compatible with the simplicial
    structure maps. This yields a commutative diagram of spans
    \begin{equation}
    \begin{tikzcd}[ampersand replacement=\&]
        {\rkan_{\gamma}(\cF)([1])} \& {\rkan_{\gamma}(\cF)([0])} \&
        {\rkan_{\gamma}(\cF)([1])} \& \cdots \&
        {\rkan_{\gamma}(\cF)([0])} \& {\rkan_{\gamma}(\cF)([1])} \\
        \\
        {\cF\big(\sU(A_0,A_1)\big)} \& {\cF\big(\sU(A_1)\big)} \&
        {\cF\big(\sU(A_1,A_2)\big)} \& \cdots \&
        {\cF\big(\sU(A_{n-1})\big)} \&
        {\cF\big(\sU(A_{n-1},A_{n})\big)}
        \arrow[from=1-1, to=1-2]
        \arrow[from=1-1, to=3-1]
        \arrow[from=1-2, to=3-2]
        \arrow[from=1-3, to=1-2]
        \arrow[from=1-3, to=3-3]
        \arrow[from=1-4, to=1-3]
        \arrow[from=1-4, to=1-5]
        \arrow[from=1-5, to=3-5]
        \arrow[from=1-6, to=1-5]
        \arrow[from=1-6, to=3-6]
        \arrow[from=3-1, to=3-2]
        \arrow[from=3-3, to=3-2]
        \arrow[from=3-4, to=3-3]
        \arrow[from=3-4, to=3-5]
        \arrow[from=3-6, to=3-5]
    \end{tikzcd}
    \end{equation}
    in which all vertical arrows are weak equivalences. Passing to the
    corresponding homotopy pullbacks, we obtain a weak equivalence between the
    homotopy pullback computing the Segal map for $\rkan_\gamma(\cF)$ and the
    homotopy pullback on the right-hand side of \eqref{eq segal}. Since the
    former is a weak equivalence by the previous argument, the latter is as
    well. Thus $\cF$ satisfies the interval Segal condition.

    Conversely, if $\cF$ is fibrant and satisfies the interval Segal condition,
    the same reasoning can be run backwards to show that the Segal maps of
    $\rkan_\gamma(\cF)$ are weak equivalences, which is equivalent to
    $\mathbb{L}\gamma^*(\varphi)$-locality.
\end{proof}

\begin{proposition}\label{ch2:: prop: model structure interval Segal space}
    There exists a model structure on the category of pre-interval spaces
    $\fun(\Int,\sset)$, which we will denote by $\mathsf{ISSp}$, with the
    following properties:
    \begin{enumerate}
        \item The cofibrations are the injective cofibrations.
        \item The fibrant objects are interval Segal spaces, i.e.\ functors
              $\cF\colon \Int \to \sset$ satisfying the interval Segal condition \Cref{ch2:: def: interval Segal condition}.
        \item The weak equivalences are the $\mathbb{L}\gamma^{*}(\varphi)$-local
              equivalences. In particular, a weak equivalence between interval
              Segal spaces is a levelwise weak equivalence.
        \item The model structure is left proper, combinatorial, and simplicial.
    \end{enumerate}
    We call this the model category of \textbf{interval Segal spaces}. 
\end{proposition}

\begin{proof}
    Existence of the left Bousfield localization of
    $\Fun(\Int,\sset)_{\mathsf{inj},h}$ at $\mathbb{L}\gamma^{*}(\varphi)$ gives
    (1) and (4), since cofibrations are unchanged and left properness,
    combinatoriality, and simpliciality are preserved. Assertion (2) follows
    from \Cref{ch2:: lem: fibrant Segal condition}, which identifies the
    fibrant objects with interval Segal spaces. Finally, (3) is the general
    description of weak equivalences between fibrant objects in a left
    Bousfield localization: between interval Segal spaces,
    $\mathbb{L}\gamma^{*}(\varphi)$-local equivalences coincide with levelwise
    weak equivalences.
\end{proof}

\begin{theorem}\label{ch2:: thm: Quillen-equivalence-SSp-ISSp}
    The Quillen equivalence previously established between simplicial spaces and
    interval spaces
    \begin{equation}
        \begin{tikzcd}[ampersand replacement=\&]
            {\fun(\bDelta^{\op},\sset)_{\mathsf{inj}}} \&\&
            {\fun(\Int,\sset)_{\mathsf{inj},h}}
            \arrow[""{name=0, anchor=center, inner sep=0}, "{\gamma^*}",
                   curve={height=-16pt}, from=1-1, to=1-3]
            \arrow[""{name=1, anchor=center, inner sep=0}, "{\rkan_{\gamma}}",
                   curve={height=-16pt}, from=1-3, to=1-1]
            \arrow["\dashv"{anchor=center, rotate=90}, draw=none, from=1, to=0]
        \end{tikzcd}
    \end{equation}
    induces a Quillen equivalence
    \begin{equation}
        \begin{tikzcd}[ampersand replacement=\&]
            {\mathsf{SSp}} \&\& {\mathsf{ISSp}}
            \arrow[""{name=0, anchor=center, inner sep=0}, "{\gamma^*}",
                   curve={height=-18pt}, from=1-1, to=1-3]
            \arrow[""{name=1, anchor=center, inner sep=0}, "{\rkan_{\gamma}}",
                   curve={height=-18pt}, from=1-3, to=1-1]
            \arrow["\dashv"{anchor=center, rotate=90}, draw=none, from=1, to=0]
        \end{tikzcd}
    \end{equation}
    between the model category of Segal spaces and the model category of
    interval Segal spaces, upon further localizing along $\varphi$ and
    $\mathbb{L}\gamma^*(\varphi)$, respectively. 
\end{theorem}

\begin{proof}
    This follows directly from \Cref{ch2:: prop: Segal-model-structure} and
    \Cref{ch2:: prop: model structure interval Segal space}, together with
    \Cref{ap:: thm: transfer of localizations}, applied to the Quillen
    equivalence $\gamma^* \dashv \rkan_{\gamma}$ between
    $\Fun(\Delta^{\op},\sset)_{\mathsf{inj}}$ and $\Fun(\Int,\sset)_{\mathsf{inj},h}$.
\end{proof}

\begin{observation}\label{ch2:: obs: complete Segal localization}
If we left Bousfield localize this Quillen equivalence along
$\{J \longra \bDelta[0]\}$ in $\mathsf{SSp}$ and along
$\LL \gamma^* (J \longra \bDelta[0])$ in $\mathsf{ISSp}$, we get a new pair of
model structures whose fibrant objects are complete Segal spaces in
$\mathsf{SSp}$ and $\LL \gamma^* (J \longra \bDelta[0])$-local objects in
$\mathsf{ISSp}$. By the derived adjunction, we have that if $\cF$ is
$\LL \gamma^* (J \longra \bDelta[0])$-local then the induced map
\begin{equation}
    \Map^h(\LL\gamma^* \bDelta[0], \cF) \xra{\,\,\,\simeq \,\,\,} \Map^{h}(\LL\gamma^* J, \cF)
\end{equation} 
is a weak equivalence, which translates to $\Rrkan_\gamma(\cF)$ being a
complete Segal space. We call the resulting model structure $\mathsf{CISSp}$.
From \Cref{ap:: thm: transfer of localizations} we have the following Quillen
equivalence 
\begin{equation}
    \begin{tikzcd}[ampersand replacement=\&]
        {\mathsf{CSSp}} \&\& {\mathsf{CISSp}}
        \arrow[""{name=0, anchor=center, inner sep=0}, "{\gamma^*}", curve={height=-12pt}, from=1-1, to=1-3]
        \arrow[""{name=1, anchor=center, inner sep=0}, "{\rkan_{\gamma}}", curve={height=-12pt}, from=1-3, to=1-1]
        \arrow["\dashv"{anchor=center, rotate=90}, draw=none, from=1, to=0]
    \end{tikzcd}
\end{equation}
\end{observation}

\begin{theorem}[Complete interval Segal spaces]\label{ch2:: thm: complete interval model}
    There exists a model structure on the category of pre-interval spaces
    $\fun(\Int,\sset)$, which we will denote by $\mathsf{CISSp}$, with the
    following properties:
    \begin{enumerate}
        \item The cofibrations are the injective cofibrations.
        \item The fibrant objects are complete interval Segal spaces.
        \item The weak equivalences between complete interval
              Segal spaces are the levelwise weak equivalences.
        \item The model structure is left proper, combinatorial, and simplicial.
    \end{enumerate}
    We call this the model category of \textbf{complete interval Segal spaces}.
    Moreover, restriction and right Kan extension along $\gamma$ induce a
    Quillen equivalence
    \begin{equation}
        \mathsf{CSSp}
        \ \underset{\rkan_{\gamma}}{\overset{\gamma^{*}}{\rightleftarrows}}\
        \mathsf{CISSp}.
    \end{equation}
    In particular, $\mathsf{CISSp}$ presents the $(\oo,1)$-category of
    $(\oo,1)$-categories.
\end{theorem}
\begin{proof}
    The theorem follows from \Cref{ch2:: obs: complete Segal localization}.
\end{proof}

\begin{remark}
    One of the most important results in Rezk's original paper \cite{Rezk01} is the fact that $\mathsf{CSSp}$ is compatible with its Cartesian closed structure, which ensures that for a pair of complete Segal spaces $X$ and $Y$, the internal-hom $X^Y$ is again a complete Segal space. In order to ensure that $\mathsf{CISSp}$ satisfies a similar property, a little more work would be needed, since this is not a property transferred through a Quillen equivalence.
\end{remark}

\begin{notation}
    We denote the full subcategory of $\mathsf{ISSp}$ spanned by the \emph{complete} Segal interval spaces as $\cC\mathsf{ISSp}$.
\end{notation}

\section{Homotopy theory of interval Segal spaces}
\label{ch2:: sec: homotopy theory int}

The theory of $(\oo,1)$-categories' primary goal is to smoothly accommodate category theory with homotopy theory, so before we end this chapter we should articulate the sense in which interval Segal spaces achieve this. To do this, we construct a fundamental invariant of every $(\oo,1)$-category, its \emph{homotopy category}. 

We begin with the canonical definition, obtained by passing from an interval Segal space to its associated Segal space via the homotopy right Kan extension. We then unpack this definition in terms of explicit homotopy limits over $\Int$.

\begin{observation}\label{ch2:: obs: RRkan-gamma-pointwise}
Let \(\cF \colon \Int \longra \sset\) be fibrant in
\(\Fun(\Int,\sset)_{\inj,h}\).
By the construction of the homotopy right Kan extension, for every
\([n]\in \bDelta\) we have a natural weak equivalence
\begin{equation}
  \Rrkan_{\gamma} \cF([n])
  \;\simeq\;
  \hoLimit[([n]\xra{\alpha}\gamma(I))\in [n]/\Int] \cF(I).
\end{equation}
Moreover, by \Cref{ch2:: prop: hom initial} the inclusion
\(\cJ_n \hookrightarrow [n]/\Int\) is homotopy initial, so there is a
natural weak equivalence
\begin{equation}
  \Rrkan_{\gamma} \cF([n])
  \;\simeq\;
  \hoLimit[([n]\xra{\id}\gamma(I))\in \cJ_n] \cF(I).
\end{equation}
For \(n=0,1\) this expresses the \(0\)- and \(1\)-simplices of the
associated Segal space in terms of explicit homotopy limits over \(\Int\).
\end{observation}

\begin{definition}\label{ch2:: def: homotopy category ISSP}
  Let \(\cF\colon \Int \longra \sset\) be an interval Segal space, i.e.\ a
  fibrant object of \(\mathsf{ISSp}\).
  We define its associated \define{homotopy category}
  \begin{equation}
    \ho(\cF) \;\overset{\mathsf{def}}{=}\;
    \ho\big(\Rrkan_{\gamma}\cF\big)
  \end{equation}
  to be the homotopy category of the induced Segal space
  \(\Rrkan_{\gamma}\cF\colon \bDelta^{\op}\longra \sset\), as in
  \Cref{ch2:: def: homotopy category segal space}.
\end{definition}

\begin{remark}\label{ch2:: rem: objects-morphisms-as-holims}
Let \(\cF\colon \Int \longra \sset\) be an interval Segal space.
By Observation~\ref{ch2:: obs: RRkan-gamma-pointwise}, the simplicial set
of objects of the associated Segal space is
\begin{equation}
  \big(\Rrkan_{\gamma}\cF\big)_0
  \;\simeq\;
  \hoLimit[I \in \cJ_0] \cF(I),
\end{equation}
and the simplicial set of \(1\)-simplices is
\begin{equation}
  \big(\Rrkan_{\gamma}\cF\big)_1
  \;\simeq\;
  \hoLimit[J \in \cJ_1] \cF(J).
\end{equation}
The source and target maps
\begin{equation}
  s,t \colon
  \big(\Rrkan_{\gamma}\cF\big)_1 \longra
  \big(\Rrkan_{\gamma}\cF\big)_0
\end{equation}
are induced by precomposition with the maps
\(\langle 0\rangle,\langle 1\rangle \colon [0]\ra [1]\) in~\(\bDelta\), and they can be understood as taking a morphism over $\sU(A,B)$ and then restricting to $\sU(A)$ and $\sU(B)$, respectively — where both maps are given by gluing two of the opens or by covering one of the closed intervals.

\nid Given \(x,y \in \Ob\big(\ho(\cF)\big)\), the mapping space
\(\Map_{\Rrkan_{\gamma}\cF}(x,y)\) is the homotopy fibre of
\begin{equation}
  \Map^{h}_{\Rrkan_{\gamma}\cF}(x,y) \,\, \overset{\sf def}{=}\,\, \mathsf{hofiber}_{(x,y)}\left(
  \big(\Rrkan_{\gamma}\cF\big)_1 \longra
  \big(\Rrkan_{\gamma}\cF\big)_0 \times
  \big(\Rrkan_{\gamma}\cF\big)_0\right)
\end{equation}
over the point \((x,y)\), exactly as in \Cref{ch2:: cons: objects mapping space}.
\end{remark}

\section{Symmetric monoidal structures}
\label{ch2:: sec: symmetric monoidal structure}

Symmetric monoidal structures are key to the study of field theories and it is thus imperative that we discuss how they will be incorporated into our constructions. There is an entire taxonomy of ways to encode symmetric monoidal structures into homotopical or higher categorical constructions — a profound study of this topic can be found in Lurie's book \cite{LurieHA}. The two main ones, dual via straightening-unstraightening, are as algebras for the $\cE_{\oo}$-operad and as a Segal $\Gamma$-object — we will, given the future applications we have in mind, work only with the latter. 

\subsection{As \texorpdfstring{$\Gamma$}{Gamma}-objects}

Lurie defined what it means to be a commutative monoid in an arbitrary $(\oo,1)$-category, based on ideas dating back to Segal. We expand on that definition in this section and we construct a model structure on a certain category of diagrams in $\cC\mathsf{ISSp}$ whose fibrant objects are the symmetric monoidal ones. Our sources for this section are \cite{CS19}, \cite[Appendix~A]{JFS15}, and of course \cite{LurieHA}.

\begin{definition}
    Let $\Gamma$ be a category with the following data and structure:
    \begin{enumerate}
        \item \textbf{Objects:} Finite sets 
        \begin{equation}
            \langle m \rangle = \{0\,\cdots, m\}, \quad m\geq 0
        \end{equation}
        which we consider to be \emph{based} at 0. 
        \item \textbf{Morphisms:} Morphisms are basepoint-preserving functions 
        \begin{equation}
            f\maps \langle m \rangle \longra \langle k \rangle 
        \end{equation} 
        such that $f(0) = 0$.
    \end{enumerate}
    This category is called the \define{Segal category}. For each $n\geq 0$, let
      \begin{equation}
        i_k \colon \langle n\rangle \longra \langle 1\rangle,
        \qquad 1\leq k\leq n,
      \end{equation}
      be the maps picking out the $k$-th non-basepoint element. 
\end{definition}

\begin{observation}
    $\Gamma$ is the skeleton of $\Fin_*$, the category of based finite sets and basepoint-preserving maps among them.
\end{observation}

In \cite{CS19}, Calaque--Scheimbauer construct a model for the
$(\oo,1)$-category of symmetric monoidal $(\oo,1)$-categories using
$\Gamma$-objects in complete Segal spaces. Concretely, they consider the
projective model structure on
\begin{equation}
  \fun(\Gamma,\mathsf{CSSp})_{\mathsf{proj}}
\end{equation}
and perform a left Bousfield localization at a suitable set of maps
\begin{equation}
  S_{\Gamma} \subset \fun(\Gamma,\mathsf{CSSp})_{\mathsf{proj}}
\end{equation}
encoding the Segal and unit conditions in arity $m\geq 0$; see
\cite[Definition~3.6]{CS19} for an explicit description in terms of the
(co)Yoneda embedding. The fibrant objects in the localized model
structure
\begin{equation}
  L_{S_{\Gamma}}\fun(\Gamma,\mathsf{CSSp})_{\mathsf{proj}}
\end{equation}
are precisely the \emph{symmetric monoidal complete Segal spaces} in the
sense of \cite{CS19}, and the underlying $(\oo,1)$-category of this
model category is the $(\oo,1)$-category of symmetric monoidal
$(\oo,1)$-categories and symmetric monoidal functors.

\smallskip

On the other hand, \Cref{ch2:: thm: Quillen-equivalence-SSp-ISSp}
constructed a simplicial Quillen equivalence
\begin{equation}
  \begin{tikzcd}[ampersand replacement=\&]
    {\mathsf{CSSp}} \&\& {\mathsf{CISSp}}
    \arrow[""{name=0, anchor=center, inner sep=0}, "{\gamma^*}", curve={height=-12pt}, from=1-1, to=1-3]
    \arrow[""{name=1, anchor=center, inner sep=0}, "{\rkan_{\gamma}}", curve={height=-12pt}, from=1-3, to=1-1]
    \arrow["\dashv"{anchor=center, rotate=90}, draw=none, from=1, to=0]
  \end{tikzcd}
\end{equation}
between complete Segal spaces and complete interval Segal spaces. By
\cite[Proposition~A.3.3.6]{HTT}, this induces a simplicial Quillen
equivalence on $\Gamma$-diagrams
\begin{equation}\label{ch2:: eq: Gamma-diagram-Quillen}
  \begin{tikzcd}[ampersand replacement=\&]
    {\fun(\Gamma, \mathsf{CSSp})_{\mathsf{proj}}} \&\& {\fun(\Gamma, \mathsf{CISSp})_{\mathsf{proj}}}
    \arrow[""{name=0, anchor=center, inner sep=0}, "{(\gamma^*)^{\Gamma}}", curve={height=-12pt}, from=1-1, to=1-3]
    \arrow[""{name=1, anchor=center, inner sep=0}, "{\rkan_{\gamma}^{\Gamma}}", curve={height=-12pt}, from=1-3, to=1-1]
    \arrow["\dashv"{anchor=center, rotate=90}, draw=none, from=1, to=0]
  \end{tikzcd}
\end{equation}

Let $\widetilde{S}_{\Gamma} \subset \fun(\Gamma,\mathsf{CISSp})_{\mathsf{proj}}$ be any
set of maps weakly equivalent to the image of $S_{\Gamma}$ under
$(\gamma^*)^{\Gamma}$. By the general
existence theorem for left Bousfield localizations
(Theorem~\ref{ap:: thm: transfer of localizations}) and the invariance
of left Bousfield localization under Quillen equivalence, the adjunction
\eqref{ch2:: eq: Gamma-diagram-Quillen} descends to a simplicial Quillen
equivalence
\begin{equation}\label{ch2:: eq: Gamma-localized-Quillen}
  \begin{tikzcd}[ampersand replacement=\&]
    {L_{S_{\Gamma}}\fun(\Gamma, \mathsf{CSSp})_{\mathsf{proj}}} \&\& {L_{\widetilde{S}_\Gamma}\fun(\Gamma, \mathsf{CISSp})_{\mathsf{proj}}}
    \arrow[""{name=0, anchor=center, inner sep=0}, "{(\gamma^*)^{\Gamma}}", curve={height=-12pt}, from=1-1, to=1-3]
    \arrow[""{name=1, anchor=center, inner sep=0}, "{\rkan_{\gamma}^{\Gamma}}", curve={height=-12pt}, from=1-3, to=1-1]
    \arrow["\dashv"{anchor=center, rotate=90}, draw=none, from=1, to=0]
  \end{tikzcd}.
\end{equation}

Thus the localized model category
$L_{\widetilde{S}_\Gamma}\fun(\Gamma, \mathsf{CISSp})_{\mathsf{proj}}$ presents a
Quillen-equivalent
$(\oo,1)$-category of symmetric monoidal $(\oo,1)$-categories as
$L_{S_{\Gamma}}\fun(\Gamma,\mathsf{CSSp})_{\mathsf{proj}}$.

\begin{proposition}\label{ch2:: prop: fibrant-sym-monoidal-CISSp}
  A diagram
  \begin{equation}
    \cF^{\langle \bullet \rangle} \colon \Gamma \longra \mathsf{CISSp}
  \end{equation}
  is fibrant in the localized model category
  $L_{\widetilde{S}_\Gamma}\fun(\Gamma, \mathsf{CISSp})_{\mathsf{proj}}$ if and only if
  the following conditions hold:
  \begin{enumerate}
    \item For every $m\geq 0$, the value
      $\cF^{\langle m\rangle}$ is a complete interval Segal space (i.e.\
      fibrant in $\mathsf{CISSp}$).
    \item For every $m\geq 0$, the induced map
      \begin{equation}\label{ch2:: eq: Gamma-Segal-map}
        \cF^{\langle m\rangle}
        \longra
        \big(\cF^{\langle 1\rangle}\big)^{\times m}
      \end{equation}
      induced by the product of $\cF^{\lag \bullet \rag}$ evaluated on the inert maps, is a weak equivalence in $\mathsf{CISSp}$.
    \item The value at $\langle 0\rangle$ is essentially terminal: the
      terminal map
      \begin{equation}
        \cF^{\langle 0\rangle} \longra \ast
      \end{equation}
      is a weak equivalence in $\mathsf{CISSp}$.
  \end{enumerate}
\end{proposition}

\begin{proof}
  On the $\mathsf{CSSp}$-side this is exactly the characterization of
  fibrant objects in $L_{S_{\Gamma}}\fun(\Gamma,\mathsf{CSSp})_{\mathsf{proj}}$
  given in \cite[Definition~3.6]{CS19}. The Quillen equivalence
  \eqref{ch2:: eq: Gamma-localized-Quillen} transports this description
  to the $\mathsf{CISSp}$-side.
\end{proof}

We package this as our working definition.

\begin{definition}\label{ch2:: def: symmetric-monoidal-CISSp}
  A \textbf{symmetric monoidal (complete) interval Segal space} is a
  fibrant object of the localized model category
  $L_{\widetilde{S}_\Gamma}\fun(\Gamma, \mathsf{CISSp})_{\mathsf{proj}}$, equivalently a
  functor
  \begin{equation}
    \cF^{\langle \bullet \rangle} \colon \Gamma \longra \mathsf{CISSp}
  \end{equation}
  satisfying the conditions of
  Proposition~\textup{\ref{ch2:: prop: fibrant-sym-monoidal-CISSp}}. We denote this model structure on $\fun(\Gamma, \mathsf{CISSp})_{\mathsf{proj}}$ as $\mathsf{CISSp}^{\otimes}$.
\end{definition}

\begin{remark}
  Via the Quillen equivalence
  \eqref{ch2:: eq: Gamma-localized-Quillen}, the homotopy theory of
  symmetric monoidal (complete) interval Segal spaces is equivalent to
  the homotopy theory of symmetric monoidal complete Segal spaces of
  \cite{CS19}, and hence presents the $(\oo,1)$-category of symmetric
  monoidal $(\oo,1)$-categories. When no confusion can arise, we will
  simply speak of \emph{symmetric monoidal $(\oo,1)$-categories} and
  implicitly model them by symmetric monoidal (complete) interval Segal
  spaces.
\end{remark}

\subsection{Other models for \texorpdfstring{$(\oo,1)$}{(infty,1)}-categories}

It is important to make a note of a handful of different models for
$(\oo,1)$-categories we will be using throughout this thesis.  This entire
chapter is dedicated to proposing a variant model of $(\oo,1)$-categories that
is better suited to construct differential-topological objects --- yet,
throughout the rest of this thesis, other models will be employed and
comparison results between different models will be assumed.  To keep the text
as self-contained as possible, we briefly recall here some of the models we
use and how they are related.

\subsubsection{Quasicategories}

As we have discussed in the introduction to this chapter, an
$(\oo,1)$-category can be understood as a category enriched in
$(\oo,0)$-categories, or $\oo$-groupoids.  By the \emph{Homotopy Hypothesis},
the theory of $\oo$-groupoids, taken up to weak homotopy equivalence, agrees
with the theory of Kan complexes up to homotopy equivalence --- which, by
\Cref{ch2:: thm: homotopy hyp}, is equivalent to the ordinary homotopy theory
of topological spaces.

\begin{definition}\label{ch2:: def: quasicategory}
  A \textbf{quasicategory} $X$ is a simplicial set
  $X_{\bullet} \colon \bDelta^{\op} \to \set$ which satisfies the following
  lifting condition: for every $n \geq 1$, for every $0 < p < n$, and for
  every $v \colon \Lambda^{n}_{p} \to X$ there exists (not necessarily
  uniquely) a \emph{filler} for the commutative diagram
  \[
  \begin{tikzcd}[ampersand replacement=\&]
  	{\Lambda^{n}_{p}} \&\& X \\
  	\\
  	{\Delta^{n}} \&\& {.}
  	\arrow["v", from=1-1, to=1-3]
  	\arrow[hook, from=1-1, to=3-1]
  	\arrow[dashed, from=3-1, to=1-3]
  \end{tikzcd}
  \]
  In other words, $X$ satisfies the \emph{inner Kan condition}.
\end{definition}

\begin{observation}
  A quasicategory $X$ encodes the homotopy theory of its mapping spaces in its
  higher simplices.  Given a pair $x,y \in X_0$ of objects, one defines the
  \emph{mapping space} $\Map_{X}(x,y)$ as the homotopy pullback
  \[
  \begin{tikzcd}[ampersand replacement=\&]
    {\Map_{X}(x,y)} \&\& {X^{\Delta^1}} \\
    \\
    {*} \&\& {X^{\Delta^0}\times X^{\Delta^0} \quad ,}
    \arrow[from=1-1, to=1-3]
    \arrow[from=1-1, to=3-1]
    \arrow["{(d_1,d_0)}", from=1-3, to=3-3]
    \arrow["{(x,y)}"', from=3-1, to=3-3]
  \end{tikzcd}
  \]
  where $X^{K} \coloneqq \Hom_{\sset}(K,X)$ and $d_0,d_1\colon X_1\to X_0$
  are the face maps.  It is a theorem of Lurie that $\Map_X(x,y)$ is always a
  Kan complex; see, for instance, \cite[§2.2.5]{HTT}.
\end{observation}

\begin{observation}
  The theory of quasicategories is especially useful because it offers a very
  compact way of presenting an $(\oo,1)$-category: one constructs a simplicial
  set and checks a collection of lifting conditions.  Depending on the
  $(\oo,1)$-categories one wants to use, this is often tractable in practice.
  Another strong motivation to work with quasicategories is that they admit a
  good notion of $(\oo,1)$-category of functors: for quasicategories $X$ and
  $Y$ there is an internal hom $Y^{X}$, again a quasicategory, which models
  the $(\oo,1)$-category of functors $X \to Y$; see \cite[§1.2.7]{HTT}.
\end{observation}

\begin{remark}
  For every (ordinary) category $\cC$, the (ordinary) simplicial nerve
  $\sN(\cC)$ is a quasicategory.  In this sense, ordinary 1-categories embed
  fully faithfully into the $(\oo,1)$-category of all quasicategories.  We
  will almost always tacitly conflate an ordinary category with its nerve.
\end{remark}

\begin{warning}
  In some parts of the literature this ambiguity is avoided at great cost,
  with authors carefully notating whenever a category is viewed through its
  embedding into quasicategories.  We will not follow this practice except
  when it is genuinely helpful to avoid confusion.
\end{warning}

\medskip

The theory of quasicategories as avatars for $(\oo,1)$-categories is developed
systematically in \cite{HTT}.  In particular, $\oo$-versions of foundational
constructions and theorems of category theory such as the Yoneda lemma,
adjoint functor theorems, and Kan extensions are defined and suitably
generalized in that framework.

\subsubsection{Topological categories}

\begin{definition}
  A \textbf{topological category} is a category $\cC$ enriched in (say)
  compactly generated weak Hausdorff topological spaces.  Thus $\Ob(\cC)$ is
  a set, and for every pair of objects $x,y \in \Ob(\cC)$ there is a mapping
  space
  \[
    \Map_{\cC}(x,y) \in \Top,
  \]
  together with continuous composition maps and units satisfying the usual
  associativity and unitality axioms.
\end{definition}

Starting from a topological category $\cC$, one obtains a simplicial category
by applying the singular complex functor levelwise to the mapping spaces:
\[
  \Sing(\cC)(x,y) \coloneqq \Sing\bigl(\Map_{\cC}(x,y)\bigr) \in \sset.
\]
Lurie defines the \emph{topological nerve} $\sN^{t}(\cC)$ as the coherent
simplicial nerve of this simplicial category:
\[
  \sN^{t}(\cC) \,\coloneqq\, \sN\bigl(\Sing(\cC)\bigr) \;\in\; \sset.
\]

\begin{proposition}
  Let $\cC$ be a (small) topological category.  Then the topological nerve
  $\sN^{t}(\cC)$ is a quasicategory.
\end{proposition}

\begin{proof}[Sketch of proof]
  The singular construction $\Sing(\cC)$ is a simplicial category whose
  mapping spaces are Kan complexes.  The coherent nerve of such a simplicial
  category is a quasicategory; see \cite[§1.1.5]{HTT}.  By definition
  $\sN^t(\cC)$ is exactly this coherent nerve.
\end{proof}

\begin{example}
  Let $\cC$ be the topological category whose objects are CW-complexes and
  whose mapping spaces are the compact-open mapping spaces
  $\Map_{\Top}(X,Y)$.  Then $\sN^{t}(\cC)$ is an $(\oo,1)$-category of spaces,
  equivalent (via a suitable chain of Quillen equivalences) to the
  $(\oo,1)$-category presented by the Kan--Quillen model structure
  $\sset_{\mathrm{QK}}$.
\end{example}

\subsubsection{Simplicially-enriched categories}

\begin{definition}
  A \textbf{simplicial category} is a category $\cC$ enriched in the cartesian
  monoidal category of simplicial sets.  Concretely, this means that for each
  pair of objects $x,y \in \Ob(\cC)$ there is a simplicial set
  $\Map_{\cC}(x,y) \in \sset$, together with simplicial composition maps
  \[
    \Map_{\cC}(y,z) \times \Map_{\cC}(x,y)
    \longrightarrow
    \Map_{\cC}(x,z)
  \]
  and units $*\to\Map_{\cC}(x,x)$ satisfying the usual associativity and unit
  axioms.  A functor between simplicial categories
  $F \colon \cC \to \cD$ is \emph{simplicial} if it induces morphisms of
  simplicial sets
  \[
    \Map_{\cC}(x,y) \longrightarrow \Map_{\cD}\bigl(F(x),F(y)\bigr)
  \]
  compatible with composition and identities.  We write $\mathsf{SC}$ for the
  (ordinary) category of simplicial categories and simplicial functors.
\end{definition}

\begin{remark}
  Applying the singular complex functor to the mapping spaces of a topological
  category $\cC$ produces a simplicial category $\Sing(\cC)$, and conversely
  geometric realization of the mapping spaces of a simplicial category
  produces a topological category.  These two constructions are inverse up to
  homotopy, and they induce a Quillen equivalence between suitable model
  structures on topological and simplicial categories; see, for instance,
  \cite[Ch.~2]{Ber18}.  In practice, we will freely pass between topological
  and simplicial categories without comment.
\end{remark}

A key notion of weak equivalence between simplicial categories is the
following.

\begin{definition}
  A simplicial functor $F \colon \cC \to \cD$ is a
  \textbf{Dwyer--Kan equivalence} if:
  \begin{enumerate}
    \item for every pair of objects $x,y \in \Ob(\cC)$, the induced map
      \[
        \Map_{\cC}(x,y) \longrightarrow \Map_{\cD}\bigl(F(x),F(y)\bigr)
      \]
      is a weak equivalence of simplicial sets; and
    \item the induced functor on homotopy categories
      \[
        \ho(\cC) \longrightarrow \ho(\cD),
        \qquad
        \ho(\cC)(x,y) \coloneqq \pi_0\Map_{\cC}(x,y),
      \]
      is essentially surjective.
  \end{enumerate}
\end{definition}

\begin{theorem}[Bergner]
  There is a cofibrantly generated model structure
  $\mathsf{SC}_{\kan}$ on $\mathsf{SC}$ whose weak equivalences are the Dwyer--Kan
  equivalences.  In this model structure, fibrant objects are precisely the
  simplicial categories whose mapping spaces are Kan complexes (i.e.~categories
  enriched in Kan complexes), and we refer to $\mathsf{SC}_{\kan}$ as the
  \emph{$\kan$-model structure} on simplicial categories.
\end{theorem}

There is a close relationship between this model structure and the Joyal model
structure on quasicategories.

\begin{theorem}
  There exists a Quillen equivalence
  \[
    \mathfrak{C} \colon \sset_{\mathrm{Joyal}}
    \rightleftarrows
    \mathsf{SC}_{\kan} \colon \sN
  \]
  between the Joyal model structure on simplicial sets and the Bergner model
  structure on simplicial categories, where the right adjoint $\sN$ is the
  \emph{coherent nerve} functor.  In particular, both model categories
  present the same underlying $(\oo,1)$-category of $(\oo,1)$-categories.
\end{theorem}

\begin{remark}
  Via Rezk's model structure on complete Segal spaces 
   $\mathsf{CSSp}$ recalled earlier, there is in fact a zigzag
  of Quillen equivalences
  \[
    \sset_{\mathrm{Joyal}}
    \;\simeq_Q\;
    \mathsf{CSSp}
    \;\simeq_Q\;
    \mathsf{SC}_{\kan}
  \]
  presenting the same homotopy theory of $(\oo,1)$-categories; see for
  instance \cite[Ch.~2]{Ber18}.  We will use these comparison results
  implicitly when moving between quasicategories, complete Segal spaces,
  simplicial categories, and the interval Segal spaces introduced in this
  chapter.
\end{remark}

\begin{remark}[Notational convention]
  We will sometimes write a map $\cC \to \cD$ from a topological category
  $\cC$ directly into an $(\oo,1)$-category modelled as a quasicategory or a
  complete Segal space.  By this we always mean the composite
  \[
    \sN\bigl(\Sing(\cC)\bigr) \longrightarrow \cD,
  \]
  where $\Sing$ takes each topological mapping space $\Map_{\cC}(x,y)$ to its
  singular complex and $\sN$ is the coherent simplicial nerve.  This
  convention allows us to suppress a number of intermediate functors from the
  notation.
\end{remark}

%% file: figures/Ch2/INT10.tikz
\begin{tikzpicture}
	\begin{pgfonlayer}{nodelayer}
		\node [style=0 strat] (0) at (-3, 0) {};
		\node [style=0 strat] (1) at (0, 0) {};
		\node [style=0 strat] (2) at (3, 0) {};
		\node [style=none] (3) at (-3, -0.75) {$0$};
		\node [style=none] (4) at (0, -0.75) {$1$};
		\node [style=none] (5) at (3, -0.75) {$2$};
	\end{pgfonlayer}
\end{tikzpicture}

%% file: Chapters/Ch3.tex
\chapter{The \texorpdfstring{$(\oo,1)$} --category of cobordisms}

The main goal of this thesis is to produce an $(\oo,1)$-category of stratified cobordisms with defects, and to sketch a higher version thereof. Successfully achieving this will require exploring a toolkit which can feel quite ornate and unfamiliar to a large swath of practicing mathematicians. To improve on this current state of affairs, introduce important ideas and results, and build a consistent narrative, this section is dedicated to introducing the two constructions we aim to generalize: the \emph{unstratified} $(\oo,n)$-category of cobordisms and the ordinary symmetric monoidal category of cobordisms with defects. 

\medskip

\section{A new construction of the \texorpdfstring{$(\oo,1)$}{(infty,1)}-bordism category}

In this section we present a new construction of the $d$-dimensional bordism
$(\oo,1)$-category. This construction is a spiritual cousin of earlier work:
most notably that of Galatius--Madsen--Tillman--Weiss \cite{GMTW09}, later
developed by Lurie \cite{Lur09}, and more recently by Scheimbauer--Calaque
\cite{CS19}. 

The pioneering work of Galatius--Madsen--Tillman--Weiss offered the first tools to construct an $(\oo,1)$-category of cobordisms. In their paper on the homotopy type of the cobordism category, they construct a topological $d$-dimensional bordism category by using classifying spaces. Following their work, Lurie, in his sketch of a proof of the cobordism hypothesis \cite{Lur09}, put forward a first construction of a model for an $(\oo,n)$-category of cobordisms using $n$-fold complete Segal spaces. This contained a small mistake preventing it from satisfying the \emph{essential constancy} condition. In order to address this, in her PhD thesis and in a subsequent paper with Damien Calaque, Scheimbauer--Calaque propose an alternative definition, which fixes Lurie's while retaining the use of $n$-fold complete Segal spaces. Despite this almost successful attempt at a \emph{corrigendum}, the work of Scheimbauer--Calaque contains a separate flaw when it claims to prove that the simplicial sets which model the levels of its $n$-fold Segal spaces are Kan complexes---which translates to a loss of fibrancy. We will explore this in the next chapter when we discuss the higher $n$ case.

Our approach will be to build $(\oo,n)$-categories using an alternative model, related to that of $n$-fold complete Segal spaces but more adapted to the construction of higher categories of cobordisms. For $n=1$, the case considered in this section, this means constructing an interval Segal space, whose theory we developed in \Cref{ch2:: title}. In the course of this construction we will pick up some tools from the theory of sheaves on manifolds, in particular as they are used and developed in the works of Berwick-Evans--Boavida de Brito--Pavlov \cite{BBEP24} and Ayala--Francis--Rozenblyum \cite{AFR1}. 

\subsection{A recollection of previous constructions}

Before getting into the deeper technical aspects of the mathematical work in this chapter, it is useful to first place it in its historical context. This serves both as a pedagogical aid for future users of higher cobordism categories and as a warm-up to the ideas of Calaque--Scheimbauer following Lurie, whose interpretation of cobordisms is, in spirit, different from the one found in earlier work on the subject. 

Instead of introducing the work of Lurie, as a means to point out what the work of Calaque--Scheimbauer fixes, we will instead present the main ideas driving the work and point out how it differs from Lurie and what fails in the original paper by the same author on the \emph{Cobordism Hypothesis}.

We are used to formulating cobordism categories from the \emph{bottom up}. In other words, we often define the cobordism category as a category whose $k$-morphisms are $(d-k)$-dimensional manifolds with corners and whose composition is given by gluing along the ``right'' boundary components. While this interpretation is retained in the works of Lurie and Calaque--Scheimbauer, their approach can be described as \emph{top down}. One of the most contrived aspects of defining higher categories of bordisms is composition, since it comes with a compounding pile of choices whose compatibility is easy to check, but very hard to write down.

This chapter treats only the $(\oo,1)$ case. We begin by discussing the work of
Lurie and Calaque--Scheimbauer, both to explain the background and to motivate
our choice of interval Segal spaces. 

For every $k\geq 0$, the space $(\bord_d^{\sf un})_k$ should encode the classifying space for composable chains of $k$ unoriented bordisms:
\begin{equation}
    \Sigma_0 \xra{M_1} \Sigma_1 \xra{M_2} \cdots \xra{M_k} \Sigma_k.
\end{equation}
Lurie's proposal is to define this space by picking auxiliary data and then showing it does not alter the relevant homotopy type. In this case, the auxiliary data is a finite-dimensional real vector space $V$ and a smooth embedding
\begin{equation}
    M \hookra V\times \RR
\end{equation}  
of the composed $d$-dimensional bordism $M = M_1 \sqcup_{\Sigma_1} \cdots \sqcup_{\Sigma_{k-1}} M_k$. Once we fix this embedding of $M$, we can encode the gluing regions i.e. the $\Sigma_i$'s, as certain intersections or ``cuts'' by transverse hyperplanes $V\times \{t_i\}$. To keep this visually grounded, let us consider the following example
\begin{equation}
    \tikzfig{figures/Ch3/BOR6}{\input{figures/Ch3/BOR6.tikz}}
\end{equation} 
where the arrows identify the levels $t_i$ and the corresponding cuts. Lurie's idea is thus to define the classifying space for composable chains of bordisms of length $k$ by considering the set of all smooth embeddings of $d$-dimensional unoriented manifolds $M \hookra V \times \RR$ together with a choice $\{t_0\leq \cdots \leq t_k\} \subset \RR$ of $k$ points and such that the restricted projection  $\pi_{|M}\maps M \fib \RR$ is proper, a necessary technical condition we will unpack afterwards, and whose set of critical values $\mathsf{Crit}(\pi_{|M})$ is disjoint from $\{t_0\leq \cdots \leq t_k\}$. Importantly, this set can be endowed with a topology which encodes the correct homotopy type --- and although details matter here, we will leave a deeper discussion for the end of this chapter.

Lurie's definition \cite[Def 2.2.9]{Lur09}, despite requiring a lot of legwork to set up properly, offers a perfectly valid construction of a Segal space, which can be further completed to encode the $(\oo,1)$-category of unoriented cobordisms. Despite this, the Segal condition proves hard to check, and it does so for reasons which are unsurprising at this stage: cutting manifolds across a single submanifold, or gluing them if you see it as the reverse operation, is easier to do if you first specify collars. Failure to do so begets issues of coherence and pseudo-isotopy — which, although not intractable, increase the burden of working with such models. This issue can be ameliorated by placing the legwork elsewhere, and this is at the core of the variant introduced in the work of Calaque--Scheimbauer. 

To build collars into the definition, Calaque--Scheimbauer introduce an auxiliary $(\oo,1)$-category of closed intervals $\cI$. Calaque--Scheimbauer denote this category by $\mathrm{Int}$; we use $\cI$ to distinguish it from our interval category $\Int$. As a complete Segal space, its $k$-levels are defined as $(k+1)$-tuples 
\begin{equation}
    \cI_k = \{([a_0,b_0]\leq \cdots \leq [a_k,b_k])\,:\, [a_i,b_i]\subset \RR \}
\end{equation}
topologized via the standard Euclidean topology on $\RR^{2k+2}$. With these intervals, they improve on Lurie's definition by pairing the embedding $M \hookra V\times \RR$ with a choice of a $k$-simplex $(I_0\leq \cdots \leq I_k)\in \cI_k$ such that $\mathsf{Crit}_{|\pi}$ is instead disjoint from $I_0 \cup \cdots \cup I_k$. Visually, we can depict a $k$-simplex as follows 
\begin{equation}
    \scalebox{0.8}{\tikzfig{figures/Ch3/BOR7}{\input{figures/Ch3/BOR7.tikz}}}
\end{equation}
where the red shaded regions in $M$ are manifestly the pre-images of $\pi_{|M}$ over the intervals $I_i$, and correspond to the collars along which the bordism decomposes. This definition has a slightly larger amount of overhead than Lurie's, but it offers a clearer way of checking the Segal condition and it fixes a failure of essential constancy. 

Our definition tries to reduce the extra machinery in Calaque--Scheimbauer's interval $(\oo,1)$-category $\cI$ by, once again, placing it elsewhere. We use interval Segal spaces, and with this model intervals are part and parcel of the definition of $(\oo,1)$-categories, rather than an auxiliary tool which merely supports the construction of bordism categories. This approach will prove especially useful when we decide to construct the $(\oo,n)$-category of conically smooth stratified cobordisms --- offering almost a machine for churning out higher categories of topological objects embedded in Euclidean space.

In our construction we no longer talk about $k$-simplices as direct encoders of chains of composable morphisms, instead we use a disjoint union of open intervals $I \subset [0,1]$ containing the endpoints. An $I$-level element will be, similarly to the above constructions, an embedding $M \hookra V\times [0,1]$, which we depict as
\begin{equation}
    \scalebox{0.8}{\tikzfig{figures/Ch3/BOR8}{\input{figures/Ch3/BOR8.tikz}}}
\end{equation}
where the restricted projection $\pi_{|M}\maps M \fib [0,1]$ must be proper and such that $\mathsf{Crit}(\pi_{|M})$ is disjoint from $[0,1]\setminus I$. Note that the red shaded regions of $M$ correspond to the pre-image of $I$ under the restricted projection $\pi_{|M}$, and they contain the critical points of the map --- it is over these regions that the cobordism will undergo meaningful changes in topology as it may contain any finite number of critical points with respect to the projection --- on the other hand, in its complementary regions the property outlined above implies that the embedding must be diffeomorphic to a cylinder of a $(d-1)$-dimensional manifold. We will explore this model and its subtle distinctions to the ones discussed above in the next sections.

\subsection{Sheaves on manifolds}

One of the distinctions in our methodology when approaching higher categories of cobordisms, as compared to Lurie \cite{Lur09} and Calaque--Scheimbauer \cite{CS19}, has to do with how we set up the spatial structure. Before any technical detail it is useful to record the aim of this section via the slogan: \emph{offering tools to convert simplicial homotopy theory into smooth homotopy theory}. To properly set up the importance of this issue, we should explain what fails in aforementioned works.

In the work of Calaque--Scheimbauer, the spatial structure of the $n$-fold Segal spaces is constructed with recourse to the smooth singular complex. More concretely, Lemma 5.7 of \cite{CS19} states that the smooth $l$-simplices of a certain topological space can be neatly described as a pair of an embedding $(M \hookra V \times \RR)$ and a sequence of intervals $I_0 \leq \cdots \leq I_k$ satisfying three conditions. From \cite[Corollary 4.47]{CW}, it is understood that a manifold $M$ with boundary or corners is not fibrant as a diffeological space; meaning its associated smooth singular complex $\mathsf{Sing}^{\sf sm}(M)$ is not, in general, a Kan complex. Tracking the reasoning in \cite{CS19}, it is unclear how to ensure that the topological space in question has the fibrancy properties required for the Segal-space argument.

To avoid dealing with the cumbersome point-set topologies on spaces of submanifolds, one can instead work with simplicial sets. These, however, are rarely Kan complexes. This failure can be detected by testing whether it is possible to glue manifolds along closed subsets. Such gluing procedures, much like gluing smooth functions, are generally impossible along closed sets. This line of thinking is what motivates the appearance of sheaves on the site of smooth manifolds, and in particular the framework developed by Berwick-Evans, Boavida de Brito, and Pavlov. Their theory provides precisely the homotopical tools needed to perform such “gluing” operations along closed subsets, up to smooth homotopy.

\subsection{Sites and sheaves}

As noted above, in contrast to Calaque--Scheimbauer we do not model the spatial
structure of our $(\oo,1)$-categories via smooth singular simplices. Instead,
we work with a homotopical theory of ``spaces'' that has superior gluing
properties — namely, ($\oo$-)sheaves on the site of smooth manifolds.  This sheaf-theoretic viewpoint clarifies the necessary gluing arguments. Our main references for these methods are
\cite{BBEP24} and \cite{MW07}. 

\medskip

We begin by recalling the basic sheaf-theoretic formalism. Let $\cC$ be a small category. The presheaf category $\psh(\cC)$ is the free cocompletion of $\cC$: any diagram in $\cC$ that may fail to admit a colimit does admit one in $\psh(\cC)$ after applying the Yoneda embedding
\begin{equation}
  \cC \hookra \psh(\cC).
\end{equation}
In particular, if $\man$ denotes the category whose objects are
finite-dimensional smooth manifolds without boundary and whose morphisms are
smooth maps, then $\man$ does not admit all colimits (e.g. the span
\[
* \longleftarrow S^0 \xra{\nssphere} S^1
\]
has no colimit in $\man$), whereas $\psh(\man)$ does. However, the colimits in
$\psh(\man)$ are generally \emph{not} the ``correct'' geometric ones.

A simple indication of this failure is disjoint union. Given smooth manifolds
$M_1,M_2$, the coproduct $M_1\sqcup M_2$ exists in $\man$ and is given by
disjoint union. A map
\[
S^0 \longra M_1 \sqcup M_2
\]
amounts to choosing either two points in a single component or one point in
each component. On the other hand, in presheaves we have
\[
\Hom_{\psh(\man)}(rS^0,rM_1 \sqcup rM_2)
\simeq \man(S^0,M_1)\,\sqcup\,\man(S^0,M_2),
\]
so a map $rS^0 \to rM_1 \sqcup rM_2$ is equivalent to choosing two points in
$M_1$ or two points in $M_2$. In other words, we are gluing incorrectly: the
geometric context encoded in $\man$ has been lost. Sheaves repair exactly this
loss by forcing certain ``geometric'' cones to become colimits.

\medskip

Abstractly, the geometric content is encoded by specifying which families of
maps $\{c_i\to c\}$ should behave like covers, and then requiring the
presheaf to reconstruct its value on $c$ as the limit of its values on the
$c_i$ and their intersections. This is the notion of a Grothendieck topology.

\begin{definition}\label{ch3:: def: site}
   Let $\cC$ be a category. A \textbf{Grothendieck topology} on $\cC$ assigns
   to every object $c\in \cC$ a collection $\tau(c)$ of families
   $\{f_i\maps c_i \longra c\}$, called \textbf{covering families}, such that:
   \begin{enumerate}
     \item If $f\maps d \longra c$ is an isomorphism, then
     $\{f\maps d \longra c\}$ is a covering family.
     \item For a covering family $\{f_i \maps c_i \longra c\}$ and any map
     $g\maps d \longra c$, the fiber product
     \[
        \begin{tikzcd}[ampersand replacement=\&]
            {d \times_{c}c_i} \&\& d \\
            \\
            {c_i} \&\& c
            \arrow[from=1-1, to=1-3]
            \arrow[from=1-1, to=3-1]
            \arrow["\lrcorner"{anchor=center, pos=0.125}, draw=none, from=1-1, to=3-3]
            \arrow[from=1-3, to=3-3]
            \arrow[from=3-1, to=3-3]
        \end{tikzcd}
     \]
     exists in $\cC$ and the family of projections
     $\{\pi_i \maps d\times_c c_{i} \longra d\}$ is a covering family.
     \item If $\{f_i \maps c_i \longra c\}$ is a covering family and, for each
     $i$, $\{g_{ji}\maps c_{ji} \longra c_{i}\}$ is a covering family, then
     the composites
     \[
        \{f_i \circ g_{ji}\maps c_{ji} \longra c\}
     \]
     form a covering family of $c$.
   \end{enumerate}
   A covering family with a single map $f\maps d \longra c$ is called a
   \textbf{covering map}. A category $\cC$ endowed with a Grothendieck
   topology $\tau$ is called a \textbf{site}.
\end{definition}

We now introduce the only example of a Grothendieck topology we will need.

\begin{definition}\label{ch3:: def: smooth site}
    Let $\man$ denote the category whose objects are all finite-dimensional
    smooth manifolds without boundary and smooth maps thereamong. The site
    structure is given by \emph{open covers}. For $M\in\man$, a
    \textbf{covering family} is a collection of open embeddings
    \[
      \{U_i \hookra M\}_{i\in I}
    \]
    such that
    \[
      M \,=\, \bigcup_{i \in I} U_i,
    \]
    where we conflate each open embedding with its image.
\end{definition}

\begin{definition}\label{ch3:: def: sheaf}
    Let $\cC$ be a site with Grothendieck topology $\tau$. A functor
    $\cF\maps \cC^{\op} \longra \set$ is called a \textbf{presheaf}. A
    presheaf $\cF$ is a \textbf{sheaf} (with respect to $\tau$) if for each
    covering family $\{f_i\maps c_i \longra c\}_{i\in I}$ the dashed arrow in
    \[
      \cF(c) \dashrightarrow \prod_{i\in I}\cF(c_i)
      \rightrightarrows \prod_{i,j\in I}\cF(c_i\times_c c_j)
    \]
    exhibits $\cF(c)$ as the equalizer of the two parallel arrows.
\end{definition}

In what follows we are primarily interested in sheaves on the smooth site $\man$ of Definition~\ref{ch3:: def: smooth site}.
\bigskip

\subsubsection{Concordance after Madsen--Weiss}

We now recall the notion of concordance from Madsen--Weiss as it is presented in \cite{MW07} --- these ideas can be found in earlier papers of Waldhausen \cite{Wa} and Weiss--Williams \cite{MW95}. This is the basic equivalence relation that will later be promoted to a homotopical condition for $\oo$-sheaves.

\begin{definition}[Concordance]\label{ch3:: def: concordance}
  {\cite[Definition~2.1.3]{MW07}}
  Let $\cF\maps \man^{\op} \longra \set$ be a presheaf and let
  $\pi \maps X \times \RR \longra \RR$ be the projection. Two sections
  $\sigma_0,\sigma_1 \in \cF(X)$ are \textbf{concordant} if there exists a
  section $\Sigma \in \cF(X\times\RR)$ such that:
  \begin{enumerate}
    \item There is an open neighborhood
    $X \times (-\oo,0] \subset U_{-}$ with
    $\Sigma_{|U_{-}} = (\pi^* \sigma_{0})_{|U_{-}}$.
    \item There is an open neighborhood
    $X \times [0,+\oo) \subset U_{+}$ with
    $\Sigma_{|U_{+}} = (\pi^* \sigma_{1})_{|U_{+}}$.
  \end{enumerate}
  The section $\Sigma$ is called a \textbf{concordance} from $\sigma_0$ to
  $\sigma_1$, and we write $\sigma_0 \sim_c \sigma_1$.
\end{definition}

\begin{example}
  Let $N$ be a smooth manifold and let
  \[
    \cF = C^{\oo}(-,N)\maps \man^{\op} \longra \set
  \]
  be the presheaf $M \mapsto C^{\oo}(M,N)$. Given $f,g\in C^{\oo}(M,N)$, a
  concordance $H\in C^{\oo}(M\times\RR,N)$ is precisely a smooth homotopy
  $H\maps M\times\RR\to N$ from $f$ to $g$ that is eventually constant in
  each end, i.e.\ $H(x,t)=f(x)$ for $t<0$ and $H(x,t)=g(x)$ for $t>0$.
\end{example}

\begin{lemma}\label{ch3:: lem: concordance is equiv}
    Let $\cF\maps \man^\op \longra \set$ be a sheaf. Then the relation
    $\sim_c$ on $\cF(X)$ is an equivalence relation.
\end{lemma}

\begin{proof}
  Reflexivity is witnessed by $\pi^*\sigma\in\cF(X\times\RR)$ for any
  $\sigma\in\cF(X)$. Symmetry is obtained by precomposing a concordance
  $\Sigma$ with $-\id_{\RR}$. For transitivity, suppose
  $\sigma_0\sim_c\sigma_1$ and $\sigma_1\sim_c\sigma_2$ with concordances
  $\Sigma_{01},\Sigma_{12}\in\cF(X\times\RR)$. Choosing a diffeomorphism
  $\phi\maps\RR\to\RR$ that sends an open neighborhood of $+\infty$ to all of
  $\RR$ and reparameterizing $\Sigma_{12}$ along $\phi$, one can glue the two
  concordances along an open neighborhood of $X\times\{0\}$ using the sheaf
  property, and then reparameterize again to obtain a concordance from
  $\sigma_0$ to $\sigma_2$. Details follow the argument in \cite[Lemma~3.1]{BBEP24}.
\end{proof}

\begin{remark}\label{ch3:: rmk: MW-representability}
  Given a sheaf $\cF\maps\man^\op\to\set$ and $X\in\man$, consider the set
  \[
    \cF[X] \;\overset{\sf def}{=}\; \cF(X)/\sim_c.
  \]
  This assignment is functorial in $X$ and defines a presheaf
  $\cF[\bullet]\maps\man^\op\to\set$. It is rarely a sheaf again. Madsen--Weiss
  show that this presheaf is nevertheless representable in a weak sense:
  there exists a space $|\cF|$ such that for every manifold $X$ there is a
  natural bijection
  \[
    \pi_0\Map(X,|\cF|)\;\simeq\;\cF[X].
  \]
  Moreover, $|\cF|$ can be constructed as the geometric realization of the
  simplicial set $[n]\mapsto \cF(\Delta^n_e)$, where $\Delta^n_e$ is the
  extended simplex in the sense of \cite[Definition~2.4.2]{MW07}, and this
  representability extends to relative concordance classes and higher
  homotopy groups.
\end{remark}

Madsen--Weiss work with set-valued sheaves. In our application we need the
fully homotopical version developed by Berwick-Evans--Boavida de Brito--Pavlov
\cite{BBEP24}, where one replaces sheaves of sets by $\oo$-sheaves of spaces.
\bigskip

\subsubsection{\texorpdfstring{$\oo$}{infty}-sheaves and the concordance resolution}

We now upgrade the previous discussion to $\oo$-sheaves. Following
\cite{BBEP24}, by a \emph{space} we mean a simplicial set, and by a
\emph{simplicial space} a simplicial object in spaces.

\begin{definition}\label{ch3:: def: infinity-sheaf}
  {\cite[Definition~2.15]{BBEP24}}
  A presheaf $\cF\maps \man^{\op}\to\sset$ is an \textbf{$\oo$-sheaf} if for every manifold $M$ and every open cover \[\{U_i\hookra M\}_{i\in I},\] the  canonical map
  \begin{equation}
    \cF(M)\;\longra\;
    \hoLimit[\bDelta]\Big(
      \prod_{i_0}\cF(U_{i_0})
      \doubleleftarrow
      \prod_{i_0,i_1}\cF(U_{i_0}\cap U_{i_1})
      \tripleleftarrow\cdots
    \Big)
   \end{equation}
  is a weak equivalence of spaces.
\end{definition}

\begin{remark}
    The homotopy limit in the definition of an $\oo$-sheaf for a given open cover $\cU = \{U_\alpha \hookra M\}$ of $M$, is taken over the cosimplicial space 
    \begin{equation}
        \bDelta \xlonghookra{\,\,\,\,\, (\check{\cU}_{\bullet})^{\op}\,\,\,\,\,\,} \man^\op \xlonghookra{\,\, \cF \,\,} \sset,
    \end{equation} 
    where $\check{\cU}_{\bullet}$ is the Cech complex
    \begin{align}
        \check{\cU}_{\bullet} \maps \bDelta^{\op} &\longra \man \\
        [n] &\longmapsto \coprod_{i_{0},\cdots,i_{n}} U_{i_{0}}\cap U_{i_{1}} \cap \cdots \cap U_{i_{n}}.
    \end{align}
\end{remark}

\begin{theorem}[Binary descent criterion for \(\oo\)-sheaves]\label{ch3:: thm: BBEP-sheaf-criterion}
{\cite[Theorem~5.1]{BBEP24}}
A presheaf \(\cF\colon \man^{\op}\to\sset\) is an \(\oo\)-sheaf (i.e.\ satisfies descent for open covers)
if and only if \(\cF(\emptyset)\simeq *\) and the following two conditions hold:
\begin{enumerate}
  \item \label{binary cover} \textbf{Binary open covers}: For every manifold \(M\) and every open cover \(M=U\cup V\),
  the commutative square
  \begin{equation}
    \begin{tikzcd}[ampersand replacement=\&]
        {\cF(M)} \&\& {\cF(V)} \\
        \\
        {\cF(U)} \&\& {\cF(U\cap V)}
        \arrow[from=1-1, to=1-3]
        \arrow[from=1-1, to=3-1]
        \arrow[from=1-3, to=3-3]
        \arrow[from=3-1, to=3-3]
    \end{tikzcd}
  \end{equation}
  is a homotopy pullback square in \(\sset\).
  \item \label{disj union} \textbf{Disjoint unions}: If 
  \begin{equation}
    M \cong \bigsqcup_{i\in I} U_i
  \end{equation} is a (possibly infinite) disjoint union of
  submanifolds \(U_i\), then the canonical map
  \begin{equation}
    \cF(M)\;\longra\;\prod\nolimits^{h}_{i\in I} \cF(U_i)
  \end{equation}
  is a weak equivalence.
\end{enumerate}
\end{theorem}

Set-valued sheaves and groupoid-valued stacks are precisely those presheaves that are $\oo$-sheaves with values in discrete spaces and groupoids, respectively.

\begin{definition}\label{ch3:: def: concordance-invariant}
  {\cite[Definition~2.4]{BBEP24}}
  A presheaf $\cF\maps\man^{\op}\to\sset$ is \textbf{concordance-invariant} if
  for every manifold $M$ the map
  \begin{equation}
    \cF(M)\longra \cF(M\times\RR)
  \end{equation}
  induced by the projection $M\times\RR\to M$ is a weak equivalence.
\end{definition}

Let $\AA^n\subset\RR^{n+1}$ denote the smooth extended $n$-simplex
\begin{equation}
  \AA^n
  \;=\;
  \{(x_0,\dots,x_n)\in\RR^{n+1}\mid x_0+\cdots+x_n=1\},
\end{equation}
which defines a cosimplicial object $\AA^\bullet\maps\bDelta\to\man$. This cosimplicial object endows $\man$ with an enrichment in spaces by declaring 
\begin{equation}
    \Hom_{\man}(M,N)_{\bullet} \; \overset{\sf def}{=} \; \Hom_{\man}(M \times \AA^{\bullet}, N).
\end{equation}
A presheaf $\cF\maps \man^{\op} \longra \sset$ that can be enriched is therefore one which sends smooth homotopies to simplicial homotopies, and smooth homotopy equivalences to simplicial homotopy equivalences. In particular, an enriched presheaf is always concordance-invariant.

Informally, this construction turns a moduli sheaf into a space by retaining
its smooth families: a family over $\AA^1$ gives a concordance, while families
over higher affine simplices encode the corresponding higher homotopies.
For the bordism sheaves below, these are families of embedded bordisms, with
their cutting regions retained throughout the family.

\begin{definition}\label{ch3:: def: concordance-resolution}
  {\cite[Definition~2.5]{BBEP24}}
  Let $\cF\maps\man^{\op}\to\sset$ be a presheaf. Its \textbf{concordance
  resolution} is the presheaf
  \begin{equation}
    \cF(-\times \AA^{\bullet})\maps \man^{\op} \longra \fun(\bDelta^{\op},\sset).
  \end{equation}
  We denote the post-composition of $\cF(-\times \AA^\bullet)$ with the homotopy colimit (of simplicial spaces) as
  \begin{equation}
    \mathbf{B}\cF\maps\man^{\op}\longra\sset,\qquad
    \mathbf{B}\cF(M)\;=\;
    \left\vert\,[k]\longmapsto\cF(M\times \AA^k)\,\right\vert,
  \end{equation}
   We use the diagonal of the associated bisimplicial set as a model here.
\end{definition}

\begin{remark}
    The concordance resolution $\cF(-\times \AA^\bullet)$ lifts to an enriched functor defined on the hom-sets by 
    \begin{equation}
        \Hom_{\man}(M \times \AA^{n}, N) \longra \Hom_{\fun(\bDelta^\op,\sset)}(\cF(N \times \AA^{\bullet})\times \bDelta[n], \cF(M \times \AA^{\bullet})),
    \end{equation}
    for all $n\geq 0$. Post-composing this with the map to its homotopy colimit functor translates to a simplicial map
    \begin{equation}
        \Hom_{\man}(M \times \AA^{n}, N) \longra \Hom_{\sset}\left(\bB\cF(N)\times \Delta^n, \bB\cF(M )\right).
    \end{equation}
    Concluding that $\bB\cF$ is both enriched in spaces, and therefore concordance-invariant — hence the naming choice.
\end{remark}

Berwick-Evans, Boavida de Brito, and Pavlov prove that if $\cF$ is already an $\oo$-sheaf, then $\mathbf{B}\cF$ remains an $\oo$-sheaf and in fact represents concordance classes of sections \cite{BBEP24}.

\begin{theorem}\label{ch3:: thm: BBEP-classifying}
  {\cite[Theorems~1.1 and~1.2]{BBEP24}}
  Let $\cF$ be an $\oo$-sheaf on $\man$. Define the \emph{classifying space} of $\cF$ to be a Kan fibrant replacement $B\cF$ of the simplicial set $\mathbf{B}\cF(*)$. Then:
  \begin{enumerate}
    \item For every manifold $M$ there is a natural weak equivalence
      \begin{equation}
        (\mathbf{B}\cF)(M)\;\simeq\;
        \Map^{h}(\Sing M, B\cF),
      \end{equation}
      where $\Sing M$ is the singular complex and $\Map^{h}$ denotes the
      derived mapping space.
    \item In particular, $\mathbf{B}\cF$ is an $\oo$-sheaf.
    \item The canonical map
      \[
        \pi_0(\mathbf{B}\cF(M))\;\longrightarrow\;
        \pi_0\big(\cF(M)\big)/\sim_c
      \]
      is a bijection, so $B\cF$ represents concordance classes of sections of
      $\cF$ in the same way that $|\cF|$ represents concordance classes for
      set-valued sheaves in \cite{MW07}. Moreover, for sheaves of sets one
      obtains $\pi_n(B\cF,x)\cong \pi_n$ of the Madsen--Weiss representing
      space.
  \end{enumerate}
\end{theorem}

\subsubsection{Weak Kan fibrations and realization fibrations}

The proof of Theorem~\ref{ch3:: thm: BBEP-classifying} relies on a careful
comparison between the homotopy limit condition defining $\oo$-sheaves and
the homotopy colimits appearing in the concordance resolution. The technical
bridge is the notion of a \emph{weak Kan fibration} of simplicial spaces.

We briefly recall the key definitions from \cite[Section~3]{BBEP24}.

\begin{definition}\label{ch3:: def: weak-RLP}\cite[Definition~3.3]{BBEP24}
  Let $f\maps X\to Y$ be a Reedy fibration between Reedy fibrant simplicial
  spaces, and let $i\maps A\to B$ be a map of simplicial spaces. We say that
  $f$ has the \textbf{weak right lifting property} with respect to $i$ if for
  every commutative square
  \begin{equation}
       \begin{tikzcd}[ampersand replacement=\&]
            A \&\& X \\
            \\
            B \&\& Y
            \arrow[""{name=0, anchor=center, inner sep=0}, "\beta", from=1-1, to=1-3]
            \arrow["i"', from=1-1, to=3-1]
            \arrow["f", from=1-3, to=3-3]
            \arrow["{\tilde{\alpha}}"{description}, dashed, from=3-1, to=1-3]
            \arrow["\alpha"', from=3-1, to=3-3]
            \arrow["H"{description}, color={rgb,255:red,214;green,92;blue,92}, between={0.3}{0.7}, Rightarrow, from=0, to=3-1]
        \end{tikzcd}
    \end{equation}
  there exists a map $\tilde\alpha\colon B\to X$ and a \textbf{vertical
  homotopy} \[H\colon A\times\Delta[1]\to X\] satisfying the following conditions:
      \begin{enumerate}
        \item $A\times \Delta[\{0\}] \hookra A \times \Delta[1] \xra{H} X= \beta$.
        \item $A\times \Delta[\{1\}] \hookra A \times \Delta[1] \xra{H} X = \widetilde{\alpha}\circ i$.
        \item The following diagram commutes:
        \begin{equation}
        \begin{tikzcd}[ampersand replacement=\&]
            {A\times \Delta[1]} \&\& X \&\& Y \\
            \\
            A \&\& B
            \arrow["H", from=1-1, to=1-3]
            \arrow["\pi"', from=1-1, to=3-1]
            \arrow["f", from=1-3, to=1-5]
            \arrow["i"', hook, from=3-1, to=3-3]
            \arrow["\alpha"', from=3-3, to=1-5]
            \arrow["\tilde{\alpha}", dashed, from=3-3, to=1-3]
        \end{tikzcd}
        \end{equation}
    \end{enumerate}
\end{definition}

Informally, weak lifting translates to the ability of solving lifting problems for $i$ up to a specified homotopy that is \emph{vertical} over $Y$. For a more illustrative explication of the motivation behind this definition see \cite[Remark~3.4]{BBEP24}

\begin{definition}\label{ch3:: def: weak-kan-fibration}
  {\cite[Definition~3.8]{BBEP24}}
  A Reedy fibration $f\maps X\to Y$ of simplicial spaces is a
  \textbf{weak Kan fibration} if it has the weak right lifting property with
  respect to the maps
  \[
    h_{i,n,k}\colon
    \mathsf{Sd}^{i}\big(\Lambda_k[n]\big)\longrightarrow
    \mathsf{Sd}^{i}\big(\Delta[n]\big)
  \]
  for all $i\geq 0$, $n\geq1$ and $0\leq k\leq n$, where
  $\mathsf{Sd}^{i}$ denotes the $i$-fold subdivision functor on simplicial
  spaces.
\end{definition}

This notion interpolates between genuine Kan fibrations and the geometric fibrations needed for concordance resolutions. The key property is that weak Kan fibrations control the interaction between homotopy pullbacks and geometric realization.

\begin{definition}\label{ch3:: def: realization-fibration}
  A map $f\maps X\to Y$ of simplicial spaces is a
  \textbf{realization fibration} if for every $Z\to Y$ the induced map
  \begin{equation}
    \big|X\times_Y^{h} Z\big|
      \longrightarrow
    |X|\times_{|Y|}^{h}|Z|
  \end{equation}
  is a weak equivalence of spaces.
\end{definition}

\begin{theorem}\label{ch3:: thm: weak-kan-realization}
  {\cite[Theorem~3.17]{BBEP24}}
  A weak Kan fibration is a realization fibration.
\end{theorem}

In particular, geometric realization commutes with homotopy pullbacks along
weak Kan fibrations. This is the mechanism allowing BBEP to commute the
homotopy limits in the $\oo$-sheaf condition with the homotopy colimits
appearing in the concordance resolution $M\longmapsto\big|[k]\mapsto\cF(M\times \AA^k)\big|$, and thus to prove Theorem~\ref{ch3:: thm: BBEP-classifying}.

\begin{example}
    For exposition purposes, let us unpack the example given in \cite[Example~3.18]{BBEP24}, since it highlights the need for the lifts against all simplicial subdivisions when defining weak Kan fibrations. Let $X$ be the union of three non-degenerate 1-simplices and the projection $f\maps X \ra \Delta[1]$ as in the diagram below:
    \begin{equation}
        \begin{tikzcd}[ampersand replacement=\&]
            \bullet \& \bullet \\
            \bullet \& \bullet \\
            \\
            \bullet \& \bullet
            \arrow[from=1-1, to=1-2]
            \arrow[from=1-1, to=2-1]
            \arrow[""{name=0, anchor=center, inner sep=0}, from=2-1, to=2-2]
            \arrow[""{name=1, anchor=center, inner sep=0}, from=4-1, to=4-2]
            \arrow["f", between={0.2}{0.8}, from=0, to=1]
        \end{tikzcd}
    \end{equation} 
    A key couple of observations follow:
    \begin{enumerate}
        \item The map $f$ is not a Kan fibration. The following outer 2-horn (indicated in a separate color)
        \begin{equation}
            \begin{tikzcd}[ampersand replacement=\&]
                \bullet \& \bullet \\
                \bullet \& \bullet
                \arrow[color={rgb,255:red,214;green,92;blue,92}, from=1-1, to=1-2]
                \arrow[color={rgb,255:red,214;green,92;blue,92}, from=1-1, to=2-1]
                \arrow[from=2-1, to=2-2]
            \end{tikzcd}
        \end{equation}  
        over the degenerate 2-simplex, by inspection, does not admit a filler.
        \item Despite the point above, the map $f$ admits weak right lift against $h_0$; on the other hand, it does not admit a weak lift against $h_2$. To understand this example, consider the following map $\mathsf{Sd}^{2}(\Lambda^{2}[0]) \longra X$:
        \begin{equation} 
        \begin{tikzcd}[ampersand replacement=\&]
            \& \bullet \& \bullet \&\&\& \bullet \\
            {} \& \bullet \\
            \& \bullet \& \bullet \& \bullet \& \bullet \& \bullet
            \arrow[from=1-2, to=1-3]
            \arrow[from=1-2, to=2-2]
            \arrow[from=1-6, to=1-3]
            \arrow[from=3-2, to=2-2]
            \arrow[from=3-2, to=3-3]
            \arrow[from=3-4, to=3-3]
            \arrow[from=3-4, to=3-5]
            \arrow[from=3-6, to=3-5]
        \end{tikzcd}
        \end{equation} 
        then any homotopy lifting $\Delta[2] \longra \Delta[1]$ would fail to be vertical.
        \item This example indicates the importance of the weak lifts against repeated subdivisions. Indeed, this map is not a realization fibration: homotopy fiber over $0$ is contractible whilst the homotopy fiber over $1$ is weak homotopy equivalent to a pair of points — therefore $f$ is not a quasi-fibration, thereby not a realization fibration.
    \end{enumerate}
\end{example}

The specific weak Kan fibrations we will use arise from \emph{endpoint
evaluation} for $\oo$-sheaves. Consider the inclusion of the boundary
$\partial \AA^1 = \AA^0\sqcup \AA^0 \hookra \AA^1$, and form the induced map of
simplicial spaces
\begin{equation}
  \cF(\AA^{\bullet}\times \AA^{1}) \longra
  \cF(\AA^{\bullet}\times \partial \AA^{1})
\end{equation}
for an $\oo$-sheaf $\cF\maps\man^{\op}\to\sset$.

\begin{proposition}\label{ch3:: prop: weak-kan-endpoints}\cite[Proposition~4.2]{BBEP24}
  Let $\cF\maps\man^{\op}\to\sset$ be an $\oo$-sheaf. Then the map of
  simplicial spaces
  \begin{equation}
    \cF(\AA^{\bullet}\times \AA^{1}) \longra
    \cF(\AA^{\bullet}\times \partial \AA^{1})
  \end{equation}
  is a weak Kan fibration.
\end{proposition}

This proposition is the crucial glue between the concordance-theoretic framework of Berwick-Evans--Boavida de Brito--Pavlov and the Segal-type constructions which we will entertain in the last section: it allows us to treat endpoint evaluation for concordance data as a type of quasifibration (see \cite{May90}) whose behaviour is well-controlled under realization.

The strategy of using this endpoint weak Kan fibration to transport the concordance-resolution argument of \cite[Theorem~5.1]{BBEP24} into the proof of the interval Segal condition was suggested to the author by Pedro Boavida de Brito.

\medskip

We will apply these general results in the concrete situation where $\cF$ is an $\oo$-sheaf of $d$-dimensional manifolds (and later of bordisms and stratified bordisms). The upshot is that the concordance classes of such geometric data are represented by explicitly described classifying spaces, constructed functorially from $\cF$ via the concordance resolution.

\section{An interval Segal space of \texorpdfstring{$d$}{d}-dimensional bordisms}

We will now construct an interval Segal space whose associated $(\oo,1)$-category is that of $d$-dimensional unoriented bordisms. To construct this, we will make use of the theory of smooth $\oo$-sheaves, i.e.\ $\oo$-sheaves on the site of smooth manifolds, as a conduit for constructing appropriate spaces of smooth manifolds. This adds a distinct flavor to our approach, and clarifies the minor pathologies identified in previous constructions --- as explained in the introduction to this chapter. 

Before proceeding with definitions, let us first sketch the path we will take in this section. We start by defining a general construction
\[
\widetilde{\bord}_{d,V}(-,-)
\]
which is a bifunctor on $\man^{\op}\times \Int$ that, when fixing $T\in \Int$, produces a sheaf of sets on $\man^{\op}$ and captures a type of moduli-space construction for $d$-dimensional smooth manifolds embedded in a fixed finite-dimensional real vector space $V$, fibered over smooth manifolds of varying dimension. It is designed so that it also encompasses our desired classifying spaces of bordisms. Upon taking its concordance resolution and further restricting to $\bDelta^{\op}$, we get a simplicial set $\bord_{d,V}$, or a space, that we will show is the classifying space of $\gamma(T)$-composites of $d$-dimensional bordisms.

Lastly, the most technical part of this section is showing that this composite
construction, initially only a pre-interval space, is in fact an interval Segal
space whose homotopy category is equivalent to the $d$-dimensional bordism
category $\mathsf{Cob}_d$.

\begin{notation}
Given an object $T\in \Int$, we write $T^{\mathsf{inn}}$ for the union of connected open intervals which do not contain the endpoints. For objects $T$ with $0$ connected inner components we have $T^{\mathsf{inn}} = \emptyset$. 
\end{notation}

\begin{notation}
Given an object $I \in \Int$, we call the union
\[
(I^{\mathsf{c}}\cup I^{\mathsf{inn}}) \subset [0,1]
\]
the \textbf{core} of $I$ and denote it by $\mathsf{cor}(I)$. For every $I\in \Int$ the core $\core(I)$ is a connected closed interval of $[0,1]$.
\end{notation}

\begin{definition}\label{ch3:: def: bord_d}
    Let $V$ be a finite-dimensional real vector space. We define a functor 
    \begin{equation}
    \widetilde{\bord}^{B}_{d,V} \maps \mathsf{Man}^{\op}\times \Int \longrightarrow \set    
    \end{equation}
    sending a pair $(M,T)$, where $M\in \man$ and $T\in \Int$, to the set of bounded closed subsets 
    \begin{equation}
    W\subseteq M\times \mathsf{cor}(T) \times V
    \end{equation} 
    satisfying the following three conditions:
    \begin{enumerate}
    \item \label{id: dim} $W$ is a smooth $(d+\mathsf{dim}(M))$-dimensional submanifold with boundary $\partial W$ such that 
    \begin{equation}
        \partial W \subset \partial \big(M \times \mathsf{cor}(T) \times V\big) = M \times \partial (\mathsf{cor}(T)) \times V.   
    \end{equation}
    In other words, the boundary of $W$ is contained in the boundary of the ambient space.
    \item \label{id: categorical} The restricted projection 
    \begin{equation}
        \pi\maps \big(W \cap (M\times T^{\mathsf c} \times V) \big) \hookra  M \times T^{\mathsf c} \times V \fib  M \times T^{\mathsf c}
    \end{equation} 
    is a submersion. 
    \item \label{id: spatial} The projection 
    \begin{equation}
    \pi_W\maps W \hookra M\times \mathsf{cor}(T) \times V \fib M 
    \end{equation}
    is a submersion.
    \end{enumerate} 

    \noindent The functor $\widetilde{\bord}^{B}_{d,V}$ sends:
    \begin{itemize}
    \item[($\man$):] A morphism 
    \[
    (f,\mathsf{id}_{T})\maps (N,T) \longra (M,T)
    \]
    induced by a smooth map $f\maps N \longra M$ to the pullback submanifold 
    \[
    (f \times \id)^{*} W \subseteq N \times \mathsf{cor}(T) \times V,
    \]
    where $W\subseteq M\times \mathsf{cor}(T) \times V$ is a submanifold with boundary as above, and $(f\times \id)^\ast W$ is defined by the pullback diagram
    \begin{equation}
        \begin{tikzcd}[ampersand replacement=\&]
            {(f\times \id)^\ast W} \& {W} \\
            \\
            N\times \mathsf{cor}(T) \times V \& M\times \mathsf{cor}(T) \times V,
            \arrow[from=1-1, to=1-2]
            \arrow[hook, "(f\times \id)^\ast \pi_{|W}"',from=1-1, to=3-1]
            \arrow["\lrcorner"{anchor=center, pos=0.125}, draw=none, from=1-1, to=3-2]
            \arrow[hook, from=1-2, to=3-2]
            \arrow["f\times \id"', from=3-1, to=3-2]
        \end{tikzcd}
    \end{equation}
    which is a $\big(d + \mathsf{dim}(N)\big)$-dimensional submanifold with boundary
    \begin{equation}
        (f\times \id)^* W \subseteq N \times \mathsf{cor}(T) \times V,
    \end{equation}
    since the pullback preserves the fibers over the boundary $\partial (N \times \mathsf{cor}(T) \times V)$.

    \item[($\Int$):] A morphism $(\id_M, T\subset S)$ is defined by \emph{forgetting} connected components and restricting the submanifold with boundary 
    \begin{equation}
        W \subseteq M \times \mathsf{cor}(S) \times V
    \end{equation}
    along the inclusion $M \times \mathsf{cor}(T)\times V \hookra M \times \mathsf{cor}(S)\times V$. We denote the restriction by $W_{|T}$, which comes equipped with a smooth embedding 
    \begin{equation}
        W_{|T} \hookra W.
    \end{equation}
    \end{itemize}  
\end{definition}

\begin{remark}
    Two points are worth recording about this definition. First, regarding well-definedness of the functor, the pullback $(f\times \id)^*W$ admits a smooth $(d+\mathsf{dim}(N))$-dimensional embedded submanifold structure (with boundary) in $N\times \mathsf{cor}(T) \times V$ because, in the following diagram
    \begin{equation}
        \begin{tikzcd}[ampersand replacement=\&]
            {(f\times \id)^\ast W} \&\& W \\
            \\
            {N\times \mathsf{cor}(T) \times \RR^{\oo}} \&\& {M\times \mathsf{cor}(T) \times \RR^{\oo}} \\
            \\
            N \&\& M
            \arrow[from=1-1, to=1-3]
            \arrow[hook, from=1-1, to=3-1]
            \arrow["\lrcorner"{anchor=center, pos=0.125}, draw=none, from=1-1, to=3-3]
            \arrow[hook, from=1-3, to=3-3]
            \arrow["{f\times \id}"', from=3-1, to=3-3]
            \arrow[two heads, from=3-1, to=5-1]
            \arrow["\lrcorner"{anchor=center, pos=0.125}, draw=none, from=3-1, to=5-3]
            \arrow[two heads, from=3-3, to=5-3]
            \arrow["f"', from=5-1, to=5-3]
        \end{tikzcd}
    \end{equation}
    when we vertically paste the two pullback squares and inspect them, we note that $f^\ast W = (f\times \mathsf{id})^\ast W$. The dimension condition is satisfied by inspecting the submersivity conditions.

    \medskip

    Additionally, we must check that the pullback submanifold $f^*W$ satisfies the submersivity conditions~(\ref{id: categorical}) and~(\ref{id: spatial}). This follows because the projection $\pi_{|W}$ is a submersion, and submersions are stable under pullback along smooth maps. 
    
    \medskip

    \nid More explicitly, we have the following commutative diagram:
    \begin{equation}
    \begin{tikzcd}[ampersand replacement=\&]
        \& {(f\times \id)^\ast W_{|}} \&\& {W_{|}} \\
        {(f\times \id)^\ast W} \&\& W \\
        \& {N\times T^{\mathsf{c}}\times V} \&\& {M\times T^{\mathsf{c}}\times V} \\
        {N\times\mathsf{cor}(T)\times V} \&\& {M\times\mathsf{cor}(T)\times V} \\
        \& {N\times T^{\mathsf{c}}} \&\& {M\times T^{\mathsf{c}}}
        \arrow[from=1-2, to=1-4]
        \arrow[from=1-2, to=2-1]
        \arrow[color={rgb,255:red,214;green,153;blue,92}, hook, from=1-2, to=3-2]
        \arrow["\lrcorner"{anchor=center, pos=0.125}, draw=none, from=1-2, to=3-4]
        \arrow["\lrcorner"{anchor=center, pos=0.125, rotate=-90}, draw=none, from=1-2, to=4-1]
        \arrow[from=1-4, to=2-3]
        \arrow[hook, from=1-4, to=3-4]
        \arrow["\lrcorner"{anchor=center, pos=0.125, rotate=-90}, draw=none, from=1-4, to=4-3]
        \arrow[from=2-1, to=2-3]
        \arrow[hook, from=2-1, to=4-1]
        \arrow["\lrcorner"{anchor=center, pos=0.125}, draw=none, from=2-1, to=4-3]
        \arrow[hook, from=2-3, to=4-3]
        \arrow[from=3-2, to=3-4]
        \arrow[from=3-2, to=4-1]
        \arrow[color={rgb,255:red,214;green,153;blue,92}, from=3-2, to=5-2]
        \arrow["\lrcorner"{anchor=center, pos=0.125}, draw=none, from=3-2, to=5-4]
        \arrow[hook', from=3-4, to=4-3]
        \arrow[from=3-4, to=5-4]
        \arrow[from=4-1, to=4-3]
        \arrow["{f\times \id}"', from=5-2, to=5-4]
    \end{tikzcd}
\end{equation}
    This diagram renders manifest the fact that the composite of the two orange arrows is a submersion, since it is the pullback of a submersion along a smooth map, as desired.   
\end{remark}

\subsubsection{Digression on interpretation}

Understanding this definition as a precursor to defining an $(\oo,1)$-category of bordisms requires further unpacking. In \Cref{ch2:: subsec: interval_spaces_model_infinity_cats} we explained that, for interval Segal spaces, the space of $k$-composable morphisms is encoded in the functor's value at a choice of interval $T$ with $\gamma(T)=[k]$. So far, as we explained in the introduction to this section, we do not yet have a pre-interval space, let alone an interval Segal space. Instead, we have a co-presheaf in $\Int$ with values in presheaves on $\Man$ --- taking the transpose under the power-mapping adjunction.

Nevertheless, we can already outline the blueprint for how we would encode a sequence of $k$-composable bordisms
\begin{equation}
    M_1 \xra{\Sigma_1} M_2 \xra{\Sigma_2} \cdots \xra{\Sigma_{k-1}} M_k
\end{equation} 
under this setup. For simplicity and grounding, let us examine the case $d=1$. Choose $T = \sU(A,B,C)\in \Int$ so that $\gamma(T) = [2]$; an element of the set $\widetilde{\bord}_{d,V}(\ast, \sU(A,B,C))$ is a $1$-dimensional submanifold with boundary
\[
W\subset \ast\times \core(T)\times V,
\]
which could be represented as 
\begin{equation}
    \scalebox{0.8}{\tikzfig{figures/Ch3/BOR1}{\input{figures/Ch3/BOR1.tikz}}}
\end{equation}
\smallskip

\noindent In which case we have a $1$-dimensional submanifold $W\subset \ast \times \sU(A,B,C) \times V$ where 
\begin{equation}
\sU(A,B,C) = \scalebox{0.6}{\tikzfig{figures/Ch3/BOR2}{\input{figures/Ch3/BOR2.tikz}}}.
\end{equation}
The only non-trivial submersivity condition is (\ref{id: categorical}), and it translates to the projection
\begin{equation}
    \pi \maps \big(W\cap (\ast \times (A\cup B \cup C)\times V)\big) \fib \ast \times (A\cup B \cup C)
\end{equation}
onto the closed intervals $A$, $B$, and $C$ being a submersion. Recall from the boundary version of Ehresmann's fibration theorem, \Cref{ap:: prop: ehresmann boundary}, that proper submersions are smooth fiber bundles. Moreover, in this case the targets are closed intervals which are contractible, so from \Cref{ap:: prop: triviality contractibility} we have that 
\begin{equation}
    W_{|A_1\sqcup A_2 \sqcup A_3} \simeq (W_1 \sqcup W_2 \sqcup W_3) \times [0,1]. 
\end{equation}    
On the other hand, over the intervals $A_i$ (depicted in red) there are no submersivity conditions so the projection restricted to these open subsets can admit any assortment of critical points. 

These red regions, where the submanifold is allowed to change its topology, store the content of morphisms. In this chapter, these data are $d$-dimensional unoriented bordisms. Later in this chapter and, more importantly, in Chapter~5, we will generalize this construction to include additional structure, i.e.\ fiberwise tangential structures and stratifications. 

A few more words are still deserved. We note that the non-shaded regions are effectively storing data for the sources and targets of our morphisms, already with a ``baked-in'' choice of collar.  

We end this digression with a discussion on spatial structure. More concretely, we address how and where this construction encodes the idea that higher morphisms, or invertible morphisms in the $(\oo,1)$-category $\bord_d$, should be diffeomorphisms of $d$-dimensional bordisms, and homotopies should be isotopies between such diffeomorphisms, etc. Paths will correspond to concordances of the point $\ast$, more concretely to sections for $M=\AA^1$. A path is thus a $(d+1)$-dimensional smooth submanifold $\Sigma \subset \AA^1 \times \sU(A,B,C)\times V$, for example 
\begin{equation}
    \tikzfig{figures/Ch3/BOR12}{\input{figures/Ch3/BOR12.tikz}}.
\end{equation}
Condition (\ref{id: categorical}) ensures that the collars behave as expected, and condition (\ref{id: spatial}) enforces the idea that paths of composable sequences of bordisms are diffeomorphisms of the composed bordism which intertwine with the decomposition. This can be seen explicitly by inspecting the fibers over the inner components.

The next few results establish the locality properties of $\widetilde{\bord}_{d,V}$ needed for the construction of our main object. In turn, these are what allow the resulting interval Segal space to have the correct symmetric monoidal homotopy category. 

\begin{proposition}\label{ch3:: prop: sheaf condition}
    The adjoint functor 
    \begin{equation}
    \widetilde{\bord}_{d,V}(-,-)\maps \Int \longra \Fun(\man^{\op},\sset)
    \end{equation}
    factors through $\shv(\man)$, the category of $\oo$-sheaves on $\man$ relative to the Grothendieck topology of open covers. In other words, for every $T\in \Int$ the functor
    \[
    \widetilde{\bord}_{d,V}(-,T)\maps \man^{\op}\longra \set
    \]
    is a smooth $\oo$-sheaf of (discrete) simplicial sets.
\end{proposition}

\begin{proof}
    We will check conditions (\ref{binary cover}) and (\ref{disj union}) from \Cref{ch3:: thm: BBEP-sheaf-criterion}. We begin with (\ref{binary cover}): let $M$ be a smooth manifold and let $\{U_1,U_2\}$ be an open cover of $M$. We must show that the \emph{decomposition} map
    \begin{equation}
        \label{eq1}
        \widetilde{\bord}_{d,V}(M,T) \longrightarrow \widetilde{\bord}_{d,V}(U_1,T)\underset{\widetilde{\bord}_{d,V}(U_1\cap \,U_2,T)}{\times}\widetilde{\bord}_{d,V}(U_2,T)
    \end{equation}
    is an isomorphism. Injectivity is almost immediate: given a pair of different submanifolds with boundary 
    \[
    W,W' \subset \big(M\times \core(T)\times V\big),
    \]
    their images under the decomposition map \eqref{eq1} must differ, manifestly, in at least one of the factors of the pullback.
    
     The surjectivity of the decomposition map is downstream of the fact that our construction works with submanifolds embedded into some larger universe, as opposed to a set (or even space) of manifolds --- in this case $M \times \core(T) \times V$. An element of the codomain of \eqref{eq1} is thus equivalent to the datum of
    \begin{equation}
    \big(W_1 \subset U_1 \times \mathsf{cor}(T)\times V\big), \quad \big(W_2 \subset U_2 \times \mathsf{cor}(T)\times V \big),
    \end{equation} 
    a pair of submanifolds with boundary satisfying conditions (\ref{id: categorical}) and (\ref{id: spatial}) of \Cref{ch3:: def: bord_d}, whose restrictions along $\big(U_1 \cap U_2\big)\times \core(T)\times V$ agree pointwise --- both equal to the submanifold
    \begin{equation}
    W_{|U_{1}\cap U_{2}}\subset \big( U_1\cap U_2\big) \times \mathsf{cor}(T)\times V.
    \end{equation}
    To construct an element in $\widetilde{\bord}_{d,V}(M,T)$ which decomposes to this data, we glue, i.e.\ take the following pushout
    \begin{equation}
        \begin{tikzcd}
            {W_{|U_{1}\cap\,U_{2}}} && {W_{2}} \\
            \\
            {W_{1}} && {W_{1}\,\underset{W_{|U_{1}\cap\,U_{2}}}{\cup} \, W_2}
            \arrow[hook, from=1-1, to=1-3]
            \arrow[hook, from=1-1, to=3-1]
            \arrow[hook, from=1-3, to=3-3]
            \arrow[hook, from=3-1, to=3-3]
            \arrow["\lrcorner"{anchor=center, pos=0.125, rotate=180}, draw=none, from=3-3, to=1-1]
        \end{tikzcd}
    \end{equation}
    in $\man$. This pushout restricts to $W_1$ and $W_2$ on the open subsets $U_1$ and $U_2$, respectively, and satisfies conditions (\ref{id: categorical}) and (\ref{id: spatial}). The latter follows from the fact that submersion is an open condition, and the pushout was carried out along an open cover. 

    Condition (\ref{disj union}) from \Cref{ch3:: thm: BBEP-sheaf-criterion} is checked similarly, using the fact that disjoint unions of manifolds correspond to coproducts in $\man$ and that the construction is taken componentwise.
\end{proof}
\begin{remark}
    The simplicial set $\bord_{d,V}(\AA^{\bullet},T)$ is a $0$-weak Kan complex for every $T\in \Int$. This follows from the general fact that for every concordance-invariant sheaf $\cF$, the simplicial set 
    \begin{equation}
        \big | \cF(\AA^{\bullet} \times *)\big | \in \sset
    \end{equation}
    is a $0$-weak Kan complex.
\end{remark}

Having entertained the locality of our construction with respect to the smooth manifold variable, we now assert yet another form of locality; again fixing $T \in \Int$, the next proposition proves a ``naive sheaf'' property with respect to binary open covers of the closed interval $[0,1]$ itself.

\begin{proposition}[Naive sheaf condition in the interval direction]\label{ch3:: prop: naive sheaf con}
Let $T\in \Int$, and let $\{U,V\}$ be an open cover of $\core(T)\subset [0,1]$. For each $k\geq 0$, the canonical map
\begin{equation}
\label{eq:naive sheaf interval levelwise}
\bord_d\big(\AA^{k},T\big)\longrightarrow
\bord_d\big(\AA^{k},T\big)_{|U}\underset{\bord_d\big(\AA^{k},T\big)_{|U\cap V}}{\times}
\bord_d\big(\AA^{k},T\big)_{|V}
\end{equation}
is a bijection of sets. Equivalently, the simplicial set $\bord_d(\AA^{\bullet},T)$ is the strict pullback
\begin{equation}
\label{eq:naive sheaf interval simplicial}
\bord_d(\AA^{\bullet},T)\;\cong\;
\bord_d(\AA^{\bullet},T)_{|U}\underset{\bord_d(\AA^{\bullet},T)_{|U\cap V}}{\times}
\bord_d(\AA^{\bullet},T)_{|V},
\end{equation}
where
\begin{equation}
\label{eq:bord-restriction-U}
\bord_d(\AA^{k},T)_{|U}
=
\Big\{
W\subseteq \AA^{k}\times\big(U\cap \core(T)\big)\times V
\Big\}.
\end{equation}
Here $\bord_d(\AA^{k},T)_{|U}$ denotes the set of all $(d+k)$-dimensional
submanifolds with boundary
\[
W\subseteq \AA^{k}\times\big(U\cap \core(T)\big)\times V
\]
which satisfy conditions (\ref{id: dim}),(\ref{id: categorical}), and (\ref{id: spatial}) from \Cref{ch3:: def: bord_d}, and the restriction maps are given by intersecting with
\[
\AA^{k}\times\big(U\cap V\cap \core(T)\big)\times V\subseteq 
\AA^{k}\times\big(U\cap\core(T)\big)\times V,
\]
and similarly for $V$.
\end{proposition}

\begin{proof}
Fix $k\geq 0$ and abbreviate
\begin{equation}
Z = \AA^{k}\times \core(T)\times V.
\end{equation}
An element of $\bord_d(\AA^{k},T)$ is a $(d+k)$–dimensional submanifold with boundary
\begin{equation}
W\subseteq Z
\end{equation}
satisfying conditions (\ref{id: dim})--(\ref{id: spatial}) of \Cref{ch3:: def: bord_d}. Its images in
$\bord_d(\AA^{k},T)_{|U}$ and $\bord_d(\AA^{k},T)_{|V}$ are given by intersection with
$\AA^{k}\times\big(U\cap \core(T)\big)\times V$ and
$\AA^{k}\times\big(V\cap \core(T)\big)\times V$, respectively, so injectivity of
\eqref{eq:naive sheaf interval levelwise} is immediate.

For surjectivity, start with a pair of submanifolds
\begin{equation}
W_U\subseteq \AA^{k}\times\big(U\cap\core(T)\big)\times V,\qquad
W_V\subseteq \AA^{k}\times\big(V\cap\core(T)\big)\times V,
\end{equation}
each satisfying (\ref{id: dim})--(\ref{id: spatial}), whose restrictions to
\begin{equation}
\AA^{k}\times\big(U\cap V\cap\core(T)\big)\times V
\end{equation}
coincide as subsets. Define
\begin{equation}
W := W_U \cup W_V \subseteq Z.
\end{equation}
Since $U$ and $V$ cover $\core(T)$ and $W_U$, $W_V$ agree on the overlap, $W$ is a well-defined subset of $Z$.

Standard differential-topological arguments ensure that $W$ is a smooth
submanifold with boundary: all the structure is determined locally, and at each point of $W$ either $W$ coincides with $W_U$, or with $W_V$, or with their common restriction in the overlap. Condition (\ref{id: dim}) is local, so it also holds for $W$.

The submersivity conditions (\ref{id: categorical}) and (\ref{id: spatial}) are likewise local on the source: they hold on $W_U$ and $W_V$ by assumption, hence they hold on $W$, since submersion is an open condition. This shows that $W$ lies in $\bord_d(\AA^{k},T)$ and restricts to $(W_U,W_V)$, proving surjectivity and hence the bijectivity of \eqref{eq:naive sheaf interval levelwise}. The statement for the simplicial set follows levelwise.
\end{proof}

We now come to the main idea of this section. We will use the concordance resolution of sheaves, restricted to the smooth cosimplicial object, in order to construct an interval Segal space. For each interval object $I$, the functor
\begin{equation}
    \widetilde{\bord}_{d,V}(-,I) \maps \man^{\op}\longra \mathsf{Set}
\end{equation}
is set-valued, so evaluating it on the smooth cosimplicial simplex already produces a simplicial set
\begin{equation}
    \widetilde{\bord}_{d,V}(\AA^\bullet,I).
\end{equation}
Thus, in the present case, the concordance-resolution construction introduces no additional simplicial direction beyond this one. The next result is both a definition and a theorem.

\medskip

\begin{theorem}\label{ch3:: thm: interval segal bordism cat}
    Let $d \geq 0$ be a non-negative integer and let $V$ be a finite-dimensional vector space of sufficiently high dimension. The composite functor
    \begin{align}
        \bord_{d,V} \maps \Int &\longra \shv(\man)  \longra \sset  \\
         I &\longmapsto \widetilde{\bord}_{d,V}(-,I) \longmapsto \bB\widetilde{\bord}_{d,V}(*,I) = \big|\,\widetilde{\bord}_{d,V}(\AA^\bullet,I)\,\big| 
    \end{align}
    is an interval Segal space.
\end{theorem}

\smallskip

\begin{proof}
    \textbf{Homotopical property.}
    We begin by checking that the functor sends $\pi_0$-equivalences to weak equivalences in $\sset$. Concretely, for every open embedding $T\hookra S$ such that $\pi_0(T) \simeq \pi_0(S)$, the associated map of simplicial sets 
    \begin{equation}
        \bord_{d,V}(T) \xra{\,\,\,\,\mathsf{forget}_{T \ra S}\,\,\,\,} \bord_{d,V}(S),
    \end{equation}   
    is a weak equivalence. There is a map
    \[
    \mathsf{ext}\maps \bord_{d,V}(S) \longra \bord_{d,V}(T)
    \]
    which ``squeezes'' the underlying submanifold with boundary, in particular bunching up critical points. Let 
    \begin{equation}
        T = U_1 \cup \cdots \cup U_k, \quad S = \widetilde{U}_1 \cup \cdots \cup \widetilde{U}_k,
    \end{equation} 
    with $0\in U_1,\widetilde{U}_1$ and $1\in U_k,\widetilde{U}_k$ --- note that $U_i \subseteq \widetilde{U}_i$. Let $s\maps [0,1]\longra [0,1]$ be a diffeomorphism which maps $\widetilde{U}_i \mapsto U_i$ and is the identity elsewhere; essentially a rescaling map obtained by smoothing a piecewise-linear map. We now define
    \begin{align*}
        \bord_{d,V}(S) &\xra{\,\,\,\mathsf{ext}_{S \longra T}\,\,\,} \bord_{d,V}(T) \\
        \big( W \subset \AA^\bullet \times [0,1]\times V \big) &\longmapsto \big((\id_{\AA^\bullet} \times\, s \times \id_{V})^{*}W\subset \AA^\bullet \times [0,1]\times V \big).
    \end{align*}
    The map is manifestly well-defined since $(\id_{\AA^\bullet} \times\, s \times \id_{V})$ is a diffeomorphism, i.e.\ it preserves condition (\ref{id: dim}) and the submersivity conditions (\ref{id: categorical}) and (\ref{id: spatial}) over $T$. 

    Neither of the two composites is strictly the identity, so we must articulate appropriate simplicial homotopies witnessing the desired equivalence. To do so, we observe that a simplicial homotopy 
    \begin{equation}
        \bord_{d,V}(S)\times \Delta^1 \longra \bord_{d,V}(S)
    \end{equation} 
    amounts, via adjunction, to a simplicial map 
    \begin{equation}
        \widetilde{\bord}_{d,V}(\AA^\bullet,S) \longra \widetilde{\bord}_{d,V}(\AA^\bullet \times \AA^1, S)
    \end{equation}
    whose further restriction to $\{0\}$ and $\{1\}$ yields maps identical to $\big(\mathsf{forget}\circ\mathsf{ext}\big)$ and $\id$, respectively. We construct this map by pulling back along 
    \begin{align*}
        \big(\AA^\bullet \times (\AA^1 \times [0,1]) \times V\big) \xra{\id_{\AA^\bullet}\times h \times \id_{V}} \big( \AA^\bullet \times  [0,1]\times V \big),
    \end{align*}
    where $h\maps \AA^1 \times [0,1] \longra [0,1]$, defined as $h(t,x) = (1-t) \cdot s(x) + t \cdot x$, is a smooth linear isotopy between the diffeomorphism $s$ and the identity --- whose graph we depict in the figure below: 
    \begin{equation}
        \scalebox{0.8}{\tikzfig{figures/Ch3/BOR14}{\input{figures/Ch3/BOR14.tikz}}}
    \end{equation}
    Conditions (\ref{id: dim}), (\ref{id: categorical}), and (\ref{id: spatial}) are preserved since the isotopy is a diffeotopy, which ensures a well-defined simplicial map. 
    
    To establish that $(\mathsf{ext}\circ \mathsf{forget}) \simeq \id_{\bord_{d,V}(T)}$ we construct the analogous simplicial homotopy 
    \begin{equation}
    \widetilde{\bord}_{d,V}(T)\times \Delta^1 \longra \widetilde{\bord}_{d,V}(T)
    \end{equation}
    this time using $s^{-1}(x)$ in place of $s(x)$ when describing the homotopy $h$.

    \medskip

    \textbf{Interval Segal condition.}
    Recall that the interval Segal condition (\Cref{ch2:: def: interval Segal condition}) states that the following diagram of simplicial sets\footnote{These happen to be $0$-weak Kan complexes, and not entirely arbitrary simplicial sets.} is a homotopy pullback square for every $\sU(A,B,C)\in \cJ_2$:
    \begin{equation}\label{eq:Segal condition for bordism}
        \begin{tikzcd}[ampersand replacement=\&]
            {\bord_d\big(U(A,B,C)\big)} \&\& {\bord_d\big(U(A,B)\big)} \\
            \\
            {\bord_d\big(U(B,C)\big)} \&\& {\bord_d\big(U(B)\big)}
            \arrow[from=1-1, to=1-3]
            \arrow[from=1-1, to=3-1]
            \arrow["h"{description, pos=0.3}, "\lrcorner"{anchor=center, pos=0.125}, draw=none, from=1-1, to=3-3]
            \arrow[from=1-3, to=3-3]
            \arrow[from=3-1, to=3-3]
        \end{tikzcd}
    \end{equation}
     Proving that the diagram above is a homotopy pullback square is equivalent to showing that the following diagram is a homotopy pullback square in $\sset$: 
    \begin{equation}\label{eq pullback 1}
        \begin{tikzcd}[ampersand replacement=\&]
            {\bord_d\big(U(A,B,C)\big)} \&\& {\bord_d\big(U(A,B)\big)\times \bord_d\big(U(B,C)\big)} \\
            \\
            {\Map_{\sset}\Big(\Delta^1,\bord_d\big(U(B)\big)\Big)} \&\& {\Map_{\sset}\Big(\partial\Delta^1,\bord_d\big(U(B)\big)\Big)}
            \arrow[from=1-1, to=1-3]
            \arrow[from=1-1, to=3-1]
            \arrow["h"{description, pos=0.3}, "\lrcorner"{anchor=center, pos=0.125}, draw=none, from=1-1, to=3-3]
            \arrow[from=1-3, to=3-3]
            \arrow[from=3-1, to=3-3]
        \end{tikzcd}
    \end{equation}
    Note that all we did was apply a Kan fibration replacement. We now take advantage of the geometric nature of our construction to re-express the homotopy-theoretic path fibration using concordances: namely, we use the fact that, for all $T\in \Int$, the $\oo$-sheaf 
    \begin{equation}
    |\widetilde{\bord}_{d,V}(-\times \AA^{\bullet},T)| \maps \man^{\op} \longra \sset 
    \end{equation}
    is concordance-invariant, i.e.\ for every $M$ the projection
    \[
    |\widetilde{\bord}_{d,V}(M\times \AA^{\bullet},T)| \xra{\,\,\simeq\,\,} |\widetilde{\bord}_{d,V}(M\times \AA^{\bullet}\times \AA^1,T)|
    \]
    induces a weak equivalence of simplicial sets. Thereby, paths in $\bord_d(T)$ are equivalent to concordances in $\bord_d(T)$. We can extend the original diagram~(\ref{eq pullback 1}) by a homotopy equivalent one as follows:
    \begin{equation}
        \begin{tikzcd}[ampersand replacement=\&]
            {\bord_d\big(U(A,B,C)\big)} \&\& {\bord_d\big(U(A,B)\big)\times \bord_d\big(U(B,C)\big)} \\
            \\
            {\mathsf{Map}_{\sset}\Big(\Delta^1,\bord_d\big(U(B)\big)\Big)} \&\& {\mathsf{Map}_{\sset}\Big(\partial\Delta^1,\bord_d\big(U(B)\big)\Big)} \\
            \\
            {|\widetilde{\bord}_{d,V}\big(\AA^{\bullet}\times \AA^1,U(B)\big)|} \&\& {\bord_d\big(U(B)\big)\times \bord_d\big(U(B)\big)}
            \arrow[from=1-1, to=1-3]
            \arrow[from=1-1, to=3-1]
            \arrow["\lrcorner"{anchor=center, pos=0.125, rotate=45}, draw=none, from=1-1, to=3-3]
            \arrow[from=1-3, to=3-3]
            \arrow[from=3-1, to=3-3]
            \arrow["\simeq"', from=3-1, to=5-1]
            \arrow[equals, from=3-3, to=5-3]
            \arrow["{\text {weak Kan fibration}}"', from=5-1, to=5-3]
        \end{tikzcd}
    \end{equation} 
    where the left bottom vertical arrow is a weak equivalence precisely due to the concordance-invariance of $\bord_d\left(\sU(B)\right)$. This means that the homotopy limit of the outer span is weakly equivalent to the homotopy limit of the inner span
    \begin{equation}
        \begin{tikzcd}[ampersand replacement=\&]
            {\mathsf{holim}\Big[\bord_d\big(U(A,B)\big)\times \bord_d\big(U(B,C)\big) \ra \bord_d (U(B))^{\times 2} \la \bord_d\big(\AA^{\bullet}\times \AA^1,U(B)\big)\Big]} \\
            \\
            {\mathsf{holim}\Big[\bord_d\big(U(A,B)\big)\times \bord_d\big(U(B,C)\big) \ra \bord_d (U(B))^{\times 2}\la \mathsf{Map}_{\sset}\Big(\Delta^1, \bord_d\big(U(B)\big)\Big)\Big]}.
            \arrow[from=1-1, to=3-1]
        \end{tikzcd}
    \end{equation}
    Therefore, the interval Segal condition (\ref{eq:Segal condition for bordism}) holds if and only if the following map is a weak equivalence:
    \begin{equation}\label{ch3:: eq: segal condition 1}
        \bord_{d}\big(U(A,B,C)\big) \xra{\,\, \simeq \,\,} \Big(\bord_d\big(U(A,B)\big)\times \bord_d\big(U(B,C)\big)\Big) \underset{\bord_d\big(U(B)\big)\times \bord_d\big(U(B)\big)}{\times^h} \bord_d\big(\AA^{\bullet}\times \AA^1,U(B)\big).
    \end{equation}

    Using the fact (proved in \Cref{ch3:: prop: weak-kan-endpoints}) that for concordance-invariant $\oo$-sheaves, the endpoint map 
    \begin{equation}
    |\widetilde{\bord}_{d,V}(\AA^\bullet\times \AA^1,U(B))| \longra |\widetilde{\bord}_{d,V}(\AA^\bullet \times \partial \AA^1, U(B))|
    \end{equation} 
    is a weak Kan fibration, together with the fact that weak Kan fibrations are realization fibrations (\Cref{ch3:: thm: weak-kan-realization}), we see that showing the interval Segal condition is equivalent to showing that the following map is a weak equivalence of simplicial spaces with the injective model structure:
    \begin{equation}
    \widetilde{\bord}_{d,V}\big(\AA^{\bullet}, U(A,B,C)\big) \xra{\,\, \simeq \,\,} \Big(\widetilde{\bord}_{d,V}\big(\AA^{\bullet},U(A,B)\big)\times \widetilde{\bord}_{d,V}\big(\AA^{\bullet},U(B,C)\big)\Big) \underset{\widetilde{\bord}_{d,V}\big(\AA^{\bullet},U(B)\big)^{\times 2}}{\times^h} \widetilde{\bord}_{d,V}\big(\AA^{\bullet}\times \AA^1,U(B)\big).
    \end{equation}  

    \noindent Alternatively, the right-hand side of the above diagram can be expressed as the following iterated (strict) fiber product:
    \begin{equation}
        \label{eq:iterated fiber product for Segal condition 1}
        \widetilde{\bord}_{d,V}\big(\AA^{\bullet},U(A,B)\big) \underset{\widetilde{\bord}_{d,V}\big(\AA^\bullet,U(B)\big)}{\times} \widetilde{\bord}_{d,V}\big(\AA^{\bullet}\times \AA^1,U(B)\big) \underset{\widetilde{\bord}_{d,V}\big(\AA^\bullet,U(B)\big)}{\times} \widetilde{\bord}_{d,V}\big(\AA^{\bullet},U(B,C)\big).
    \end{equation}
    
    \medskip

    \nid As observed in \Cref{ch2:: obs: discreteness}, weak equivalences between discrete simplicial spaces amount to levelwise set bijections. To show that \eqref{eq:iterated fiber product for Segal condition 1} is weakly equivalent to $\widetilde{\bord}_{d,V}\big(\AA^\bullet,U(A,B,C)\big)$, we adapt the proof of \cite[Theorem~5.4]{BBEP24} to our setting, using the naive sheaf condition from \Cref{ch3:: prop: naive sheaf con}. Concretely, given an open cover $\{U,V\}$ of $[0,1]$, the simplicial set $\bord_d(\AA^\bullet,T)$ is the strict pullback of the diagram
    \begin{equation}
        \bord_d(\AA^\bullet,T)_{|U} \longrightarrow \bord_d(\AA^\bullet,T)_{|U\cap V} \longleftarrow \bord_d(\AA^\bullet,T)_{|V},
    \end{equation} 
    where the vertical subscript bars denote the levelwise restriction of the simplicial set to the corresponding open subset of the interval $[0,1]$. 

    With this in mind, let us pick a cover $\cU = \{U,V\}$ of $[0,1]$ such that $U$ covers $(A\cup B)$ and $V$ covers $(B\cup C)$, which we depict as
    \begin{equation}
        \tikzfig{figures/Ch3/BOR9}{\input{figures/Ch3/BOR9.tikz}}
    \end{equation}
    such that the red shaded region and the blue shaded region cover the opens $U$ and $V$, respectively. Additionally, we define a partition of unity $\{f_{U},f_{V}\}$ subordinate to $\cU$. Define smaller opens $U' = f_{U}^{-1}[0,1/3)$ and $V' = f_V^{-1}(2/3,1]$, noting that, since $f_U + f_V \equiv 1$, they do not overlap, $U' \cap V' = \emptyset$, and thus the cover $\cU' = \{U', U\cap V, V'\}$ refines $\cU$. There is an induced map of simplicial sets
    \begin{equation}
    \begin{tikzcd}
        \bord_d\big(\AA^{\bullet},U(A,B)\big)     
        \underset{\bord_d (\AA^\bullet,U(B))}{\times^{h}} 
        \bord_d\big(\AA^{\bullet}\times \AA^1,U(B)\big) 
        \underset{\bord_d (\AA^\bullet,U(B))}{\times^{h}} 
        \bord_d\big(\AA^{\bullet},U(B,C)\big) \ar[dd, "\mathsf{res}"] \\ \\
        \bord_d\big(\AA^{\bullet},U(A,B)\big)_{|U'}     
        \underset{\bord_d (\AA^\bullet,U(B))_{|U'}}{\times^{h}} 
        \bord_d\big(\AA^{\bullet}\times \AA^1,U(B)\big) 
        \underset{\bord_d (\AA^\bullet,U(B))_{|V'}}{\times^{h}} 
        \bord_d\big(\AA^{\bullet},U(B,C)\big)_{|V'} 
    \end{tikzcd}
    \end{equation}
    which is a weak equivalence, by applying \Cref{ch3:: prop: naive sheaf con} to the covers $\{U',U\cap V\}$ of $U$ and $\{V',U\cap V\}$ of $V$.

    We now construct the following additional pair of simplicial maps:
    \begin{equation}
    \begin{tikzcd}
        \bord_d\big(\AA^{\bullet},U(A,B)\big)     
        \underset{\bord_d\big(\AA^\bullet,U(B)\big)}{\times} 
        \bord_d\big(\AA^{\bullet}\times \AA^1,U(B)\big) 
        \underset{\bord_d\big(\AA^\bullet,U(B)\big)}{\times} 
        \bord_d\big(\AA^{\bullet},U(B,C)\big) \ar[dd, "\mathsf{res}"] \\ \\
        \bord_d\big(\AA^{\bullet},U(A,B)\big)_{|U'}     
        \underset{\bord_d\big(\AA^\bullet,U(B)\big)_{|U'}}{\times} 
        \bord_d\big(\AA^{\bullet}\times \AA^1,U(B)\big) 
        \underset{\bord_d\big(\AA^\bullet,U(B)\big)_{|V'}}{\times} 
        \bord_d\big(\AA^{\bullet},U(B,C)\big)_{|V'} \ar[dd, bend right=25, "F"', shorten <=2pt, shorten >=2pt]   \\ \\ 
        \bord_d\big(\AA^{\bullet},U(A,B)\big)_{|U'}     
        \underset{\bord_d\big(\AA^\bullet,U(B)\big)_{|U'}}{\times} 
        \bord_d\big(\AA^{\bullet},U(B)\big) 
        \underset{\bord_d\big(\AA^\bullet,U(B)\big)_{|V'}}{\times} 
        \bord_d\big(\AA^{\bullet},U(B,C)\big)_{|V'} \ar[uu, bend right=25, "\pr^\ast"']      
    \end{tikzcd}
    \end{equation}
    about which we now expand:
    \begin{enumerate}
        \item As mentioned above, the map $\mathsf{res}$ is an isomorphism of simplicial sets, by \Cref{ch3:: prop: naive sheaf con}.
        \item Crucially, the induced map from $\bord_d\big(\AA^\bullet, U(A,B,C)\big)$ to the fiber product in the bottom of the diagram is also an isomorphism of simplicial sets, again using \Cref{ch3:: prop: naive sheaf con} with the refined cover $\{U',U\cap V,V'\}$.
        \item Lastly, we must articulate why $F$ is a homotopy inverse to $\pr^\ast$. There is a smooth map $c\maps \AA^1 \longra \AA^1$ which is strictly increasing and satisfies the following additional conditions: 
        \begin{enumerate}
        \item $c_{|(-\oo,1/3)} \equiv 0$, and 
        \item $c_{|(2/3,+\oo)} \equiv 1$. 
        \end{enumerate}
        Composing $c \circ (f_{V})_{|U\cap V}\maps U\cap V \longra \AA^1$ we get a function that vanishes on $U'\subset U\cap V$ and is constantly equal to $1$ on $V'\subset U\cap V$; we call this map $f$, which we depict in the following blue graph, as an aid to the reader:
        \begin{equation}
            \scalebox{0.8}{\tikzfig{figures/Ch3/BOR10}{\input{figures/Ch3/BOR10.tikz}}}
        \end{equation}
        The map $F$ sends a smooth submanifold 
        \[
        W\subset (\AA^\bullet \times \AA^1)\times U(B)\times V
        \]
        to the pullback along 
        \begin{equation}
        \id_{\AA^{\bullet}}\times f \times \id_{V}\maps \AA^{\bullet}\times U(B) \times V \longra \AA^{\bullet}\times \AA^1 \times U(B) \times V,
        \end{equation}
        obtaining a smooth $(d+\bullet + 1)$-dimensional submanifold $\big(W_f \subset \AA^\bullet\times U(B)\times V\big)$. For intuition, we sketch its action in dimension $d=1$:
        \begin{equation}
            \!\!\!\!\scalebox{0.8}{\tikzfig{figures/Ch3/BOR11}{\input{figures/Ch3/BOR11.tikz}}} \rightsquigarrow \scalebox{0.75}{\tikzfig{figures/Ch3/BOR12}{\input{figures/Ch3/BOR12.tikz}}}
        \end{equation}
        On the left we drew a section, i.e.\ a $(d+1+\bullet)$-dimensional manifold admitting a submersion onto $\AA^{1}\times U(B)$. On the right we depict the interpolating submanifold induced by $f$. Lastly, the next figure captures $(\id_{\AA^\bullet} \times f \times \id_V)$ with the orange graph; we leave the fibers over $U' \cap U(B)$ and $U(B)\cap V'$ in the original green so as to highlight how $W_f$ interpolates the two:
        \begin{equation}
            \scalebox{0.8}{\tikzfig{figures/Ch3/BOR13}{\input{figures/Ch3/BOR13.tikz}}}
        \end{equation}
        This figure is meant to elucidate this map's well-definedness --- concretely, the fact that it preserves the source at $U'$ and the target at $V'$:
         \begin{equation}
        (W_f)_{|U'} = (\mathsf{source})^\ast(W)_{|U'}, \quad
        (W_f)_{|V'} = (\mathsf{target})^\ast(W)_{|V'}, 
        \end{equation} 
        which is a direct consequence of how $f$ is defined as a constant map on $0$ and $1$ over $U'$ and $V'$, respectively. Observe that $(\pr^\ast W)_f = (W \times \AA^1)_f =  W$, stemming from the fact that 
        \begin{equation}
            \AA^\bullet \times U(B)\times V \xra{\,\id \times f \times \id\,} \AA^\bullet \times \AA^1 \times U(B)\times V \xra{\,\pr\,} \AA^\bullet \times U(B) \times V = \id,
        \end{equation}
        thus proving that $F \circ \pr^\ast = \id$. Intuitively, this map sends a $d$-dimensional submanifold $W$ to its corresponding cylinder over $\AA^1$ and then applies $F$, recovering $W$. 
        \item Finally, we must prove that $\pr^* \circ F \simeq \id$. The composite $\pr^* \circ F$ sends a submanifold 
        \begin{equation}
            \big(W \subset (\AA^\bullet \times \AA^1) \times U(B) \times V \big) \longmapsto \big((W_f \times \AA^1) \subset (\AA^\bullet \times \AA^1) \times U(B) \times V \big).
        \end{equation} 
        We are interested in finding a simplicial homotopy 
        \begin{equation}
            \eta \maps \bord_d(\AA^\bullet \times \AA^{1}, U(B)) \times \Delta^1 \longra \bord_d(\AA^\bullet \times \AA^1 , U(B))
        \end{equation}
        from $\pr^\ast \circ F$ to $\id$. 
        We first apply the adjunction to reframe the homotopy as a simplicial map to the path space  
        \begin{equation}
        \begin{tikzcd}[ampersand replacement=\&]
            {\bord_d(\AA^{\bullet}\times \AA^{1},U(B))} \&\& {\mathsf{Map}_{\sset}\big(\Delta^1,\bord_d(\AA^{\bullet}\times \AA^{1},U(B))\big)} \\
            \\
            \&\& {\bord_d\big((\AA^{\bullet}\times \AA^{1})\times \AA^1,U(B)\big)}
            \arrow["\eta", from=1-1, to=1-3]
            \arrow["(*)"', curve={height=18pt}, from=1-1, to=3-3]
            \arrow["\simeq", from=1-3, to=3-3]
        \end{tikzcd}
        \end{equation}
        and then, using concordance-invariance, we get that the homotopy is ultimately equivalent to a simplicial map to $\bord_d\big((\AA^{\bullet}\times \AA^{1})\times \AA^1,U(B)\big)$. Let us consider the smooth map
        \begin{equation}
            \id_{\AA^\bullet} \times \,H \times \id_V \maps \big((\AA^\bullet \times \AA^1) \times \AA^1 \times U(B) \times V \big) \longra \big((\AA^\bullet \times \AA^1)  \times U(B) \times V \big) 
        \end{equation}
        where $H\maps (U(B)\times \AA^1) \times \AA^1 \longra U(B) \times \AA^1$ is defined as the smooth homotopy 
        \begin{equation}
            H\big((x,t),s\big) = \big(x, h\big((x,t),s)\big)
        \end{equation}
        via the smooth linear interpolation
        \begin{align}
        h\big((x,t),s\big) &= (1-s)\cdot f \circ \pr_{U(B)}(x,t) + s \cdot \pr_{\AA^1}(x,t) \\
        &= (1-s)\cdot f(x) + s\cdot t.
        \end{align}
        The smooth homotopy $H$ starts at $H\big((x,t),0\big) = (x,f(x))$ and ends at $H\big((x,t),1\big) = (x,t)$, i.e.\ the identity. Now we check that the diagonal map $(*)$ satisfies the properties we are after: when we ``evaluate'' the homotopy at $\{i\}$, by pulling back along $\AA^0 \xra{\langle i \rangle^*} \AA^1$, we have that:
        \begin{equation}
            \big(\id\times H \times \id\big)((u,i),t,x,v) = \begin{cases}
                (u,f(x),x,v) & \text{ if } i=0, \\
                (u,t,x,v) & \text{ if } i = 1.
            \end{cases}
        \end{equation}
        For $i=1$ we manifestly get the identity, as desired, and for $i=0$ we have 
        \begin{equation}
        (\id \times \,\pr\times \id) \circ(\id_{\AA^\bullet}\times \,f\times \id_V),
        \end{equation}
        thus finishing the proof.  
    \end{enumerate}
\end{proof}

\begin{observation}
We pause here to reflect on the difference between this proof of the Segal condition and that of Calaque--Scheimbauer \cite[Proposition 5.19]{CS19}. Note first that we dispense with the need to pass through various subspaces of the homotopy fiber product, and, more importantly, we get to use the intrinsic geometric nature of our construction as the conduit for simplicial homotopy theory. This feature, we believe, may pay dividends in further applications to manifold topology and topological field theories. In a way, we see our construction as a means to discretize their $(\oo,1)$-category of intervals, thus simplifying it further.
\end{observation}

The following result characterizes the structural conditions that any model of the bordism category must satisfy. To make sense of these, we must first introduce some less well-studied Lie groups first. As far as we are aware, they were first defined in \cite{GMTW09}. 

\begin{definition}\label{ch3:: def: iso bordism}
 Let $(M,h_0,h_1)$ be a $d$-dimensional smooth unoriented bordism from $\Sigma_0$ to $\Sigma_1$ where we have collars
 \begin{align}
    h_0 \maps \Sigma_0 \times [0,1) &\hookra M \\
    h_1 \maps \Sigma_1 \times (0,1] &\hookra M
 \end{align}
 where $\partial M$ is the image of $h_{i}(\{i\}\times \Sigma_i)$ for $i=0,1$. Let $\Diff(\{M,h_0,h_1\})$ be the subgroup of diffeomorphisms of $M$ which must satisfy the following additional condition: 
 \begin{itemize}
    \item[($*$)] There exists a diffeomorphism $[0,1) \longra [0,1)$ making the following diagram commute:
    \begin{equation}
        \begin{tikzcd}
            {h_0(\Sigma_0\times [0,1))} && {\varphi\big(h_0(\Sigma_0\times [0,1))\big)} \\
            \\
            {[0,1)} && {[0,1)}
            \arrow["\varphi", from=1-1, to=1-3]
            \arrow[from=1-1, to=3-1]
            \arrow[from=1-3, to=3-3]
            \arrow[dashed, from=3-1, to=3-3]
        \end{tikzcd}
    \end{equation} 
    with an analogous condition for $h_1$.
 \end{itemize}
\end{definition}

\begin{observation}
    It is instructive to think of the group $\Diff(\{M,h_0,h_1\})$ as the subgroup of diffeomorphisms of $M$ which are a product of diffeomorphisms near the boundary $\partial M$.
\end{observation}

We want to not only consider bordisms, as they are interpreted as 1-morphisms, but also chains of composable bordisms. Moreover, we want to study the space of such composable chains, and so we must, in some fashion, generalize the diffeomorphism group in order to accommodate these additional blocks of data.

\begin{definition}\label{ch3:: def: diffeo group}
    Let $(M,\{h_i\}_{i=0}^{k})$ be a composed chain of $k$-composable $d$-dimensional unoriented bordisms 
    \begin{equation}
        \Sigma_0 \xra{M_1} \Sigma_1 \xra{M_2} \cdots \xra{M_{k}} \Sigma_k 
    \end{equation}
    together with the accompanying collaring data 
    \begin{align}
        h_0 \maps [0,1) \times &\Sigma_0 \hookra M \\
        h_i \maps (0,1) \times &\Sigma_i \hookra M \\
        h_k \maps (0,1]  \times &\Sigma_k \hookra M.
    \end{align}
    Let $\Diff(M,\{h_i\}_{i=0}^{k})$ denote the subgroup of diffeomorphisms of $M$ which satisfy the additional set of conditions:
    \begin{itemize}
        \item[$(*)$] There exist diffeomorphisms $(0,1) \longra (0,1)$ making the following diagram commute:
        \begin{equation}
            \begin{tikzcd}
                {h_i(M_{i}\times(0,1))} && {\varphi\big(h_i(M_{i}\times(0,1))\big)} \\
                \\
                {(0,1)} && {(0,1)}
                \arrow["\varphi", from=1-1, to=1-3]
                \arrow[from=1-1, to=3-1]
                \arrow[from=1-3, to=3-3]
                \arrow[dashed, from=3-1, to=3-3]
            \end{tikzcd}
        \end{equation} 
        with analogous conditions for $h_0$ and $h_1$.
    \end{itemize}
\end{definition}

We now record the fundamental properties an $(\oo,1)$-category $X$ must have in order to be deemed an $(\oo,1)$-category of $d$-dimensional unoriented bordisms: 

\begin{definition}
Let $\cF$ be an $(\oo,1)$-category. We say that $\cF$ \emph{is} the $(\oo,1)$-category of $d$-dimensional bordisms if it satisfies the following three conditions:
\begin{enumerate}
    \item \label{ch3:: cond: bordism cat objects} The underlying space of objects $\Obj(\cF)$ must be weakly equivalent to 
    \begin{equation}
        \coprod_{[\Sigma]}\mathbf{B}\!\Diff(\Sigma),
    \end{equation}
    where the coproduct ranges over all diffeomorphism classes of \emph{closed} $(d-1)$-dimensional unoriented manifolds.  
    \item \label{ch3:: cond: bordism cat morphisms} The space $\Hom(\cF)$ of morphisms should be weakly equivalent to the space 
    \begin{equation}
        \coprod_{[M,\partial_{\nu}M_i]}\mathbf{B}\!\Diff(M,\{\partial_{\nu}M\}),
    \end{equation} 
    where the coproduct ranges over all diffeomorphisms classes of unoriented $d$-dimensional cobordisms. Here $\Diff(M,\{\partial_\nu M\})$ is the subgroup of diffeomorphisms of $M$ that restrict to diffeomorphisms of collars on $\partial_0 M$ and $\partial_1 M$. 
    \item \label{ch3:: cond: bordism cat composable} The space of $k$-composable morphisms
    \begin{equation}
        \underbrace{ \Hom(\cF) \underset{\Obj(\cF)}{\times} \Hom(\cF) \underset{\Obj(\cF)}{\times} \cdots \underset{\Obj(\cF)}{\times} \Hom(\cF) \underset{\Obj(\cF)}{\times} \Hom(\cF) }_{k \text{ times}}
    \end{equation} 
    should be weakly equivalent to the space 
    \begin{equation}
         \coprod_{[M_0 \longra \cdots \longra M_k]}\mathbf{B}\!\Diff(M,\{h_i\}_{i=0}^{k}),
    \end{equation}
    where the coproduct is considered over all the diffeomorphism classes of $k$-composable chains of unoriented $d$-dimensional bordisms. This diffeomorphism group contains all the diffeomorphisms of $M$, the composed bordism, which additionally restrict to diffeomorphisms on a tubular neighborhood of $M_i \hookra M$ for all $i=1,\ldots,k$.
\end{enumerate}

\end{definition}

\begin{observation}
We acknowledge the redundancy in the above definition. Condition~(\ref{ch3:: cond: bordism cat morphisms}) is implied by condition~(\ref{ch3:: cond: bordism cat composable}), yet we chose to separate them for clarity.
\end{observation}

\begin{theorem}\label{ch3:: thm: bordism category characterization}
The $(\oo,1)$-category $\bord_{d}$ is \emph{the} $d$-dimensional bordism category.
\end{theorem}

\begin{proof}
    We verify the three conditions from the previous definition. Since $\bord_d$ is an interval Segal space, its values on any two interval objects with the same image under $\gamma$ are weakly equivalent. Thus, for the three conditions, it is enough to work with fixed representatives in $\cJ_0$, $\cJ_1$, and more generally in $\cJ_k$.

    For the object space, fix $\sU(I)\in \cJ_0$. A $k$-simplex of $\bord_d\big(\sU(I)\big)$ is a smooth embedded submanifold
    \[
        W\subset \AA^k\times \core\big(\sU(I)\big)\times \RR^\infty
    \]
    such that the projection
    \[
        \pi_W \maps W\longra \AA^k\times \core\big(\sU(I)\big)
    \]
    is a proper submersion. By \Cref{ap:: prop: Ehresmann fibration no boundary}, $\pi_W$ is therefore a smooth fiber bundle with closed fiber some $(d-1)$-dimensional manifold $\Sigma$. Since the base $\AA^k\times \core\big(\sU(I)\big)$ is contractible, \Cref{ap:: prop: triviality contractibility} identifies this bundle with the trivial bundle
    \[
        \Sigma\times \AA^k\times \core\big(\sU(I)\big).
    \]
    Write
    \[
        B = \AA^k\times \core\big(\sU(I)\big).
    \]
    Choosing a trivialization of $\pi_W$ identifies $W$ with a smoothly varying family of embedded copies of $\Sigma$ inside $\RR^\infty$, parametrized by $B$, hence with a smooth map
    \[
        B\longra \Emb(\Sigma,\RR^\infty).
    \]
    A second trivialization differs from the first by an automorphism of the trivial bundle $\Sigma\times B\to B$, and therefore induces a smooth map $B\to \Diff(\Sigma)$. Thus the resulting parametrized family of embeddings is well-defined only up to the right action of $\Diff(\Sigma)$. In other words, the connected component of $\bord_d\big(\sU(I)\big)$ corresponding to the diffeomorphism type of $\Sigma$ is the usual simplicial model for the moduli of embeddings of $\Sigma$ into $\RR^\infty$ modulo reparametrization by $\Diff(\Sigma)$. This is precisely the embedded-bundle picture used in \cite[\S2.1, p.~4]{GMTW09}: there the orbit map
    \[
        \Emb(\Sigma,\RR^\infty)\longra \Emb(\Sigma,\RR^\infty)/\Diff(\Sigma)
    \]
    is identified as a principal $\Diff(\Sigma)$-bundle, and the quotient is denoted $B_\infty(\Sigma)$. Since $\Emb(\Sigma,\RR^\infty)$ is contractible by Whitney's embedding theorem in the stable Euclidean space, again as recalled in \cite[\S2.1, p.~4]{GMTW09}, this component is a model for $\mathbf{B}\!\Diff(\Sigma)$. Taking the coproduct over all diffeomorphism classes of closed $(d-1)$-manifolds proves condition~(\ref{ch3:: cond: bordism cat objects}).

    The proof of condition~(\ref{ch3:: cond: bordism cat morphisms}) is analogous, but the relevant family direction is now the simplicial variable $\AA^k$, not the interval direction. Fix $\sU(A,B)\in \cJ_1$, and for a $k$-simplex set
    \[
        W_{A,B}=W\cap \Big(\AA^k\times \big[\mathsf{mid}(A),\mathsf{mid}(B)\big]\times \RR^\infty\Big).
    \]
    Condition (\ref{id: spatial}) from \Cref{ch3:: def: bord_d} implies that the projection
    \[
        q_W \maps W_{A,B}\longra \AA^k
    \]
    is a submersion, and it is proper because $W$ is bounded and the interval $\big[\mathsf{mid}(A),\mathsf{mid}(B)\big]$ is compact. Therefore $W_{A,B}$ should be regarded as a smooth family, over $\AA^k$, of compact $d$-manifolds with boundary. Condition (\ref{id: categorical}), restricted to the distinguished closed pieces $A$ and $B$, identifies smoothly varying incoming and outgoing collars in this family. Hence each fiber of $q_W$ is naturally a compact bordism, and the whole family is a smoothly varying family of bordisms with fixed source and target.

    Since the base $\AA^k$ is contractible, this family is trivial up to isomorphism. After choosing one fiber $M$, the residual datum is again a smoothly varying family of neat embeddings of $M$ into $\big[\mathsf{mid}(A),\mathsf{mid}(B)\big]\times \RR^\infty$, well-defined only up to smooth maps $\AA^k\to \Diff(M,\{\partial_\nu M\})$. Thus each connected component is modeled on the quotient of a contractible neat embedding space by the collar-preserving diffeomorphism group, hence is a model for $\mathbf{B}\!\Diff(M,\{\partial_\nu M\})$. Taking the coproduct over diffeomorphism classes of bordisms proves condition~(\ref{ch3:: cond: bordism cat morphisms}).

    Lastly, let $\sU(I_0,\ldots,I_k)\in \cJ_k$. Here the same point applies: the relevant smooth family is again parametrized by the simplicial variable. A simplex of $\bord_d\big(\sU(I_0,\ldots,I_k)\big)$ determines, by restricting to the closed interval spanning the chosen midpoints, a proper submersion to $\AA^\bullet$ whose fibers are compact $d$-manifolds with boundary. Condition (\ref{id: categorical}) records, along the distinguished closed pieces of $\core\big(\sU(I_0,\ldots,I_k)\big)$, the collar neighborhoods of the intermediate cutting hypersurfaces. Repeating the same trivialization argument over the contractible base $\AA^\bullet$ shows that every connected component is determined by a fixed composed bordism
    \[
        \Sigma_0 \xra{M_1} \Sigma_1 \xra{M_2} \cdots \xra{M_k} \Sigma_k,
    \]
    together with the chosen collar neighborhoods of the intermediate hypersurfaces. The automorphism group of this data is precisely the group $\Diff(M,\{h_i\}_{i=0}^{k})$ from \Cref{ch3:: def: diffeo group}. As above, the remaining embedding data are contractible in the stable Euclidean ambient space, so the corresponding connected component is a model for $\mathbf{B}\!\Diff(M,\{h_i\}_{i=0}^{k})$. Taking the coproduct over all diffeomorphism classes of $k$-composable chains proves condition~(\ref{ch3:: cond: bordism cat composable}).

    We conclude that $\bord_d$ satisfies the three defining properties above, and hence it is the $(\oo,1)$-category of $d$-dimensional bordisms.
\end{proof}

\begin{remark}
    The proof above isolates the two geometric inputs that will need stratified analogues later on: Ehresmann's fibration theorem, and the triviality of smooth bundles over contractible bases. Together, these reduce each simplex to a smoothly varying family of embeddings of a fixed manifold, or bordism, into $\RR^\infty$.
\end{remark}

\section{Symmetric monoidal structure}

One of the defining features of the cobordism categories studied and used in topological field theories is that it is symmetric monoidal relative to disjoint union. Constructing a symmetric monoidal structure for an $(\oo,1)$-category can be done in different flavors, as discussed in \Cref{ch2:: sec: symmetric monoidal structure}. We opt for using $\Gamma$-spaces because it offers a clearer route to formalizing models of $(\oo,1)$-categories of symmetric monoidal $(\oo,1)$-categories.

\begin{construction}\label{ch3:: cons: Gamma-space}
  For every $m\geq 0$, let $\bord_d^{\langle m \rangle}$ denote the following interval Segal space:
\begin{align}
    \bord_d^{\langle m \rangle} \maps \Int &\longra \sset \\
    I &\longmapsto \Big(\,[k] \mapsto \{M \subset \AA^{k}\times \mathsf{cor}(I)\times \RR^{\oo}\}\,\Big)
\end{align} 
where $M$ is given by a disjoint union $M = M_{1} \sqcup \cdots \sqcup M_{m}$ where each $M_{i}\subset \AA^{k} \times \mathsf{cor}(I)\times \RR^{\oo}$ is a $k$-simplex of $\bord_d(I)$. Given a pointed map of pointed finite sets, $f\maps \langle m \rangle \longra \langle n \rangle$, we construct a map of interval Segal spaces. This is the natural transformation whose component on $I \in \Int$ is
\begin{equation}
    \begin{tikzcd}
        {\bord_d^{\langle m \rangle }(I)} & {M = M_1 \sqcup \cdots \sqcup M_m} \\
        \\
        {\bord_d^{\langle n \rangle }(I)} & {(f_*)_I = \Big(\coprod_{f(i) = 1}M_i\Big) \sqcup \cdots \sqcup \Big(\coprod_{f(i) = n}M_i\Big)}
        \arrow[from=1-1, to=3-1]
        \arrow["{{(f_*)_I}}"', maps to, from=1-2, to=3-2]
    \end{tikzcd}
\end{equation}
where we omit the ambient data of the embeddings for simplicity.
\end{construction}

\nid The structure outlined in \Cref{ch3:: cons: Gamma-space} can be assembled so as to endow $\bord_d$ with a structure of a symmetric monoidal interval Segal space, as defined in \Cref{ch2:: sec: symmetric monoidal structure}.

\begin{theorem}\label{ch3:: thm: symmetric monoidal}
    The assignment 
    \begin{align}
        \bord_{d}^{\langle \bullet \rangle} \maps \Gamma &\longra \mathsf{ISSp} \\
        \langle m \rangle & \longmapsto \bord_d^{\langle m \rangle}
    \end{align}
    is functorial and endows the interval Segal space $\bord_d$ of unoriented $d$-dimensional bordisms with a symmetric monoidal structure.
\end{theorem}

\begin{proof}

    Verification that $\bord_{d}^{\langle \bullet \rangle}$ is a symmetric monoidal structure reduces to proving that the induced $\Gamma$-Segal map
    \begin{equation}
        \prod_{1 \leq k \leq m} i_k \maps \bord_{d}^{\langle m \rangle} \longra \underbrace{\bord_{d}^{\langle 1 \rangle} \times \cdots \times \bord_{d}^{\langle 1 \rangle}}_{m \text{ times}}
    \end{equation}
    is an equivalence of interval Segal spaces, i.e., a level-wise equivalence — which we will henceforth denote as $S_m$, to avoid cluttering notation. To show this, fix embeddings
    \begin{equation}
        \iota_k \maps \RR^{\oo} \hookra \RR^{\oo}, \quad \iota_k(\RR^{\oo}) \cap \iota_l(\RR^{\oo}) = \emptyset \mbox{ for } k, l \in \{ 1, 2, \ldots, m \}, \, \, k \neq l .
    \end{equation}
    For every $I\in \Int$, we construct an explicit homotopy inverse
    \begin{equation}
        r_{m} \maps \bord_{d}^{\langle 1 \rangle} \times \cdots \times \bord_{d}^{\langle 1 \rangle} \xra{\,\,\,\,} \bord_{d}^{\langle m \rangle}
    \end{equation}
    where we first consider the disjoint union of all $M_i$ sitting inside $\AA^{k}\times \mathsf{cor}(I)\times \big(\RR^{\oo}\big)^{\sqcup m}$ and then take its image under the action of the coproduct of the $\iota_k$ for all $k=0,\ldots,m$; this is subsumed by the formula
    \begin{equation}
        (M_1,\ldots, M_m) \longmapsto \bigcup_{k = 1}^{m}(\id_{\AA^k\times (I \cup I)}\times \iota_k)(M_1\sqcup \cdots \sqcup M_m) \subset \AA^k \times \mathsf{cor}(I)\times \RR^{\oo}
    \end{equation}
    
    Both composites are shown to be homotopy equivalent to the identity using as argument the weak contractibility of the space $\Emb(\RR^{\oo},\RR^{\oo})$. Let us start with 
    \begin{equation}
         \Bigg(\prod_{0\leq k \leq m} i_k \Bigg)  \circ r_{m} \maps \bord_{d}^{\langle 1 \rangle}(I) \times \cdots \times \bord_{d}^{\langle 1 \rangle}(I) \longra \bord_{d}^{\langle 1 \rangle}(I) \times \cdots \times \bord_{d}^{\langle 1 \rangle}(I)  
    \end{equation}
    which sends a tuple of $k$-simplices $(M_1,\ldots,M_m)$ in $\bord_d(I)$ to the disjoint union post-composed with the embedding $M = M_1 \sqcup \cdots \sqcup M_m \subset \AA^{k}\times \mathsf{cor}(I)\times \RR^{\oo}$, before finally sending each $M_i$ to its corresponding $k$-simplex in $\bord_d^{\langle 1 \rangle}$ obtained by forgetting all the others. There is a homotopy equivalence from this map to the identity on $\bord_{d}^{\langle 1 \rangle} \times \cdots \times \bord_{d}^{\langle 1 \rangle}$ defined by constructing an isotopy $H\maps \RR^{\oo}\times \AA^1 \longra \RR^{\oo}$ which for the $k$-th entry sends $\iota_k$ to $\id{\RR^{\oo}}$. The isotopy is well-defined precisely because it acts solely on the $\RR^{\oo}$ factor, whilst preserving all submersivity conditions.

    The other composite
    \begin{equation}
          r_{m} \circ \Bigg(\prod_{0\leq k \leq m} i_k \Bigg) \maps \bord_{d}^{\langle m \rangle}(I) \longra \bord_{d}^{\langle m \rangle}(I)
    \end{equation}
    acts on a disjoint union $M = M_1 \sqcup \cdots \sqcup M_m$,  where $\big(M \subset \AA^{k} \times \mathsf{cor}(I)\times \RR^{\oo}\big)$ first by sending it to $(M_1,\ldots, M_m)$ where each $M_i$ is just $M$ with the other $(m-1)$ components deleted, followed by a re-embedding of each component into $\RR^{\oo}$ using the fixed choices of embeddings $\bigcup_{k = 0}^{m}\iota_k \maps (\RR^{\oo})^{\sqcup k} \hookra \RR^{\oo}$. This admits a homotopy from it to the identity on $\bord_{d}^{\langle m \rangle}(I)$ defined through an isotopy $H\maps \RR^{\oo}\times \AA^{1} \longra \RR^{\oo}$ chosen so that $H_1 \circ \iota_k = e_k$, where $e_k$ encodes the embedding of each $M_k$ into $\AA^{k}\times \mathsf{cor}(I)\times \RR^{\oo}$.   

    \begin{remark}[Operadic viewpoint]\label{ch43:: rem :operadic-symm-monoidal}
        It is worth noting that the symmetric monoidal structure on $\bord_d$ could also be packaged by exhibiting $\bord_d$ as an algebra over an $\cE_\infty$-operad, rather than via the Segal $\Gamma$-object formalism used here.

        Recall that one convenient model for $\cE_n$ is the \emph{little $n$-disks} (or little $n$-cubes) operad, whose $k$-ary operations are rectilinear embeddings
        \begin{equation}
        \bigsqcup_{i=1}^k (0,1)^n \hookrightarrow (0,1)^n
        \end{equation}
        (with operadic composition given by substitution of embeddings); see \cite[\S5.1.1]{LurieHA}.
        In particular, one may model $E_\infty$ by ``little cubes in $\RR^\infty$'', and algebras over such an operad are precisely commutative monoid objects up to coherent homotopy; cf.\ again \cite[\S5.1.1]{LurieHA}.

        In our geometric setting, for each rectilinear embedding
        \begin{equation}
        \iota:\;\bigsqcup_{i=1}^k \RR^\infty \hookrightarrow \RR^\infty
        \end{equation}
        disjoint union of $k$ families of bordisms (inside $M\times \core(T)\times \RR^\infty$) followed by postcomposition with $\iota$ defines a $k$-ary operation
        \begin{equation}
        \mu_\iota:\;\bord_{d}^{\times k}\longrightarrow \bord_{d}.
        \end{equation}
        The operadic compatibilities are inherited from functoriality of postcomposition with embeddings, and the contractibility of the relevant embedding spaces in $\RR^\infty$ forces the required higher coherences.

        We nevertheless work with the equivalent Segal $\Gamma$-object presentation, since it interfaces more directly with the model-categorical framework we have used so far.
    \end{remark}
\end{proof}

\section{Homotopy category as the bordism 1-category}

Before studying $(\oo,1)$-categories we studied ordinary 1-categories, and the first instances of the bordism category date back to at least Atiyah's 1988 seminal paper \cite{At88}. To confirm that we have successfully lifted the ordinary bordism category to the $(\oo,1)$-categorical setting, we should compute its associated symmetric monoidal  homotopy 1-category, and verify that they are equivalent as symmetric monoidal categories. To achieve this, we will use the construction presented in \Cref{ch2:: def: homotopy category ISSP}.

\begin{remark}\label{ch3:: rmk: homotopy ehresmann}
    In the next theorem we will work with the homotopy category of $\bord_{d}$. Recall from \Cref{ch2:: def: homotopy category ISSP} that this category is defined as the homotopy category of the Segal space 
    \[
    X  \,\overset{\sf def}{=}\, \mathsf{RKan}_{\gamma}(\bord_d).
    \]

    \nid The objects of this homotopy category are thus the 0-simplices of $X$. A 0-simplex of this homotopy limit can be understood as a coherent family $\big(\Sigma_{I}, [0] \xra{\id} \gamma(\sU(I))\big)$ indexed by objects of $\cJ_0$, where $\Sigma_{I}\in \bord_d\big(\sU(I))_0$. The coherence comes from the fact that for every $f\maps \sU(I) \ra \sU(J)$, there is a 1-simplex in $\bord_d\big(\sU(J)\big)$ from $(\Sigma_I)_{|\sU(J)}$; which amounts to a concordance $\widetilde{\Sigma}_{I,J}\in \bord_{d}(\sU(J))_1$ relating $\bord_{d}\big(\sU(J)\big)$ with $\Sigma_J$.
    
    \medskip

    \nid  We illustrate this with the following diagram, where we depict the pair 
    \medskip
    \begin{equation}
        \tikzfig{figures/Ch3/BOR15}{\input{figures/Ch3/BOR15.tikz}}
    \end{equation} 

    \medskip
    
    \nid The bottom line is that, for every $0$-simplex and its associated coherent family, we can associate a unique diffeomorphism class of $(d-1)$-dimensional unoriented smooth manifolds. This class is given by considering the diffeomorphism class of the fiber of 
    \begin{equation}
    \Sigma_I \longra \sU(I) \times \RR^{\oo}
    \end{equation}
    over any $t\in \sU(I)$. This fiber is constant, up to diffeomorphism, since it is a smooth fiber bundle as a consequence of \Cref{ap:: prop: ehresmann boundary}.
\end{remark}

\begin{theorem}\label{ch3:: thm: homotopy category bordism}
    There is an equivalence of symmetric monoidal 1-categories 
    \begin{equation}
        \ho(\bord_d) \,\, \simeq \,\, \mathsf{Cob}_d
    \end{equation}
    between the homotopy category of $\bord_{d}$ and the unoriented $d$-dimensional bordism category.
\end{theorem}

\begin{proof}
    We construct functor 
    \begin{equation}
        h \maps \ho(\bord_d) \longra \mathsf{Cob}_d.
    \end{equation}
    Let $\underline{\Sigma} = (\Sigma_{I}, [0] \xra{\id} \gamma(I),\pi_{I})$ denote an object of $\ho(\bord_d)$, i.e., a coherent family of $(d-1)$-dimensional smooth unoriented manifolds. We define
    \begin{equation}
        h \big(\underline{\Sigma}\big) \, \overset{\mathsf{def}}{=} \, \pi_{I}^{-1}(\mathsf{mid}(I)), \quad \text{for some } I\in \cJ_0.
    \end{equation}
    By \Cref{ch3:: rmk: homotopy ehresmann}, this is well-defined up to diffeomorphism.

    Now fix objects $\underline{\Sigma}$ and $\underline{\Sigma'}$ of $\ho(\bord_d)$. A morphism in
    \[
        \Hom_{\ho(\bord_d)}\big(\underline{\Sigma},\underline{\Sigma'}\big)
    \]
    is represented by a coherent family of $d$-dimensional unoriented bordisms compatible with the given source and target data. Choosing some $\sU(A,B)\in \cJ_1$ and a representative
    \[
        M \hookra \sU(A,B)\times \RR^{\oo},
    \]
    we define
    \begin{equation}
        M_{A,B} = \pi^{-1}\big([\mathsf{mid}(A),\mathsf{mid}(B)]\big), \quad h\big([M]\big) \, \overset{\sf def}{=} \, \big[M_{A,B}\big],
    \end{equation}
    where we take the diffeomorphism class of the extracted bordism $M_{A,B}$. This is well-defined: changing the choice of $\sU(A,B)$ only changes the extracted fiber by a diffeomorphism, again by \Cref{ch3:: rmk: homotopy ehresmann}, and concordant representatives determine the same class in $\mathsf{Cob}_d$.

    \medskip

    To prove that this functor is essentially surjective, let $\Sigma$ be a closed $(d-1)$-dimensional smooth manifold. The cylinder $\Sigma\times [0,1]$ admits a neat embedding into $\sU(A)\times \RR^{\oo}$ for every $\sU(A)\in \cJ_{0}$ by the manifold-with-boundary version of Whitney's embedding theorem, for instance \cite[Theorem 8.9]{CS19}. Choosing these embeddings compatibly, via linear diffeomorphisms $\sU(A)\cong [0,1]$, they assemble into a coherent family and hence define an object of $\ho(\bord_d)$ whose image under $h$ is the diffeomorphism class of $\Sigma$.

    \medskip
    
    It remains to show that $h$ is fully faithful. Fix objects $\underline{\Sigma}$ and $\underline{\Sigma'}$ of $\ho(\bord_d)$, and write
    \[
        \Sigma = h(\underline{\Sigma}), \qquad \Sigma' = h(\underline{\Sigma'}).
    \]
    By \Cref{ch3:: thm: bordism category characterization}, the total space of $1$-morphisms in $\bord_d$ is weakly equivalent to a coproduct of classifying spaces indexed by diffeomorphism classes of bordisms. Restricting the source-target map to the fiber over the pair $(\underline{\Sigma},\underline{\Sigma'})$ therefore selects precisely the sub-coproduct indexed by bordisms from $\Sigma$ to $\Sigma'$. Hence the corresponding mapping space is weakly equivalent to
    \begin{equation}
        \coprod_{[M \maps \Sigma \longra \Sigma']} \mathbf{B}\!\Diff(M,\{\partial_{\nu}M\}),
    \end{equation}
    where the coproduct ranges over diffeomorphism classes of compact $d$-dimensional bordisms from $\Sigma$ to $\Sigma'$. Passing to path components, and using that each classifying space $\mathbf{B}\!\Diff(M,\{\partial_{\nu}M\})$ is path connected, we obtain a bijection
    \begin{equation}
        \Hom_{\ho(\bord_d)}\big(\underline{\Sigma},\underline{\Sigma'}\big)
        \cong
        \big\{[M \maps \Sigma \longra \Sigma']\big\}
        =
        \Hom_{\mathsf{Cob}_d}\big(h(\underline{\Sigma}),h(\underline{\Sigma'})\big).
    \end{equation}
    Under this identification, the induced map on hom-sets is exactly the map defined above by extracting the bordism over the interval between the chosen midpoints. Hence $h$ is fully faithful. Since it is also essentially surjective, it is an equivalence of categories. Compatibility with disjoint union is immediate from the construction, so this is an equivalence of symmetric monoidal $1$-categories.
\end{proof}

\section{Tangential structures}
\label{ch3:: sec: tangential structures}

When studying manifold topology, it is natural to ask what additional structures manifolds may be endowed with. Generally speaking, the structures we care about are the simplest ones: \emph{tangential structures}. These structures can be thought of as families of linear data, such as group representations.

\begin{definition}
    An $n$-dimensional tangential structure is a topological space $B$ together with a fibration $B \longra \BO(n)$. A \define{stable tangential structure} is a topological space $\fB$ together with a fibration $\fB \longra \BO$. 
\end{definition}

\begin{remark}
    A stable tangential structure $\fB \longra \BO$ induces, for every $n\geq 0$, an $n$-dimensional tangential structure via pullback 
    \begin{equation}
        \begin{tikzcd}[ampersand replacement=\&]
            {B(n)} \&\& {\BO(n)} \\
            \\
            \fB \&\& \BO
            \arrow[from=1-1, to=1-3]
            \arrow[from=1-1, to=3-1]
            \arrow["\lrcorner"{anchor=center, pos=0.125}, draw=none, from=1-1, to=3-3]
            \arrow[hook, from=1-3, to=3-3]
            \arrow[from=3-1, to=3-3]
        \end{tikzcd}
    \end{equation}
    where $\BO(n) \hookra \BO$ denotes the canonical inclusion.
\end{remark}

\begin{remark}
    Since we can always factor a map of topological spaces $B \longra \BO(n)$ as \[B \xra{\,\,\,\sim\,\,\,} B^{\sf fib} \longra \BO(n),\] where the second map is a fibration, we might as well define an $n$-dimensional tangential structure as a map $B \longra \BO(n)$.
\end{remark}

\begin{definition}
    Let $M$ be a $k$-dimensional smooth manifold and let $B\longra \BO(n)$ be an $n$-dimensional tangential structure. A $B$\define{-framing} on $M$ is a map $\alpha \maps M \longra B$, together with a homotopy making the following diagram commute up to the specified homotopy:
    \begin{equation}
    \begin{tikzcd}[ampersand replacement=\&]
        \&\& B \\
        M \&\& {\BO(n).}
        \arrow[from=1-3, to=2-3]
        \arrow["\alpha", from=2-1, to=1-3]
        \arrow[""{name=0, anchor=center, inner sep=0}, "{\widetilde{\mathsf{\tau}_M}}"', from=2-1, to=2-3]
        \arrow[between={0.4}{0.6}, Rightarrow, from=1-3, to=0]
    \end{tikzcd}
    \end{equation}
    The map $\widetilde{\tau}_M$ is the map classifying the $n$-dimensional stabilized tangent bundle 
    \begin{equation}
        \widetilde{\sT M} \,\, \overset{\sf def}{=} \,\, \epsilon^{n-k}_{M} \oplus \sT M.
    \end{equation}   
\end{definition}

\begin{remark}
    Instead of defining a $B$-framing as a lift together with a homotopy, we can equivalently state that it is a lift in $\Top$ as an $(\oo,1)$-category — this makes the choice of a homotopy implicit.
\end{remark}

\begin{example}
    Let $G$ be a Lie group together with a linear representation $\rho \maps G \ra \sG \sL(n$). Recall that the Gram-Schmidt orthogonalization offers a strong deformation retract of $\sO(n) \xhookrightarrow{\simeq} \sG\sL(n)$. The representation $\rho$ gives rise to a map of spaces
    \begin{equation}
        \rho \maps \sB G \longra \sB G\sL(n) \xra{\simeq} \sB\sO(n).
    \end{equation}
    This representation classifies a smooth rank $n$ bundle $E$ on $\sB G$. To endow a smooth $n$-dimensional manifold $M$ with a $(\sB G,\rho)$-framing is to give a map $\varphi\maps M \longra \sB\sG$ together with an isomorphism of rank $n$ vector bundles over $M$:
    \begin{equation}
        \varphi^{*}E \xra{\,\,\,\,\cong\,\,\,\,} \sT M.
    \end{equation}
    Ultimately, a $\sB G$-framing amounts to a reduction of the structure group of $\sT M$ from $\sO(n)$ to $\SO(n)$. Equivalently, a manifold admits a $\sB G$-framing if and only if its tangent bundle's transition functions can be chosen so as to take values in $G$.
\end{example}

\begin{example}
    The short exact sequence 
    \begin{equation}
        1 \longra \SO (n) \longra \sO(n) \longra \sO(1) \longra 1
    \end{equation}
    induces a group homomorphism $\SO(n) \longra \sO (n)$ for every $n \geq 1$. A smooth $n$-manifold with a $\sB\!\SO (n)$-framing is thus a smooth manifolds endowed with an additional choice of an orientation.   
\end{example}

\begin{example}
    Not every tangential structure is obtained from a smooth linear representation of a Lie group. Let $X$ be a space, an $X$-framing is defined via the projection 
    \begin{equation}
        X \times \sB\sO(n) \xra{\sf pr} \sB\sO(n).
    \end{equation}
    An $X$-framed smooth manifold $M$ is a smooth $n$-dimensional manifold together with a continuous map $M \longra X \times \sB\sO(n)$ over $\sB\sO(n)$; this data is determined by a continuous map $M \longra X$. In the case where $X = \cD_n$ is a discrete space, i.e., a set, a continuous map $M \longra \cD_n$ factors through $M \ra \pi_{0}(M)$, which is equivalent to a \emph{coloring} or a labeling of the connected components of $M$ by colors, the elements of $\cD_n$. This tangential structure admits immediate generalizations: we can combine a smooth linear representation of a Lie group $G$, together with a coloring by $\cD_n$:
    \begin{equation}
        \big(B \longra \sB\sO(n)\big) \simeq \big(\cD_{n}\times \sB G \xra{\,\,\, \id \times \rho\,\,\,} \cD_{n}\times \sB \sO(n) \xra{\,\,\,\pr\,\,\,} \sB\sO(n)\big).
    \end{equation}
    An $(\cD_n\times G,\rho)$-framing on $M$ amounts to a coloring of its connected components by colors in $\cD_n$ together with a $G$-reduction of its tangent bundle.
\end{example}

When we have a family of smooth manifolds, i.e. a smooth bundle of manifolds with fiber a smooth manifold $M$, it is natural to ask how might families of smooth manifolds with a given tangential structure behave. To study such families, we introduce the \emph{fiberwise tangent bundle of a smooth bundle $\pi_W\maps W \longra M$} as the kernel of the derivative of $\pi_W$:

\begin{equation}
    \mathsf{Ker}(\sD\pi_{l}) \longra M.
\end{equation}
This bundle is also often called the vertical tangent bundle.

Using the fiberwise tangent bundle we can formulate what we mean by a family of $B$-framed smooth manifolds. Roughly, a family of $B$-framed smooth manifolds is a smooth fiber bundle $\pi_W \maps W \longra M$ together with a $B$-trivialization of its fiberwise tangent bundle.

\begin{definition}
    Let $\pi_W \maps W \longra M$ be a smooth fiber bundle. We say that it is a smooth family of $B$-framed manifolds if there is a map $\alpha\maps M \longra B$ together with a homotopy making the following diagram commute:
    \begin{equation}
    \begin{tikzcd}[ampersand replacement=\&]
        \&\& B \\
        M \&\& {\BO(n).}
        \arrow[from=1-3, to=2-3]
        \arrow["\alpha", from=2-1, to=1-3]
        \arrow[""{name=0, anchor=center, inner sep=0}, "{\widetilde{\mathsf{\tau}_{W/M}}}"', from=2-1, to=2-3]
        \arrow[between={0.4}{0.6}, Rightarrow, from=1-3, to=0]
    \end{tikzcd}
    \end{equation}
    The map $\widetilde{\mathsf{\tau}_{W/M}}\maps M \longra \BO(n)$ is the classifying map of the fiberwise tangent bundle over $M$.
\end{definition}

\begin{definition}\label{ch3:: def: bord_d tangential}
    Let $V$ be a finite-dimensional real vector space and let $e\maps B \longra \sB\sO(d)$ be a map of spaces. Let us define a functor 
    \begin{equation}
    \widetilde{\bord}^{B}_{d,V} \maps \mathsf{Man^{\op}}\times \Int \longrightarrow \set    
    \end{equation}
    sending a pair $(M,T)$, where $M\in \man^{\op}$ and $T\in \Int$, to the set of bounded closed subsets 
    \begin{equation}
    W\subseteq M\times \mathsf{cor}(T) \times V
    \end{equation} 
    satisfying the following four conditions:
    \begin{enumerate}
    \item $W$ is a smooth $(d+\mathsf{dim}(M))$-dimensional submanifold with boundary $\partial W$ such that 
    \begin{equation}
        \partial W \subset \partial \big(M \times \mathsf{cor}(T) \times V\big) = M \times \partial (\mathsf{cor}(T)) \times V   
    \end{equation}
    the boundary of $W$ must be contained in the boundary of the ambient space.
    \item The restricted projection 
    \begin{equation}
        \pi\maps \big(W \cap (M\times T^{\mathsf c} \times V) \big) \hookra  M \times T^{\mathsf c} \times V \fib  M \times T^{\mathsf c}
    \end{equation} 
    is a submersion. 
    \item The restricted projection 
    \begin{equation}
    \pi_W\maps W \hookra M\times \mathsf{cor}(T) \times V \fib M 
    \end{equation}
    is a submersion.
    \item A map of spaces $\alpha \maps W \longra B$ making the following diagram commute up to homotopy:
    \begin{equation}
        \begin{tikzcd}[ampersand replacement=\&]
            \&\& B \\
            W \&\& {\BO(d)}
            \arrow["e", from=1-3, to=2-3]
            \arrow["\alpha", dashed, from=2-1, to=1-3]
            \arrow["{\widetilde{\tau_{W/M}}}"', from=2-1, to=2-3]
        \end{tikzcd}
    \end{equation}   
    This amounts to a choice of a map together with a homotopy from $e\circ \alpha$ to $\widetilde{\tau}_{W/M}$, where the latter classifies the fiberwise tangent bundle of the smooth submersion $\pi_W\maps W \longra M$. 
    \end{enumerate} 

    \noindent The functor $\widetilde{\bord}^{B}_{d,V}$ sends:
    \begin{itemize}
    \item[($\man$):] The morphism $(f,\mathsf{id}_{T})\maps (N,T) \longra (M,T)$ induced by the smooth map $f\maps N \longra M$, which sends the submanifold with boundary
    \begin{equation}
         W\subseteq M\times \mathsf{cor}(T) \times V 
    \end{equation}
    to $(f \times \id)^{*} W$ — the pullback submanifold 
    \begin{equation}
        \begin{tikzcd}[ampersand replacement=\&]
            {(f\times \id)^\ast W} \& {W} \\
            \\
            N\times \mathsf{cor}(T) \times V \& M\times \mathsf{cor}(T) \times V,
            \arrow[from=1-1, to=1-2]
            \arrow[hook, "(f\times \id)^\ast \pi_{|W}"',from=1-1, to=3-1]
            \arrow["\lrcorner"{anchor=center, pos=0.125}, draw=none, from=1-1, to=3-2]
            \arrow[hook, from=1-2, to=3-2]
            \arrow["f\times \id"', from=3-1, to=3-2]
        \end{tikzcd}
    \end{equation}
    which is a $\big(d + \mathsf{dim}(N)\big)$-dimensional submanifold with boundary
    \begin{equation}
        (f\times \id)^* W \subseteq N \times \mathsf{cor}(T) \times V.
    \end{equation}
Which is well-defined since the pullback preserves the fibers over the boundary $\partial (N \times \mathsf{cor}(T) \times V)$. Moreover, we have that $(f\times \id)^* W$ comes equipped with a map of spaces to $\sB\sO(d)$ and a corresponding homotopy, both obtained through pre-composition with the top map of the pullback diagram
    \begin{equation}
        \begin{tikzcd}
            &&& B \\
            \\
            {(f\times \id)^*W} & W && {\sB\sO(n)}
            \arrow[from=1-4, to=3-4]
            \arrow[from=3-1, to=3-2]
            \arrow[dashed, from=3-2, to=1-4]
            \arrow["{\sT_{W/M}}"', from=3-2, to=3-4]
        \end{tikzcd}
    \end{equation}

    \item[($\Int$):] The pair $(\id_M, T\subset S)$ by \emph{forgetting} connected components and restricting the submanifold with boundary 
    \begin{equation}
        W \subseteq N \times \mathsf{cor}(S) \times V
    \end{equation}
    along $N \times \mathsf{cor}(T)\times V \hookra N \times \mathsf{cor}(S)\times V$. The restriction $W_{|T}$ comes equipped with a smooth embedding 
    \begin{equation}
        W_{|T} \hookra W
    \end{equation}
    which induces a $B$-framing on the restricted fiberwise tangent bundle, as desired.
    \end{itemize}  
\end{definition}

\begin{remark}
    The proof of the interval Segal condition rests entirely on the two locality statements \Cref{ch3:: prop: sheaf condition} and \Cref{ch3:: prop: naive sheaf con}. These also hold for $\widetilde{\bord}_{d}^B$ since $B$-framed manifolds enjoy the same locality properties. Therefore, we skip the proof of the next proposition.
\end{remark}

\begin{proposition}\label{ch3:: prop: Segal condition with structrue}
    Let $d \geq 0$ be a non-negative integer, let $B \longra \BO(d)$ be a tangential structure, and let $V$ be a finite-dimensional vector space of sufficiently high dimension. The composite functor
    \begin{align}
        \bord_{d,V}^{B} \maps \Int &\longra \shv(\man)  \longra \sset  \\
         I &\longmapsto \widetilde{\bord}_{d,V}^{B}(-,I) \longmapsto \bB\widetilde{\bord}_{d,V}^{B}(*,I) = \left|\,\widetilde{\bord}_{d,V}^{B}(\AA^\bullet,I)\,\right| 
    \end{align}
    is an interval Segal space. 
\end{proposition}

%% file: Chapters/Ch4.tex
\chapter{Conically smooth stratified spaces}

We now turn to the theory of \emph{conically smooth stratified spaces}. These were first introduced in the work of Ayala--Francis--Tanaka and further elaborated in the subsequent works of Ayala--Francis--Rozenblyum \cite{AFT1, AFR1, AFR2}. They will form the geometric foundations for our stratified higher bordism categories, and their adoption for that role is undoubtedly one of the more novel constructions contained in this thesis. Part of this novelty lies in the recency of their introduction, and in the bespoke categorical toolkit required to master them well enough so as to use them effectively. 

It is a stated goal of this chapter to not only introduce the reader to the main definitions and results on conically smooth stratified spaces, but also expand, elaborate and even refine the exposition of these geometric objects. The overarching goal being to encourage their broader adoption in contexts where stratified spaces naturally arise.

\section{Stratified spaces}

When we introduce smooth manifolds formally, it is customary to first introduce the notion of topological space followed by the definition of topological manifold, after which smooth manifolds are defined as topological manifolds with an additional notion of smooth atlas. This is the trajectory we will emulate for conically smooth stratified spaces; the first section will introduce the equivalent to topological manifolds for stratified spaces.       

\begin{definition}\label{ch4:: def: upward closed}
    Let $P$ be a \textit{partially ordered set}, henceforth referred to by its customary shorthand poset. We say that a subset $U$ of $P$ is \textbf{upward closed} if for every $x \in U$, its upward closure, i.e., the set of elements $x\leq y$, is contained in $U$. 
\end{definition} 

\begin{observation}
Every poset is naturally endowed with the structure of category, where for every $p,q \in P$ there is a unique morphism $q \longra p$ if and only if $q\leq p$. Composition derives from the transitivity of $\leq$ while associativity is trivial.
\end{observation}

\begin{observation}
Given a poset $P$, there is an upward closed subset $P_{\geq p}$ associated to every $p \in P$ which is defined as 
\begin{equation}
    P_{\geq p} = \{ q \in P\,:\, p \leq q \}.
\end{equation}
Interpreting $P$ as a category, the subset $P_{\geq p}$, which is itself a sub-poset of $P$, can also be interpreted as the overcategory $p/P$.  
\end{observation}

\begin{definition}\label{ch4:: def: strat_topology}
    Let $P$ be a poset. We define a topology for $P$ for which an open subset $U \subset P$ is \emph{open} if it is upward closed. We will call this topology the Alexandroff topology. 
\end{definition}

\begin{example}
    The finite totally ordered set $[n]$ is a poset for every $n\geq -1$. The open subsets of $[n]$ are given by the image of the maps 
    \begin{equation}
    \langle n-k, n-k +1, \ldots, n \rangle \maps [k] \hookra [n]
    \end{equation} for every $0\leq k\leq n$. Alternatively, we can write these as $[n]_{\geq k}$.
\end{example}

\begin{observation}
We check that this is in fact a topology. Let $\{U_{i}\}_{i=1}^{n}$ be a finite list of open subsets of $P$, if $x \in U_1 \cap \cdots \cap U_n$ then if $y\geq x$ we have $y\in U_1,\cdots, U_n$ it thus belongs to the finite intersection.       
\end{observation}

A map $P \xra{f} Q$ of posets is continuous if and only if it is order-preserving. Indeed, the preimage of an upward closed subset under an order-preserving map is again upward closed, and conversely continuity forces $x\leq y$ to imply $f(x)\leq f(y)$. 

\begin{observation}\label{ch4: obs:: posets as tops}
    This means that there exists a fully faithful functor 
    \begin{equation}
        \mathsf{Posets} \hookra \Top
    \end{equation}
    from the category of posets with the Alexandroff topology and continuous maps thereamong, to the category of topological spaces.
\end{observation}

As noted by Ayala--Francis--Tanaka, the next definition first appears in Lurie's \emph{Higher Algebra} \cite[Appendix~A.5]{LurieHA}. 

\begin{definition}\label{ch4:: def: p-stratified space}
    Let $P$ be a poset, we define a $P$-stratified space to be a pair 
    \begin{equation}
        (X,S\maps X \longra P)
    \end{equation} where $X$ is a topological space and $S\maps X \longra P$ is a continuous map of topological spaces. 
\end{definition}

For $p\in P$, we write $X_p=S^{-1}(p)$ and call it the $p$-stratum of $X$.

\begin{observation}
There are many different definitions of stratified space in the literature, but all coalesce around the notion of $P$-stratified space. 
\end{observation}

\begin{observation}
    Ordinary topological spaces, or unstratified topological spaces, are encoded as stratified spaces with a trivial stratification, or a stratification by the trivial poset $P\simeq *$.
\end{observation}

As with any good definition of a new mathematical object, we must give a definition of stratified continuous map, which we will henceforth just call \emph{stratified map}. 
\begin{definition}\label{ch4:: def: stratified map}
    Let $(X, S\maps X \longra P)$ and $(Y, S'\maps Y \longra Q)$ be $P$-stratified and $Q$-stratified spaces, respectively. A \textbf{continuous stratified map} $f \maps (X,S\maps X \longra P ) \longra (Y,S'\maps Y \longra Q)$ is a pair $(f \maps X \longra Y, \alpha\maps P \longra Q)$ where $f$ and $\alpha$ are continuous maps such that: 
    \begin{equation}
        \begin{tikzcd}[ampersand replacement=\&]
            X \&\& Y \\
            \\
            P \&\& Q
            \arrow["f"{description}, from=1-1, to=1-3]
            \arrow["S"', from=1-1, to=3-1]
            \arrow["{S'}", from=1-3, to=3-3]
            \arrow["\alpha"{description}, from=3-1, to=3-3]
        \end{tikzcd}
    \end{equation}
    A pair of morphisms $(f,\alpha)$ and $(g,\beta)$ is composed by horizontally stacking the respective commutative diagrams. The composition is well-defined on the nose.
\end{definition}

\begin{remark}
    One compact and categorical way of thinking about how the theory of stratified topological spaces might be defined is as the following pullback category:
    \begin{equation}
        \mathsf{Ar}(\Top) \underset{\Top}{\times} \mathsf{Posets},
    \end{equation} 
    where the arrow $\mathsf{Ar}(\Top) \longra \Top$ sends a continuous map $X \longra Y$ to its target $Y$.
\end{remark}

\nid The next example indicates how $P$-stratifications connect to the typically adopted notion of stratification as a decomposition into locally closed subspaces. 

\begin{example}
    Let $X$ be a topological space, and let $P = \{a < b\}$ be a two-element poset. The Alexandroff topology is such that the only open subsets are $P$, $\{b\} \subset P$, and the $\emptyset$. A continuous map $S\maps X \longra P$ is specified entirely by the choice of a closed subset $X_0 \subseteq X$, corresponding to the pre-image of $a\in P$. Therefore, a $P$-stratified space in this case is equivalent to a topological space $X$ and a choice of a closed subset $X_0 \subseteq X$ — the continuous map $S$ sends points of $X_0$ to $a$ and points of $X\setminus X_0$ to $b$.
\end{example}

\begin{example}
    Posets do not need to be linearly ordered. As an example, let $P$ denote the \emph{span-shaped} poset
    \begin{equation}
        P = \big(l\, \la * \ra r \, \big) 
    \end{equation}  
    A $P$-stratified space $X$ consists of two open strata $X_l$ and $X_r$, together with a closed stratum $X_*$ contained in the complement of $X_l \cup X_r$. 
    Let us illustrate with a concrete example how a $P$-stratified torus looks like:
    \begin{equation}
    \tikzfig{figures/Strat/STR1}{\input{figures/Strat/STR1.tikz}}
    \end{equation}
    The two strata $X_{l}$ and $X_{r}$ are open cylinders while $X_*$ is the closed subspace which is contained in the complement of $X_{l} \cup X_{r} \subsetneq X$.
\end{example}

\begin{example}
To hone in on the intuition for what a stratified map is, we note the following pair of examples 
\begin{equation}
    \tikzfig{figures/Strat/STR10}{\input{figures/Strat/STR10.tikz}} \qquad \qquad \tikzfig{figures/Strat/STR11}{\input{figures/Strat/STR11.tikz}}
\end{equation}
\smallskip

\noindent where on the left we depict as a source of the map a stratified circle which has eight different strata: it has four 0-dimensional strata (the black dots) and four 1-dimensional open arcs; the target is a stratified closed interval with five different strata: it has three 0-dimensional strata and two 1-dimensional strata given by the open arcs. The map is the vertical projection induced by the framing of the page. On the right, the spaces are the same but for the absence of a single 0-dimensional stratum. The rightmost map fails to be stratified because the lower arc of the source stratified circle gets smeared across different strata.     
\end{example}

\begin{example}
    Consider the two stratified circles of the example before. The underlying spaces are the same, so we can ask whether the identity homeomorphism is stratified or not. The answer is yes, but only in one direction — the maps 
    \begin{equation}
        \tikzfig{figures/Strat/STR12}{\input{figures/Strat/STR12.tikz}}
    \end{equation}
    are not both stratified. Only the top map is, and we call it a \textbf{refinement} precisely because it is witnessed by a homeomorphism of the underlying topological spaces.
\end{example}

\medskip
We have a working definition of \emph{topological} stratified space and continuous stratified map, now we would like to articulate what it means for a stratified space to be openly embedded into another stratified space. This is a fundamental touchstone in the path towards defining an appropriate site of topological stratified manifolds — which, as we have previously mentioned, will be the scaffolding for a theory of \emph{smooth} stratified spaces. 

Before introducing the definition of stratified open embedding, let us give the definition for the unstratified case first as a warm-up and a refresher.

\begin{definition}
    Let $X$ and $Y$ be topological spaces and $f\maps X \longra Y$ a continuous map. We say that $f$ is an \textbf{open embedding} if $f$ is \emph{open}, i.e., sends open subsets of $X$ to open subsets of $Y$, and an embedding i.e., the image $f(X) \subseteq Y$ is homeomorphic to $X$. Alternatively, we say that $f$ is a \textbf{closed embedding} if it is \emph{closed}, i.e., if it sends closed subsets of $X$ to closed subsets of $Y$ and also an embedding.
\end{definition}

\begin{definition}\label{ch4:: def: strat open emb}
    A continuous stratified map $(f,\alpha)\maps (X\longra P) \longra (Y \longra Q)$ is a \textbf{stratified open embedding} if:
    \begin{itemize}
        \item the underlying map $f\maps X \longra Y$ is an open embedding of topological spaces;
        \item for every $p \in P$, the restriction
        \begin{equation}
            f_{|X_{p}} \maps X_{p} \longra Y_{\alpha(p)}
        \end{equation}
        is an open embedding of topological manifolds.
    \end{itemize}
\end{definition}

Now we turn to two ways of producing stratified spaces from old ones, which will be especially important for our theory — \emph{cones} and \emph{products}. Together they will generate local models for our singular spaces.

\begin{definition}
    Let $(S \maps X \longra P)$ and $(S' \maps Y \longra Q)$ be a pair of stratified topological spaces. We define the following continuous map
    \begin{equation}
    \big(S\times S' \maps (X \times Y) \longra (P \times Q) \big),
    \end{equation}
    where $S\times S'$ is the product of the stratification maps. We call this the \textbf{product stratification}.
\end{definition}

\begin{definition}\label{ch4:: def: cones}
    Let $(X\longra P)$ be a stratified topological space. We define its associated \textbf{open cone} $\sC(X \longra P)$ to be the stratified topological space obtained as follows: 
    \begin{itemize}
        \item $\sC(X)$: the open cone on $X$ is the pushout of topological spaces
        \begin{equation}
            \sC(X) \, \overset{\sf def}{=} \,  \big(X\times \RR_{\geq 0}\big) \coprod_{X\times \{0\}} *.
        \end{equation}
        \item $\sC(P)$: the cone on $P$ is obtained by formally adjoining a minimal element $*$. Equivalently as the pushout of posets
        \begin{equation}
            \sC(P) \, \overset{\sf def}{=} \, (P \times [1]) \coprod_{P \times \{0\}} *.
        \end{equation}
    \end{itemize}
    The stratification $\sC(X) \longra \sC(P)$ sends $x \in X_{p} \times \RR_{> 0}$ to $p$ and sends the cone point $* \in \sC(X)$ to $*\in \sC(P)$. The cone point is thus sent to the cone point.
\end{definition}

\begin{remark}
    The poset $P \times [1]$ has two canonically isomorphic copies of $P$ embedded in it, call them $P_0$ and $P_1$, together with an arrow connecting every element of $P_0$ to its image in $P_1$. By taking the cone we are essentially collapsing $P_0$ and thus getting a new minimal element called $*$.  
\end{remark}

Since we are attempting to define a notion of topological stratified manifold, we must first isolate the correct topological class of stratified spaces.
\begin{definition}\cite[Definition~2.1.15]{AFT1}
    A \textbf{$C^0$ stratified space} is a paracompact stratified topological space belonging to the smallest full subcategory of stratified topological spaces such that:
    \begin{enumerate}
        \item it contains the empty stratified space;
        \item if $(X \ra P)$ is compact and belongs to it, then so does its open cone $\sC(X \ra P)$;
        \item if $(X \ra P)$ belongs to it, then so does $\RR \times (X \ra P)$, where $\RR$ carries the trivial stratification;
        \item it is closed under stratified open embeddings;
        \item if a stratified topological space admits an open cover by objects in the subcategory, then it also belongs to the subcategory.
    \end{enumerate}
\end{definition}

Since spaces of the form $\RR^i\times \sC(Z)$ are the local models singled out by this definition, we now isolate them.
\begin{definition}\label{ch4:: def: C0 basics}
    A \textbf{$C^0$ basic} is a $C^0$ stratified space of the form 
    \begin{equation}
        \RR^{i} \times \sC(Z)
    \end{equation}
    where $Z$ is a compact $C^0$ stratified space, $i\geq 0$ and $\RR^{i}$ is taken as the usual topological manifold with the trivial stratification. 
\end{definition}

\begin{lemma}\cite[Lemma~2.2.2]{AFT1}
    Let $X \longra P$ be a stratified topological space whose underlying space is second countable and Hausdorff. Consider the collection 
    \begin{equation}
        \{f_\alpha\maps U_\alpha \hookra X\mkern1mu | \, f_\alpha \text{ an open embedding} \}
    \end{equation}
    where $U_\alpha$ ranges over all $C^0$ basics. This collection forms a basis for the topology of $X$ if and only if $X$ is a $C^0$ stratified space.
\end{lemma}

\subsection{Dimension and depth}

When dealing with topological manifolds, there is a fundamental local invariant --- \emph{dimension}. For $C^0$ stratified spaces, this notion demands revisiting, since not every point admits a neighborhood homeomorphic to an open $\RR^{n}$ for some $n \geq 0$. 

\begin{definition}
    Let $(X \longra P)$ be a $C^0$ non-empty stratified space. For $x\in X$, we define its \textbf{local dimension} $\dim_x(X)$ to be the \emph{covering dimension} of $X$ at $x$, and we define the dimension of $(X\longra P)$ to be 
    \begin{equation}
        \mathsf{dim}(X) \mkern1mu \overset{\sf def}{=} \mkern1mu \underset{x \in X}{\mathsf{sup}} \mkern1mu\mkern1mu \mathsf{dim}_{x}(X).
    \end{equation}
    We say that $X$ has \textbf{pure dimension $n$} if it has constant local dimension, i.e., if $\mathsf{dim}_{x}(X) = n$ for all $x\in X$. We adopt the convention $\mathsf{dim}(\emptyset) = -1$.
\end{definition}

\begin{definition}
    Let $(X \longra P)$ be a non-empty $C^0$ stratified space and let $x\in X_{p}$ belong to the stratum indexed by $p\in P$. We define the \textbf{local depth} of $X$ at $x$ by
    \begin{equation}
        \depth_{x}(X) \mkern1mu \overset{\sf def}{=} \mkern1mu \mathsf{dim}_{x}(X) - \mathsf{dim}(X_{p}),
    \end{equation}
    and we define the \textbf{depth} of $X$ by
    \begin{equation}
        \depth(X) \mkern1mu \overset{\sf def}{=} \mkern1mu \underset{x \in X}{\mathsf{sup}} \mkern1mu \depth_x(X).
    \end{equation}
    We adopt the convention $\depth(\emptyset)=-1$.
\end{definition}

\begin{observation}
Before wrapping up the fundamentals of the theory of $C^0$-stratified spaces, let us end with the following constructions and results.
\end{observation}

\begin{construction}
    Let $\big(X \ra P \big)$ be a $C^0$ stratified space. There is a stratification $\big(X \ra P^{\sf con}\big)$ of $X$, where the $P^{\sf con}$ has as a set of objects the set of path components 
    \begin{equation}
        \pi_0\Big(\coprod_{p \in P} X_p\Big).
    \end{equation}
    Let $q$ and $q'$ denote a pair of connected strata in $X$, this receives a partial order from the relation 
    \begin{equation}
        q \leq q' \quad \Leftrightarrow \quad X_q \subset \overline{X}_{q'}.  
    \end{equation}
    The poset $P$ receives a natural surjective poset map from $P^{\sf con}$, defined by sending a connected component $q$ to its respective stratum. This map of posets also has the special property of being \emph{discrete} — the fiber contains solely identity arrows. There is a sense in which this factorization of $X \ra P$ 
    \begin{equation}
        X \longra P^{\sf con} \longra P
    \end{equation}
    is initial.
\end{construction}

\begin{remark}
    If a stratification $X \ra Q$ is such that all of its strata are connected, we call it a \define{characteristic stratification} or the continuous map itself \define{characteristic}.
\end{remark}

\begin{lemma}
    The category of factorizations of $\big(X \ra P\big)$ as a characteristic map followed by a discrete map of posets, is contractible.
\end{lemma}
\begin{proof}
  This is proved in \cite[Lemma B.37]{DD21}.   
\end{proof}
\medskip

\nid This in particular means that whether we work with a $P$-stratification or its associated stratification by connected components, is essentially the same — we will often use this fact to stipulate stratifications merely by highlighting their individual connected components, without a direct mention of the continuous map to a poset $P$, or even the poset $P$.

\section{Conical differentiability and smoothness}

Unburdened by the desire or the need to be succinct, we will now introduce a thorough discussion on conical smoothness. The reason for doing this is twofold: the first is premised on our appreciation for the subtle, yet dizzying, inductive nature of its definition; the second on our perception of their usefulness as a toolkit to the future studies of topological field theories — potentially decorated with any assortment of other adjectives e.g., with \emph{topological defects}.  

\subsection{Conically smooth along the cone locus}
One might argue that the single most important contribution of Ayala--Francis--Tanaka in their foundational paper on stratified spaces (see \cite{AFT1}), was their introduction of the notion of conical smoothness along cone loci — a mild regularity condition that calls for an adaptation of the usual derivative for smooth manifolds — defined by the existence of a limit --, without appealing to the secant $k$-planes that come with working within an ambient manifold, as is the case for Whitney stratifications. This adaptation chooses instead to efficiently work with what is left of a vector space around every point in a $C^{0}$ stratified space. 

To initiate this discussion, let $(Z\ra P)$ be a compact topological stratified space and let $U = \RR^{i} \times \sC(Z)$ be a topological basic. So as to easily distinguish points in $U$, let us fix coordinates
\begin{equation}
    \big(u,[r,z]\big) \in U
\end{equation}
where $\big( u,r,z \big)\in \RR^{i} \times \RR_{\geq 0} \times Z$, and $[s,z]$ is the image of $(s,z)$ under the quotient map $\pi_Z \maps \RR_{\geq 0} \times Z \longra \sC(Z)$ collapsing the zero fiber onto the cone point. 

\begin{definition}
    Let $U = \RR^{i} \times \sC(Z)$ be a topological basic. We define its \textbf{tangent cone} to be the space
    \begin{equation}
        \sT U \,\overset{\sf def}{=}\, \sT\RR^{i} \times \sC(Z).
    \end{equation}
\end{definition}

Since $\RR^{i}$ is framed, we can re-write the factors in the tangent cone as follows:
\begin{equation}
    \sT \RR^{i} \times \sC(Z) \simeq (\RR^{i}_{\vec{v}} \times \RR^{i}) \times \sC(Z) \simeq \RR^{i}_{\vec{v}} \times \big(\RR^{i} \times \sC(Z)\big) = \RR^{i}_{\vec{v}} \times U.
\end{equation} 
\begin{observation}
This identification suggests we can use coordinates 
\begin{equation}
    \big(\vec{v},u,[s,z]\big),\quad \text{ with } \vec{v} \in \RR^{i}_{\vec{v}},\big(u,[s,z]\big)\in U.
\end{equation}
Now we explicate the sense in which the tangent cone captures some of what is left of a vector space in a basic neighborhood. In a genuine real vector space — the case, if $Z=\emptyset$ — one is allowed to translate vectors along any direction and to scale vectors by any scalar amount $\lambda \in \RR$. We can re-package this observation formally into a statement about how the group
\begin{equation}
    \mathsf{Aff}^{+}(\RR^{i}) = \big(\RR^{i} \rtimes \RR_{>0}\big)
\end{equation} 
freely acts on $\RR^{i}$, with the action determined by 
\begin{align}
    \alpha \maps (\RR^{i} \rtimes \RR_{>0}) \times \RR^{i} &\longra \RR^{i} \\
    (\vec{v},\lambda,u) &\longmapsto (\vec{v} + \lambda u).
\end{align}
Given an arbitrary continuous map $f\maps \RR^{i} \longra \RR^{j}$ and a choice of some point $u\in \RR^{i}$, this action is instrumental in defining the derivative of $f$ at $u$:
\begin{align}
    \label{ch4:: eq: quotient derivative}
    \sD f\maps \sT\RR^{i} \simeq \RR^{i}_{\vec{v}} \times \RR^{i} &\longra \sT\RR^{j} \simeq \RR^{j}_{\vec{v}} \times \RR^{j} \\
    \label{ch4:: eq: quotient derivative 2}
    (\vec{v},u) &\longmapsto \underset{\lambda \ra 0^+}{\mathsf{lim}} \Bigg(\frac{~~f\big(\alpha_{\lambda,\vec{v}}(u)\big)-f(u)~~}{\lambda},\,\, f(u)\Bigg)
\end{align}
If the limit does not exist for $(\vec{v},u)$, we say the derivative of $f$ at $u$ does not exist. We can do better — the action map $\alpha_{\lambda,\vec{v}}$ admits an inverse induced by the inverse in the group
\begin{equation}
    \alpha_{\lambda,\vec{v}}^{-1}(u) = \alpha_{\frac{1}{\lambda},-\frac{\vec{v}}{\lambda}}(u) = \frac{1}{\lambda}(u-\vec{v})
\end{equation}
and the limit in \ref{ch4:: eq: quotient derivative 2} can be re-written as 
\begin{equation}
     \sD f(\vec{v},u) \, = \, \bigg(\underset{\lambda \ra 0^+}{\mathsf{lim}}\,\, \alpha_{\lambda,f(u)}^{-1}\circ f \circ \alpha_{\lambda,\vec{v}},\,\,f(u)\bigg).
\end{equation}
Thus the derivative may be viewed as the continuous extension, along the inclusion $\RR^{i}\times \RR_{>0} \hookra \RR^{i}\times \RR_{\geq 0}$, of the \emph{conjugated action} on $f$.

The basic neighborhood $U=\RR^{i}\times \sC(Z)$, assuming that $Z$ is not-empty retains an action by $\mathsf{Aff}^{+}(\RR^{i})$. This is the key geometric observation permitting an immediate bootstrap of the conventional definition of differentiability into a context where maps require preservation of strata. The action we are indicating here is as follows:
\begin{align}
    \gamma \maps  \RR_{>0} \times \sT \RR^{i} \times \sC(Z) &\longra \RR_{>0} \times \sT \RR^{i} \times \sC(Z) \\
    (\lambda,u,\vec{v},[s,z]) &\longmapsto (\lambda, \lambda\vec{v} + u, u, [\lambda s,z]). 
\end{align}
This action captures the infinitesimal displacements which are permissible in the vicinity of a point sitting at the origin of a basic, in partial analogy with how the affine group acts on the smooth tangent bundle. 
\end{observation}

We illustrate the action $\gamma$ with two examples, with the intent to consolidate the exposition.\footnote{In both cases, color is used for mere visual aid, and is not in any sense fundamental or intrinsic to the data of a stratified space. Apologies in advance for the readers who might otherwise not have the ability to distinguish all colors.}

\begin{example}\label{ch4:: ex: scaling and translation}
    In the first example, we consider the basic $U=\RR \times \sC(S^{0}) \simeq \RR \times \sC(*_+ \sqcup *_{-})$ which we choose to represent as follows
    \begin{equation}
        \tikzfig{figures/Strat/STR2}{\input{figures/Strat/STR2.tikz}} \quad \simeq \quad \RR \,\, \times \,\, \tikzfig{figures/Strat/STR3}{\input{figures/Strat/STR3.tikz}}
    \end{equation}
    with a choice of identification $*_+$ with orange (the left patch) and $*_{-}$ with pink (the right patch). The \emph{cone loci} is identified with the green strand. Note of course that we implicitly \emph{rounded} the corners to obtain a disk. We can attempt to visualize the action $\gamma$ of $\RR^{i} \rtimes \RR_{>0}$: 
    \begin{enumerate}
        \item The action of scaling $\gamma_{\lambda, 0}$ along the radial direction can be depicted as follows
        \begin{equation}
         \tikzfig{figures/Strat/STR2}{\input{figures/Strat/STR2.tikz}}\,\,\,\, \xra{\,\,\,\,\gamma_{\lambda,0}\,\,\,\,}\,\,\,\, \tikzfig{figures/Strat/STR4}{\input{figures/Strat/STR4.tikz}} 
        \end{equation}
        where we have chosen $0< \lambda < 1$. Note also that this map is just the identity on the $\RR$ i.e., the Euclidean coordinate. 
        \item The action of translation $\gamma_{1,\vec{v}}$ along the cone locus can be depicted as follows
        \begin{equation}
            \tikzfig{figures/Strat/STR6}{\input{figures/Strat/STR6.tikz}}\,\,\,\, \xra{\,\,\,\,\gamma_{1,\vec{v}}\,\,\,\,} \,\,\,\, \tikzfig{figures/Strat/STR5}{\input{figures/Strat/STR5.tikz}} 
        \end{equation}
        where we have made a light rescaling to better indicate the image of the mapping.
    \end{enumerate}
    Noticeably, both of these maps are stratified and homeomorphisms. Note also that, despite the fact that the underlying space is homeomorphic to $\RR^2$, we dispense of any translation along a direction not in the cone locus — as that action would budge the strata \emph{off} the cone locus, and circumvent any interpretation as a conically smooth infinitesimal action. Finally, observe that neither action interacts with the $z$ coordinate. That is because we cannot scale or translate points of an arbitrary compact topological stratified space. 
\end{example}

\begin{example}
    As a second example, let us consider the open cone on a circle stratified by its north and south poles $U = \sC(\SS^{1})$ where $\SS^1 = (q \maps S^{1} \ra [1])$ and $q$ sends ${N,S}$ to $\{0\}$ and the remaining two open arcs to $\{1\}$, which we depict as follows
    \begin{equation}
        \tikzfig{figures/Strat/STR7}{\input{figures/Strat/STR7.tikz}} \,\,\,\, \simeq \,\,\,\, \text{open cone on } \,\,\tikzfig{figures/Strat/STR8}{\input{figures/Strat/STR8.tikz}}.
    \end{equation}
    Despite diverging from \Cref{ch4:: ex: scaling and translation} only in the fact that it contains a point stratum at the origin, the implication is that it no longer admits any non-trivial action by translation. The only action is by scaling along the radial coordinate.
    \begin{equation}
        \tikzfig{figures/Strat/STR7}{\input{figures/Strat/STR7.tikz}}\,\,\,\, \xra{\,\,\,\,\gamma_{\lambda,0}\,\,\,\,}\,\,\,\, \tikzfig{figures/Strat/STR9}{\input{figures/Strat/STR9.tikz}} 
    \end{equation}
    The cone locus has changed from being 1- to 0-dimensional. This severely restricts the kind of infinitesimal automorphisms allowed, and as described in \Cref{ch4:: obs: con smooth circles} it has strong implications for the comparison between the ordinary notion of $C^{1}$ and the weaker $C^{1}$ along the cone locus.
\end{example}

Immediately, we are led to consider a notion of differentiability for basics modeled on the continuous extension of the conjugated action. Let us define this explicitly.

\begin{definition}\cite[Definition 3.1.4]{AFT1}
    \label{ch4:: def: conically smooth along}
    Given a pair of compact stratified topological spaces $\big(Z\ra P\big)$ and $\big(W\ra Q\big)$ and a continuous stratified map 
    \begin{equation}
        f \maps \RR^{i} \times \sC(Z) \longra \RR^{j} \times \sC(W),
    \end{equation} 
    we say that \textbf{$f$ is $C^1$ along $\RR^{i}$} if\footnote{Often this is instead stated as \emph{along the cone locus}, since $\RR^{i}\times \{*\}$ contains all the points in a basic which are in the deepest stratum.}:
    \begin{itemize}
        \item \textbf{(Cone):} The underlying map of posets $\sC(P) \longra \sC(Q)$ maps the cone point $*_P$ to the cone point $*_Q$.
        \item \textbf{($C^1$ along):} If the solid diagram 
        \begin{equation}
            \begin{tikzcd}
                {\RR_{>0}\times \sT\RR^{i}\times \sC(Z)} &&& {\RR_{>0}\times \sT\RR^{j}\times \sC(W)} \\
                \\
                {\RR_{\geq 0}\times \sT\RR^{i}\times \sC(Z)} &&& {\RR_{\geq 0}\times \sT\RR^{j}\times \sC(W)}
                \arrow["{\gamma^{-1}\circ f_{\Delta}\circ \gamma}", from=1-1, to=1-4]
                \arrow[hook, from=1-1, to=3-1]
                \arrow[hook, from=1-4, to=3-4]
                \arrow["{\sD f}"', dashed, from=3-1, to=3-4]
            \end{tikzcd}
        \end{equation}
        admits a continuous extension $\sD f$.  
    \end{itemize}
    For $r>1$, we say that $f$ is \textbf{$C^{r}$ along $\RR^{i}$} if it is $C^1$ along $\RR^{i}$ and, under the identification
    \begin{equation}
        \sT\RR^{i}\times \sC(Z) \simeq \RR^{i}\times \RR^{i}\times \sC(Z),
    \end{equation}
    the derivative
    \begin{equation}
        \sD f \maps \RR^{i}\times \RR^{i}\times \sC(Z) \longra \RR^{j}\times \RR^{j}\times \sC(W)
    \end{equation}
    is $C^{r-1}$ along $\RR^{i}\times \RR^{i}$. We say that $f$ is \textbf{conically smooth along $\RR^i$} if it is $C^r$ along $\RR^i$ for every $r\geq 1$.
\end{definition}

\begin{observation}\label{ch4:: obs: con smooth circles}
Conical smoothness along $\RR^{i}$ for a given basic, is a weaker notion than that of smoothness, even in cases where the cone point admits a neighborhood homeomorphic to an open disk. To see this, we observe the example outlined in \cite[Example 3.1.9]{AFT1}: choosing polar coordinates $\sC(S^{n-1})\cong \RR^{n}$ we note that any arbitrary homeomorphism $f\maps S^{n-1} \ra S^{n-1}$ will induce a conically smooth map (but generally not smooth) $\sC(f)\maps \sC(S^{n-1}) \ra \sC(S^{n-1})$ at the cone point $*$. 

\medskip

\nid This is sensible, since all we are demanding is that the \emph{radial directional derivatives} exist and are continuous — we are taking the limit of the conjugation of $f$ exclusively by the radial action, since no other is available. Yet, this is not enough to ensure $\sC(f)$ is $C^1$ at the origin — this example is an indicator of the type of regularity we will work with as we venture into conically smooth stratified geometry.   
\end{observation}

\subsection{Stratified maps}

The full inductive construction of the category of basics and the sheaf of maximal atlases in \cite[Section~3.2]{AFT1} is somewhat elaborate. Since our aim is to proceed to the geometric applications, we record only the streamlined formulation that we shall use: first the recursive definition on basics, and then the chartwise definitions for atlases and arbitrary maps.

\begin{definition}[Conically smooth maps between basics]\cite[Definition~3.3.1]{AFT1}
    \label{ch4:: def: conically smooth map}
    We declare the unique map $\emptyset \to \emptyset$ to be conically smooth. Assume inductively that conically smooth maps have already been defined for basics of strictly smaller depth. Let
    \begin{equation}
        f \maps \RR^{i}\times \sC(Z) \longra \RR^{j}\times \sC(Y)
    \end{equation}
    be a continuous stratified map between basics. Then $f$ is \define{conically smooth} if:
    \begin{enumerate}
        \item \textbf{Hits the cone-point:} if, on the level of posets, $f$ sends the cone-point to the cone-point, then
        \begin{itemize}
            \item $f$ is conically smooth along $\RR^{i}$ in the sense of \Cref{ch4:: def: conically smooth along};
            \item the restriction
            \begin{equation}
                f_{|f^{-1}(\RR^{j}\times Y \times \RR_{>0})} \maps f^{-1}(\RR^{j}\times Y \times \RR_{>0}) \longra \RR^{j}\times Y \times \RR_{>0}
            \end{equation}
            is conically smooth.
        \end{itemize}
        \item \textbf{Misses the cone-point:} if, on the level of posets, $f$ misses the cone-point, then it factors as
        \begin{equation}
            \begin{tikzcd}[ampersand replacement=\&]
                {\RR^{i}\times \sC(Z)} \&\& {\RR^{j}\times \sC(Y)} \\
                \\
                \& {\RR^{j}\times Y \times \RR_{>0}} \& {.}
                \arrow["f", from=1-1, to=1-3]
                \arrow["{f_0}"', dashed, from=1-1, to=3-2]
                \arrow[hook, from=3-2, to=1-3]
            \end{tikzcd}
        \end{equation}
        and we say that $f$ is conically smooth if $f_0$ is conically smooth.
    \end{enumerate}
\end{definition}

\begin{definition}\cite[Definition~3.2.10]{AFT1}
    A \textbf{conically smooth atlas} on a $C^0$ stratified space $M$ is a collection of charts
    \begin{equation}
        \big\{\varphi_\alpha \maps U_\alpha \hookra M\big\}_\alpha
    \end{equation}
    where each $U_\alpha$ is a basic, the images $\{\varphi_\alpha(U_\alpha)\}_\alpha$ cover $M$, and whenever
    \begin{equation}
        x \in \varphi_\alpha(U_\alpha)\cap \varphi_\beta(U_\beta),
    \end{equation}
    there exists a basic $W$ and conically smooth open embeddings
    \begin{equation}
        W \hookra U_\alpha, \qquad W \hookra U_\beta
    \end{equation}
    whose composites with $\varphi_\alpha$ and $\varphi_\beta$ agree on their images, and this common image contains $x$. Here, a \textbf{conically smooth open embedding between basics} means a conically smooth map between basics whose underlying stratified map is a stratified open embedding. A \textbf{conically smooth stratified space} is a $C^0$ stratified space equipped with a maximal conically smooth atlas.
\end{definition}

\begin{definition}
    Let $f \maps \big(M,\mathcal{A}\big) \longra \big(N,\mathcal{B}\big)$ be a continuous stratified map between conically smooth stratified spaces. We say that $f$ is \define{conically smooth} if for every pair of charts
    \begin{equation}
        (\varphi\maps U \hookra M) \in \mathcal{A}, \qquad (\psi \maps V \hookra N) \in \mathcal{B}
    \end{equation}
    such that
    \begin{equation}
       f\big(\varphi(U)\big) \subset \psi(V),
    \end{equation}
    the induced map between basics
    \begin{equation}
        \psi^{-1} \circ f \circ \varphi \maps U \longra V
    \end{equation}
    is conically smooth in the sense of \Cref{ch4:: def: conically smooth map}.
\end{definition}

\begin{proposition}
    If $M\xra{\,\, f \,\,} N \xra{\,\, g \,\,} Q$ is a pair of composable conically smooth maps, then the composite is also a conically smooth map. 
\end{proposition}
\begin{proof}
    It suffices to check this result again at the level of basics.
\end{proof}

\begin{observation}
    Let $f \maps M \longra N$ be a conically smooth map, and let $N$ be trivially stratified, i.e., a smooth manifold. Per \Cref{ch4:: def: conically smooth map} we have that since, locally, $f$ (as a map of posets) preserves the cone-point $ U = \RR^{i} \times \sC(Z) \longra \RR^{j} = V$ for some basics $U$ of $M$ and $V$ of $N$, it suffices to check that $f$ is conically smooth along $\RR^{i}$.
\end{observation}

\begin{observation}
    A map of smooth manifolds, equivalently of depth $0$ conically smooth stratified spaces regarded as trivially stratified, is conically smooth. 
\end{observation}

\begin{definition}
    Let $\strat$ denote the ordinary category whose objects are conically smooth stratified spaces and whose morphisms are conically smooth maps. Let $\snglr \subset \strat$ denote the wide subcategory whose morphisms are the conically smooth maps whose underlying stratified maps are stratified open embeddings. Let $\bsc \subset \snglr$ denote the full subcategory spanned by basics.
\end{definition}

\begin{proposition}\cite[Corollary~3.4.10]{AFT1}
    The category $\strat$ admits products.
\end{proposition}

Finally, we note the following useful result.

\begin{lemma}\cite[Lemma 3.4.5]{AFT1}
    Let $\big(M \ra P, \mathcal{A}\big)$ be a conically smooth stratified space, and let $Q \hookra P$ be a consecutive inclusion of posets. The pullback $C^0$ stratified space 
    \begin{equation}
        \label{eq 2}
        \big(M_{|Q} \ra Q\big) \longra (M \ra P)
    \end{equation}
    inherits a conically smooth structure such that the map \Cref{eq 2} is also conically smooth.
\end{lemma}

\begin{cor}
    Let $\big(M \ra P,\mathcal{A}\big)$ be a conically smooth stratified space. Then we have that for every $p\in P$:
    \begin{enumerate}
        \item $M_{\leq p}$ is a conically smooth stratified subspace.
        \item $M_{p}$ is a smooth manifold.
        \item $M_{\nless p}$ is a conically smooth stratified subspace.
    \end{enumerate}
\end{cor}
\begin{proof}
    Note that for $(1)$ we consider $P^{\leq p} \hookra P$; which is consecutive because: it is a full subcategory of $P$ (spanned by all elements $x\leq p$.), and because given a pair $(x\leq p), (y\leq p) \in P_{\leq p}$, we have that every $z\in P$ such that $(x \leq z \leq y)$ also satisfies $ z\leq p$. For $(2)$ we have that $\{p\} \hookra P$ is trivially consecutive. Lastly, $(3)$ is consecutive because it is a full subcategory — it is the full subcategory spanned by all the objects $x \in P$ such that $x \geq p$ or $x$ is not related to $p$. Moreover, given a pair $x,y \in P_{\nless}$, if $(x \leq z \leq y)$, surely $z\in P_{\nless p}$.      
\end{proof}

\subsubsection{Classes of maps}

Conically smooth stratified maps, by virtue of the potential for differing behaviour along distinct strata, accommodate a strictly larger wealth of sub-classes when compared to the smooth unstratified ones.

\begin{definition}[Conically smooth submersion]
    \label{ch4:: def: con sm subm}
    Let $f\maps X \longra Y$ be a conically smooth map. We call it a \define{conically smooth submersion} if there is a collection of images $\{\varphi_\alpha(O_\alpha)\subseteq Y\}_\alpha$, indexed by commutative diagrams 
    \begin{equation}\label{eq:def:con sm submersion}
    \begin{tikzcd}[ampersand replacement=\&]
	{U_\alpha\times O_\alpha } \&\& X \\
	\\
	{O_\alpha} \&\& Y
	\arrow["{\psi_{\alpha}}", hook, from=1-1, to=1-3]
	\arrow["{\sf pr}"', from=1-1, to=3-1]
	\arrow["\lrcorner"{anchor=center, pos=0.125}, draw=none, from=1-1, to=3-3]
	\arrow["f", from=1-3, to=3-3]
	\arrow["{\varphi_\alpha}"', hook, from=3-1, to=3-3]
    \end{tikzcd}
    \end{equation}
    where the horizontal arrows are conically smooth open embeddings, and this collection constitutes a basis for the topology of $Y$.
\end{definition}

\begin{remark}\label{ch4:: rmk: con sm submersions}
    Note that \Cref{ch4:: def: con sm subm} looks different from the usual definition for submersive smooth maps. Namely, for categories of smooth manifolds and smooth maps therein, we record submersivity as a condition on the rank of the derivative map. It is a standard fact from differential topology that submersions are locally isomorphic to the standard projection
    \begin{equation}
    \RR^{n+k} \simeq \RR^n \times \RR^k\fib \RR^n.
    \end{equation}
    \Cref{ch4:: def: con sm subm} takes this characterization of submersions as their defining feature.
\end{remark}
\begin{observation}
    An implication of this definition is that, restricting to any stratum $X_p \subset X$, the induced map of smooth manifolds
    \begin{equation}
        f_{|X_p}\maps X_p \longra Y_{f(p)}
    \end{equation}
    is a smooth submersion. In particular, a conically smooth submersion among depth $0$ stratified spaces is an ordinary smooth submersion.
\end{observation}

\begin{definition}[Conically smooth fiber bundles]\label{ch4:: def: con sm fiber bundles}
    Let $f\maps X \longra Y$ be a conically smooth map. We call it a \define{conically smooth fiber bundle} if the collection of images $\varphi_\alpha(O_\alpha)\subseteq Y$, indexed by pullback diagrams 
    \begin{equation}\label{eq:def:fiber bundle}
    \begin{tikzcd}[ampersand replacement=\&]
	{F_\alpha\times O_\alpha } \&\& X \\
	\\
	{O_\alpha} \&\& Y
	\arrow["{\psi_{\alpha}}", hook, from=1-1, to=1-3]
	\arrow["{\sf pr}"', from=1-1, to=3-1]
	\arrow["\lrcorner"{anchor=center, pos=0.125}, draw=none, from=1-1, to=3-3]
	\arrow["f", from=1-3, to=3-3]
	\arrow["{\varphi_\alpha}"', hook, from=3-1, to=3-3]
    \end{tikzcd}
    \end{equation}
    where the horizontal arrows are conically smooth open embeddings, and this collection constitutes a basis for the topology of $Y$.
\end{definition}

\begin{example}
    In the figure below, we depict an example of a conically smooth fiber bundle $\pi \maps X \longra Y$ 
    \begin{equation}
        \tikzfig{figures/Strat/STR15}{\input{figures/Strat/STR15.tikz}}
    \end{equation}
    \medskip

    \nid For every point $y \in Y$ of the base, there is a neighborhood $y \in U \subset Y$ such that there is a conically smooth isomorphism
    \begin{equation}
        F \times U \, \simeq  \, \pi^{-1}(U) \simeq \,  U \times_{Y} X,
    \end{equation} where $F$ is the discrete space with two points.
\end{example}

\begin{example}
    Perhaps just as useful to the introduction of conically smooth fiber bundles, is the indication of a good counter-example. The projection map $\pi \maps X \longra Y$, depicted in the diagram below
    \begin{equation}
        \tikzfig{figures/Strat/STR16}{\input{figures/Strat/STR16.tikz}}
    \end{equation}
    \medskip 

    \nid is not a conically smooth fiber bundle. Any conically smooth open embedding $U \hookra Y$ which contains the 0-stratum in the base, will not be in the basis for a topology of $Y$, since restricting to $U$ does not yield a stratified projection of a fixed stratified fiber $F$:
    \begin{equation}
        \tikzfig{figures/Strat/STR17}{\input{figures/Strat/STR17.tikz}}
    \end{equation}
    \medskip
    
    \nid Modifying the total space   
    \begin{equation}
        \tikzfig{figures/Strat/STR17}{\input{figures/Strat/STR17.tikz}} \quad \squigra \quad \tikzfig{figures/Strat/STR18}{\input{figures/Strat/STR18.tikz}} 
    \end{equation}
    \medskip

    \nid turns it into a conically smooth fiber bundle with fiber 
    \begin{equation}
     \tikzfig{figures/Strat/STR19}{\input{figures/Strat/STR19.tikz}} \quad .   
    \end{equation}
\end{example}

\begin{definition}[(Weakly) Constructible bundles]\label{ch4:: def: (weak) cbl}
    Let $f \maps \big(X \ra P \big) \longra \big(Y \ra Q \big)$ be a conically smooth stratified map. We call it a \define{weakly constructible bundle} if, for each $q\in Q$, the vertical map in the pullback
    \begin{equation}
        \begin{tikzcd}[ampersand replacement=\&]
            {X_q} \&\& X \\
            \\
            {Y_q} \&\& Y
            \arrow[hook, from=1-1, to=1-3]
            \arrow[from=1-1, to=3-1]
            \arrow["\lrcorner"{anchor=center, pos=0.125}, draw=none, from=1-1, to=3-3]
            \arrow["f", from=1-3, to=3-3]
            \arrow[hook, from=3-1, to=3-3]
        \end{tikzcd}
    \end{equation} 
    is a conically smooth fiber bundle (\Cref{ch4:: def: con sm fiber bundles}). The map $f$ is said to be a \define{constructible bundle} if:
    \begin{itemize}
        \item For $\depth(X) = 0$, the map $f$ is a smooth fiber bundle of smooth manifolds \footnote{If $\depth(X)=0$ even if $\depth(Y)>0$ the map must factor through a single stratum of $Y$.}
        \item For $\depth(X)>0$, then:
        \begin{itemize}
            \item $f$ is a weakly constructible bundle.
            \item For each $q\in Q$ the map 
            \begin{equation}
                \link_{X_{q}}(X) \longra X_q \underset{Y_{q}}{\times} \link_{Y_{q}}(Y)
            \end{equation}
            is a constructible bundle.
        \end{itemize}
    \end{itemize}
\end{definition}

\begin{observation}
Note that for any given sub-stratified space $Y \hookra X$, the conically smooth stratified space $\link_{Y}(X)$ has strictly less depth than $X$, which validates the inductive step.
\end{observation}

\begin{example}
    Let $f\maps X \longra Y$ be the projection 
    \begin{equation}
    \sC(\scalebox{0.3}{\tikzfig{figures/Strat/STR28}{\input{figures/Strat/STR28.tikz}}}) \longra \RR_{\geq 0}
    \end{equation} 
    which we diagramatically depict as 
    \begin{equation}
        \tikzfig{figures/Strat/STR29}{\input{figures/Strat/STR29.tikz}}
    \end{equation}
    \medskip

    \nid This is a weakly constructible bundle since for every stratum of the base, the pullback is a conically smooth fiber bundle. The fiber over $\{0\}$ is the cone point and the fiber over any point in the interior stratum $\RR_{>0}$ is an interval stratified by the boundary and a single $0$-stratum in the interior.
\end{example}

A key fact, and a hard result to prove, is that weakly constructible bundles, unlike constructible bundles, do not in general pull back along arbitrary conically smooth stratified maps.

\begin{theorem}[Pullback stability]\label{ch4:: thm: cbl pullback}\cite[Lemma~6.1.11]{AFR1}
    Let $f \maps X \longra Y$ be a constructible bundle and $g\maps Y' \longra Y$ an arbitrary conically smooth stratified map, the pullback 
    \begin{equation}
        \begin{tikzcd}[ampersand replacement=\&]
            {X \underset{Y}{\times} Y'} \&\& X \\
            \\
            {Y'} \&\& Y
            \arrow[from=1-1, to=1-3]
            \arrow["{f'}"', from=1-1, to=3-1]
            \arrow["\lrcorner"{anchor=center, pos=0.125}, draw=none, from=1-1, to=3-3]
            \arrow["f", from=1-3, to=3-3]
            \arrow["g"', from=3-1, to=3-3]
        \end{tikzcd}
    \end{equation}
    exists in the category $\strat$. The map $f'\maps X\underset{Y}{\times} Y' \longra Y'$ is again a constructible bundle. It is also proper if $f$ is proper.
\end{theorem}

\begin{definition}[Refinement]\label{ch4:: def: refinement}
    A map $f\maps X \longra Y$ is a \define{refinement} if it is a homeomorphism of underlying topological spaces, and for each $p\in P$
    \begin{equation}
        f_{|X_p}\maps X_p \longra Y_{f(p)}
    \end{equation} 
    is an isomorphism onto its image as smooth manifolds.
\end{definition}

\begin{example}
    This diagram depicts a refinement:
    \begin{equation}
        \tikzfig{figures/Strat/STR25}{\input{figures/Strat/STR25.tikz}} \;\;\xra{\,\,\, f \,\,\,}\;\; \tikzfig{figures/Strat/STR26}{\input{figures/Strat/STR26.tikz}}
    \end{equation}
    \medskip

    \nid which is the identity on the underlying topological space. Since it is stratified, it implies that $f_{|X_{p}}$ is an isomorphism onto its image. 
\end{example}

\subsection{Differential topology of conically smooth stratified spaces}
In this subsection we extend classical results from differential topology from smooth manifolds to conically smooth stratified spaces. The key constructions here will be: 
\begin{enumerate}
    \item The $(\oo,1)$-category $\Strat$, a suitable  $\mathsf{Kan}$-enrichment of $\strat$ .
    \item Categories of basics.
    \item Conically smooth stratified spaces with corners.
    \item Unzipping: a suitably functorial way of resolving singularities, akin to \emph{blow-up} in algebraic geometry.
    \item Tubular neighborhoods around strata and constructible sub-stratifications.
\end{enumerate}

One of the features of smooth manifolds which sets them apart from algebraic varieties, for example, is the ability to patch together local data and constructions into global ones using partitions of unity. These are likewise available for conically smooth stratified spaces.

\begin{lemma}[Conically smooth partitions of unity]\label{lem:partitions of unity}\cite[Lemma 7.1.1]{AFT1}
    Every conically smooth stratified space admits a conically smooth partition of unity.
\end{lemma}

\subsubsection{\texorpdfstring{Topologizing $\strat$}{Topologizing strat}}

For finer details about this construction, see Sections~2.1--2.4 in \cite{AFR1} — we will give only a schematic overview. 

We observe that the $\oo$-category of topological spaces can be defined as the Dwyer--Kan localization of $\Top$ along the set of homotopy equivalences. This is a strategy we can aptly reproduce in the case of conically smooth stratified spaces if we first define what is meant by \define{stratified homotopy equivalence}. 

\begin{definition}[Stratified homotopy equivalence]\label{ch4:: def: stratified homotopy equivalence}\cite[Definition~2.1.1]{AFR1}
    Given a pair of conically smooth stratified maps $f,g\maps X \longra Y$, a \define{stratified homotopy} 
    \begin{equation}
        H \maps X \times \RR \longra Y \times \RR
    \end{equation}
    is a conically smooth stratified map $H$ for which the following commutative diagram commutes:
    \begin{equation}
        \begin{tikzcd}[ampersand replacement=\&]
            X \&\& Y \\
            {X\times \RR} \&\& {Y \times \RR} \\
            X \&\& Y
            \arrow["f", from=1-1, to=1-3]
            \arrow["{\{0\}}"', from=1-1, to=2-1]
            \arrow["{\{0\}}", from=1-3, to=2-3]
            \arrow["H"{description}, from=2-1, to=2-3]
            \arrow["{\{1\}}", from=3-1, to=2-1]
            \arrow["g"', from=3-1, to=3-3]
            \arrow["{\{1\}}"', from=3-3, to=2-3]
        \end{tikzcd}
    \end{equation}
    The map $f$ is said to be a \define{stratified homotopy equivalence} if there exists a conically smooth stratified map $g$ together with stratified homotopies from $gf$ and $fg$ to $\id_{X}$ and $\id_{Y}$, respectively. Let us denote with $\cJ$ the collection of all stratified homotopy equivalences. 
\end{definition}

Our first pass at a definition is the following.

\begin{definition}
    Let 
    \begin{equation}
    \strat[\cJ^{-1}]
    \end{equation} 
    be the Dwyer--Kan localization of $\strat$ along the stratified homotopy equivalences. We denote the associated $(\oo,1)$-category by $\Strat$.
\end{definition}

\begin{observation}
This definition captures the correct homotopy theory of conically smooth stratified spaces, but as with most naive localizations, it does not offer a concrete model — an equivalent $(\oo,1)$-category modeled after a simplicially enriched category, for example.   
\end{observation}

A first simplification is to localize along a narrower class of morphisms.

\begin{definition}
    Let $X \times \RR \xra{\; \mathsf{pr}\;} X$ be the canonical projection. Define $\cR$ to be the collection of all projections for all $X\in \strat$. This collection is manifestly contained in the stratified homotopy equivalences, i.e.
        \begin{equation}
        \cR \subset \cJ.
        \end{equation}   
\end{definition}

\begin{lemma}\cite[Lemma~2.1.7]{AFR1}
    The canonical map between Dwyer--Kan localizations 
    \begin{equation}
        \strat[\cR^{-1}] \longra \strat[\cJ^{-1}]
    \end{equation}
    is an equivalence of $(\oo,1)$-categories.
\end{lemma}
\begin{proof}[Proof sketch]
    Let $F \maps \strat \longra \cC$ be a functor which sends every projection in $\cR$ to an equivalence in $\cC$. For each $X\in \strat$ and each endpoint inclusion $\{i\}\maps X \hookra X\times \RR$ with $i\in \{0,1\}$, the composite
    \begin{equation}
        X \xra{\;\{i\}\;} X\times \RR \xra{\;\mathsf{pr}\;} X
    \end{equation}
    is the identity. Since $F(\mathsf{pr})$ is an equivalence, the 2-out-of-3 property implies that $F(\{i\})$ is also an equivalence.

    Now let $H\maps X\times \RR \longra Y\times \RR$ be a stratified homotopy between $f,g\maps X\longra Y$. By definition of stratified homotopy, we have
    \begin{equation}
        f = \mathsf{pr}_Y \circ H \circ \{0\}, \qquad g = \mathsf{pr}_Y \circ H \circ \{1\}.
    \end{equation}
    After applying $F$, the morphisms $F(f)$ and $F(g)$ differ only by replacing the endpoint inclusion $F(\{0\})$ with $F(\{1\})$, and both of these are equivalences. Hence $F(f)$ is an equivalence if and only if $F(g)$ is.

    Finally, if $f\maps X\longra Y$ is a stratified homotopy equivalence with homotopy inverse $g\maps Y\longra X$, then the homotopies $g\circ f \simeq \id_X$ and $f\circ g \simeq \id_Y$ imply that $F(g)\circ F(f)$ and $F(f)\circ F(g)$ are equivalent to identities. Therefore $F(f)$ is itself an equivalence. This shows that every functor inverting $\cR$ also inverts $\cJ$, which is exactly the required universal property.
\end{proof}

We regard $\strat$ as a site by declaring a family of morphisms $\{U_i\hookra X\}_{i\in I}$ to be a covering family if each $U_i\hookra X$ is a conically smooth open embedding and the opens $\{U_i\}_{i\in I}$ cover the underlying topological space of $X$.
Write
\[
  \psh(\strat)\;:=\;\Fun(\strat^{\op},\Spaces)
\]
for the $(\oo,1)$-category of presheaves, and
\[
  \shv(\strat)\;\subset\;\psh(\strat)
\]
for the full subcategory of (hypercomplete) $\infty$-sheaves for this Grothendieck topology.

\begin{definition}
    The $\oo$-category of \define{constructible sheaves} on $\strat$ is the full $(\oo,1)$-subcategory 
    \begin{equation}
        \shv^{\cbl}(\strat) \; \subset \; \shv(\strat).
    \end{equation}
\end{definition}

In particular, if we restrict a constructible sheaf $\cF$ to the subsite $\man \hookra \strat$ of smooth manifolds we obtain a concordance-invariant sheaf. Contrary to the case of smooth manifolds, the concordance resolution 
\begin{equation}
    \psh(\strat) \longra \psh(\strat), \quad \cF \longmapsto \left| \cF(-\times \AA^{\bullet}) \right|
\end{equation}
does not, in general, report back a constructible sheaf. In Section 2.3 of \cite{AFR1} the authors provide a class of sheaves whose concordance resolution (alias \emph{topologizing diagram} in \cite{AFR1}) is constructible — see \cite[Theorem~2.3.14]{AFR1}. These are called isotopy sheaves, and as demonstrated in the following definition, have a very well-defined list of checkable criteria.
\begin{definition}\cite[Definition~2.3.8]{AFR1}
    The category of \define{isotopy sheaves} is the full subcategory 
    \begin{equation}
        \mathsf{Isot} \subset \Cat_{/\strat}
    \end{equation}
    consisting of right fibrations $\cI \longra \strat$ whose straightening $\cI \maps \strat^{\op} \longra \mathsf{Gpd}$ satisfies the following conditions:
    \begin{enumerate}
        \item \textbf{Sheaf:} The straightening $\cI\maps \strat^{\op} \longra \mathsf{Gpd}$ is a sheaf. 
        \item \textbf{Isotopy extension:} For each stratified space $K$, the following conditions are satisfied:
        \begin{enumerate}
            \item For each $q\geq 0$, the functor between groupoids 
            \begin{equation}
                \cI(K \times \AA^q) \longra \cI(K \times \partial \AA^q)
            \end{equation}
            is an isofibration.
            \item For each $0\leq i \leq q$, the functor between groupoids
            \begin{equation}
                \cI(K \times \AA^q) \longra \cI(K \times (\Lambda^q)_i)
            \end{equation}
            is surjective on objects and on Hom-sets.
        \end{enumerate}
    \end{enumerate}
\end{definition}

There is a limit-preserving natural inclusion $\shv^{\cbl}(\strat) \xlonghookra{\;\;} \shv(\strat)$ which admits a left adjoint:
\begin{equation}
    L\maps \shv(\strat) \longra \shv^{\cbl}(\strat).
\end{equation}
When restricted to the subcategory $\mathsf{Isot}\subset \shv(\strat)$, the value of this functor is not only a constructible sheaf, it is its associated concordance resolution.
\begin{proposition}\cite[Lemma~2.3.16]{AFR1}\label{ch4:: def: isotopy sheaf realizes to cbl}
    If $\cI$ is an isotopy sheaf, then there is a canonical natural transformation 
    \begin{equation}
        L\cI \longra \left| \cI(-\times \AA^\bullet) \right|,
    \end{equation}
    and it is an equivalence.
\end{proposition} 
We can use this result to provide a $\mathsf{Kan}$-enrichment of $\strat$ as follows:
\begin{enumerate}
    \item For every $K\in \strat$, the forgetful functor 
    \begin{equation}
        \strat_{/K} \longra \strat
    \end{equation} 
    which straightens to the representable $\strat(-,K)$, is an isotopy sheaf. The verification uses the existence of regular neighborhoods for constructible subspaces, a fact we establish later in the chapter.
    \item  The functor
    \begin{equation}
        \strat \longra \shv(\strat), \quad K \longmapsto \strat(-,K)
    \end{equation}
    lands in $\mathsf{Isot}$. Applying $L$, produces a functor 
    \begin{equation}
        \strat \longra \shv^{\cbl}(\strat), \quad K \longmapsto L\strat(-,K).
    \end{equation}
\end{enumerate}
These observations lead us to following definition.
\begin{definition}\cite[Definition~2.4.2]{AFR1}
    The $(\oo,1)$-category $\Strat$ of conically smooth stratified spaces is the essential image of 
    \begin{equation}
    \strat \longra \shv^{\cbl}(\strat), \qquad K \longmapsto L\strat(-,K).
    \end{equation}
\end{definition}
The $\mathsf{Kan}$-enrichment is then derived from \Cref{ch4:: def: isotopy sheaf realizes to cbl} as follows: given a pair of conically smooth stratified spaces $X,Y\in \strat$, the associated simplicial set of morphisms from $X$ to $Y$ is computed by 
\begin{equation}
    \Hom_{\shv^{\cbl}}\left(L\strat(-,X),L\strat(-,Y)\right).
\end{equation} 
Using the adjunction $L \dashv i$, we have that 
\begin{align}
    \Hom_{\shv^{\cbl}}\left(L\strat(-,X),L\strat(-,Y)\right) &\simeq \Hom_{\shv}\left(\strat(-,X),iL\strat(-,Y)\right) \\
    &\simeq iL\strat(-,Y)(X) \\
    &= L\strat(-,Y)(X).
\end{align}
From \Cref{ch4:: def: isotopy sheaf realizes to cbl} we have that 
\begin{equation}
    L\strat(-,Y) \simeq \left|\strat(-\times \AA^\bullet,Y)\right|,
\end{equation}
which implies that 
\begin{equation}
    \Hom_{\shv^{\cbl}}(X,Y) \,\, \simeq \,\, L\strat(-,Y)(X) \,\, \simeq \,\, \left|\strat(X\times \AA^{\bullet},Y)\right|,
\end{equation} where we are interpreting $X$ and $Y$ as their images under the localization.

\begin{theorem}\cite[Theorem 2.4.5]{AFR1}\label{ch4:: thm: Strat enrichment}
    The $(\oo,1)$-category obtained through the localization $\strat[\cJ^{-1}]$ is equivalent to the $(\oo,1)$-category $\Strat$ obtained as the $\mathsf{Kan}$-enrichment of $\strat$ via:
    \begin{equation}\label{eq: st enrich}
        \Hom_{\Strat}(X,Y) \,\, \overset{\sf def}{=} \,\, \Strat(X\times \AA^{\bullet},Y).
    \end{equation}
\end{theorem}

\begin{observation}
    We could equivalently describe the $\mathsf{Kan}$-enrichment~(\ref{eq: st enrich}) of $\strat$ by 
\begin{equation}
    \Strat(X,Y) \,\, \overset{\sf def}{=}\,\, \strat_{/\AA^{\bullet}}(X\times \AA^{\bullet},Y\times \AA^{\bullet}).
\end{equation}
\end{observation}

\begin{notation}
    We denote by $\Bsc$ and $\Snglr$ as the $\mathsf{Kan}$-enriched subcategories of $\Strat$ whose objects are those of $\bsc$ and $\snglr$, respectively. In particular, we will refer to these not as $\mathsf{Kan}$-enrichments but as $(\oo,1)$-categories.
\end{notation}

\medskip

\medskip

\subsubsection{Categories of basics}

In \Cref{ch3:: sec: tangential structures} we discussed the notion of tangential structures for smooth manifolds. In this section we aim to introduce an analog of $\BO(d)$ and the tangent classifier $M \longra \BO(d)$, i.e.\ the map which classifies the tangent bundle of the $d$-dimensional smooth manifold $M$. The goal being to introduce tangential structures and singularity types.

As is observed in the introduction of \cite{AFT1}, a space over $\BO(d)$, which we can assume to be a Kan fibration, is equivalent data to that of a functor $\BO(d)^{\op} \longra \Top$. Through the marvelous observation that 
\begin{equation}
    \sO(d) \longra \Emb_{\BO(d)}(\RR^{d},\RR^{d})
\end{equation} 
is a weak homotopy equivalence — a combination of the Gram--Schmidt orthogonalization and taking derivatives at the origin — we have that 
\begin{equation}
    \psh(\sD_{d}) \, \simeq \, \Top_{/\BO(d)}
\end{equation}
between presheaves on $\sD_d$ and the spaces mapping into $\BO(d)$. Here we are defining $\sD_d$ as $\man^{\emb}$ restricted to the $d$-dimensional Euclidean spaces. 

Recall that every $p$-stratum of a conically smooth stratified space $(X \ra P)$ is a $d_{p}$-dimensional smooth manifold. Therefore, a suitable proposed generalization of $\sD_d$ is $\Bsc$, the $(\oo,1)$-category of all basics and conically smooth open embeddings thereamong.

\begin{definition}\cite[Definition~1.1.3]{AFT1}
    The \define{tangent classifier} is the composite functor
    \begin{equation}
        \tau \maps \Snglr \xra{\,\,\cY\,\,} \psh(\Snglr) \xra{\,\,\iota^*\,\,} \psh(\Bsc).
    \end{equation}
\end{definition}

\begin{remark}
To a conically smooth stratified space $M$, the functor $\tau$ associates to it the presheaf
\begin{equation}
    \Snglr(-,M)\maps \Bsc^{\op} \longra \Top
\end{equation}
\end{remark}

\begin{observation}
Restricting the tangent classifier $\tau$ to $\Man_{d}^{\emb}$ it is proved in \cite[Corollary~4.4.9]{AFT1} that its image lands in $\psh(\sD_d) \subset \psh(\Bsc)$. More concretely, it is the usual tangent classifier for smooth manifolds which classifies the tangent bundle. 
\end{observation}

\begin{remark}
    Let $\entr(M \ra P)$ denote the $(\oo,1)$-category which is the unstraightening of $\tau(M)$. Unpacking this we note that since 
    \begin{equation}
        \tau(M\ra P) = \Snglr(\iota(-),M)\maps \Bsc^{\op} \longra \Top,
    \end{equation}
    then its unstraightening is the $(\oo,1)$-overcategory 
    \begin{equation}
        \Bsc_{/M}.
    \end{equation}
    We denote this $(\oo,1)$-category as the \define{enter-path category} of $M$. The enter-path category comes with a natural projection functor 
    \begin{equation}
        \Bsc_{/M} \longra \Bsc
    \end{equation}
    which is a right fibration.
\end{remark}

\begin{example}
    Let us consider a simple example, let $M= \sC(* \sqcup *)$
    \begin{equation}
       M \; = \; \tikzfig{figures/Strat/STR30}{\input{figures/Strat/STR30.tikz}}\, .
    \end{equation}
    Objects in $\entr(M)$ are of essentially three types:
    \begin{itemize}
        \item Open embeddings of open intervals into the left of the origin
        \begin{equation}
            \tikzfig{figures/Strat/STR31}{\input{figures/Strat/STR31.tikz}}.
        \end{equation} 
        \medskip
        \item Open embeddings of open intervals into the right of the origin
        \begin{equation}
            \tikzfig{figures/Strat/STR32}{\input{figures/Strat/STR32.tikz}}.
        \end{equation} 
        \medskip
        \item Open embeddings of $M$ into itself
        \begin{equation}
            \tikzfig{figures/Strat/STR33}{\input{figures/Strat/STR33.tikz}}.
        \end{equation} 
        \medskip
    \end{itemize}
    Morphisms are also simple. There are no morphisms from objects of the first type into the second type, and vice-versa. There is essentially one morphism from an object of the first or second type into one of the third type, depicted below as follows:
\end{example}

This lead us naturally to the following generalization of $B$-framed manifolds.

\begin{definition} \cite[Definition~5.0.2]{AFT1} \label{ch4:: def: categories of basics}
    An $(\oo,1)$\define{-category of basics} is a right fibration $\cB \longra \Bsc$. Given such an $(\oo,1)$-categoy of basics $\cB$, the $(\oo,1)$-category of $\cB$-manifolds is the pullback $(\oo,1)$-category 
    \begin{equation}
        \begin{tikzcd}[ampersand replacement=\&]
            {\Mfld(\cB)} \&\& {\left( \rfib_{\Bsc} \right)_{/\cB}} \\
            \\
            \Snglr \&\& {\rfib_{\Bsc}}
            \arrow[from=1-1, to=1-3]
            \arrow[from=1-1, to=3-1]
            \arrow["\lrcorner"{anchor=center, pos=0.125}, draw=none, from=1-1, to=3-3]
            \arrow[from=1-3, to=3-3]
            \arrow[from=3-1, to=3-3]
        \end{tikzcd}
    \end{equation} 
\end{definition}

\begin{remark}
    Unpacking definition \Cref{ch4:: def: categories of basics} we note the analogies with the smooth case. Instead of a space $B$ and a continuous map $B \longra \BO(d)$ we now have an $(\oo,1)$-category $\cB$ and a right fibration $\Bsc$. If $\cB = B$ is an $\oo$-groupoid, then a right fibration $B \longra \Bsc$ must factor through an $\oo$-subgroupoid of $\Bsc$. Inspecting such a subgroupoid, we observe that it must be a copy of $\BO(d)$ for some $d$, thereby demonstrating that it is a stratified augmentation of the smooth unstratified case.

    A conically smooth stratified space $M$ with a $\cB$-manifold structure amounts to a morphism — a natural transformation to be more precise — in $\rfib_{\Bsc}$ with source and target $\entr(M) \longra \Bsc$ and $\cB \longra \Bsc$, respectively
    \begin{equation}
        \begin{tikzcd}[ampersand replacement=\&]
            \&\& \cB \\
            \\
            {\entr(M)} \& {} \& \Bsc
            \arrow[between={0.5}{0.9}, Rightarrow, from=1-3, to=3-2]
            \arrow[from=1-3, to=3-3]
            \arrow["g", dashed, from=3-1, to=1-3]
            \arrow[from=3-1, to=3-3]
        \end{tikzcd}
    \end{equation}
\end{remark}

\begin{observation}
As we will see, the notion of $\cB$-manifold incubates not just the more typical notion of tangential structure, like orientations and spin structures, but also \emph{singularity types} — the latter being a specification of the kinds of links allowed for any given stratum.
\end{observation}

Specifying an $(\oo,1)$-category of basics can be split into two sequential steps:

\begin{enumerate}
    \item Specifying \emph{singularity types} which amounts to specifying a fully-faithful functor 
    \begin{equation}
        \cL \xlonghookra{\, \,\,\,} \Bsc.
    \end{equation}
    \item Specifying an essentially sujective functor 
    \begin{equation}
        \cK \longra \cL.
    \end{equation}
\end{enumerate}
This sequence is sustained by the fact a right fibration (any functor between $(\oo,1)$-categories, in fact) (see \cite[Lemma~2.8.3]{Run25}) admits a factorization into a fully-faithful functor composed with an essentially surjective functor.

Let us begin with singularity types.
\begin{definition}
    Let $\cC$ be an $(\oo,1)$-category and $\cL \subseteq \cC$ an $(\oo,1)$-subcategory. $\cL$ is a \define{sieve} if $\cL \hookra \cC$ is a right fibration. This is equivalent to the following checkable condition:
    \begin{itemize}
        \item \textbf{(Sieve condition)}] \label{con: sieve condition} For each $c \in \cL$, the space 
        \begin{equation}
            \Hom_{\cC}(c',c) \neq \emptyset
        \end{equation}
        implies that $c' \in \cL$.
    \end{itemize}
\end{definition}

\begin{remark}
    Checking a $\mathsf{Kan}$-enriched full subcategory is a right fibration can be done by checking the sieve condition.  
\end{remark}

\begin{example}
    If $\cB = \sD_d$, then $\cD_d \longra \Bsc$ is simply a sieve, hence a right fibration. A $\sD_d$-manifold is a smooth $d$-dimensional manifold — we can interpret this as saying that only a single ``singularity type'' is available, the smooth $d$-dimensional one, or the absence of a singularity.
\end{example}

\begin{example}
    The full subcategory $\Bsc_{\leq k,n} \subset \Bsc$ is a sieve. More concretely, given a basic $U$ of depth $l \leq k$ and pure dimension $n$, then any basic $V$ admitting an open embedding $V \hookra U$ must be of depth $l'$ and pure dimension and $l'\leq l \leq k$. 
\end{example}

\begin{example}
    We now present an example of a subcategory which is not a sieve. Let $\Bsc_1$ be the subcategory $\Bsc_1 \subset \Bsc$  dimension 1 basics. Considering the full subacategory spanned by $\cL = \{\sC(\{n\})\}_{n \geq 1}$, we fail to construct a sieve since 
    \begin{equation}
        \Hom_{\Bsc_1}(\RR, \sC(\{n\})) \neq \emptyset
    \end{equation}
    and $\RR \notin \cL$, as an example
    \begin{equation}
        \tikzfig{figures/Strat/STR38}{\input{figures/Strat/STR38.tikz}}
    \end{equation}
    we indicate in blue the image of a conically smooth open embedding $\RR \hookra \sC({1,2,3})$.
\end{example}

\begin{example}
    A typical example considered in mathematical physics and knot theory is that of a single $d$-dimensional proper submanifold $L\subset M$ of an $n$-dimensional smooth manifold $M$. We would like to consider $L$ as a singularity locus, so we would like it to form a class of $\cB$-manifolds for some $\cB$. To do so, we focus on the basic neighborhoods covering points in $L\subset M$ and $M$, respectively. 
    
    Let us define $\sD_{d\subset n}$ by choosing only two objects of $\Bsc$, namely:
    \begin{equation}
        \mathsf{Ob}(\sD_{d\subset n}) = \{\RR^{n}, \RR^{n-d}\times \sC(S^{n-d-1})\}.
    \end{equation}
    The morphism Kan complexes are as follows:
    \begin{itemize}
        \item $\sD_{d\subset n}(\RR^{n},\RR^{n}) = \mathsf{Emb}(\RR^{n},\RR^{n})$.
        \item $\sD_{d\subset n}(\RR^{n},\RR^{n-d}\times \sC(S^{n-d-1})) = \mathsf{Emb}(\RR^{n},\RR^{n}\setminus \RR^{d}) \simeq \sO(n) \times S^{n-d-1}$.
        \item $\sD_{d\subset n}(\RR^{n-d}\times \sC(S^{n-d-1}),\RR^{n}) = \emptyset.$
        \item $\sD_{d\subset n}(\RR^{n-d}\times \sC(S^{n-d-1}),\RR^{n-d}\times \sC(S^{n-d-1}) \, \simeq \, \sO(d) \times \Aut(S^{n-d-1})$.
    \end{itemize}
    A $\sD_{d\subset n}$-manifold is thus the data of an smooth $n$-manifold together with an embedded proper $d$-submanifold $L \subset M$. 
\end{example}

Finally, we look into the specifiation of tangential structures. 

\begin{example}
Let us define $\sD_{n}^{\fr} = (\sD_{n})_{/\RR^{n}} \longra \sD_{n}$. This functor unstraightens to $\sD_{n}^{\op} \longra \mathsf{Kan}$ given by $\sD_{n}(-,\RR^{n})$ — a fiber over $\RR^{n}$ is therefore equivalent to a point of $\mathsf{Emb}(\RR^n, \RR^n)$ or a framing on $\RR^{n}$. A $\sD_n^{\sf fr}$-manifold is therefore a smooth $n$-manifold together with a framing or a trivialization of its tangent bundle. 
\end{example}

\begin{example}
    We define $\sD_{d\subset n}^{\fr}$ as the pullback $(\oo,1)$-category 
    \begin{equation}
        \begin{tikzcd}[ampersand replacement=\&]
            {\sD_{d\subset n}^{\fr}} \&\& {\sD_{ n}^{\fr}} \\
            \\
            {\sD_{d\subset n}} \&\& {\sD_{n}}
            \arrow[from=1-1, to=1-3]
            \arrow[from=1-1, to=3-1]
            \arrow["\lrcorner"{anchor=center, pos=0.125}, draw=none, from=1-1, to=3-3]
            \arrow[from=1-3, to=3-3]
            \arrow[from=3-1, to=3-3]
        \end{tikzcd}
    \end{equation}
    For further unpacking of this construction see \cite[Example~5.2.12]{AFT1}. A $\sD_{d\subset n}^{\fr}$-manifold is comprised of:
    \begin{enumerate}
        \item A framed $n$-manifold $M$.
        \item A properly embedded smooth $d$-submanifold $L$ equipped with a null-homotopy of the map
        \begin{equation}
            L \longra \mathsf{Gr}_{d}(\RR^n)
        \end{equation}
        classifying the tangent bundle as its embedded in $\sT M$ — also called its \emph{Gauss map}. To specify such a null-homotopy one can stipulate a trivialization of both the tangent and normal bundle of $L$ in a way compatible with the underlying trivialization of the tangent bundle of $M$. 
    \end{enumerate}

    As a simple example of a $\sD_{1\subset 2}^{\fr}$ one might consider $\RR^{2}$ with the standard framing and the $x$-axis $L = \{y = 0\}$ together with its standard framing (orientation, equivalently). A non-example in a similar vein is given by $L = S^1 \hookra \RR^2$. Note that the Grassmannian of $1$-planes in $\RR^{2}$ is homeomorphic to $S^1$. Therefore, the Gauss map $L = S^1 \xra{\,\,l\,\,} \mathsf{Gr}_{1}(\RR^2)\simeq S^1$ has a degree and it is $+1$ — it winds once around the origin. As such, it is a standard algebraic topology fact that $l$ does not admit a null-homotopy.  
\end{example}

\subsubsection{Conically smooth stratified spaces with corners}

Now we turn to the goal of incorporating boundaries and corners into the class of conically smooth stratified spaces. To do so, we begin by introducing a special poset, designed with the purpose of capturing the local models of manifolds with $\lag n \rag$-corners \emph{à la} Jänich--Joyce--Laures--Melrose (see \cite{Jan68}\cite{Joy10}\cite{Lau00}\cite{Mel96} for definitions and key results and details). We will follow \cite[Section~7.2]{AFT1}. 

\begin{notation}
    Let $I$ be a finite set. We denote $\cP(I)$ as the poset of subsets of $I$, with the caveat that the ordering is given by reverse inclusion, i.e.\ given a pair of subsets $U, U' \subseteq I$, we have $U\leq U'$ if $U' \subset U$.    
\end{notation}

\begin{example}
    Let $I= \{1,\ldots,n\}$, then we can define a stratification 
    \begin{equation}
        (\RR_{\geq 0})^{n} \xra{\;\;q\;\;} \cP(I), \quad (x_{i})_{i \in I} \longmapsto \{j\,:\, x_{j} = 0\} \subseteq I. 
    \end{equation}
    For $n=1$ and $n=2$ we have the following depictions of this stratification:
    \smallskip
    \begin{equation}
        \tikzfig{figures/Strat/STR20}{\input{figures/Strat/STR20.tikz}} \qquad \qquad \tikzfig{figures/Strat/STR21}{\input{figures/Strat/STR21.tikz}}
    \end{equation}
    \medskip

    \nid Note that in the $n=2$ case, the integers which label the stratum describe which coordinates are vanishing, so they change contravariantly to their standard ordering — that is why $\cP(I)$ is ordered with the ``reverse'' inclusion. 
\end{example}

\begin{observation}
    If $I = \{1\}$, then $\cP(I) = \{\emptyset, \{1\}\}$.
\end{observation}

\begin{definition}[Topologically coCartesian]\cite[Definition 7.2.2]{AFT1}\label{ch4:: def: top coCartesian}
    Let $X = (X \ra P)$ be a conically smooth stratified space. We define a continuous map 
    \begin{equation}
        X \xra{c} \cP(I)
    \end{equation}
    to be \define{topologically coCartesian} if it factors as a map of posets as follows:
    \begin{equation}
        \begin{tikzcd}[ampersand replacement=\&]
            X \&\& P \\
            \& {\cP(I)}
            \arrow[from=1-1, to=1-3]
            \arrow["c"', from=1-1, to=2-2]
            \arrow["{\overline{c}}", dashed, from=1-3, to=2-2]
        \end{tikzcd}
    \end{equation} 
    and if, for every $T \subset I$ i.e.\ for every $T \in \cP(I)$, the inclusion of the $T$-stratum $ X_T \xlonghookra X$ into $X$ extends to a conically smooth open embedding over the poset inclusion:
    \begin{equation}
        \begin{tikzcd}[ampersand replacement=\&]
            {(\RR_{\geq 0})^{T} \times X_T} \&\& X \\
            \\
            {\cP(T)} \&\& {\cP(I)}
            \arrow[hook, from=1-1, to=1-3]
            \arrow["{q\times c_{|\cP(T)}}"', from=1-1, to=3-1]
            \arrow["c", from=1-3, to=3-3]
            \arrow[hook, from=3-1, to=3-3]
        \end{tikzcd}
    \end{equation}
\end{definition} 

\begin{example}
    If $M \longra *$ is a smooth manifold and $I = \{*\}$, then a topologically coCartesian map
    \begin{equation}
        M \xra{\;\; c \;\;} \{\emptyset, \{*\}\}
    \end{equation}  
    determines a pair of subspaces $M_{\emptyset}$ and $M_{*}$ together with a (conically) smooth map 
    \begin{equation}
        (\RR_{\geq 0}) \times M_{*} \xra{\;\;\;} M
    \end{equation}
    such that $\{0\}\times M_\ast \xra{\;\;\;} M$ is the stratum inclusion. 
\end{example}

\begin{example}
    Let the following diagram
    \smallskip 
    \begin{equation}
        M \;\; = \;\; \tikzfig{figures/Strat/STR22}{\input{figures/Strat/STR22.tikz}}
    \end{equation}
    \medskip

    \nid denote a stratification $M$ of the unit square $[0,1]^2$. This decomposes as:
    \smallskip
    \begin{equation}
        M_{\emptyset}\;\; = \;\; \tikzfig{figures/Strat/STR23}{\input{figures/Strat/STR23.tikz}}, \qquad \qquad  M \setminus M_{\emptyset} \;\; = \;\; \tikzfig{figures/Strat/STR24}{\input{figures/Strat/STR24.tikz}}
    \end{equation}
    \medskip
    with $M\setminus M_{\emptyset}$ further decomposing in the obvious way. 
\end{example}

\begin{definition}
    Fix a finite set $I$. The category $\Bsc^{\lag I \rag}$ is defined as follows:
    \begin{enumerate}
        \item \textbf{Objects:} A basic $U\in \Bsc$ endowed with a topologically coCartesian map $U \xra{\; c \;} \cP(I)$. 
        \item \textbf{Morphisms:} Given a pair of objects $\left(c \maps U \ra P \xra{\overline{c}} \cP(I)\right)$ and $\left(d \maps V \ra Q \xra{\overline{d}} \cP(I)\right)$, the set of morphisms is defined as the pullback:
        \begin{equation}
             \begin{tikzcd}[ampersand replacement=\&]
                {\Bsc^{\lag I \rag}\left(\left(U,c\right),\left(V,d\right)\right)} \&\& {\Bsc(U,V)} \\
                \\
                {\Hom_{\Top_{/\cP(I)}}(P,Q)} \&\& {\Hom_{\Top}(P,Q).}
                \arrow[from=1-1, to=1-3]
                \arrow[from=1-1, to=3-1]
                \arrow["\lrcorner"{anchor=center, pos=0.125}, draw=none, from=1-1, to=3-3]
                \arrow[from=1-3, to=3-3]
                \arrow[from=3-1, to=3-3]
            \end{tikzcd}
        \end{equation}
        Since the lower horizontal map is injective so then the top horizontal must be an inclusion on $\pi_0$. In particular, this  implies that $\Bsc^{\lag I \rag}$ is a sieve. We can thus define \begin{equation}
            \Snglr^{\lag I \rag} \, \overset{\sf def}{=} \, \Mfld(\Bsc^{\lag I \rag})
        \end{equation}
        which we will call the $(\oo,1)$\define{-category of conically smooth stratified spaces with} $\lag I \rag$\define{-corners}.   
    \end{enumerate}
\end{definition}

\subsubsection{Unzipping and vector fields}

We now turn to a key construction, \emph{unzipping}, which offers a bridge between the world of conically smooth stratified spaces and that of smooth manifolds with corners, thereby opening up a fruitful avenue for translating statements, definitions, and theorems in one into definitions, statements, and theorems in the other.   

In particular, unzipping will lead us to the definition of parallel and inward-pointing vector fields. These are the correct notions of conically smooth infinitesimal automorphisms for conically smooth stratified spaces with $\lag I \rag$-corners. In turn, these lead to vector field flows and ultimately to a generalization of the Ehresmann fibration theorem (see \Cref{ap:: prop: Ehresmann fibration no boundary}) for conically smooth stratified spaces.

In a slogan: unzipping associates to every conically smooth stratified spaces $(M \ra P)$ of interior depth $k$ a smooth manifold $\unzip(M)$ with $\lag k \rag$-corners (no singularities) and it does so functorially. This procedure works by morally mirroring the notion of algebraic blow-up for algebraic varieties: it removes the singularity locus and glues in its associated spherical bundle. A simple example is shown below:
\begin{equation}
    \tikzfig{figures/Strat/STR39}{\input{figures/Strat/STR39.tikz}} \;\;\; \xra{\,\, \unzip \,\,} \;\;\; \tikzfig{figures/Strat/STR40}{\input{figures/Strat/STR40.tikz}}
\end{equation}
We remove the only singular locus, the 0-dimensional stratum at the origin, and we replace it using the spherical bundle associated to its normal bundle.

\begin{definition}[Unzip at depth $k$]\cite[Definition~7.3.2]{AFT1}
    Fix $k \geq 0$ and let $\overline{X} = (X,c_X)$ be a conically smooth stratified space with $\lag I \rag$-corners, whose interior depth is $k$. We define the following category $\mathsf{Res}(\overline{X},k)$
    \begin{itemize}
        \item \textbf{Objects:} A conically smooth stratified space with $\lag * \sqcup I \rag$-corners 
        \begin{equation}
            \overline{Y} \xra{c_{Y}} \cP(*\sqcup I)
        \end{equation} 
        together with a conically smooth map 
        \begin{equation}
            p \maps \overline{Y} \longra X
        \end{equation}
        such that:
        \begin{itemize}
            \item The following diagram commutes:
            \begin{equation}
                \begin{tikzcd}[ampersand replacement=\&]
                    {\overline{Y}} \&\& X \\
                    \\
                    {\cP(*)\times \cP(I)} \&\& {\cP(I)}
                    \arrow[from=1-1, to=1-3]
                    \arrow[from=1-1, to=3-1]
                    \arrow[from=1-3, to=3-3]
                    \arrow[from=3-1, to=3-3]
                \end{tikzcd}
            \end{equation}
            where we have used the identification $\cP(*\sqcup I) \simeq \cP(*) \times \cP(I)$.
            \item The subspaces of $Y$ \[\overline{\partial}_* \overline{Y} = p^{-1}\overline{X}_{k}\]
            agree.
        \end{itemize}
        \item A morphism from $(Y,c_{Y},p)$ to $(Y',c_{Y'},p')$ is a conically smooth map 
        \begin{equation}
            \begin{tikzcd}[ampersand replacement=\&]
                {\overline{Y}} \&\& {\overline{Y}'} \\
                \\
                \& {\cP(*\sqcup I)\times X}
                \arrow[from=1-1, to=1-3]
                \arrow[from=1-1, to=3-2]
                \arrow[from=1-3, to=3-2]
            \end{tikzcd}
        \end{equation} 
    \end{itemize}
    We define a $k$-unzip of $\overline{X}$ to be 
    \begin{equation}
        \pi_k \maps \unzip_{k}(X) \longra X 
    \end{equation}
    a \emph{final} object of $\mathsf{Res}(\overline{X},k)$. If we restrict $\unzip_k(X)$ to the boundary $\overline{\partial}_{*}\unzip_{k}(X)$ we get the \define{link}:
    \begin{equation}
        \link_{k}(X) \xra{\, \pi_k \,} X_{k}.
    \end{equation} 
\end{definition}

\begin{example}
    Let us explore unzipping through a simple example. Let 
    \begin{equation}
        M \;\; = \;\; \tikzfig{figures/Strat/STR41}{\input{figures/Strat/STR41.tikz}}
    \end{equation}
    \medskip
    
    \nid be $\RR^2$ stratified by the origin and the $y$-axis minus the origin. The 0-dimensional stratum has depth 2, therefore we have that 
    \begin{equation}
        \unzip_2(M) \;\; = \;\; \tikzfig{figures/Strat/STR42}{\input{figures/Strat/STR42.tikz}}
    \end{equation}
    \medskip

    \nid The resulting conically smooth stratified space $\unzip_2(M)$ is a conically smooth stratified space with $\lag * \rag$-corners, i.e.\ a boundary — the link is the blue boundary circle stratified with the 0-dimensional dots. The interior depth of $\unzip_2(M)$ is 1, so we can apply $\unzip_1$, and this would remove the red stratum and replace it by the sphere bundle associated to its normal bundle 
    \begin{equation}
        \unzip_{1}\left(\unzip_{2}(M)\right) \;\; = \;\; \tikzfig{figures/Strat/STR43}{\input{figures/Strat/STR43.tikz}}
    \end{equation} 
\end{example}

\begin{remark}\cite[Example~7.3.4]{AFT1}
    Let us consider $\unzip_2(M)$ from the example above. A possible object of $\mathsf{Res}\left(\unzip_2(M),1\right)$ is a curve $\gamma \maps (0,1] \longra \unzip_2(M)$ which must be transverse to $W$ (the depth $1$ stratum) at $1$ and only at $1$ — this ensures it is a manifold with boundary in $\unzip_2(M)$ and its boundary is contained in $W$. For example, let the yellow curve $\gamma$ be such an object
    \begin{equation}
        \tikzfig{figures/Strat/STR44}{\input{figures/Strat/STR44.tikz}}
    \end{equation}
    \medskip

    \nid Moreover, there is a clear morphism in $\mathsf{Res}\left(\unzip_2(M),1\right)$ from $\gamma$ to $\unzip_2(M)$, described below 
    \begin{equation}
        \tikzfig{figures/Strat/STR44}{\input{figures/Strat/STR44.tikz}} \;\; \longra \;\; \tikzfig{figures/Strat/STR45}{\input{figures/Strat/STR45.tikz}}
    \end{equation}
\end{remark}
\medskip

Unzipping allows us to discuss the ``right'' notion of vector fields for conically smooth stratified spaces. Before doing so, we introduce a few key components and useful notation.

Given a smooth manifold $M$, let us denote $\vartheta(M)$ for the vector space $\Gamma(\sT M \ra M)$ of smooth sections or smooth vector fields on $M$. Fixing a smooth map $f\maps M \longra N$, we have the following linear maps
\begin{equation}
    \vartheta(N) \longra \Gamma(f^* \sT N \ra M), \qquad (N \xra{X} \sT N) \longmapsto (M \xra{f} N \xra{X} \sT N)
\end{equation}
which we will call restriction, and 
\begin{equation}
    \vartheta(M) \xra{\sD f} \Gamma(f^* \sT N \ra M), \qquad (M \xra{X} \sT M) \longmapsto (M \xra{X} \sT M \xra{\sD f} \sT N)
\end{equation}
which we will call the derivative map. We will be considering smooth vector fields on smooth manifolds with corners, so let us provide a gentle reminder of how these are defined. Smooth vector fields are defined locally, so it suffices to determine what these mean on open subsets of the form 
\begin{equation}
    \RR^{i} \times \RR^{j}_{\geq 0}, \qquad \vartheta(\RR^{i} \times \RR^{j}_{\geq 0}) = \Gamma\left( \sT \RR^{i+j}_{|\RR^{i}\times \RR^{j}_{\geq 0}} \longra \RR^{i}\times \RR^{j}_{\geq 0} \right).
\end{equation}
In particular, if $M$ is a smooth manifold with boundary, say $\RR_{\geq 0}$, then a possible vector field might point outward. We are interested in defining sub-sheaves of vector fields which point-inward or tangently to the faces of the boundary.

\begin{definition}
    Let $\overline{M}$ be a smooth manifold with corners. We define the sub-sheaf of \define{inward-pointing vector fields} as
    \begin{equation}
        \vartheta^{\sf in}(\overline{M}) \, = \, \{X \maps \RR^{i}\times \RR^{j}_{\geq 0} \ra \RR^{i+j}\,:\,  \mathsf{pr}_{i}\left(X_{|\RR^{i}\times \RR^{S} \times \{0\}}\right) > 0, \,\,\,\, \forall\, i \in \{1,\ldots,j\}\setminus S \,\,\,\, \forall\, S \in \cP(\{1,\ldots,j\})\}
    \end{equation}
    and that of \define{parallel vector fields} as
    \begin{equation}
        \vartheta^{\sf ||}(M) \, = \, \{X \maps \RR^{i}\times \RR^{j}_{\geq 0} \ra \RR^{i+j}\,:\,  X_{|\RR^{i}\times \RR^{S} \times \{0\}} \subset \RR^{i}\times \RR^{S} \times \{0\} \,\,\,\, \forall\, S \in \cP(\{1,\ldots,j\})\}.
    \end{equation} 
\end{definition}

\begin{example}
    Consider the case $i=1$ and $j=1$, i.e.,\ $\RR \times \RR_{\geq 0}$.
    \begin{equation}
        \tikzfig{figures/Strat/STR46}{\input{figures/Strat/STR46.tikz}} 
    \end{equation} 
    \medskip

    \nid The vectors are inward-pointing, parallel, and non-inward-pointing from left to right, respectively.
\end{example}

Lastly, we note that given a conically smooth stratified space $M\ra P$ and a linear subposet $S \subset P$ we can articulate the corresponding face of its unzip:
\begin{equation}
    \partial_{S} \unzip(M) \, = \, \cap_{s \in S} \partial_{s} \unzip(M).
\end{equation}
It comes with a constructible bundle (smooth fiber bundle) $\pi_S \maps \partial_{S} \unzip(M) \longra M_{\mathsf{Min}(S)}$.

\begin{definition}[Parallel vector fields]\cite[Definition~8.1.2]{AFT1}
    Let $M = (M\ra P)$ be a conically smooth stratified space. We define the vector space $\vartheta(M)$ of \define{parallel vector fields} on $M$ as the limit 
    \begin{equation}
        \begin{tikzcd}[ampersand replacement=\&]
            {\vartheta(M)} \&\&\& {\vartheta^{||}\left(\unzip(M)\right)} \\
            \\
            \&\&\& {\prod_{S\underset{\text{linear}}{\subset}P}\vartheta\left(\partial_{S}\unzip(M)\right)} \\
            \\
            {\prod_{p\in P} \vartheta(M_p)} \&\&\& {\prod_{S\underset{\text{linear}}{\subset}P} \Gamma\left(\pi_{S}^*\sT M_{\mathsf{Min}(S)}\ra \partial_{S}\unzip(M)\right)}
            \arrow["{\widetilde{(-)}}", from=1-1, to=1-4]
            \arrow["{(-)_{|M_{-}}}"', from=1-1, to=5-1]
            \arrow["{(-)_{|\partial_{-}}}", from=1-4, to=3-4]
            \arrow["{\sD\pi_{-}}", from=3-4, to=5-4]
            \arrow["{(\pi_{-})^{*}}"', from=5-1, to=5-4]
        \end{tikzcd}
    \end{equation}
\end{definition}

\begin{example}
    The notion of parallel vector field is best understood through an example. Let us consider the following stratification of the open disk
    \begin{equation}
        M \;\; = \;\; \tikzfig{figures/Strat/STR47}{\input{figures/Strat/STR47.tikz}}.
    \end{equation}
    \medskip
    
    \nid We might change it by both breaking into strata or taking its total unzip:
    \begin{equation}
        \tikzfig{figures/Strat/STR48}{\input{figures/Strat/STR48.tikz}} \qquad , \qquad \tikzfig{figures/Strat/STR49}{\input{figures/Strat/STR49.tikz}}
    \end{equation}
    \medskip 

    \nid Any parallel vector field has to be such that it restricts to a vector field on every stratum; moreover, they must satisfy a compatibility condition stipulated by their lift to $\partial_S \unzip(M)$ for varying $S$.     
\end{example}

\begin{remark}
    Fixing a conically smooth stratified space $M$, the collection of parallel vector fields organizes into a sheaf of vector spaces:
    \begin{equation}
        \vartheta_{M} \maps \mathsf{Opens}(M)^{\op} \longra \mathsf{Vect}_{\RR}, \quad \vartheta_M(U) \, = \, \vartheta(U).
    \end{equation} 
\end{remark}

\begin{proposition}[Flows]\cite[Lemma~8.1.3]{AFT1}\label{ch4:: def:flows}
For a conically smooth stratified space $M = (M \longra P)$, and for $V$ a parallel vector field on $M$, there exists a conically smooth map 
\begin{equation}\label{eq:def:flows}
    \gamma^V \maps M \times \RR \dashrightarrow M
\end{equation}
defined on an open neighborhood of $M\times \{0\}$, satisfying the following conditions:
\begin{enumerate}
    \item For each conically smooth map $\epsilon \maps M \longra \RR$, the composition 
    \begin{equation}\label{eq:def:flows:iso}
        M \xra{~ \id_M \times \epsilon ~} M \times \RR \dashrightarrow{\gamma^V} M
    \end{equation}
    is a conically smooth isomorphism, provided it is defined. 
    \item \label{additivity for flow} There is an equality 
    \begin{equation}\label{eq:def:flows:action}
        \gamma^V_{t}\big(\gamma^V_{s}(x)\big) = \gamma^V_{s+t}(x)
    \end{equation}
    provided the maps are defined.
    \item For each point $x\in M$ belonging to a stratum $x\in M_p \subset M$,
    \begin{equation}\label{eq:def:flows:tangent}
        \frac{\sf d}{{\sf d}t}\gamma_t^{V}(x)|_{t=0} = V_{|M_p}(x).
    \end{equation}
\end{enumerate}
We call $\gamma^V$ the \define{conically smooth flow of }$V$.
\end{proposition}

\begin{definition}[Inward-pointing vector fields]\cite[Definition~8.1.2]{AFT1}
    Let $\overline{M} = (\overline{M}\ra P,c)$ be a conically smooth stratified space with corners. We define the vector space $\vartheta^{\sf in}(M)$ of \define{inward vector fields} on $\overline{M}$ as the limit 
    \begin{equation}
        \begin{tikzcd}[ampersand replacement=\&]
            {\vartheta^{\sf in}(M)} \&\&\& {\vartheta^{||}\left(\unzip(\overline{M})\right)} \\
            \\
            \&\& {\prod_{S\underset{\text{linear}}{\subset}P}\vartheta\left(\partial_{S}\unzip(M)\right)} \& {\prod_{S\underset{\text{linear}}{\subset}P}\Gamma\left(\sT \unzip(\overline{M})_{|\partial_S}\ra \partial_S \right)} \\
            \\
            {\prod_{p\in P} \vartheta^{\sf in}(M_p)} \&\& {\prod_{S\underset{\text{linear}}{\subset}P} \Gamma\left(\pi_{S}^*\sT \overline{M}_{\mathsf{Min}(S)}\ra \partial_{S}\unzip(\overline{M})\right)}
            \arrow["{\widetilde{(-)}}", from=1-1, to=1-4]
            \arrow[from=1-1, to=3-3]
            \arrow["{(-)_{|\overline{M_{-}}}}"', from=1-1, to=5-1]
            \arrow["{(-)_{|\partial_{-}}}", from=1-4, to=3-4]
            \arrow[from=3-3, to=3-4]
            \arrow["{\sD\pi_{-}}", from=3-3, to=5-3]
            \arrow["{(\pi_{-})^{*}}"', from=5-1, to=5-3]
        \end{tikzcd}
    \end{equation}
\end{definition}

\begin{proposition}[Inward flows]\cite[Lemma~8.1.3]{AFT1}\label{ch4:: def: in flows}
For a conically smooth stratified space with corners $\overline{M} = (M \longra P,c)$, and for $V$ an inward vector field on $\overline{M}$, there exists a conically smooth map 
\begin{equation}\label{eq:def:in flows}
    \gamma^V \maps M \times \RR \dashrightarrow M
\end{equation}
defined on an open neighborhood of $M\times \{0\}$, satisfying the following conditions:
\begin{enumerate}
    \item For each conically smooth map $\epsilon \maps \overline{M} \longra \RR$, and for each linear subposet $S\subset P$, the composition 
    \begin{equation}\label{eq:def:in flows:iso}
        \partial_{S}\unzip(\overline{M}) \xra{\id \times \epsilon_{|\partial_{S}}} \partial_{S}\unzip(\overline{M})\times \RR_{\geq 0} \dashrightarrow \unzip(\overline{M})
    \end{equation}
    is a refinement onto its image, provided it is defined. 
    \item There is an equality 
    \begin{equation}\label{eq:def:in flows:action}
        \gamma^V_{t}\big(\gamma^V_{s}(x)\big) = \gamma^V_{s+t}(x)
    \end{equation}
    provided the maps are defined.
    \item For each point $x\in M$ belonging to a stratum $x\in M_p \subset M$,
    \begin{equation}\label{eq:def:in flows:tangent}
        \frac{\sf d}{{\sf d}t}\gamma_t^{V}(x)|_{t=0} = V_{|M_p}(x).
    \end{equation}
    
\end{enumerate}
\end{proposition}

From the existence of flows we can establish the existence of regular neighborhoods for arbitrary proper constructible embedding $Y\hookra X$.

\begin{proposition}\cite[Proposition~8.2.5]{AFT1}
    Let $Y\hookra X$ be a proper constructible embedding. Then there is a conically smooth map 
    \begin{equation}
        \sC(\pi_Y) \longra X
    \end{equation}
    under $Y\hookra X$ such that it is a refinement onto its image, and its image is open. The conically smooth stratified space $\sC(\pi_Y)$ is built as the fiberwise open cone of $\link_{Y}(X) \xra{\pi_Y} Y$.
\end{proposition}

The existence of regular neighborhoods allows us to establish the following important structural fact about $\Strat$.

\begin{proposition}\cite[Lemma~2.5.2]{AFR1}
    Given a conically smooth stratified space $M$ together with its deepest stratum $M_k \subset M$, the following is a pushout square in the $(\oo,1)$-category $\Strat$:
    \begin{equation}
        \begin{tikzcd}[ampersand replacement=\&]
            {\link_{M_{k}}(M)} \&\& {\unzip(M)} \\
            \\
            {M_k} \&\& M
            \arrow[hook, from=1-1, to=1-3]
            \arrow[from=1-1, to=3-1]
            \arrow[from=1-3, to=3-3]
            \arrow[hook, from=3-1, to=3-3]
        \end{tikzcd}
    \end{equation} 
\end{proposition}

Conically smooth stratified spaces can thus be obtained from their unzip, by iteratively gluing strata to the shallower faces. More concretely, assuming $k$ is the interior depth of $M$:
\begin{equation}
    \unzip_{k}(M) = X_{k-1} \coprod_{\link_{X_{k-1}}(\unzip_{k-1}(M))} \unzip_{k-1}(M)
\end{equation}

\subsubsection{Ehresmann fibration theorem}
A key result in building the infrastructure around higher categories of cobordisms is the Ehresmann fibration theorem, in its versions with boundary and without boundary (see \Cref{ap:: prop: Ehresmann fibration no boundary}). We now turn to proving its conically smooth analogue — in essence, a conically smooth version of \emph{Thom's first isotopy lemma} (see Mather's seminal notes \cite[Proposition~11.1]{Mat70}). 

We begin by establishing the fact that, akin to smooth submersions, conically smooth submersions admit lifts of parallel vector fields.

\begin{proposition}\label{prop: lifts weak}
    Let $\pi\maps X \longra K$ be a weakly constructible bundle. Then every parallel vector field $V \in \vartheta(K)$ lifts to a parallel vector field $\widetilde{V}\in \vartheta(X)$, i.e., $\sD \pi (\widetilde{V}) = V$.
\end{proposition}
\begin{proof}
    The map $\pi\maps X \longra K$ determines a morphism between sheaves of vector spaces on $\unzip(X)$:
\begin{equation}
    \sD \tilde{\pi} \colon \vartheta_{\unzip(X)} \ra \tilde{\pi}^{\ast}\vartheta_{\unzip(K)}
\end{equation}
because the differential of a face-preserving map again preserves the corresponding boundary faces. Locally on $\unzip(K)$, the map $\tilde{\pi}$ is isomorphic to a projection $U\times F \longra U$, and for a projection every parallel vector field on $U$ lifts to $(V,0)$ on $U\times F$. Hence $\sD\tilde{\pi}$ is surjective as a map of sheaves. We therefore have the commutative diagram
\begin{equation}
\begin{tikzcd}[ampersand replacement=\&]
	{(q_X)_\ast\vartheta_{\unzip(X)}} \&\& {(q_X)_\ast (\tilde{\pi}^\ast\vartheta_{\unzip(K)})} \\
	\&\& {\pi^\ast (q_K)_{\ast}(\vartheta_{\unzip(K)})} \\
	{\vartheta_X} \&\& {\pi^\ast\vartheta_K}
	\arrow["{(q_X)_\ast(\sD \tilde{\pi})}", from=1-1, to=1-3]
	\arrow[from=1-1, to=3-1]
	\arrow[from=1-3, to=2-3]
	\arrow[from=2-3, to=3-3]
	\arrow[dashed, from=3-1, to=3-3]
\end{tikzcd}
\end{equation}
The vertical maps identify parallel vector fields on $X$ and $K$ with the corresponding pushforwards from their unzips, so the bottom horizontal map is induced from the top one. Since the top map is surjective on sections over every open set, the same is true for $\sD\pi$. Thus every parallel vector field on $K$ lifts locally to $X$, and because $\vartheta_X$ is a sheaf of vector spaces these local lifts glue to a global lift.
\end{proof}

\begin{proposition}\label{prop:lifts of vector fields}
    Let $f\maps X \longra Y$ be a conically smooth submersion, and let $V$ be a parallel vector field on $Y$. There exists a parallel vector field $W$, defined on $X$, such that $\sD f(W) = V.$  
\end{proposition}
\begin{proof}
     Because $f$ is a submersion, we can find a basis $\cU_{Y}= \{\varphi_\alpha(O_\alpha)\}_\alpha$ for the topology of $Y$ together with conically smooth open embeddings $\psi_\alpha \maps U_\alpha \times O_\alpha \hookrightarrow X$ for which 
    \begin{equation}
        \big(f\circ \psi_\alpha \big)(u,o) = \varphi_\alpha(o).
    \end{equation} 
    On the local model $U_\alpha \times O_\alpha \longra O_\alpha$, the restriction $V_\alpha := V_{|O_\alpha}$ lifts to the parallel vector field $(0,V_\alpha)$. Transporting this field across $\psi_\alpha$ gives a parallel lift $W_\alpha$ on the open subset $\psi_\alpha(U_\alpha\times O_\alpha)\subset X$.

    Choose a conically smooth partition of unity $\{\eta_\alpha\}$ on $Y$ subordinate to the cover $\cU_Y$. Then
    \begin{equation}
        W:=\sum_\alpha (f^\ast\eta_\alpha)\,W_\alpha
    \end{equation}
    is a well-defined parallel vector field on $X$. Since each $W_\alpha$ lifts $V_{|O_\alpha}$ and the coefficients come from the base, the same partition-of-unity argument as in the smooth case shows that $W$ lifts $V$. 
\end{proof}

These results lead us closer to the proof of the Ehresmann fibration theorem. We are only missing a standard fact about flow extensions over compact subsets. 

\begin{lemma}\label{lem:compact flows}
    Let $K\subset X$ be a compact subset of a conically smooth stratified space $Z =(X\longra P)$, and let $V$ be a parallel vector field on $X$ with maximal flow given by $\gamma^V \maps A \longra X$, with $A\subset X\times \RR$ an open neighborhood of $X\times \{0\}$. Let $x\in X_p \subset X$, and let 
    \begin{equation}\label{eq:lem:compact flows}
        \{x\}\times J_x  = A \cap (\{x\}\times \RR)
    \end{equation}
    such that $J_x = (a_x,b_x)$ and contains $0$. If $b_x<\infty$, then there exists $\epsilon > 0$ such that $\gamma^V_t(x)\notin K$ for $b_x - \epsilon < t < b_x$. Likewise, if $a_x>-\infty$, then there exists $\epsilon > 0$ such that $\gamma^V_t(x)\notin K$ for $a_x < t < a_x+\epsilon$.
\end{lemma}
\begin{proof}
    Since $A$ is an open neighborhood of $K\times \{0\}$, compactness gives $\epsilon>0$ such that
    \begin{equation}\label{eq:lem:prf:compact flows cover}
        K\times [-\epsilon,\epsilon] \subset A.
    \end{equation}
    In other words, every point of $K$ can be flowed for at least time $\epsilon$ in both directions.

    Suppose $b_x<\infty$ and that there exists a sequence $t_n \uparrow b_x$ with $\gamma^V_{t_n}(x)\in K$. For each $n$, the point $\gamma^V_{t_n}(x)$ lies in $K$, so by \Cref{eq:lem:prf:compact flows cover} the flow is defined for an additional time interval $[-\epsilon,\epsilon]$ starting at $\gamma^V_{t_n}(x)$. Therefore $\gamma^V_{t_n+\epsilon}(x)$ is defined. For $n$ large enough, however, $t_n+\epsilon>b_x$, contradicting the maximality of $J_x$. Hence $\gamma^V_t(x)\notin K$ for all $t$ sufficiently close to $b_x$ from the left.

    The argument for the left endpoint is identical, applied to a sequence $t_n \downarrow a_x$.
\end{proof}

\begin{remark}\label{rmk:flow out of compact}
\Cref{lem:compact flows} shows that a maximal trajectory cannot terminate while remaining inside a fixed compact subset. In particular, if the image of a trajectory is contained in a compact subset for every bounded time interval, then the corresponding flow is defined for all time.
\end{remark}

We can now finally state and prove the main theorem of this section — the conically smooth stratified analogue of the Ehresmann fibration theorem. 

\begin{theorem}[Thom/Ehresmann/Morse for conically smooth stratified spaces]\label{ch4:: thm : stratified ehresmann}
    Let $\pi\maps X \longra Y$ be a conically smooth proper submersion. Then $\pi$ is a conically smooth fiber bundle.
\end{theorem}
\begin{proof}
Recall that $\pi$ is a conically smooth fiber bundle if for every $y\in Y$, there exists a neighborhood $U\subset Y$, which we can assume centered at $y$, and a conically smooth isomorphism
\begin{equation}
    h_U\maps \pi^{-1}(U) \longra \pi^{-1}(\{y\})\times U
\end{equation}
participating in the following commutative diagram: 
\begin{equation}\label{eq:thm:Ehresmann}
\begin{tikzcd}[ampersand replacement=\&]
	{\pi^{-1}(U)} \&\& {\pi^{-1}(\{y\})\times U} \\
	\\
	\& U \& {.}
	\arrow["{h_U}", from=1-1, to=1-3]
	\arrow["{\pi_{|\pi^{-1}(U)}}"', from=1-1, to=3-2]
	\arrow["{{\sf pr}}", from=1-3, to=3-2]
\end{tikzcd}
\end{equation}
Because the definition is local on the base, after shrinking around $y$ we may suppose that $Y=\RR^k$ and that $y=0\in \RR^k$. Let $\partial_1,\ldots, \partial_k$ be the coordinate vector fields on $\RR^k$. By \Cref{prop:lifts of vector fields}, there exists a parallel vector field $\partial_{i}^X$ on $X$ for each $i\in \{1,\cdots, k\}$ such that 
\begin{equation}\label{eq:thm:str_ehrs:prf:lift}
    (\sD_{x}\pi)\,\big(\partial^X_{i}(x)\big)= \partial_{i}(\pi(x)) \quad \forall x\in X. 
\end{equation}
By \Cref{eq:def:flows}, each of the lifted parallel vector fields $\partial^X_{i}$ on $X$ generates a local flow 
\begin{equation}\label{eq:thm:prf:flow}
\gamma^{\partial^X_{i}}\maps X \times \RR \dashrightarrow X,
\end{equation} defined on an open neighborhood $A_i \subset X\times \RR$ of $X\times \{0\}.$ Projecting this flow down to $\RR^k$, we note that it is given by the simple explicit formula
\begin{equation}
    \pi\big(\gamma^{\partial^{X}_{i}}_{t}(x)\big) = \pi(x) + (0,\ldots,0,t,0,\ldots,0),
\end{equation} where the only non-vanishing coordinate is $i$. 

Fix $x\in X$ and $i\in \{1,\ldots,k\}$, and write
\begin{equation}
    \{x\}\times J_{x,i}=A_i\cap (\{x\}\times \RR)
\end{equation}
for the maximal interval of definition of the flow through $x$. Suppose $b_{x,i}<\infty$. Choose numbers $a<0<b_{x,i}<b$. The formula above shows that
\begin{equation}
    \pi\big(\gamma^{\partial^{X}_i}_{t}(x)\big)\in K:=\pi(x)+\{(0,\ldots,0,t,0,\ldots,0)\mid t\in [a,b]\}\subset \RR^k
\end{equation}
for every $t\in J_{x,i}\cap (a,b)$. The set $K$ is compact, so properness of $\pi$ implies that $\pi^{-1}(K)$ is compact. Hence the image of the trajectory on $J_{x,i}\cap (a,b)$ lies in a fixed compact subset of $X$. By \Cref{lem:compact flows}, this is incompatible with a finite right endpoint. Therefore $b_{x,i}=+\infty$. The same argument applied to the left endpoint shows $a_{x,i}=-\infty$. Thus each $\partial_i^X$ is complete.

We introduce the notation $t_{i}(x):=\pi_{i}\big(\pi(x)\big)$, where $\pi_i$ is the projection onto the $i$-th factor. We define
\begin{equation}
    h\maps X \longrightarrow \pi^{-1}(\{0\})\times \RR^k
\end{equation}
as follows:
\begin{align}
    h(x) &= \bigg(\Big(\gamma^{\partial^{X}_{1}}_{-t_1(x)}\circ \cdots \circ \gamma^{\partial^{X}_{k}}_{-t_k(x)}\Big)(x), \pi(x)\bigg).
\end{align}
Projecting the first component to $\RR^k$ and using the displayed formula for the base flow, each negative-time flow subtracts the corresponding coordinate of $\pi(x)$. Hence the first component indeed lands in $\pi^{-1}(\{0\})$, so $h$ is well-defined. The map is conically smooth because it is built from conically smooth flows and projections. To show it is an isomorphism, define
\begin{equation}
    \bar{h}\big(x_0,(t_1,\cdots,t_k)\big) = \Big(\gamma^{\partial^{X}_{k}}_{t_k}\circ \cdots \circ \gamma^{\partial^{X}_{1}}_{t_1}\Big)(x_0).
\end{equation} 
This map is again conically smooth. The identities $\bar{h}\circ h = \id_X$ and $h\circ \bar{h}=\id_{\pi^{-1}(0)\times \RR^k}$ follow by repeatedly applying the additivity relation from \Cref{eq:def:flows}; no commutativity between the different flows is needed because the two compositions are written in opposite orders. Therefore $\pi$ is locally isomorphic to the projection $\pi^{-1}(\{0\})\times \RR^k \longra \RR^k$, as required.
\end{proof}

\begin{remark}
    We will need a variant of the stratified Ehresmann lemma for conically smooth stratified spaces \emph{with corners}, for our applications to stratified bordisms. This is the content of the next theorem.
\end{remark}

\begin{cor}\label{ch4:: cor: stratified ehresmann boundary}
    Let $\big(M \ra P\big)$ and $K$ be a conically smooth stratified space and a smooth manifold, respectively. Assume $\big(M \ra P\big)$ and $K$ both carry boundary structures. If $\pi\maps M \ra K$ is a conically smooth proper submersion such that $\pi(\partial M) \subseteq \partial K$, then $\pi$ is a conically smooth fiber bundle. 
\end{cor}
\begin{proof}
    The argument of the proof follows that of the unstratified case, exhibited in \Cref{ap:: prop: Ehresmann fibration no boundary}. For boundary points, the vector field lift must necessarily be inward-pointing, and the existence of the corresponding flow is ensured by \Cref{ch4:: def: in flows}.
\end{proof}

\subsubsection{Exit-path \texorpdfstring{$(\oo,1)$}{(infty,1)}-category}

The singular complex $\sing(M)$ of a smooth manifold $M$ is a fundamental construction in homotopy theory and manifold topology. Any conically smooth stratified space $\big(M \ra P\big)$ also admits a singular complex $\sing(M)$. Yet given, for example, a refinement of $\big(M\ra P\big) \longra (M' \ra Q)$, where the underlying spaces are homeomorphic, this induces an evident equivalence between singular complexes $\sing(M)\longra \sing(M')$ — the caveat is that $M$ and $M'$ are not in general stratified homotopy equivalent.

The construction of the exit-path category (or its 2-truncation bicategory) dates back to David Treumann's work on exit-path 2-categories and constructible stacks \cite{Tre09} together with Lurie's work (see \cite{LurieHA}) on the $(\oo,1)$-category version — the latter is the one which we will be discussing here. Before proceeding with the definition, we must first introduce an important stratification of the topological standard $n$-simplex $\Delta^n$.

\begin{definition}\label{ch4:: def: standard stratification}
There is a cosimplicial object in $\strat$, hence also in $\Strat$,
\begin{align}
    \mathsf{st} \maps \bDelta &\longra \strat \\
    [p] &\longmapsto \big(\,\Delta^{p} \xra{m} [p] \,\big)
\end{align}
where $m\maps (t_0,\cdots, t_p) \longmapsto \mathsf{Max}\{i\, |\, t_i \neq 0\}$.
\end{definition}

This cosimplicial object is well-defined, since giving a morphism $f \maps [p] \ra [q]$ in $\bDelta$, the usual map of spaces $|f| \maps \bDelta^{p} \longra \bDelta^{q}$ extends to a map of stratified spaces 
\begin{equation}
    \begin{tikzcd}
        {\bDelta^{p}} && {\bDelta^{q}} \\
        \\
        {[p]} && {[q]}
        \arrow["{|f|}", from=1-1, to=1-3]
        \arrow["m"', from=1-1, to=3-1]
        \arrow["m", from=1-3, to=3-3]
        \arrow["f"', from=3-1, to=3-3]
    \end{tikzcd}
\end{equation}
since $f$ is monotonic. 

\begin{observation}
We grant that this is a very unorthodox stratification of the topological $p$-simplex. Other, more pedestrian stratifications one could conceive of include, for example, a stratification by dimension of the face of the simplex or, as a substratification, just two strata: the interior and the boundary. To help make sense of this stratification we draw the stratifications until $p=3$.

\smallskip
\begin{equation}
    \tikzfig{figures/Strat/STR13}{\input{figures/Strat/STR13.tikz}}
\end{equation}
\medskip

\noindent There is a more natural way of understanding this asymmetric stratification through the lens of cones. Namely, remember the closed cone construction $\overline{\sC}(-)$ which produces a stratified space with boundary. If one takes the closed cone on the empty set 
\begin{equation}
\overline{\sC}(\emptyset \ra [0])
\end{equation} which is precisely $\Delta^1 \ra [1]$. Moreover, iterating this $(p-1)$ times we get 
\begin{equation}
    \overline{\sC}^{p}(\emptyset \ra [0]) = \big(\,\Delta^{p} \ra [p]\,\big).
\end{equation} 
So that for every $p\geq 0$, we have the stratified identification 
\begin{equation}
\overline{\sC}\big(\bDelta^{p-1}\xra{m}[p-1]\big) \, \simeq \, \big(\bDelta^{p} \xra{m} [p] \big).
\end{equation}
\end{observation}

The next theorem is both a definition and a theorem.

\begin{theorem}\cite[Corollary 3.3.6]{AFR1}
    The restricted Yoneda embedding associated to $\mathsf{st}\maps \bDelta \ra \Strat$ 
    \begin{equation}
    \exit \maps \Strat \xra{\,\,\mathsf{Yoneda}\,\,} \psh(\Strat) \xra{\,\,\,\mathsf{st}^{*}\,\,\,} \psh(\bDelta)
    \end{equation}
    lands in complete Segal spaces. We call the value of $\exit$ on a conically smooth stratified space its \textbf{exit-path $(\oo,1)$-category}. 
\end{theorem}

\begin{observation}
Unpacking this definition we see that given a conically smooth stratified space $\big(M \ra P\big)$, the space of objects of its exit-path $(\oo,1)$-category is:
\begin{equation}
    \exit(M)_{0} = \Strat(*,M)
\end{equation}
which as a Kan complex has as $q$-simplices the set $\Strat(*\times \AA^q, M)$, an element of which, as a map out of a smooth manifold, must factor through a single smooth stratum $M_p$ for some $p\in P$. The space of objects can thus be identified with
\begin{equation}
    \exit(M)_{0} \,\, \simeq \,\, \coprod_{p \in P} M_p
\end{equation} 
the space which is the disjoint union of the individual strata.
\end{observation}

\begin{observation} Observe that exit-paths do not compose strictly. Exit-paths are the $0$-simplices $\big(\gamma\maps \bDelta^1 \ra \big(M \ra P\big)\big)\in \Strat(\bDelta^{1}\times \AA^0, M)$ — these are stratified maps with an important feature: as any map out of $[0,1]$, it traces a path in $M$, but as a stratified map it forces the path to start at a deeper stratum than the one it ends. Moreover, once the path leaves its starting stratum — if it does — it never returns to it. In a very precise sense, exit-paths are extremely restricted.

Their restricted nature invites deeper reflection on how they compose — given the fact that they form an $(\oo,1)$-category. At a distance, they do not compose strictly, as is exemplified by the next figure
\begin{equation}
    \tikzfig{figures/Strat/STR14}{\input{figures/Strat/STR14.tikz}}
\end{equation}
where the underlying stratified space is $\underline{\RR}^2 =\big(\RR^2 \ra [2]\big)$. The two exit-paths $\gamma_1,\gamma_2\maps \bDelta^{1} \ra \RR^2$ are such that the source of $\gamma_2$ is the target of $\gamma_1$, which would allow for a composite continuous map $|\bDelta^1| \ra \RR^2$ -- yet, as stratified maps i.e., exit-paths, they cannot compose strictly. Of concern is whether the composition is well-defined up to homotopy. To see why this is true we make use of the homotopy colimit associated to the \emph{double closed cones} in $\Strat$

\begin{equation}
    \begin{tikzcd}
        {\overline{\sC}(\emptyset)} && {\overline{\sC}(L)} \\
        \\
        {\overline{\sC}^2(\emptyset)} && {\overline{\sC}^2(L)}
        \arrow["{\overline{\sC}(\emptyset \hookra L)}", hook, from=1-1, to=1-3]
        \arrow["{\{1\}}"', hook, from=1-1, to=3-1]
        \arrow["{\{1\}}", hook, from=1-3, to=3-3]
        \arrow["{\overline{\sC}^2(\emptyset \hookra L)}"', hook, from=3-1, to=3-3]
        \arrow["\lrcorner"{anchor=center, pos=0.125, rotate=180}, draw=none, from=3-3, to=1-1]
    \end{tikzcd}
\end{equation}
where $L$ is any compact conically smooth stratified space — note that in $\strat$ this fails to be a colimit, indicating the need to work in the suitably homotopical context of $(\oo,1)$-categories. An implication of this diagram being a pushout in $\Strat$, is that for every other conically smooth stratified space $Z$, the diagram 

\begin{equation}
    \label{ch4:: eq: exit-paths compose}
    \begin{tikzcd}
        {\Strat\big(\overline{\sC}(\emptyset),Z\big)} && {\Strat\big(\overline{\sC}(L),Z\big)} \\
        \\
        {\Strat\big(\overline{\sC}^2(\emptyset),Z\big)} && {\Strat\big(\overline{\sC}^2(L),Z\big)}
        \arrow["{\overline{\sC}(\emptyset \hookra L)}"', from=1-3, to=1-1]
        \arrow["{\{1\}}", from=3-1, to=1-1]
        \arrow["{\{1\}}"', from=3-3, to=1-3]
        \arrow["{\overline{\sC}^2(\emptyset \hookra L)}", from=3-3, to=3-1]
    \end{tikzcd}
\end{equation}
is a pullback of Kan complexes in $\sset$. If $Z$ is $(M\ra P)$ and if $L=\Delta^{\{0\}}$, then we have that 
\begin{equation}
     \Strat\big(\overline{\sC}(\Delta^{\{0\}}),Z\big)\underset{\Strat\big(\overline{\sC}(\emptyset),Z\big)}{\times} \Strat\big(\overline{\sC}^2(\emptyset),Z\big)
\end{equation}
records the space of composable exit-paths. The pullback condition then says that given a pair of composable exit-paths $\gamma_1,\gamma_2$, up to homotopy, there is a unique exit 2-simplex $\sigma\maps \Delta^{\{0<1<2\}} \ra M$ whose $\Delta^{\{0<1\}}$ and $\Delta^{\{1<2\}}$ faces are $\gamma_1$ and $\gamma_2$, respectively.
\end{observation}

\subsubsection{Embedding into Euclidean space}

Studying spaces of embeddings of smooth manifolds has turned out almost as important and fruitful as studying the manifolds themselves. In this section, we initiate the exploration of the homotopical properties of embedding spaces of conically smooth stratified spaces into Euclidean space. We begin by posing a conjecture that every conically smooth stratified space $M$ of dimension $n$ admits a conically smooth embedding into $\RR^{2n + 1}$ — we prove it for the special case of depth 1, and indicate a possible path to the proof for general depth —, we prove that every conically smooth stratified space $M$ admits a conically smooth embedding into some Euclidean space, that the space of conically smooth embeddings of $M$ into $\RR^{\oo}$ is weakly contractible, and finally that every conically smooth stratified space $M$ with boundary admits a neat embedding into $\RR_{\geq 0} \times \RR^{m}$.  

One of the most important results in differential topology is the \emph{Whitney embedding theorem}. In its weakest form it asserts that every compact smooth manifold admits a proper embedding into Euclidean space.

\begin{theorem}[Whitney embedding theorem — weak version]\label{ch4:: thm: Whitney weak}
    Let $M$ be a compact smooth manifold. There exists a smooth embedding $M \hookra \RR^{N}$ for a sufficiently large $N$
\end{theorem}

\noindent The proof of this result rests almost entirely on the fact that smooth manifolds admit partitions of unity, and Whitney showed that it can be improved significantly by using transversality arguments. 

\begin{theorem}[Whitney embedding theorem — strong version]\label{ch4:: them: Whitney strong}
    \label{thm:ch7::whitney strong}
Let $M$ be a compact smooth $n$-manifold. There exists an immersion of $M$ into $\RR^{2n-1}$ and an embedding into $\RR^{2n}$.     
\end{theorem}

In working with conically smooth stratified spaces, which are defined abstractly — as opposed to Whitney stratified spaces, which are always explicitly embedded into a smooth manifold — in order to understand the topology of their embedding spaces, we must first show a result analogous to the Whitney embedding theorem, ideally in its strongest form. Before doing so, we introduce a quick technical result.

\begin{lemma}\cite[Lemma~1.3.1]{AFR1}\label{ch4:: lem: shrinking} For every open cover $\cU$ of a conically smooth stratified space $K$, there exists a refining open cover $\cV$ with the property that for each $V\in \cV$, the image of $V \hookra K$ factors as
    \begin{equation}
        \begin{tikzcd}
            V && {\RR^{i}\times \sC(L)} \\
            \\
            && {\DD^{i}\times \overline{\sC}(L)}
            \arrow[from=1-1, to=1-3]
            \arrow["i"', hook, from=1-1, to=3-3]
            \arrow["{\overline{\phi}}"', hook, from=3-3, to=1-3]
        \end{tikzcd}
    \end{equation}
    where the vertical map is a conically smooth embedding and the diagonal map is an open embedding i.e., an isomorphism into the interior $\RR^{i}\times \sC(L) \subset \DD^{i}\times \overline{\sC}(L)$.
\end{lemma}

Using it, we can prove the stratified analogue of Whitney's weak embedding theorem \Cref{ch4:: thm: Whitney weak}.

\begin{proposition}\label{ch4:: prop: weak Whit strat}
    Let $\big(M \ra P\big)$ be a finite conically smooth stratified space. There is a proper conically smooth embedding $\big(M \ra P\big) \hookra \RR^{N}$ for some $N$. 
\end{proposition}
\begin{proof}
    We argue by induction on depth. If $\mathsf{depth}(M)=0$, the result is the usual weak Whitney embedding theorem for smooth manifolds.

    Assume $\mathsf{depth}(M)>0$. By \Cref{ch4:: lem: shrinking} and finiteness of $M$, we can choose a finite conically smooth atlas $\cA=\{U_i\}_{i=1}^{k}$ such that each chart admits a conically smooth identification
    \begin{equation}
        U_i \simeq \RR^{l_i}\times \sC(Z_i),
    \end{equation}
    where $Z_i$ is compact and satisfies $\mathsf{depth}(Z_i)<\mathsf{depth}(M)$. Choose a conically smooth partition of unity $\{\psi_i\}_{i=1}^{k}$ subordinate to this cover.

    By the inductive hypothesis, for each $i$ there exists a proper conically smooth embedding
    \begin{equation}
        \varphi_i \maps Z_i \hookra \RR^{N_i} \hookra S^{N_i}.
    \end{equation}
    Taking open cones and then the product with the identity on $\RR^{l_i}$ gives a conically smooth embedding
    \begin{equation}
        f_i \maps U_i \hookra \RR^{l_i}\times \sC(S^{N_i}) \simeq \RR^{l_i+N_i+1}.
    \end{equation}

    Because $\psi_i$ is supported in $U_i$, both maps $x\mapsto \psi_i(x)$ and $x\mapsto \psi_i(x)\,f_i(x)$ extend by zero to global conically smooth maps on $M$. Define
    \begin{equation}
        F := \big(\psi_1,\psi_1f_1,\ldots,\psi_k,\psi_kf_k\big)\maps M \longra \RR^{N}
    \end{equation}
    for a sufficiently large $N$.

    The same argument as in the smooth case now shows that $F$ is a conically smooth embedding: the scalar components $\psi_i$ detect an index for which $\psi_i>0$, and on the open set $\{\psi_i>0\}$ the pair $\big(\psi_i,\psi_i f_i\big)$ recovers the local embedding $f_i$. Hence $F$ is injective and locally a conically smooth embedding, therefore a global conically smooth embedding.

    Finally, only finitely many compactified charts occur in the construction, so $F$ is proper. This completes the inductive step.
\end{proof}

With this result, we can now prove the following important structural proposition about the space of conically smooth embeddings into Euclidean spaces of sufficiently high dimension.

\begin{proposition} \label{ch4:: prop: contractibility of embeddings}
    Let $\big(M \ra P\big)$ be a finite conically smooth stratified space. The space 
    \begin{equation}
        \Strat^{\sf emb}\big(M, \RR^{\oo}\big) \, \overset{\sf def}{=} \, \hoColimit[\RR^{k}\hookra \RR^{k+1}] \Strat^{\sf emb}(M, \RR^{k})
    \end{equation} 
    is weakly contractible.
\end{proposition}

\begin{proof}
We first observe that the vector space $\RR^{\oo}$ admits a short exact sequence 
\begin{equation}
    0 \longra \RR^{\oo} \xlongra{S^{\sf ev}} \RR^{\oo} \xra{P^{\sf ev}} \RR^{\oo} \longra 0.
\end{equation}

\noindent The two maps are defined as follows 
\begin{equation}
    S^{\sf ev}\big(x_i\big) \, = \, (0,x_1, 0,x_2, 0, \cdots), \quad P^{\sf ev}\big(x_i\big) \, = \, (x_1, x_3, x_5, \cdots).
\end{equation}
Equivalently, there is a linear isomorphism 
\begin{equation}
    \RR^{\oo} \oplus \RR^{\oo} \xra{\,\,\,\,\,S^{\sf ev} + S^{\sf od}\,\,\,\,\,} \RR^{\oo},
\end{equation}
where $S^{\sf od}\big(x_i\big) = (x_1, 0, x_2, 0, \cdots)$.  

Let us consider the smooth homotopy 
\begin{equation}
    H \maps \RR^{\oo} \times \AA^1 \longra \RR^{\oo}\times \AA^1, \quad (x,t) \longmapsto \big((1-t)\cdot x + t \cdot S^{\mathsf{ev}}(x), t \big)
\end{equation}
Assume that $H(-,t)$ fails to be injective for some $t'\in \AA^1$. Since the case $t'=0$ is just the identity, we may assume $t'\neq 0$, in which case there exist $x-x' \neq 0$ satisfying 
\begin{equation}
    (1-t')x + t' S^{\mathsf{ev}}(x) = (1-t')x' + t' S^{\mathsf{ev}}(x').
\end{equation}
This implies that 
\begin{equation}
S^{\sf ev}(x-x') = -\frac{1-t'}{t'}(x-x')
\end{equation}
or that $(x-x')$ is an eigenvector for $S^{\sf ev}$ — a contradiction.
\smallskip
\noindent We will now use this (conically) smooth isotopy of (conically) smooth embeddings to build a simplicial homotopy  
\begin{equation}
    \Strat^{\sf emb}(M, \RR^{\oo})\times \Delta^{1} \longra \Strat^{\sf emb}(M, \RR^{\oo})
\end{equation}
from $\id_{\Strat^{\sf emb}(M,\RR^{\oo})}$ to $(S^{\sf ev})_{*}$. First, observe that such a simplicial homotopy amounts to a simplicial map 
\begin{equation}
    \Strat^{\sf emb}(M, \RR^{\oo}) \longra \Strat^{\sf emb}(M \times \AA^1, \RR^{\oo})
\end{equation}   
owing to the constructibility of $\Strat^{\sf emb}$. We define this map as follows:
\begin{align}    
    \big(\, \varphi_p \maps M \times \AA^p \ra \RR^{\oo} \times \AA^p \,\big) &\longmapsto 
    \big(\, M \times \AA^p \times \AA^{1} \xra{\varphi_p \times \id} \RR^{\oo} \times \AA^p \times \AA^1  \\
    &\xra{\,\,\id \times \mathsf{swap}\,\,} \RR^{\oo} \times \AA^1 \times \AA^p \\
    &\xra{\,\, H \times \id \,\,} \RR^{\oo} \times \AA^{1} \times \AA^p \\ 
    &\xra{\,\, \sf pr \,\,} \RR^{\oo} \times \AA^p \, \big).
\end{align}
The map is well-defined since $H\big(\varphi_p(x,s),t\big)$ is a conically smooth embedding for every $s \in \AA^p$. Moreover, we have that 
\begin{equation}
    H(\varphi_p(x,s),0) = \big(\varphi_{p}(x,s),0\big), \quad H(\varphi_{p}(x,s),1) = \big(S^{\sf ev}(\varphi_{p}(x,s)),1\big),
\end{equation}
further observing the correct start- and end-points for the conically smooth isotopy.

From \Cref{ch4:: prop: weak Whit strat} we know there exists a conically smooth embedding $e \maps M \hookra \RR^{\oo}$. Since $M$ is finite, after postcomposing with a standard embedding of $\RR^{N}$ into a sphere in $\RR^{N+1}$ and then into $\RR^{\oo}$, we may assume the image of $e$ lies in the unit sphere. We now articulate a contraction of $\Strat^{\sf emb}(M,\RR^{\oo})$ onto the point represented by $S^{\sf od}\circ e$:

\begin{align}
    \Strat^{\sf emb}(M,\RR^{\oo})_p &\longra \Strat^{\sf emb}(M\times \AA^{1},\RR^{\oo})_p \\
    \big(\, \varphi_p \maps M \times \AA^p \ra \RR^{\oo} \times \AA^p \,\big) &\longmapsto \big[(x,t,s) \mapsto \big((1-t)\, S^{\sf ev}(\varphi_p(x,s)) + t\, S^{\sf od}(e(x)), s\big)\big].
\end{align}
The map is well-defined because the even coordinates recover $(1-t)\varphi_p(x,s)$ while the odd coordinates recover $t\,e(x)$. Since $e(x)$ lies on the unit sphere, its norm recovers $t$, and then $e$ recovers $x$; together with the final $\AA^p$-coordinate this shows injectivity. This provides a conically smooth isotopy of conically smooth embeddings starting at $(S^{\sf ev})_*$ and ending at the constant point $S^{\sf od}\circ e$. Through a suitable rescaling, both isotopies can be concatenated to connect $\id_{\Strat^{\sf emb}}$ to a constant map. 
\end{proof}

\Cref{ch4:: prop: contractibility of embeddings} will pay off when we attempt to construct the $(\oo,1)$-category of stratified bordisms using conically smooth embedded stratified spaces together with carefully chosen submersivity conditions — as we did in Chapter 3.   

We now consider the strong version of the Whitney embedding theorem for conically smooth stratified spaces. We do not prove it for arbitrary depth, instead leaving it as the following conjecture.

\begin{conjecture}\label{ch4:: conj: whitney strat}
Let $(M \ra P)$ be a compact $d$-dimensional conically smooth stratified space. There exists a conically smooth proper embedding 
\begin{equation}
    (M \ra P) \hookra \RR^{2d + 1}.
\end{equation}
\end{conjecture}

We provide a proof for depth 1, but before introducing the proof, we mention the work of Hiroko Natsume proving \Cref{ch4:: conj: whitney strat} for Whitney stratified spaces \cite{N80}. In particular, we make use of the following result in \cite{N80}.

\begin{proposition}\cite[Lemma~1]{N80}\label{ch4:: prop: fiberwise emb}
    Let $\pi\maps E \ra M$ be a fiber bundle with a compact manifold $F$ as the fiber and $h$ a diffeomorphism of $M$ onto $M'$. Let $S$ be a sphere bundle over $M'$ and $p$ the projection of $S$ to $M'$ 
    \begin{equation}
        \begin{tikzcd}[ampersand replacement=\&]
            E \&\& S \\
            \\
            M \&\& {M'}
            \arrow["\pi"', two heads, from=1-1, to=3-1]
            \arrow["p", two heads, from=1-3, to=3-3]
            \arrow["\simeq", from=3-1, to=3-3]
        \end{tikzcd}
    \end{equation}
    Assuming that 
    \begin{equation}
        \mathsf{dim}(S) \geq 2\cdot \mathsf{dim}(E) + 2,
    \end{equation}
    then there is a fiberwise smooth embedding $\tilde{h}\maps E \hookra S$ 
    \begin{equation}
        \begin{tikzcd}[ampersand replacement=\&]
            E \&\& S \\
            \\
            M \&\& {M'}
            \arrow["{\tilde{h}}", dashed, hook, from=1-1, to=1-3]
            \arrow["\pi"', two heads, from=1-1, to=3-1]
            \arrow["p", two heads, from=1-3, to=3-3]
            \arrow["\simeq", from=3-1, to=3-3]
        \end{tikzcd}
    \end{equation} 
\end{proposition}

Lastly, before proving the stratified strong embedding theorem for depth 1, we must introduce yet another theorem of Whitney.

\begin{theorem}[Whitney extension theorem]\label{ch4:: thm: whitney extension theorem}
    Let $A$ be a smooth $m$-dimensional manifold, and let $B\subset A$ be a closed subset of $A$. Let $M=A\setminus B$ be a smooth manifold, and let $f\maps B \longra \RR^{n}$. Then there is a map $F$ of $A$ into $\RR^{n}$ with the following properties:
    \begin{enumerate}
        \item $F_{|B}=f$.
        \item $F$ is smooth.
        \item If $n\geq 2m +1$ and $f$ is an embedding of $B$, then $F$ is an embedding of $A$.
    \end{enumerate}
\end{theorem}

We are now ready to state and prove the strong stratified Whitney embedding theorem in depth 1.

\begin{proposition}\label{ch4:: prop: Whitney strong depth 1}
    Let $\big(  M\longra P \big)$ be a compact conically smooth stratified space of dimension $n$ and depth 1. There exists a proper conically smooth embedding $M \hookra \RR^{2n + 1}$.
\end{proposition}

\begin{proof}
    Since $M$ has depth $1$, we may assume that the poset $P$ is isomorphic to $\{d<n\}$. Write $M_d$ for the deepest stratum and $M_n$ for the top-dimensional stratum.
    
    The stratum $M_d$ is a compact smooth $d$-manifold, so \Cref{thm:ch7::whitney strong} gives a smooth embedding
    \begin{equation}
        \varphi_d \maps M_d \hookra \RR^{2d} \hookra \RR^{2n+1}.
    \end{equation}
    The image $\varphi_d(M_d)$, together with its complement, endows $\RR^{2n+1}$ with a stratification 
    \begin{equation}
        \label{eq:ch7::image of embedding}
        \big(\RR^{2n+1} \longra \{d < (2n + 1)\}\big)
    \end{equation}    
    which is conically smooth and of type $\sD_{d\subset 2n+1}$. Our goal is to extend $\varphi_d$ to a conically smooth embedding of all of $M$. We first extend it over a regular neighborhood of $M_d$, and then over the top stratum.

    Taking the closures of the regular neighborhoods provides the following identifications of conically smooth stratified spaces with boundary:
    \begin{equation}
        \label{eq:ch7::regular neighborhood}
        \overline{\sC}(\pi_d) ~\simeq~ \Big(\Link_{d}(M) \times [0,1] \Big)\underset{\Link_d(M)}{\bigcup} M_d, \qquad \overline{\sC}(\widetilde{\pi}_d) ~\simeq~ \Big(\Link_{d}\big(\RR^{2n+1}\big) \times [0,1] \Big) \underset{\Link_d(\RR^{2n+1})}{\bigcup} \varphi_d(M_d).
    \end{equation}
    These identifications provide coordinates $[z,t]\in \overline{\sC}(\pi_d)$ such that $[z,t] = (z,t)$ for $z\in \Link_d(M)$ and $t\in (0,1]$, and $[z,0] = \pi_d(z)$ for $z\in\Link_d(M)$ and $t=0$. Using these coordinates, we begin by constructing a fiberwise embedding of $\Link_d(M)$ into $\Link_d(\RR^{2n+1})$ — this is possible due to \Cref{ch4:: prop: fiberwise emb}, since 
    \begin{equation}
        \mathsf{dim}\big(\Link_{d}(M)\big) = n - 1, \quad \text{and} \quad \mathsf{dim}\big(\Link_{d}(\RR^{2n+1})\big) = 2n.
    \end{equation} and $2(n-1)+2 \leq 2n$.
    Let us denote this fiberwise embedding by $\widetilde{\varphi_d}$. We can now extend this embedding as follows:
    \begin{equation}
        \overline{\varphi}_d([z,t]) = 
        \begin{cases}
            (\widetilde{\varphi_d}(z),t) & \text{ for } t>0 \\
            \varphi_d(\pi_d(z)) & \text{ for } t = 0.
        \end{cases}
    \end{equation} 
    This map is conically smooth, since it is the dashed arrow in the following commutative diagram in $\Strat$:
    \begin{equation}
    \begin{tikzcd}[ampersand replacement=\&]
        \& {\Link_d(\RR^{2n+1})} \&\& {\Link_d(M)\times [0,1]} \\
        {\Link_d(M)} \& {} \& {\Link_d(M)\times [0,1]} \\
        \& {\varphi_d(M_d)} \&\& {\overline{\sC}(\widetilde{\pi}_d)} \\
        {M_d} \&\& {\overline{\sC}(\pi_d)}
        \arrow[""{name=0, anchor=center, inner sep=0}, "0", hook, from=1-2, to=1-4]
        \arrow[two heads, from=1-2, to=3-2]
        \arrow[from=1-4, to=3-4]
        \arrow["{\widetilde{\varphi_d}}"{pos=0.2}, from=2-1, to=1-2]
        \arrow["0"{pos=0.7}, hook, from=2-1, to=2-3]
        \arrow[two heads, from=2-1, to=4-1]
        \arrow["{\widetilde{\varphi_d}\times \id}"'{pos=0.6}, from=2-3, to=1-4]
        \arrow[from=2-3, to=4-3]
        \arrow[from=3-2, to=3-4]
        \arrow["{\varphi_d}", from=4-1, to=3-2]
        \arrow[from=4-1, to=4-3]
        \arrow["\lrcorner"{anchor=center, pos=0.125, rotate=180}, draw=none, from=4-3, to=2-1]
        \arrow[dashed, from=4-3, to=3-4]
        \arrow["\lrcorner"{anchor=center, pos=0.125, rotate=180}, draw=none, from=3-4, to=0]
    \end{tikzcd}
    \end{equation}
    Thus far, we have embedded $\overline{\sC}(\pi_d)$ into $\RR^{2n+1}$. To extend this to the entire $(M\longra P)$, we work on the top stratum $M_n$, which is a smooth $n$-manifold. The restriction of $\overline{\varphi}_d$ to the closed subset $\overline{\sC}(\pi_d)\cap M_n \subset M_n$ is a smooth embedding into $\RR^{2n+1}$, so by \Cref{ch4:: thm: whitney extension theorem} it extends to a smooth embedding of all of $M_n$ into $\RR^{2n+1}$. Because this extension agrees with $\overline{\varphi}_d$ on the overlap $\overline{\sC}(\pi_d)\cap M_n$, the two pieces glue to a globally defined conically smooth embedding $M \hookra \RR^{2n+1}$. Since $M$ is compact, this embedding is automatically proper.     
\end{proof}

A manifold with boundary can be embedded into another manifold with boundary in rather pathological ways, and the same can be said about conically smooth stratified spaces. The next definition fixes pathological behaviour by requiring additional regularity in the neighborhood of the boundary. An embedding is \emph{neat} if it admits a product-like structure near the boundary.

\begin{proposition}
    \label{ch4:: prop: neat embedding strat}
    Let $\big(M \ra P\big)$ be a compact conically smooth stratified space with boundary $\big(M \ra P \fib \cP(*)\big)$. There exists a neat embedding 
    \begin{equation}
        \big(M \ra P\big) \hookra \RR_{\geq 0} \times \RR^{m}
    \end{equation}
    for some $m \geq 0$ --- where $\big(\RR_{\geq 0} \times \RR^{m}\big)$ admits the usual boundary structure.
\end{proposition}
\begin{proof}
    From \Cref{ch4:: prop: weak Whit strat}, we have a proper conically smooth embedding 
    \begin{equation}
        e \maps M \hookra \RR^{m}.
    \end{equation}
    Using an inward-pointing vector field and its flow, we obtain a collar
    \begin{equation}
        c\maps \partial M\times [0,\epsilon) \hookra M
    \end{equation}
    together with the associated boundary-defining function $\rho \maps M \longra \RR_{\geq 0}$. By the standard collar-straightening argument, after modifying $e$ by a conically smooth isotopy supported in the image of $c$, we may assume that on $c(\partial M\times [0,\epsilon/2))$ the embedding $e$ is independent of the collar coordinate.

    Taking the product of $e$ and $\rho$ gives a conically smooth map 
    \begin{equation}
        M \xra{\;\; \rho \times e \;\;} \RR_{\geq 0} \times \RR^{m}
    \end{equation}
    which is proper and injective. On the collar neighborhood $c(\partial M\times [0,\epsilon/2))$, it has the form
    \begin{equation}
        (p,t)\longmapsto \big(t,e_{|\partial M}(p)\big),
    \end{equation}
    so it is product-like near the boundary. This is exactly the neatness condition.
\end{proof}

\section{Constructible tangent bundle and structures}

For a conically smooth stratified space $M$, the constructible tangent bundle is a functor
\begin{equation}
    \sT_M \maps \exit(M)\longra \Vect^{\inj}.
\end{equation}
Concretely, it sends a point $x\in M_p$ to the ordinary tangent space $\sT_xM_p$ of the stratum containing it, and sends an exit-path to an injective linear map. The target $\Vect^{\inj}$ reflects the basic geometric feature that, as one exits a deeper stratum into a more generic one, tangent directions may be gained but not lost.

More precisely, $\sT_M$ is obtained as the composite
\begin{equation}
    \exit(M) \xra{\;\Delta_M\;} \Exit \xra{\;\sT^{\sf fib}\;} \Vect^{\inj}.
\end{equation}
The first functor $\Delta_M$ classifies the trivial constructible bundle
\begin{equation}
    M\times M \xra{\;\mathsf{pr}_1\;} M
\end{equation}
together with the diagonal section $\mathsf{diag}\maps M \longra M\times M$. The second functor $\sT^{\sf fib}$ should be thought of as the universal operation that linearizes such pointed constructible bundles. We now briefly introduce the $(\oo,1)$-categories $\Bun$, $\Exit$, and $\Vect^{\inj}$ that make this construction possible.

\begin{intuition}
One of the central constructions of Ayala--Francis--Rozenblyum (see \cite{AFR2}) is the $(\oo,1)$-category $\Bun$. Its basic universal property is that functors 
\begin{equation}
    \exit(K) \longra \Bun
\end{equation}
classify constructible bundles over $K$. The $(\oo,1)$-category $\Exit$ refines $\Bun$ by remembering a section: functors
\begin{equation}
    \exit(K) \longra \Exit
\end{equation}
classify constructible bundles over $K$ equipped with a section. The name of this category — the \emph{absolute exit-path $(\oo,1)$-category} — comes from the following property, recorded in \cite[Proposition~6.4.3]{AFR1}: for every constructible bundle $f\maps X \ra K$, the diagram
\begin{equation}
        \begin{tikzcd}[ampersand replacement=\&]
            \exit(X) \&\& \Exit \\
            \\
            \exit(K) \&\& {\Bun}
            \arrow["{(X \xra{f} K)}", hook, from=1-1, to=1-3]
            \arrow["\exit(f)"', two heads, from=1-1, to=3-1]
            \arrow[two heads, from=1-3, to=3-3]
            \arrow["(X \xra{f} K)", from=3-1, to=3-3]
        \end{tikzcd}
\end{equation} 
is a homotopy pullback diagram, indicating that the fibers over $\{f\maps X \ra K\}$ are given by $\exit(X)$. 
\end{intuition}

In particular, the projection $\mathsf{pr}_1\maps M\times M \longra M$ together with the diagonal section $\mathsf{diag}\maps M \longra M\times M$ determines a canonical functor
\begin{equation}
    \Delta_M \maps \exit(M)\longra \Exit.
\end{equation}

We now move on to the $(\oo,1)$-category $\Vect^{\sf inj}$.

\begin{intuition}
The $(\oo,1)$-category $\Vect^{\inj}$ is a generalization of the $\oo$-groupoid 
\begin{equation}
    \coprod_{i \geq 0} \BO(i).
\end{equation}
If $M$ is a smooth manifold, interpreted also as an $\oo$-groupoid, then functors
\begin{equation}
    M \longra \coprod_{i} \BO(i)
\end{equation}
classify vector bundles over $M$. The same role is played by $\Vect^{\inj}$ but for \emph{constructible} vector bundles \footnote{For simplicity, we avoid defining constructible vector bundles yet the careful reader might find it in \cite[Definition~2.2]{AFR2}.} over a conically smooth stratified space $M$, if we replace $M$ by its exit-path $(\oo,1)$-category $\exit(M)$. In practice, due to the $\inj$ superscript, we retain only constructible vector bundles whose transition maps are injective. Equivalently, along an exit-path the fiber dimension can stay the same or increase, but it cannot decrease; this is what lower semi-continuity means in the present context. This $(\oo,1)$-category is instantiated as a complete Segal space by defining it as follows:
\begin{equation}
    [p] \longmapsto \bigg|\bigg\{V_0 \hookra V_1 \hookra \cdots \hookra V_p \text{ over } \bDelta^{\bullet}\bigg\}\bigg|
\end{equation} 
where the right-hand-side is the space of injective composable sequences of finite rank vector bundles over $\bDelta^\bullet$, and such that for each $0\leq i \leq j \leq p$, the composition
\begin{equation}
    f_{j} \circ \cdots \circ f_{i}
\end{equation}
has locally constant rank. 
\end{intuition}

The remaining ingredient is the universal functor which extracts fiberwise linear data from pointed constructible bundles.
\begin{theorem}\label{ch4:: def: cbl tgt bundle}\cite[Proposition~2.15]{AFR2}
    There exists a unique symmetric monoidal functor under $\Vect^{\sim}$,
    \begin{equation}
        \sT^{\sf fib}\maps \Exit \longra \Vect^{\inj}
    \end{equation}
    such that:
    \begin{itemize}
        \item \label{closed-morphism} Sends closed morphisms $(\Strat)^{\cls} \hookra \Exit$ to equivalences.
        \item \label{embeddings} Sends embeddings $(\Strat)^{\emb} \hookra \Exit$ to equivalences.
    \end{itemize}
\end{theorem}

These two requirements encode the basic functorialities one expects from any object deserving the name tangent bundle. In the smooth case, if \(f\colon M \hookrightarrow N\) is an open embedding, then the derivative $\sD f$ identifies \(\sT M\) with the pullback of \(\sT N\), so tangent data is local for embeddings. In any suitably differential stratified setting, one further expects a restriction phenomenon along closed subspaces and closed strata; morally, passing to such a piece should not change the intrinsic linear data there. Requiring \(\sT^{\sf fib}\) to send these morphisms to equivalences is precisely the categorical way of encoding those expectations.

\begin{definition}
    The \define{constructible tangent bundle} of a conically smooth stratified space $M$ is the composite
    \begin{equation}
        \sT_M \maps \exit(M) \xra{\;\Delta_M\;} \Exit \xra{\;\sT^{\sf fib}\;} \Vect^{\inj}.
    \end{equation}
\end{definition}

This formal definition matches the concrete description given at the start of the section: $\sT_M$ sends a point $x\in M_p$ to the ordinary tangent space $\sT_xM_p$, and an exit-path to the corresponding injective linear map. A more conceptual and detailed discussion of how this composite is evaluated, both on objects and on exit-paths, can be found in \cite[Section~2]{AFR2}.

Having set up the constructible tangent bundle, we are naturally brought to the notion of conically smooth stratified \emph{tangential structure}.

\begin{definition}\label{ch4:: def: tangential structure}\cite[Definition~2.17]{AFR2}
    A \define{tangential structure} is an $(\oo,1)$-category $\cB$ equipped with a functor $\cB \longra \Vect^{\sf inj}$, and the $(\oo,1)$-category of tangential structures is the overcategory
    \begin{equation}
        (\Cat_{\oo})_{/\Vect^{\inj}}.
    \end{equation}
\end{definition}

This leads to the next definition.

\begin{definition}\label{ch4:: def: tangential structure 2}\cite[Definition~2.19]{AFR2}
    Fix a tangential structure $\cB \longra \Vect^{\inj}$ and a conically smooth stratified space $M$. A $\cB$\define{-framing on} $M$ is a lift $\varphi$
    \begin{equation}
    \begin{tikzcd}[ampersand replacement=\&]
        \& \cB \\
        {\exit(M)} \& {\Vect^{\inj}}.
        \arrow[from=1-2, to=2-2]
        \arrow["\varphi", dashed, from=2-1, to=1-2]
        \arrow["{\sT_M}"', from=2-1, to=2-2]
    \end{tikzcd}
    \end{equation} 
\end{definition}

\begin{remark}
    Note that, if $M$ is a smooth manifold then every functor $M \longra \cB$ factors through $\cB^{\sim} \hookra \cB$, the maximal $\oo$-subgroupoid. In particular this means that a $\cB$-framing for $M$ is a direct generalization of a $B$-framing, as discussed in Chapter 3.
\end{remark}

\begin{example}
    For each $n \geq 0$, the full subcategory $\Vect^{\sf inj}_{\leq n} \hookra \Vect^{\sf inj}$ on vector spaces of dimension at most $n$ is manifestly a tangential structure. A $\Vect^{\sf inj}_{\leq n}$-framing on a conically smooth stratified space $M$ is unique if it exists, and it exists precisely when all strata of $M$ have dimension at most $n$. This will be useful to us when discussing certain categories of cobordisms.
\end{example}

In this thesis we will focus exclusively on one especially tractable class of stratified tangential structures — \emph{solid tangential structures}.

\begin{definition}\cite[Definition~2.35]{AFR2}
    A \define{solid tangential structure} is a right fibration
    \begin{equation}
        \cB \longra \Vect^{\sf inj}.
    \end{equation}
\end{definition}

By straightening, to give a solid tangential structure is equivalent to stipulating the data of a functor
\begin{equation}
    F \maps (\Vect^{\sf inj})^{\op} \longra \Spaces.
\end{equation}
Thus, for every finite-dimensional vector space $V$, one has a space $F(V)$ of allowed structures on $V$, and for every injective linear map $V \hookra W$, one has a restriction map
\begin{equation}
    F(W) \longra F(V).
\end{equation}

In particular, a solid framing on a conically smooth stratified space $M$ is a lift of the constructible tangent bundle
\begin{equation}
    \sT_M \maps \exit(M) \longra \Vect^{\sf inj},
\end{equation}
and hence amounts to a compatible choice of such structure on each tangent space, varying contravariantly along exit-paths.

The point of this definition is that solid structures only remember how tangential data restricts along injections; they do not impose the stronger coherence data encoded by more general tangential structures, such as vari-framings. In this sense they are less sensitive to the underlying stratification. Indeed, it is proven by Ayala--Francis--Rozenblyum in \cite[Proposition~2.61]{AFR2} that solid tangential structures are invariant under refinement of the stratification.

We are especially interested in a subclass of solid framings, namely the \emph{solid $B$-framings}, for $B$ a space.

\begin{definition}
    Let $B$ be a space equipped with a map $B \longra \Vect^{\inj}$. A \define{solid $B$-framing} on a conically smooth stratified space $M$ is a lift
    \begin{equation}
    \begin{tikzcd}[ampersand replacement=\&]
        \& B \\
        {\exit(M)} \& {\Vect^{\inj}}.
        \arrow[from=1-2, to=2-2]
        \arrow["\varphi", dashed, from=2-1, to=1-2]
        \arrow["{\sT_M}"', from=2-1, to=2-2]
    \end{tikzcd}
    \end{equation}
\end{definition}

\begin{remark}
    Since $B$ is a space, a solid $B$-framing remembers only a homotopy-coherent choice of $B$-structure on each constructible tangent space, compatible under restriction along exit-paths. In particular, it is even less sensitive to the incidence data of the strata than a general tangential structure, and by the result quoted above it is unchanged under refinement of the stratification.
\end{remark}

\begin{example}
    The simplest case is $B=*$. Choosing the object $\RR^n\in \Vect^{\inj}$ determines a map $* \longra \Vect^{\inj}$, and for a smooth $n$-manifold $M$ a solid $*$-framing is exactly an ordinary framing $\sT M \cong M\times \RR^n$.
\end{example}

To treat families of conically smooth stratified framed manifolds, we need the fiberwise version of the constructible tangent bundle. Let $\pi\maps M \longra K$ be a conically smooth fiber bundle. The pointed constructible bundle
\begin{equation}
    M\times_{K} M \xra{\;\mathsf{pr}_1\;} M
\end{equation}
together with the diagonal section $\mathsf{diag}\maps M \longra M\times_{K} M$ determines a canonical functor
\begin{equation}
    \Delta_{\pi}\maps \exit(M)\longra \Exit.
\end{equation}

\begin{definition}
    The \define{fiberwise constructible tangent bundle} of $\pi$ is the composite
    \begin{equation}
        \sT^{\mathrm{fib}}_{\pi}\maps \exit(M) \xra{\;\Delta_{\pi}\;} \Exit \xra{\;\sT^{\sf fib}\;} \Vect^{\inj}.
    \end{equation}
\end{definition}

This is the stratified analogue of the ordinary vertical tangent bundle of a smooth fiber bundle. Indeed, if $\pi\maps M \longra K$ is a smooth fiber bundle, then the usual vertical tangent bundle is the kernel of $\sD\pi$,
\begin{equation}
    \sT^{\mathrm{vert}}_{\pi}M\coloneqq \ker(\sD\pi)\longra M,
\end{equation}
and $\sT^{\mathrm{fib}}_{\pi}$ plays exactly this role in the constructible setting. In the special case $K=*$, the fiber product $M\times_{*}M$ is just $M\times M$, so $\sT^{\mathrm{fib}}_{\pi}$ reduces to the constructible tangent bundle $\sT_M$.

\begin{definition}
    Let $\cB \longra \Vect^{\inj}$ be a tangential structure and let $\pi\maps M \longra K$ be a conically smooth fiber bundle. A \define{vertical $\cB$-framing} on $\pi$ is a lift
    \begin{equation}
    \begin{tikzcd}[ampersand replacement=\&]
        \& \cB \\
        {\exit(M)} \& {\Vect^{\inj}}.
        \arrow[from=1-2, to=2-2]
        \arrow["\varphi", dashed, from=2-1, to=1-2]
        \arrow["{\sT^{\mathrm{fib}}_{\pi}}"', from=2-1, to=2-2]
    \end{tikzcd}
    \end{equation}
\end{definition}

If $B$ is a space and $\pi\maps X \longra *$ is the terminal map, then $\sT^{\mathrm{fib}}_{\pi}=\sT_X$, so a vertical solid $B$-framing on $\pi$ is exactly the usual solid $B$-framing on $X$; when $X$ is a smooth manifold, this is the ordinary $B$-framing of Chapter 3. This is the notion we will need later when discussing families of conically smooth stratified framed manifolds.

%% file: figures/Strat/STR28.tikz
\begin{tikzpicture}
	\begin{pgfonlayer}{nodelayer}
		\node [style=none] (3) at (0, -0.5) {};
		\node [style={wide_dot}] (4) at (-2, 0.5) {};
		\node [style={wide_dot}] (5) at (0, 0.5) {};
		\node [style={wide_dot}] (6) at (2, 0.5) {};
	\end{pgfonlayer}
	\begin{pgfonlayer}{edgelayer}
		\draw (4) to (5);
		\draw (5) to (6);
	\end{pgfonlayer}
\end{tikzpicture}

%% file: figures/Strat/STR30.tikz
\begin{tikzpicture}
	\begin{pgfonlayer}{nodelayer}
		\node [style=none] (0) at (-3, 0) {};
		\node [style=none] (1) at (3, 0) {};
		\node [style=0 strat] (2) at (0, 0) {};
	\end{pgfonlayer}
	\begin{pgfonlayer}{edgelayer}
		\draw (0.center) to (1.center);
	\end{pgfonlayer}
\end{tikzpicture}

%% file: Chapters/Ch5.tex
\chapter{An \texorpdfstring{$(\oo,1)$}--category of conically smooth stratified bordisms}

We now construct an interval Segal space whose associated $(\oo,1)$-category models $d$-dimensional unoriented conically smooth stratified bordisms. The strategy is the direct stratified
analogue of the construction of $\bord_d$ in \S3.2: we encode composable strings of bordisms by a single \emph{embedded ``total bordism''} over the core of an interval object $T\in \Int$, and we obtain the relevant classifying spaces by applying the concordance resolution to an appropriate smooth $\oo$-sheaf.

Before proceeding with definitions, let us sketch the path. Fix a finite-dimensional real vector space $V$. We first define a bifunctor
\begin{equation}
\tilbordstr_{d,V}\colon \man^{\op}\times \Int \to \set
\end{equation}
which assigns to $(M,T)$ a moduli-type set of bounded closed subsets
$W\subset M\times \core(T)\times V$ satisfying boundary-fixing, conically smooth stratified regularity, and fiberwise submersivity conditions. For each fixed $T$, the functor $\tilbordstr_{d,V}(-,T)$ will be shown to satisfy descent for open covers of manifolds, hence to define a (discrete) smooth $\oo$-sheaf. Restricting along the smooth cosimplicial manifold $\AA^{\bullet} \longra \man$ and taking geometric realization produces a functor
\begin{equation}
\bordstr_{d,V}\colon \Int \longra \sset, \qquad
T \longmapsto \left|\tilbordstr_{d,V}(\AA^\bullet,T) \right|.
\end{equation}
The main technical work of this section is to prove that $\bordstr_{d,V}$ is an interval Segal space, and that its homotopy $1$-category recovers the ordinary category $\cS\mathsf{Cob}_d$ of \Cref{ch5:: def: strat bord ord}.

While the formalism parallels the unstratified case, there are three geometric inputs specific to the
stratified setting: 
\begin{enumerate}
    \item We work with stratified \emph{topological} $d$-manifolds equipped with conically smooth atlases. 
    \item We require a stratified analogue of the differential-topological facts used in Section 3.2, notably that proper conically smooth submersions (also) behave like fiber bundles. 
    \item We rely on the existence and connectedness of conically smooth embeddings of the stratified manifolds under consideration into sufficiently large Euclidean spaces. 
\end{enumerate}
With these in place, we define $\tilbordstr_{d,V}$.

\section{The unoriented stratified bordism category}

The principal goal of this thesis is to construct an $(\oo,1)$-category of stratified bordisms, $\bordstr_d$. It is useful, however, to first consider the ordinary 1-category of stratified bordisms. This will provide us with a model against which we will compare the homotopy category of $\bordstr_d$.

\subsection{Stratified \texorpdfstring{$d$}{d}-manifolds}

In this thesis we are primarily interested in studying bordism categories of conically smooth stratified spaces with additional regularity features. These are designed around the following premise: physically relevant spacetimes in topological quantum field theories are \emph{topological manifolds}, at the very least. Positive codimension strata encode topological defects which are, in principle, themselves governed by a lower-dimensional topological QFT, thereby necessitating similar regularity features, and so on.

Given the construction of the site $\Bsc_{\leq d,d} \hookra \Bsc$, containing all basics of dimension $d$ and of depth bounded by $d$, as outlined in Chapter 4, one might have thought that this sieve captures exactly the singularity types relevant for topological QFTs with topological defects. This intuition fails, however, because being of \emph{pure dimension} $d$ is strictly less restrictive than being of \emph{constant} topological dimension $d$.   

We will now define a pair of $(\oo,1)$-categories inductively: 
\begin{enumerate}
    \item Let $\cS_{-1}$ be the full $(\oo,1)$-subcategory of $\Snglr$ on $\emptyset$ and let $\cS_{\leq -1} = \cS_{-1}$.
    \item Let $\Bsc_{0}^{\cS_{-1}}$ be the full $(\oo,1)$-subcategory of $\Bsc$ on the basic
    \begin{equation}
        U = \RR^{i}\times \sC(L),
    \end{equation}
    such that $i + \big(\mathsf{dim}(L)+1\big) = 0$ and $L\in \cS_{\leq -1}$.
\end{enumerate}
The natural inclusion functor  
\begin{equation}
    \Bsc_{0}^{\cS_{-1}} \xhookrightarrow{\,\,\,\,\,\,} \Bsc
\end{equation}
is a right fibration, since it is a sieve. We can thus consider the $(\oo,1)$-category $\Mfd(\Bsc_{0}^{\cS_{-1}})$ of $\Bsc_{0}^{\cS_{-1}}$-manifolds; these are conically smooth stratified spaces which admit an atlas of basics in $\Bsc_{0}^{\cS_{-1}}$: disjoint unions of points, each finite subset of which is labelled by some element of a poset, together with stratum-preserving injections among them.

Assume now that $\cS_{k-1}$, $\cS_{\leq k-1}$, and $\Bsc^{\cS_{k-1}}_{k}$ have been defined for all $0\leq k \leq d$. We then make the following definitions.
\begin{definition}\label{ch5:: def: sphere sieve}
    Let $\cS_{d}$ be the full $(\oo,1)$-subcategory of $\Mfd(\Bsc_{d}^{\cS_{d-1}})$ on the objects $X\ra P$ such that $X$ is homeomorphic to a $d$-dimensional sphere. 

    We define $\cS_{\leq d}$ to be the full $(\oo,1)$-subcategory of $\Snglr$ whose objects are those of $\cS_{\leq d-1}$ together with those of $\cS_{d}$.
\end{definition}

\begin{definition}
    Let $\Bsc_{d+1}^{\cS_{d}}$ be the full $(\oo,1)$-subcategory of $\Bsc$ on the objects $U= \RR^{i}\times \sC(L)$ such that $L\in \cS_{\leq d}$ and
    \begin{equation}
        i + \big(\mathsf{dim}(L)+1\big) = d+1.
    \end{equation}
\end{definition}

\Cref{ch5:: def: sphere sieve} depends on the map $\Bsc_{d}^{\cS_{d-1}} \longra \Bsc$ being a right fibration. We record it as the following lemma.

\begin{lemma}\label{ch5:: lem: sieve of spheres}
    Let $d \geq 0$. The fully faithful inclusion  
    $\Bsc_{d}^{\cS_{d-1}} \longra \Bsc$ 
    is a right fibration of $(\oo,1)$-categories.
\end{lemma} 
\begin{proof}
To check this, it suffices to show that $\Bsc_{d}^{\cS_{d-1}} \subset \Bsc$ is a sieve: in other words, given a basic $U = \RR^{i}\times \sC(L)$ and an open embedding  
\begin{equation}
    \RR^{i}\times \sC(L) \hookra \RR^{k} \times \sC(\SS^{d-k-1})
\end{equation}
where $\SS^{d-k-1} \in \cS_{\leq d-k-1}$, we must show that $U \in \Bsc_{d}^{\cS_{d-1}}$. Since there is a (non-canonical) homeomorphism
\begin{equation}
    \RR^{k} \times \sC(\SS^{d-k-1}) \cong \RR^{d},
\end{equation}
the codomain is topologically a $d$-manifold. Because $\RR^{i}\times \sC(L)$ admits an open embedding into it, the domain must also be a topological $d$-manifold. For a basic of the form $\RR^{i}\times \sC(L)$, this forces the underlying topological space of $L$ to be homeomorphic to $S^{d-i-1}$. Equivalently, we have
\begin{equation}
L \in \cS_{d-i-1}, \text{ with } \mathsf{Ob}(\cS_{d-i-1}) \subset \mathsf{Ob}(\cS_{\leq d}),  
\end{equation}
ending the proof.
\end{proof}

\begin{observation}
Let us, for the sake of exposition, explicate what these $(\oo,1)$-categories consist of for $d$ in the range $0\leq d \leq 2$.

\begin{enumerate}
    \item[\un{\emph{$d = 0$}}:] The $(\oo,1)$-category $\Bsc_{0}^{\cS_-1}$ contains the single basic $\{\RR^{0}\}$, so an object in the $(\oo,1)$-category $\Mfd(\Bsc_{0}^{\cS_{-1}})$ consists of a finite collection of points stratified by some poset $P$, and a morphism consists of a stratum-preserving injection. An example is indicated in the figure below, where the colors encode different strata and the map is manifest from the relative positions:
    \medskip
    \begin{equation}
        \tikzfig{figures/Stcob/SCO2}{\input{figures/Stcob/SCO2.tikz}} \,\,\,\,\, \xhookrightarrow{\qquad} \,\,\,\,\,\tikzfig{figures/Stcob/SCO1}{\input{figures/Stcob/SCO1.tikz}} 
        \medskip
    \end{equation}
    \item[\un{\emph{$d = 1$}}:]  An object in $\cS_0$ is a stratified 0-sphere $\SS^0 = \big(S^{0}\ra P\big)$ and morphisms are stratum-preserving injections; likewise for $\cS_{\leq 0}$, with the addition of $\emptyset$. An object of $\Bsc_{1}^{\cS_{0}}$ is a basic of one of two types:
    \begin{itemize}
        \item $U = \RR \times \sC(\emptyset)$, where $\emptyset\in \cS_{\leq 0}$:
        \medskip
        \begin{equation}
            U \, = \, \tikzfig{figures/Stcob/SCO3}{\input{figures/Stcob/SCO3.tikz}}\quad ;
            \medskip
        \end{equation}
        \item $V= \sC(\SS^0)$, where $\SS^0$ is a stratified 0-sphere in $\cS_{\leq 0}$.
        \medskip
        \begin{equation}
            V \, = \, \tikzfig{figures/Stcob/SCO4}{\input{figures/Stcob/SCO4.tikz}}\quad ;
            \medskip
        \end{equation}
    \end{itemize}  
    Morphisms are open embeddings, so $\Hom_{\Bsc_{1}^{\cS_0}}(V,U) = \emptyset$ and $\Hom_{\Bsc_{1}^{\cS_0}}(U,V)$ is classified by two copies of $\Emb(\RR,\RR)$ one for each side of the cone point. Having only these two types of basics to construct the objects of $\Mfd(\Bsc_{1}^{\cS_{0}})$ means that all of them will be topological manifolds of dimension 1 with a conically smooth stratification. This category packs a few interesting features: it contains the diffeomorphism group of the circle, for example; it contains a copy of $\bDelta_{\sf int}$ — the full subcategory of inert morphisms\footnote{Consider the full subcategory of $\Mfd(\Bsc_{1}^{\cS_{0}})$ on the non-compact objects stratified by connected components. Every open embedding is then determined, up to a contractible choice, by the poset map which by definition must be inert (or consecutive).}.
    \medskip  
    \begin{equation}
        \tikzfig{figures/Stcob/SCO5}{\input{figures/Stcob/SCO5.tikz}}\,\,\,\,\, \xhookrightarrow{\qquad} \,\,\,\,\, \tikzfig{figures/Stcob/SCO6}{\input{figures/Stcob/SCO6.tikz}}
        \medskip
    \end{equation}
    On the left (the domain of the open embedding) we observe a pair of copies of $\RR$ stratified by three and two $0$-dimensional strata, respectively, and on the right (the codomain of the embedding) we observe another copy of $\RR$, this time stratified by five $0$-dimensional strata, together with a copy of $\SS^1$. 
    \item[\un{\emph{$d = 2$}}:]  An object of $\cS_1$ is a stratified 1-sphere $\SS^1 = \big(S^1 \ra P\big)$ and morphisms are stratum-preserving injections among them; the $(\oo,1)$-category $\cS_{\leq 1}$ contains $\emptyset$, all stratified $0$-spheres $\SS^0$, and all stratified 1-spheres $\SS^1$. Objects in $\Bsc_{2}^{\cS_1}$ are basics of one of the following types:
    \begin{itemize}
        \item $U = \RR^2 \times \sC(\emptyset)$, where $\emptyset \in \cS_{\leq 1}$:
        \medskip
        \begin{equation} 
            \tikzfig{figures/Stcob/SCO7}{\input{figures/Stcob/SCO7.tikz}} \quad ;
            \medskip
        \end{equation}
         \item $V = \RR^1 \times \sC(\SS^0)$, where $\SS^0 \in \cS_{\leq 1}$:
        \medskip
         \begin{equation} 
            \tikzfig{figures/Stcob/SCO8}{\input{figures/Stcob/SCO8.tikz}} \quad ;
            \medskip
        \end{equation}
         \item $W = \sC(\SS^1)$, where $\SS^1 \in \cS_{\leq 1}$:
        \medskip
         \begin{equation}    
            \tikzfig{figures/Stcob/SCO9}{\input{figures/Stcob/SCO9.tikz}} \quad ;
            \medskip
        \end{equation}
        We chose $\SS^1$ stratified by three $0$-dimensional strata, but there is an infinite collection of such singularity types, one for each way of stratifying the $1$-sphere. This rapid expansion in the number of singularity types when passing from $\Bsc_{1}^{\cS_0}$ to $\Bsc_{2}^{\cS_1}$ helps explain the increasing combinatorial complexity of working with conically smooth topological manifolds. 
        
        Objects in the $(\oo,1)$-category $\Mfd(\Bsc_{2}^{\cS_1})$ are thus topological 2-dimensional manifolds endowed with a maximal conically smooth atlas; the morphisms are conically smooth open embeddings. An example of such an object in this $(\oo,1)$-category is:
        \medskip
        \begin{equation}
            \tikzfig{figures/Stcob/SCO10}{\input{figures/Stcob/SCO10.tikz}}
        \end{equation}   
        \medskip
    \end{itemize} 

 \end{enumerate}

\end{observation}

\begin{remark}
    The two $(\oo,1)$-categories $\Bsc_{d}^{\cS_{d-1}}$ and $\Mfd(\Bsc_{d}^{\cS_{d-1}})$ inherit a symmetric monoidal structure given by disjoint union. 
\end{remark}

\begin{remark}
   The preceding discussion of the singularity types encoded by $\Bsc_{d}^{\cS_{d-1}}$ helps clarify the proof of \Cref{ch5:: lem: sieve of spheres}. Namely, it makes manifest that a sufficiently small neighborhood of every point of a basic $U \in \Bsc_{d}^{\cS_{d-1}}$ again lies in $\Bsc_{d}^{\cS_{d-1}}$.
\end{remark}

\begin{definition}
    The $(\oo,1)$-category of topological stratified $d$-manifolds of depth at most $d$ is defined as
    \begin{equation}
        \SMfd_d \, \overset{\sf def}{=} \Mfd(\Bsc_{d}^{\cS_{d-1}}).
    \end{equation}        
    The objects are topological $d$-manifolds endowed with a conically smooth atlas and the morphisms are open embeddings. The $(\oo,1)$-category admits a symmetric monoidal structure given by disjoint union:
    \begin{align}
        \sqcup \maps \SMfd_d \, \times \,  \SMfd_d &\longra \SMfd_d \\
        \big(M \ra P, N \ra Q \big) &\longmapsto (M \sqcup N) \ra (P \sqcup Q).
    \end{align}
\end{definition}

\begin{observation}
By varying the depth of the stratification, we obtain different classes of topological stratified $d$-manifolds. Recall the functor 
\begin{equation}
    \mathsf{depth} \maps \SMfd_d \longra [d].
\end{equation}
A stratified topological $d$-manifold $M$ has depth at most $k$ if for all $x \in M$, the local depth $\mathsf{depth}_x(M) \leq k$. Hence, if $k \leq l$ and $M$ has depth at most $k$, then it also has depth at most $l$.   
\end{observation}

Since we are interested in studying bordisms of stratified $d$-manifolds, we introduce the following broader class of stratified spaces.

\begin{definition}
    Let $\SMfd_{d}^{\partial}$ be the $(\oo,1)$-category which is the pullback 
    \begin{equation}
        \begin{tikzcd}[ampersand replacement=\&]
            {\SMfd_{d}^{\partial}} \&\& {\Snglr^{\partial}} \\
            \\
            {\SMfd_{d}} \&\& \Snglr \quad .
            \arrow[hook, from=1-1, to=1-3]
            \arrow[from=1-1, to=3-1]
            \arrow["\lrcorner"{anchor=center, pos=0.125}, draw=none, from=1-1, to=3-3]
            \arrow["{\mathsf{int}}", from=1-3, to=3-3]
            \arrow[hook, from=3-1, to=3-3]
        \end{tikzcd}
    \end{equation} 
    The $(\oo,1)$-category $\SMfd_d^{\partial}$ contains as objects those conically smooth stratified spaces which admit the structure of a stratified topological $d$-manifold with boundary. The morphisms are open embeddings which respect the boundary. The symmetric monoidal structure is again given by disjoint union.
\end{definition}

\begin{remark}
    It is important to remind the reader that there is a notion of depth and one of interior depth. They are not always in agreement. A smooth manifold with boundary seen as a conically smooth stratified space with $\lag * \rag$-corners has depth 1 and interior depth 0, for example. This highlights the sense in which we think of the boundary structure as an additional structure to that of the underlying stratified space.
\end{remark}

\subsection{Stratified bordisms}

\begin{definition}
Let $d\geq 0$ be an integer, and let $\Sigma_{\sf in}$ and $\Sigma_{\sf out}$ be closed stratified $(d-1)$-manifolds. A \define{stratified bordism} $M\colon \Sigma_{\sf in} \ra \Sigma_{\sf out}$ \define{of dimension} $d$ is a tuple 
\begin{equation}
    \big(M,\varphi_{\sf in},\varphi_{\sf out}\big)
\end{equation}
consisting of the following data:
\begin{enumerate}
    \item A compact stratified topological $d$-manifold with boundary $M$, together with a corner structure $c \maps M \ra \cP(*)$, i.e., a topologically coCartesian map.
    \item Conically smooth open embeddings between stratified $d$-manifolds with boundary:
    \begin{align}
    \varphi_{\sf in} \maps \RR_{\geq 0} \times \Sigma_{\sf in} & \hookra M \\
    \varphi_{\sf out} \maps \RR_{\leq 0} \times \Sigma_{\sf out} & \hookra M 
    \end{align} 
    \item Writing
    \begin{equation}
        \partial_{\sf in}M \,\, \overset{\sf def}{=} \,\, \varphi_{\sf in}\big(\{0\}\times \Sigma_{\sf in}\big), \qquad \partial_{\sf out}M \,\, \overset{\sf def}{=} \,\, \varphi_{\sf out}\big(\{0\}\times \Sigma_{\sf out}\big),
    \end{equation}
    we require that
    \begin{equation}
        \partial M \,=\, \partial_{\sf in}M \sqcup \partial_{\sf out}M.
    \end{equation}
\end{enumerate}    
\end{definition}

As in the smooth case (see \Cref{ch3:: def: iso bordism}), we can define isomorphism classes of stratified bordisms.

\begin{definition}\label{ch5:: def: iso bord str}
    Let 
    \begin{equation}
    (M,\varphi_{\sf in},\varphi_{\sf out}),\qquad (M',\varphi'_{\sf in},\varphi'_{\sf out})
    \end{equation}
    be $d$-dimensional stratified bordisms from $\Sigma_{\sf in}$ to $\Sigma_{\sf out}$.
    An \define{isomorphism of stratified bordisms} is a conically smooth diffeomorphism of topological stratified $d$-manifolds with boundary
    \begin{equation}
    F\colon M \xrightarrow{\;\cong\;} M'
    \end{equation}
    such that
    \begin{align*}
        F\circ \varphi_{\sf in} &= \varphi'_{\sf in},\\
        F\circ \varphi_{\sf out} &= \varphi'_{\sf out}.
    \end{align*}
    We write 
    \begin{equation}
        (M,\varphi_{\sf in},\varphi_{\sf out})\cong (M',\varphi'_{\sf in},\varphi'_{\sf out})
    \end{equation}
    to denote a pair of isomorphic stratified bordisms.
\end{definition}

This definition becomes especially important when we construct the stratified bordism category. Next, we observe that, just as smooth cobordism forms an equivalence relation, so does stratified cobordism. 

\begin{proposition}
    Let $\Sigma_0$ and $\Sigma_1$ be a pair of closed stratified $(d-1)$-manifolds. We say that $\Sigma_0$ and $\Sigma_1$ are \define{stratified cobordant} if there exists a stratified bordism $M\maps \Sigma_0 \ra \Sigma_1$ of dimension $d$. The relation of being stratified cobordant is an equivalence relation.
\end{proposition}

\begin{proof}
    We check reflexivity, symmetry, and transitivity.
    \begin{enumerate}
        \item \textbf{Reflexivity:} Let $\Sigma$ be a closed stratified $(d-1)$-manifold. Put $M = \Sigma\times[0,1]$ with the product stratification and with the evident boundary (corner) structure $c\maps M\ra \cP(*)$ classifying $\partial M = \Sigma\times\{0,1\}$. Let $\varphi_{\sf in}$ and $\varphi_{\sf out}$ be the evident collars, and define
        \begin{equation}
            \partial_{\sf in}M = \Sigma\times\{0\}, \qquad \partial_{\sf out}M = \Sigma\times\{1\}.
        \end{equation}
        Then $(M,\varphi_{\sf in},\varphi_{\sf out})$ is a stratified bordism $\Sigma\ra \Sigma$.

        \item \textbf{Symmetry:} Let $(M,\varphi_{\sf in},\varphi_{\sf out})$ be a stratified bordism $\Sigma_0\ra \Sigma_1$. Let $r\maps \RR\ra \RR$ be reflection $r(t)=-t$. Define
        \[
            \varphi^{\op}_{\sf in}= \varphi_{\sf out}\circ\big(r|_{\RR_{\ge 0}}\times \id_{\Sigma_1}\big),\qquad
            \varphi^{\op}_{\sf out}= \varphi_{\sf in}\circ\big(r|_{\RR_{\le 0}}\times \id_{\Sigma_0}\big).
        \]
        These are again conically smooth open embeddings of stratified $d$-manifolds with boundary, and by construction they satisfy
        \[
            \varphi^{\op}_{\sf in}\big(\{0\}\times\Sigma_1\big)= \partial_{\sf out}M,
            \qquad
            \varphi^{\op}_{\sf out}\big(\{0\}\times\Sigma_0\big)= \partial_{\sf in}M.
        \]
        Hence, $(M,\varphi^{\op}_{\sf in},\varphi^{\op}_{\sf out})$ is a stratified bordism $\Sigma_1\ra \Sigma_0$.

        \item \textbf{Transitivity:} Let $(M_0,\varphi^0_{\sf in},\varphi^0_{\sf out})\maps \Sigma_0\ra \Sigma_1$ and $(M_1,\varphi^1_{\sf in},\varphi^1_{\sf out})\maps \Sigma_1\ra \Sigma_2$ be stratified bordisms of dimension $d$. Using the collars, choose $\epsilon>0$ and form the glued topological space
        \[
            M_{01}= M_0\cup_{\Sigma_1} M_1,
        \]
        where we glue along the identified copies of $\Sigma_1$ in the outgoing boundary of $M_0$ and the incoming boundary of $M_1$. The collars identify neighborhoods of these boundary components with
        \[
            (-\epsilon,0]\times \Sigma_1 \qquad \text{and} \qquad [0,\epsilon)\times \Sigma_1,
        \]
        so the usual collar-gluing construction endows $M_{01}$ with the structure of a stratified $d$-manifold with boundary. Writing $q_i\maps M_i\ra M_{01}$ for the canonical maps, the incoming and outgoing collars are given by
        \[
            \varphi^{01}_{\sf in}= q_0\circ \varphi^0_{\sf in},
            \qquad
            \varphi^{01}_{\sf out}= q_1\circ \varphi^1_{\sf out}.
        \]

        Thus $(M_{01},\varphi^{01}_{\sf in},\varphi^{01}_{\sf out})$ is a stratified bordism $\Sigma_0\ra \Sigma_2$.
    \end{enumerate}
    Therefore stratified cobordism is an equivalence relation.
\end{proof}

\begin{remark}
    The key input in the transitivity argument is the same one used in the smooth unstratified case: collars turn the common boundary component into a product neighborhood, which makes gluing possible and canonical up to isomorphism. The only additional point in the stratified setting is that this product neighborhood must be regarded as a conically smooth stratified chart; once that is observed, the usual collar-gluing argument carries over without essential change.
\end{remark}

The next definition isolates the ordinary bordism category that the interval Segal space constructed later is meant to recover, after passage to the homotopy category. In particular, the choices of objects, morphisms, and composition below are designed so as to make that comparison precise.

\begin{definition}\label{ch5:: def: strat bord ord}
    The category $\cS\mathsf{Cob}_d$ of unoriented (conically smooth) stratified bordisms has:
    \begin{itemize}
        \item \textbf{Objects:} Objects are stratified closed $(d-1)$-manifolds.\footnote{By a closed stratified $d$-manifold we mean a stratified $d$-manifold whose underlying topological $d$-manifold has empty boundary and is compact.} 
        \item \textbf{Morphisms:} A morphism $M\maps \Sigma_{\sf in} \ra \Sigma_{\sf out}$ is an isomorphism class of $d$-dimensional stratified bordisms 
        \begin{equation}
            \left[\big(M, \varphi_{\sf in}, \varphi_{\sf out}\big)\right].
        \end{equation}
        \item \textbf{Composition:} Consider a pair of composable morphisms
        \begin{equation}
            \left[\big(M \maps \Sigma_0 \ra \Sigma_1,\varphi_{\sf in},\varphi_{\sf out}\big)\right], \qquad \left[\big(N \maps \Sigma_1 \ra \Sigma_2,\psi_{\sf in},\psi_{\sf out}\big)\right].
        \end{equation}
        Their composite is defined by gluing along the identified copies of $\Sigma_1$:
        \begin{equation}
            \left[\big(N,\psi_{\sf in},\psi_{\sf out}\big)\right]\circ \left[\big(M,\varphi_{\sf in},\varphi_{\sf out}\big)\right]
            \,\, \overset{\sf def}{=} \,\,
            \left[\big(M\cup_{\Sigma_1} N, q_M\circ \varphi_{\sf in}, q_N\circ \psi_{\sf out}\big)\right],
        \end{equation}
        where the gluing is performed along
        \begin{equation}
            \varphi_{\sf out}|_{\{0\}\times \Sigma_1}\maps \Sigma_1 \hookra M, \qquad \psi_{\sf in}|_{\{0\}\times \Sigma_1}\maps \Sigma_1 \hookra N,
        \end{equation}
        and $q_M$ and $q_N$ denote the canonical maps to the glued space.
        \item \textbf{Identity morphisms:} The identity on an object $\Sigma$ is represented by the cylinder bordism
        \begin{equation}
            \left[\big(\Sigma\times [0,1],\varphi_{\sf in},\varphi_{\sf out}\big)\right]\maps \Sigma \ra \Sigma,
        \end{equation}
        where $\varphi_{\sf in}$ and $\varphi_{\sf out}$ are the evident collars.
    \end{itemize}
\end{definition}

\begin{proposition}
    The above data defines a category. Moreover, $\cS\mathsf{Cob}_d$ is symmetric monoidal under disjoint union.
\end{proposition}
\begin{proof}
    The transitivity argument above shows that composition of composable morphisms again produces a stratified bordism. If
    \begin{equation}
        F\maps (M,\varphi_{\sf in},\varphi_{\sf out}) \xra{\;\cong\;} (M',\varphi'_{\sf in},\varphi'_{\sf out}), \qquad G\maps (N,\psi_{\sf in},\psi_{\sf out}) \xra{\;\cong\;} (N',\psi'_{\sf in},\psi'_{\sf out})
    \end{equation}
    are isomorphisms of composable stratified bordisms, then $F$ and $G$ glue to an isomorphism
    \begin{equation}
        M\cup_{\Sigma_1}N \xra{\;\cong\;} M'\cup_{\Sigma_1}N',
    \end{equation}
    so composition descends to isomorphism classes. The cylinder bordisms act as identities. Associativity of composition and the unit laws hold because collar-gluing is associative up to canonical isomorphism, and therefore strictly after passing to isomorphism classes. The symmetric monoidal structure is given on objects and morphisms by disjoint union.
\end{proof}

\section{An interval Segal space of \texorpdfstring{$d$}{d}-dimensional stratified bordisms}

We finally turn our attention to the construction of an interval Segal space of stratified bordisms. The approach is the same as in the smooth unstratified case: using intervals, together with appropriate conically smooth submersivity conditions, to construct a functor from $\Int$ landing in $\shv(\man)$, which, after post-composition with the concordance resolution, yields the desired stratified bordism $(\oo,1)$-category.

Additionally, we show that this stratified bordism $(\oo,1)$-category can be endowed with a symmetric monoidal structure and that its homotopy category is symmetric monoidal equivalent to \Cref{ch5:: def: strat bord ord}.

\begin{remark}
    We will want to consider stratified $d$-manifolds embedded in a non-stratified ambient space, namely a smooth manifold. In the case of interest this ambient space will be $\RR^{N}$, or more generally a finite-dimensional real vector space $V$. This requires the datum of a conically smooth stratified space $\big(M \ra P\big)$ together with a functor $\varphi \maps \exit(M) \ra \Vect^{\inj}_{|\sB\sO(d)}$ trivializing $\sT_{M}$ and a conically smooth embedding 
    \begin{equation}
        i \maps (M \ra P) \hookra X
    \end{equation}
    where $X$ is a smooth manifold, possibly with boundary, in which case one imposes the expected compatibility conditions along the boundary.
\end{remark}

\begin{definition}\label{ch5:: def: bord_d strat}
    Fix a finite-dimensional real vector space $V$. We define a functor 
    \begin{equation}
    \tilbordstr_{d,V}\maps \man^{\op}\times \Int \longrightarrow \set    
    \end{equation}
    sending a pair $(M,T)$, where $M\in \man^{\op}$ and $T\in \Int$, to the set of bounded closed subspaces 
    \begin{equation}
   \big( W\subseteq M\times \mathsf{cor}(T) \times V \big)
    \end{equation} 
    such that the following three conditions are satisfied:
    \begin{enumerate}
    \item \label{id: st 1} $W$ is endowed with the structure of a conically smooth stratification and a $\SMfd_{d}^{\partial}$-framing. In particular, it is a conically smooth stratified space with a boundary structure
    \begin{equation}
         (W \ra P),\quad (W \xra{c} \cP(*));
    \end{equation} 
    and the embedding must preserve the underlying stratifications, i.e., send the boundary of $W$ to the boundary of $M\times \core(T) \times V$.
    \item \label{id: st 2} The restricted projection 
    \begin{equation}
        \pi\maps \big(W \cap (M\times T^{\mathsf c} \times V) \big) \hookra  M \times T^{\mathsf c} \times V \fib  M \times T^{\mathsf c}
    \end{equation} 
    is a conically smooth submersion.
    \item \label{id: st 3} The restricted projection 
    \begin{equation}
    \pi_W\maps W \hookra M\times \mathsf{cor}(T) \times V \fib M 
    \end{equation}
    is a conically smooth submersion.
    \end{enumerate} 

    \noindent Now we define how $\tilbordstr_{d,V}$ acts on morphisms:
    \begin{itemize}
    \item[$(\man)$] The morphism $(f,\mathsf{id}_{T})\maps (M,T) \longra (N,T)$ maps the closed submanifold 
    \begin{equation}
    W\subseteq M\times \mathsf{cor}(T) \times V
    \end{equation}
    to the pullback along $(f\times \id)$:
    \begin{equation}
        \begin{tikzcd}[ampersand replacement=\&]
            {(f\times \id)^\ast W} \& {W} \\
            \\
            N\times \mathsf{cor}(T) \times V \& M\times \mathsf{cor}(T) \times V
            \arrow[from=1-1, to=1-2]
            \arrow[hook, "(f\times \id)^\ast \pi_{|W}"',from=1-1, to=3-1]
            \arrow["\lrcorner"{anchor=center, pos=0.125}, draw=none, from=1-1, to=3-2]
            \arrow[hook, from=1-2, to=3-2]
            \arrow["f\times \id"', from=3-1, to=3-2]
        \end{tikzcd}
    \end{equation}
    which is a conically smooth stratified $\big(d + \mathsf{dim}(N)\big)$-dimensional submanifold $f^\ast W \subseteq N \times \mathsf{cor}(T) \times V$.
    \item[$(\Int)$] A morphism $(\id_M, T \longra S)$ acts by \emph{forgetting} connected components and restricting the stratified $(d+\mathsf{dim}(M))$-manifold $(W\ra P)$.
    \end{itemize}  
\end{definition}

\begin{remark}
    This definition depends on the existence of a conically smooth embedding for every isomorphism class of compact stratified $d$-manifold with boundary. This is ensured by \Cref{ch4:: prop: neat embedding strat}.  
\end{remark}

\begin{observation}
    For every linear embedding of finite-dimensional real vector spaces $i\maps V \hookra W$, there is an induced natural transformation of bifunctors
    \begin{equation}
        \tilbordstr_{d,V} \xhookrightarrow{\,\,\,\,\,\,\,\,} \tilbordstr_{d,W}
    \end{equation}
    defined by post-composing the ambient embedding with $i$. We therefore define
    \begin{equation}
        \tilbordstr_{d} \,\, \overset{\sf def}{=} \,\, \hoColimit[V \in \vect^{\sf inj}]\,\tilbordstr_{d,V}
    \end{equation}
    as a bifunctor on $\man^{\op}\times \Int$. Applying the concordance resolution then gives
    \begin{equation}
        \bordstr_{d}(I) \,\, \overset{\sf def}{=} \,\, \big|\,\tilbordstr_{d}(\AA^\bullet,I)\,\big|
    \end{equation}
    for each $I\in \Int$. In this way we obtain a functor
    \begin{align}
        \bordstr_{d} \maps \Int &\longra \sset \\
        I &\longmapsto \bordstr_d(I).
    \end{align}
 
\end{observation}

The following is a critical observation rendering the definition of $\tilbordstr_{d,V}$ a sensible one.
\begin{observation}
    An important point for the soundness of the definition of $\tilbordstr_{d,V}$ is that if we pull back a conically smooth embedded submanifold, together with its conically smooth proper submersion over $M$, along a smooth map $f\maps M \longra N$, we again obtain a conically smooth proper submersion over $N$.
    
    From \Cref{ch4:: thm: cbl pullback}, we know that constructible bundles pull back along any conically smooth map. In particular, constructible bundles over a smooth base $K$ amount to conically smooth fiber bundles over $K$. From the stratified analogue of the Ehresmann fibration theorem, namely \Cref{ch4:: cor: stratified ehresmann boundary}, we know that 
    \begin{enumerate}
        \item[$(*)$] a proper conically smooth submersion $X \ra K$ over a smooth manifold $K$ is a conically smooth fiber bundle, 
    \end{enumerate} 
    ensuring stability under pullback for proper conically smooth submersions (over a smooth base $K$).
\end{observation}

\medskip

\begin{remark}
    To make contact with bordisms, let us depict an example and a pair of non-examples of a $0$-simplex in $\tilbordstr_2(\AA^{\bullet}, \sU(I))$, which illustrate how the construction captures objects in a stratified bordism category. The first diagram depicts a cylinder homeomorphic to $S^1\times [0,1]$ stratified in such a way that the tangent spaces to the $1$-dimensional stratum always project positively to the interval $I$, thereby yielding a valid $0$-simplex:
    \begin{equation}
        \tikzfig{figures/Stcob/SCO12}{\input{figures/Stcob/SCO12.tikz}}
    \end{equation}
    \medskip

    \nid Next, we depict again a surface homeomorphic to $S^1\times [0,1]$, but whose 1-dimensional stratum admits two critical points relative to the projection to $I$, thereby failing to satisfy condition (\ref{id: st 2}) — and (\ref{id: st 3}), by extension:
    \begin{equation}
        \tikzfig{figures/Stcob/SCO13}{\input{figures/Stcob/SCO13.tikz}}
    \end{equation} 
    \medskip

    \nid Lastly, we depict the same cylinder with a stratification containing a 0-dimensional stratum. At this point, the map fails to be a submersion for dimensionality reasons, thereby also failing to be a valid 0-simplex:
    \begin{equation}
        \tikzfig{figures/Stcob/SCO14}{\input{figures/Stcob/SCO14.tikz}}.
    \end{equation}
    \medskip

    \nid These examples highlight the fact that the second and third examples do not represent objects in a $2$-dimensional stratified bordism category, whereas the first does. Up to reparametrization, it captures the circle with two $0$-dimensional strata
    \begin{equation}
        \tikzfig{figures/Stcob/SCO15}{\input{figures/Stcob/SCO15.tikz}}.
    \end{equation}
    \medskip

    To get a sense of how bordisms are encoded, let us pick out an example of a $0$-simplex in $\tilbordstr_2(\AA^\bullet,\sU(I,J))$. Over the red-shaded region in the middle there are different singular points relative to the interval projection: the $0$-dimensional stratum in red and the saddle point of the pair-of-pants in the $1$-dimensional stratum.

    \begin{equation}
        \tikzfig{figures/Stcob/SCO16}{\input{figures/Stcob/SCO16.tikz}}
    \end{equation}
    \medskip 
\end{remark}

\begin{observation}
    The structure and intent behind \Cref{ch5:: def: bord_d strat} are manifestly analogous to those of \Cref{ch3:: def: bord_d} and \Cref{ch3:: def: bord_d tangential}. In particular, our path toward the interval-Segal model of the stratified bordism category follows the same sequence of steps: we first establish the analogous locality statements (\Cref{ch3:: prop: naive sheaf con} and \Cref{ch3:: prop: sheaf condition}) for $\tilbordstr_d$, and then prove the interval Segal condition.
\end{observation}

Just as in the unstratified smooth case, passing to the adjoint via the mapping-power adjunction produces an $\oo$-sheaf on $\man$.
\begin{proposition}\label{ch5:: prop: sheaf condition}
    The adjoint functor 
    \begin{equation}
    \tilbordstr_{d,V}(-,-)\maps \Int \longra \Fun(\man^{\op},\sset)
    \end{equation}
    factors through $\shv(\man)$, the category of $\oo$-sheaves on $\man$ relative to the Grothendieck topology of open covers. In other words, for every $T\in \Int$ the functor
    \[
    \tilbordstr_{d,V}(-,T)\maps \man^{\op}\longra \set
    \]
    is a smooth $\oo$-sheaf of (discrete) simplicial sets.
\end{proposition}
\begin{proof}
    The proof is entirely analogous to that of \Cref{ch3:: prop: sheaf condition}. The key observation in the stratified case is that being a conically smooth submersion is also an open condition — ensuring gluing satisfies conditions (\ref{id: st 1})--(\ref{id: st 3}).
\end{proof}

\begin{proposition}[Naive sheaf condition in the interval direction]\label{ch5:: prop: naive sheaf con}
Let $T\in \Int$, and let $\{U,V\}$ be an open cover of $\core(T)\subset [0,1]$. For each $k\geq 0$, the canonical map
\begin{equation}
\label{ch5:: eq: naive sheaf interval levelwise}
\tilbordstr_{d,V}\big(\AA^{k},T\big)\longrightarrow
\tilbordstr_{d,V}\big(\AA^{k},T\big)_{|U}\underset{\tilbordstr_{d,V}\big(\AA^{k},T\big)_{|U\cap V}}{\times}
\tilbordstr_{d,V}\big(\AA^{k},T\big)_{|V}
\end{equation}
is a bijection of sets. Equivalently, the simplicial set $\tilbordstr_{d,V}(\AA^{\bullet},T)$ is the strict pullback
\begin{equation}
\label{ch5:: eq: naive sheaf interval simplicial}
\tilbordstr_{d,V}(\AA^{\bullet},T)\;\cong\;
\tilbordstr_{d,V}(\AA^{\bullet},T)_{|U}\underset{\tilbordstr_{d,V}(\AA^{\bullet},T)_{|U\cap V}}{\times}
\tilbordstr_{d,V}(\AA^{\bullet},T)_{|V},
\end{equation}
where
\begin{equation}
\label{ch5:: eq: bord-restriction-U}
\tilbordstr_{d,V}(\AA^{k},T)_{|U}
=
\Big\{
W\subseteq \AA^{k}\times\big(U\cap \core(T)\big)\times V
\Big\}.
\end{equation}
Here $\tilbordstr_{d,V}(\AA^{k},T)_{|U}$ denotes the set of all $(d+k)$-dimensional stratified submanifolds with boundary
\[
W\subseteq \AA^{k}\times\big(U\cap \core(T)\big)\times V
\]
which satisfy conditions (\ref{id: st 1}), (\ref{id: st 2}), and (\ref{id: st 3}) from \Cref{ch5:: def: bord_d strat}, and the restriction maps are given by intersecting with
\[
\AA^{k}\times\big(U\cap V\cap \core(T)\big)\times V\subseteq 
\AA^{k}\times\big(U\cap\core(T)\big)\times V,
\]
and similarly for $V$.
\end{proposition}

\begin{proof}
    This is the stratified analogue of \Cref{ch3:: prop: naive sheaf con}, and the same argument applies. Fix $k\geq 0$ and write
    \begin{equation}
    Z = \AA^{k}\times \core(T)\times V.
    \end{equation}
    Injectivity of \eqref{ch5:: eq: naive sheaf interval levelwise} is immediate, since an element
    \begin{equation}
    W\subseteq Z
    \end{equation}
    of $\tilbordstr_{d,V}(\AA^{k},T)$ is recovered from its restrictions to
    $\AA^{k}\times(U\cap\core(T))\times V$ and
    $\AA^{k}\times(V\cap\core(T))\times V$.

    For surjectivity, let
    \begin{equation}
    W_U\subseteq \AA^{k}\times\big(U\cap\core(T)\big)\times V,\qquad
    W_V\subseteq \AA^{k}\times\big(V\cap\core(T)\big)\times V
    \end{equation}
    be compatible elements, and define $W := W_U\cup W_V\subseteq Z$. Since the two pieces agree on the overlap, this is a well-defined subset. Exactly as in Chapter 3, all the relevant structure is local on the source. Hence $W$ inherits the structure of a conically smooth stratified manifold with boundary by gluing the local models carried by $W_U$ and $W_V$. The boundary condition in (\ref{id: st 1}) is local, and the submersivity conditions (\ref{id: st 2}) and (\ref{id: st 3}) hold on $W$ because conically smooth submersion is an open condition. Therefore $W$ lies in $\tilbordstr_{d,V}(\AA^{k},T)$ and restricts to $(W_U,W_V)$.

    This proves surjectivity, hence bijectivity of \eqref{ch5:: eq: naive sheaf interval levelwise}. The simplicial-set statement follows levelwise.
\end{proof}

These two ingredients, in their unstratified versions, are the fundamental keys to proving the interval Segal condition in Chapter 3. For that reason, we will now introduce the analogue of $\bord_d$ for stratified $d$-manifolds, using the concordance resolution, and prove that the resulting pre-interval space is in fact an interval Segal space.

\begin{theorem}\label{ch5:: thm: stratbord interval segal}
    Let $d \geq 0$ be a non-negative integer and let $V$ be a finite-dimensional vector space of sufficiently high dimension. The composite functor
    \begin{align}
        \bordstr_{d,V} \maps \Int &\longra \shv(\man)  \longra \sset  \\
         I &\longmapsto \tilbordstr_{d,V}(-,I) \longmapsto \bB\tilbordstr_{d,V}(*,I) = \big|\,\tilbordstr_{d,V}(\AA^\bullet,I)\,\big| 
    \end{align}
    is an interval Segal space.
\end{theorem}

\begin{proof}
    The proof is broken down into two parts: 1) showing that $\bordstr_{d,V}$ sends morphisms in $\cW_{\pi_{0}}$ to weak equivalences of simplicial sets and 2) checking the interval Segal condition, which states — in one of its equivalent forms — that for every $\sU(A,B,C)\in \Int$, the following is a homotopy pullback square
     \begin{equation}\label{ch5:: eq: Segal condition for bordism}
        \begin{tikzcd}[ampersand replacement=\&]
            {\bordstr_{d,V}\big(U(A,B,C)\big)} \&\& {\bordstr_{d,V}\big(U(A,B)\big)} \\
            \\
            {\bordstr_{d,V}\big(U(B,C)\big)} \&\& {\bordstr_{d,V}\big(U(B)\big).}
            \arrow[from=1-1, to=1-3]
            \arrow[from=1-1, to=3-1]
            \arrow["h"{description, pos=0.3}, "\lrcorner"{anchor=center, pos=0.125}, draw=none, from=1-1, to=3-3]
            \arrow[from=1-3, to=3-3]
            \arrow[from=3-1, to=3-3]
        \end{tikzcd}
    \end{equation} 
    For 1), given a morphism 
    \begin{equation}
        T \hookra S
    \end{equation}
    inducing an isomorphism on $\pi_0$, the same diffeomorphism $s$ which interpolates between the two intervals $T$ and $S$ also provides an explicit simplicial homotopy equivalence between $\bordstr_{d,V}(S)$ and $\bordstr_{d,V}(T)$. For additional details see the first part of the proof of \Cref{ch3:: thm: interval segal bordism cat}.

    For 2), the only ingredient of the proof in the unstratified case which is not abstractly categorical and depends on the geometric context is the \emph{naive sheaf condition}, which is easily adapted to the stratified case in \Cref{ch5:: prop: naive sheaf con}. Once again, the details of the second part of the proof of \Cref{ch3:: thm: interval segal bordism cat} are the same.
\end{proof}

\begin{observation}
    We chose to first present this result, together with a complete proof of its unstratified variant, in Chapter 3 precisely because it adapts so readily to more general contexts, such as the present stratified one. In fact, we expect that this method for building higher bordism-type categories should apply in still broader settings, for instance to graph bordisms.
\end{observation}

Through disjoint union of conically smooth stratified $d$-manifolds, we can endow $\bordstr_d$ with a symmetric monoidal structure. This structure is also built in an entirely analogous fashion to that of the unstratified case. We recall the construction here for completion.

\begin{construction}\label{ch5:: cons: Gamma-space}
  For every $m\geq 0$, let $\bordstr_d^{\langle m \rangle}$ denote the following interval Segal space:
\begin{align}
    \bordstr_d^{\langle m \rangle} \maps \Int &\longra \sset \\
    I &\longmapsto \Big(\,[k] \mapsto \left\{M \subset \AA^{k}\times \mathsf{cor}(I)\times \RR^{\oo}\right\}\,\Big)
\end{align} 
where $M$ is given by a disjoint union $M = M_{1} \sqcup \cdots \sqcup M_{m}$ where each $M_{i}\subset \AA^{k} \times \mathsf{cor}(I)\times \RR^{\oo}$ is a $k$-simplex of $\bordstr_d(I)$. Given a pointed map of pointed finite sets, $f\maps \langle m \rangle \longra \langle n \rangle$, we construct a map of interval Segal spaces. This is the natural transformation whose component on $I \in \Int$ is
\begin{equation}
    \begin{tikzcd}
        {\bordstr_d^{\langle m \rangle }(I)} & {M = M_1 \sqcup \cdots \sqcup M_m} \\
        \\
        {\bordstr_d^{\langle n \rangle }(I)} & {(f_*)_I = \Big(\coprod_{f(i) = 1}M_i\Big) \sqcup \cdots \sqcup \Big(\coprod_{f(i) = n}M_i\Big)}
        \arrow[from=1-1, to=3-1]
        \arrow["{{(f_*)_I}}"', maps to, from=1-2, to=3-2]
    \end{tikzcd}
\end{equation}
where we omit the ambient data of the embeddings for simplicity.
\end{construction}

\Cref{ch5:: cons: Gamma-space} can again be proven to induce a symmetric monoidal structure, in this case on $\bordstr_d$.

\begin{theorem}\label{ch5:: thm: symmetric monoidal stratified}
    The assignment 
    \begin{align}
        \bordstr_{d}^{\langle \bullet \rangle} \maps \Gamma &\longra \mathsf{ISSp} \\
        \langle m \rangle & \longmapsto \bordstr_d^{\langle m \rangle}
    \end{align}
    is functorial and endows the interval Segal space $\bordstr_d$ of unoriented $d$-dimensional stratified bordisms with a symmetric monoidal structure.
\end{theorem}

\begin{proof}
    The details of the proof follow those of \Cref{ch3:: thm: symmetric monoidal}. Both proofs rely essentially on the contractibility of the space of conically smooth embeddings of a conically smooth stratified space $M$ into $\RR^{\oo}$ (see \Cref{ch4:: prop: contractibility of embeddings}). 
\end{proof}

\begin{lemma}\label{ch5:: lem: normalized neat embedding}
    Let 
    \begin{equation}
        \big(M,\varphi_{\sf in},\varphi_{\sf out}\big)\maps \Sigma_{\sf in} \ra \Sigma_{\sf out}
    \end{equation}
    be a compact stratified bordism. After possibly shrinking the given collars, there exists a neat embedding
    \begin{equation}
        E\maps M \hookra [0,1]\times \RR^{m}
    \end{equation}
    for some $m\geq 0$, such that
    \begin{equation}
        E(\partial_{\sf in}M)\subseteq \{0\}\times \RR^{m}, \qquad E(\partial_{\sf out}M)\subseteq \{1\}\times \RR^{m},
    \end{equation}
    and, on smaller collar neighborhoods,
    \begin{equation}
        E\circ \varphi_{\sf in}(t,x)=\big(t,e_{\sf in}(x)\big), \qquad E\circ \varphi_{\sf out}(s,x)=\big(1+s,e_{\sf out}(x)\big),
    \end{equation}
    for suitable conically smooth embeddings $e_{\sf in}$ and $e_{\sf out}$ into $\RR^{m}$.
\end{lemma}
\begin{proof}
    This is a minor adaptation of \Cref{ch4:: prop: neat embedding strat}. By compactness, after shrinking collars, we may assume that for some $\epsilon>0$ the images of
    \begin{equation}
        \varphi_{\sf in}\big([0,2\epsilon)\times \Sigma_{\sf in}\big), \qquad \varphi_{\sf out}\big((-2\epsilon,0]\times \Sigma_{\sf out}\big)
    \end{equation}
    are disjoint. Choose a proper conically smooth embedding $e\maps M \hookra \RR^{m}$ by \Cref{ch4:: prop: weak Whit strat}. Using \Cref{lem:partitions of unity}, we construct a conically smooth function
    \begin{equation}
        h\maps M \longra [0,1]
    \end{equation}
    such that
    \begin{equation}
        h\circ \varphi_{\sf in}(t,x)=t, \qquad h\circ \varphi_{\sf out}(s,x)=1+s
    \end{equation}
    on the smaller collars $[0,\epsilon)\times \Sigma_{\sf in}$ and $(-\epsilon,0]\times \Sigma_{\sf out}$, and such that $h(M\setminus \partial M)\subset (0,1)$.

    As in the proof of \Cref{ch4:: prop: neat embedding strat}, after modifying $e$ by a conically smooth bump function near the boundary, we may assume that $e$ is constant along the collar direction on these smaller collar neighborhoods. Therefore the product
    \begin{equation}
        E=(h,e)\maps M \longra [0,1]\times \RR^{m}
    \end{equation}
    is a conically smooth embedding, and on the chosen collar neighborhoods it takes the form
    \begin{equation}
        E\circ \varphi_{\sf in}(t,x)=\big(t,e_{\sf in}(x)\big), \qquad E\circ \varphi_{\sf out}(s,x)=\big(1+s,e_{\sf out}(x)\big).
    \end{equation}
    In particular, $E$ is neat and carries the incoming and outgoing boundary components to $\{0\}\times \RR^{m}$ and $\{1\}\times \RR^{m}$, respectively.
\end{proof}

\section{Homotopy category as the stratified bordism 1-category}

    Before studying $(\oo,1)$-categories we studied ordinary 1-categories, and instances of bordism categories date back at least to Atiyah's 1988 paper \cite{At88}. To confirm that the interval Segal space $\bordstr_d$ really lifts the ordinary stratified bordism category, we should compute its symmetric monoidal homotopy $1$-category and compare it with the category $\cS\mathsf{Cob}_d$ constructed earlier in this chapter.

\begin{remark}\label{ch5:: rmk: homotopy ehresmann strat}
    In the next theorem we work with the homotopy category of the Segal space
    \begin{equation}
        X \,\, \overset{\sf def}{=} \,\, \mathsf{RKan}_{\gamma}(\bordstr_d).
    \end{equation}
    Its objects are therefore coherent families of compact conically smooth stratified $(d-1)$-manifolds. Exactly as in the smooth case, each such family determines a unique isomorphism class of stratified manifold, obtained by taking the fiber over any point in the relevant interval. The point is that the projection to the interval is a proper conically smooth submersion, so \Cref{ch4:: cor: stratified ehresmann boundary} identifies all of its fibers up to conically smooth stratified diffeomorphism. The same observation applies to 1-simplices, where one obtains a stratified bordism by restricting to the closed subinterval between the chosen incoming and outgoing midpoints.
\end{remark}

\begin{theorem}\label{ch5:: thm: homotopy category stratified bord}
    There is an equivalence of symmetric monoidal 1-categories
    \begin{equation}
        \ho(\bordstr_d) \,\, \simeq \,\, \cS\mathsf{Cob}_d
    \end{equation}
    between the homotopy category of $\bordstr_d$ and the ordinary category of unoriented stratified bordisms.
\end{theorem}

\begin{proof}
    The proof follows the same pattern as the corresponding homotopy-category theorem in Chapter 3, with the smooth geometric inputs replaced by their stratified analogues from Chapter 4.

    We define a functor
    \begin{equation}
        h\maps \ho(\bordstr_d) \longra \cS\mathsf{Cob}_d.
    \end{equation}
    On objects, if $(\Sigma_I)_{I\in \cJ_0}$ is a coherent family representing an object of $\ho(\bordstr_d)$, then we set
    \begin{equation}
        h\big((\Sigma_I)_{I\in \cJ_0}\big) \,\, \overset{\sf def}{=} \,\, \pi_I^{-1}\big(\mathsf{mid}(I)\big)
    \end{equation}
    for any choice of interval $I\in \cJ_0$. By \Cref{ch5:: rmk: homotopy ehresmann strat}, this is well-defined up to conically smooth stratified diffeomorphism. On morphisms, represented by coherent families of embedded stratified bordisms over some $\sU(A,B)\in \cJ_1$, we restrict to the closed subinterval between the chosen midpoints and take the resulting isomorphism class in $\cS\mathsf{Cob}_d$. The same stratified Ehresmann argument shows that this is independent of all auxiliary choices.

    Essential surjectivity and fullness are proved exactly as in Chapter 3, except that one uses \Cref{ch5:: lem: normalized neat embedding} in place of the smooth neat-embedding statement. Indeed, every compact stratified bordism admits a normalized neat embedding into $[0,1]\times \RR^m$, and after choosing for each $\sU(A,B)\in \cJ_1$ a diffeomorphism
    \begin{equation}
        [0,1] \xra{\,\,\,\simeq\,\,\,} \core\big(\sU(A,B)\big)
    \end{equation}
    carrying $0$ and $1$ to the prescribed incoming and outgoing midpoints, these embeddings assemble into the required coherent family. The contractibility of the relevant space of interval reparameterizations is the same as in the smooth case.

    For faithfulness, suppose two morphisms in $\ho(\bordstr_d)$ have the same image in $\cS\mathsf{Cob}_d$. Then the extracted stratified bordisms are isomorphic, so after choosing a representative compact stratified bordism
    \begin{equation}
        \big(M,\varphi_{\sf in},\varphi_{\sf out}\big)\maps \Sigma \ra \Sigma'
    \end{equation}
    for that common isomorphism class, both morphisms may be represented by normalized neat embeddings of $M$ into $[0,1] \times \RR^{\oo}$. After fixing the same collar-compatible height function $h \maps M \ra [0,1]$ as in \Cref{ch5:: lem: normalized neat embedding}, the difference between the two representatives lies only in the Euclidean embedding component. By the weak contractibility of the space of conically smooth embeddings into $\RR^{\oo}$ (\Cref{ch4:: prop: contractibility of embeddings}), these two embedding components lie in the same path component, hence determine concordant $1$-simplices in $\bordstr_{d}$. Therefore they define the same morphism in $\ho(\bordstr_{d})$, so $h$ is faithful.

    It follows that $h$ is an equivalence of categories. Compatibility with the symmetric monoidal structures follows from the fact that both sides are defined by disjoint union, exactly as in the smooth case discussed in Chapter 3.
\end{proof}

\begin{remark}
    Before moving on, it is worth recording a possible further direction suggested by the preceding construction. Even in the stratified setting developed in this chapter, the parameter direction still ranges over smooth manifolds. More precisely, the input is again a bifunctor
    \begin{equation}
        \Int \times \man^{\op} \longra \set,
    \end{equation}
    which is then transposed and post-composed with a concordance-resolution procedure in order to obtain a space-valued interval Segal object. In this way, the interval variable is responsible for bordism composition, while the second variable is used to encode the spatial direction. The stratified character of the theory is therefore encoded in the fibers, not in the parameter category itself.

    This suggests a different possible generalization. One could try to replace the smooth parameter direction by a genuinely stratified one, and begin instead with a bifunctor
    \begin{equation}
        \Int \times \strat^{\op}\longra \set,
    \end{equation}
    in which the role played above by conically smooth submersions over a smooth base is replaced by the requirement of being a constructible bundle over a stratified base. Transposing such a bifunctor would produce
    \begin{equation}
        \Int \longra \Fun(\strat^{\op},\set),
    \end{equation}
    and one may then ask whether applying the topologizing diagram levelwise yields an assignment
    \begin{equation}
        \Int \longra \Fun(\strat^{\op},\sset)
    \end{equation}
    which lands in striation sheaves. Since striation sheaves form a model for $(\oo,1)$-categories \cite{AFR1}, this would mean that each interval object determines an $(\oo,1)$-category rather than merely a space. If, moreover, the resulting assignment satisfied the interval Segal condition internally, one would obtain an interval Segal object valued in $(\oo,1)$-categories, hence a double $(\oo,1)$-category, or equivalently a model for an $(\oo,2)$-category. This point of view has begun making appearances via symmetric monoidal double higher bordism categories in the recent work of Yang Yang \cite{YY22} and more recently in work of Krannich--Kupers \cite{Ku22}.

    From this perspective, the natural ambient candidate for the vertical direction is provided by constructible bundles, hence by $\Bun$. The point is not that one must decide in advance which particular subcategory of $\Bun$ is the correct one. On the contrary, that flexibility seems to be part of the picture: the work of Ayala--Francis--Rozenblyum suggests that $\Bun$ is itself a rich ambient object, with natural specializations to embeddings, refinements, closed morphisms, and related classes of stratified maps \cite{AFR1,AFR2}. If such a program can be carried out, then the horizontal direction would still encode bordism composition, while the vertical direction would retain geometric operations that are invisible once one passes directly to spaces.
\end{remark}

\section{Tangential structures}

Tangential structures on conically smooth stratified spaces form a rich and still largely unexplored theory. In this thesis we do not attempt to develop such a theory in full generality. Instead, the goal of the present section is more modest and more targeted: we want to isolate a flexible piece of \emph{local} data on stratified manifolds that plays, for defect theories, the same organizational role that ordinary tangential structures (orientation, framing, $G$-structure, \dots) play in the classical bordism category. In particular, we want a framework in which strata of different depths can carry \emph{labels} (or ``colors'') in a way that satisfies a kind of descent relative to the stratification. 

The guiding principle is that, for stratified spaces, local data is controlled by \emph{basics} — especially their links. This suggests that defect label data should be specified locally on the category of basics and
transported functorially along stratified embeddings between them. Concretely, we model a discrete
system of defect labels as a right fibration over the $(\infty,1)$-category of basics, and we focus on
the case of \emph{discrete} fibers. Via straightening, such an object can be interpreted as a
set-valued presheaf on the category of basics: it assigns to each local model a set of admissible
labels, together with restriction maps along morphisms of basics. From this point of view, a labeling of
a stratified manifold is naturally expressed as a functor out of its enter-path category, and the
functoriality built into the presheaf guarantees compatibility with refinement of stratifications and with gluing.

Besides the last section of Chapter 4, which discusses tangential structures on conically smooth stratified spaces, we recommend the recent survey by Lukas M\"uller \cite{Mu25}, as well as Ayala--Francis--Rozenblyum's original article \cite{AFR2}.

\subsection{Discrete defect data}

\begin{definition}
    We call a right fibration 
    \begin{equation}
        \Bsc_{d}^{\cS_{d-1},\DD} \longra \Bsc_{d}^{\cS_{d-1}}
    \end{equation}
    with discrete fibers a $d$\define{-dimensional defect datum}.
\end{definition}

\nid To understand this definition, we will unpack it in increasing degrees of complexity. We begin with a simple low-dimensional example.

\begin{observation}
    Recall that one of the significant features of conically smooth stratified spaces is that the obvious inclusion
    \begin{equation}
        \Aut\big(\sC(X)\big) \, \xla{\,\,\,\simeq\,\,\,} \, \Aut(X)
    \end{equation}
    is a homotopy equivalence for every compact conically smooth stratified space. An important consequence of this result is a fully homotopical characterization of the spaces of morphisms in a category of basics:
    \begin{equation}
        \Strat^{\sf emb}\big(\RR^{i} \times \sC(Z),\RR^{i} \times \sC(Z)\big)\,\, \simeq \,\, \sO(i) \times \Aut\big(\sC(Z)\big).
    \end{equation}
\end{observation}

\begin{example}\label{ch5:: ex: codim 1 defect}
Since we are considering right fibrations with discrete fibers, via straightening, we can also interpret these as $\set$-valued presheaves on $\Bsc_{d}^{\cS_{d-1}}$. 

We would like to articulate the sense in which a $k$-dimensional submanifold inside a bigger $d$-dimensional manifold might carry an additional notion of \emph{color} or \emph{defect label} — this being a key feature in exploring how different submanifolds may themselves carry different theories. Let us define
\begin{equation}
    \cF\maps \big(\Bsc_{d}^{\cS_{d-1}}\big)^\op \longra \set
\end{equation}
by the rule
\begin{equation}
    \cF(U)=
    \begin{cases}
        A, & U \cong \RR^{k}\times \sC(S^{d-k-1}),\\
        *, & U \cong \RR^{d},\\
        \emptyset, & \text{otherwise,}
    \end{cases}
\end{equation}
where, as before, the link $S^{d-k-1}$ is purposely left unstratified. Define the associated discrete fibration as $\sD_{k \subset d}^{A} \longra \Bsc_{d}^{\cS_{d-1}}$.

Let $(M \ra P)$ be a stratified topological $d$-manifold. A $\sD_{k\subset d}^{A}$-framing amounts to a functor 
\begin{equation}
g \maps \entr(M) \longra \sD^{A}_{k\subset d}
\end{equation} 
such that the above diagram of $(\oo,1)$-categories commutes:
\begin{equation}
    \begin{tikzcd}[ampersand replacement=\&]
        \&\& {\sD^{A}_{k \subset d}} \\
        {\entr(M)} \&\& {\Bsc_d^{\cS_{d-1}}}\,\, .
        \arrow[from=1-3, to=2-3]
        \arrow["g", from=2-1, to=1-3]
        \arrow[from=2-1, to=2-3]
    \end{tikzcd}
\end{equation}
In particular, such a framing can exist only when every point of $M$ is locally modeled either on $\RR^d$ or on $\RR^{k}\times \sC(S^{d-k-1})$. Let $\gamma$ be a morphism in the enter-path category from the stratum $M_p$ to the stratum $M_q$, where the strata are labelled by connected components. If both $M_p$ and $M_q$ are smooth $d$-dimensional strata, then $g(\gamma)$ is determined; but if $M_q$ is a connected $k$-dimensional stratum, locally modeled on $\RR^{k}\times \sC(S^{d-k-1})$, then 
\begin{equation}
    g(\gamma) \maps * \longra A,
\end{equation}  
which is equivalent to a choice of $a\in A$ for every connected $k$-dimensional defect stratum. 
\end{example}

\begin{example}
    Let us now consider the following $2$-dimensional discrete defect datum.
    \begin{equation}
        \cF_{\DD} \maps (\Bsc_{2}^{\cS_1})^{\op} \longra \set,
    \end{equation}
    which, at the level of objects, maps:
    \begin{itemize}
        \item $\RR^{2}$ is assigned the set $D_2$.
        \item $\RR\times \sC(S^0)$ is assigned the set $D_1$.
        \item For each stratified topological $1$-sphere $\SS^1$, the basic $\sC(\SS^1)$ is assigned a set $D_0(\SS^1)$. This set should naturally carry the action of $\Aut(\SS^1)$; when $\SS^1$ is not the unstratified circle, this amounts to the cyclic permutations of the connected components. We write
        \begin{equation}
            D_0 \, \overset{\sf def}{=} \, \coprod_{[\SS^1]\in \pi_0(\cS_1)} D_0(\SS^1) 
        \end{equation}
    \end{itemize}
    At the level of morphisms, it acts as follows:
    \begin{itemize}
        \item There are two maps $\RR^2 \ra \RR\times \sC(S^0)$ which we will map to 
        \begin{equation}
            s,t \maps D_1 \longra D_2
        \end{equation}
        which we can equivalently conceive of as a single map 
        \begin{equation}
            (s,t) \maps D_1 \longra D_2 \times D_2.
        \end{equation} 
        \item For every isomorphism class of stratified topological 1-spheres $[\SS^1]\in \pi_0(\cS_1)$, there is a map 
        \begin{equation}
            D_0(\SS^1) \longra D_1,
        \end{equation} 
        for each $0$-dimensional stratum in $\SS^1$, and a map 
        \begin{equation}
            D_0(\SS^1) \longra D_2
        \end{equation}
        for each $1$-dimensional stratum in $\SS^1$. We depict these maps using the following commutative diagram induced by the entrance-paths of $\sC(\SS^1)$ with $\SS^1$ stratified by three $0$-dimensional strata.
        \begin{equation}    
            \tikzfig{figures/Stcob/SCO9}{\input{figures/Stcob/SCO9.tikz}} \qquad \xrightarrow{\,\, \mathsf{enter-paths}\,\,} \qquad \tikzfig{figures/Stcob/SCO17}{\input{figures/Stcob/SCO17.tikz}} \qquad \xrightarrow{\,\, \cF_{\DD}\,\,} \qquad \tikzfig{figures/Stcob/SCO18}{\input{figures/Stcob/SCO18.tikz}}
            \medskip
        \end{equation}
        These maps organize to a single map 
        \begin{equation}
            D_0 \longra \coprod_{[\SS^1]\in \pi_0(\cS_1)} \left(D_1 \times_{D_2} D_1 \times_{D_2} \cdots \times_{D_2} D_1 \right)
        \end{equation}
        akin to the ``fold map'' defined in \cite[Equation~2.27]{CRS19}. We justify the existence of the fold map by noting that it is given by the iterated product of the following dashed maps
        \begin{equation}
            \begin{tikzcd}[ampersand replacement=\&]
                {D_0} \\
                \& {D_1\times_{D_2}D_1} \&\& {D_1} \\
                \\
                \& {D_1} \&\& {D_2}\, .
                \arrow[dashed, from=1-1, to=2-2]
                \arrow[curve={height=-12pt}, from=1-1, to=2-4]
                \arrow[curve={height=12pt}, from=1-1, to=4-2]
                \arrow[from=2-2, to=2-4]
                \arrow[from=2-2, to=4-2]
                \arrow["s", from=2-4, to=4-4]
                \arrow["t"', from=4-2, to=4-4]
            \end{tikzcd}
        \end{equation}
    \end{itemize} 

    This presentation of the data involved in defining a discrete right fibration over $\Bsc_{2}^{\cS_{1}}$ is very much the unoriented analogue of the Carqueville--Runkel--Schaumann notion of $2$-dimensional defect data. In fact, one can say more: if each stratum of the basic was further endowed with an orientation, and we demanded all the structure functions to be equivariant with respect to this additional structure, we would recover their definition verbatim. It is with that goal in mind that we initiate our study of orientations for conically smooth stratified spaces.
\end{example}

\subsection{Orientations}

We end this section by sketching the expected correspondence between the $2$-dimensional defect data outlined in \cite[Definition~2.4]{CRS19} and $2$-dimensional discrete right fibrations over a suitably oriented version of $\Bsc_2^{\cS_1}$. The point is not to prove a theorem in a fully specified comparison category, but rather to explain why the right-fibration viewpoint appears to recover exactly the same local ingredients.

\begin{construction}
We are especially interested in studying the case where each stratum of our stratified $d$-manifold is equipped with an orientation. We demand no further compatibilities between these orientations. 

\medskip

\nid We construct these objects by using an appropriate notion of $\cB$-framing on the constructible tangent bundle of a stratified $d$-manifold. Let $\sB\sS\sO$ be the following space 
\begin{equation}
    \sB\sS\sO \, \overset{\sf def}{=} \, \bigcup_{d \geq 0} \sB\sS\sO(d)
\end{equation}
where $\sB\sS\sO(d)$ is the Grassmannian of $d$-planes in $\RR^{\oo}$ equipped with an orientation. The map which, for every $d\geq 0$, forgets the orientation 
\begin{equation}
    \sB\sS\sO(d) \longra \sB\sO(d)
\end{equation}
induces a map 
\begin{equation}
    \sB\sS\sO \longra \sB\sO.
\end{equation}
This is the stable tangential structure of orientations.
\medskip

\nid Recall that lifts
\begin{equation}
    \begin{tikzcd}[ampersand replacement=\&]
        \&\& {\sB\sS\sO} \\
        M \&\& {\sB\sO}
        \arrow[from=1-3, to=2-3]
        \arrow[from=2-1, to=1-3]
        \arrow["\tau"', from=2-1, to=2-3]
    \end{tikzcd}
\end{equation}
correspond to a choice of orientation for the vector bundle over $M$ which $\tau$ classifies. We can therefore package orientations into a tangential structure $\Vect^{\inj, \BSO}$, which we define as a modification of the original complete Segal space:

\smallskip

\begin{equation}
    [p] \longmapsto \Big| \Big\{V_0 \xlonghookra{f_1} V_1 \xlonghookra{f_2}\ldots \xlonghookra{f_p} V_p \text{ over } \AA^{\bullet} \Big\} \Big|
\end{equation}
\smallskip

\nid where each $V_{i} \ra \AA^{\bullet}$ is equipped with a map $\AA^{\bullet} \longra \BSO$ endowing it with an orientation. \emph{Crucially}, the only requirement is that fiberwise isomorphisms are orientation-preserving. The maximal $\oo$-subgroupoid of this $(\oo,1)$-category can be identified with  
\begin{equation}
    \bigsqcup_{k \geq 0} \BSO(k).
\end{equation}

\medskip

\nid A slightly more \emph{naive}, but useful, perspective on this $(\oo,1)$-category can be seen through its homotopy category. Its set of objects contains all finite-dimensional real vector spaces $V$, together with an $\sS\sO(k)$-torsor $X_k$ and an equivariant map $X_k \longra \mathsf{Fr}(V)$ to the space of frames, which amounts to a choice of orientation for $V$. Morphisms are injections of linear maps between them. That this is a complete Segal space follows arguments similar to those used to prove $\Vect^{\inj}$ is a complete Segal space.

\medskip

To bound the possible dimensions of the spaces we are working with, consider the pullback $(\oo,1)$-category
\begin{equation}
    \begin{tikzcd}[ampersand replacement=\&]
        {\Vect^{\inj,\BSO}_{\leq d}} \&\& {\Vect^{\inj,\BSO}} \\
        \\
        {\Vect^{\inj}_{\leq d}} \&\& {\Vect^{\inj}} \, .
        \arrow[from=1-1, to=1-3]
        \arrow[from=1-1, to=3-1]
        \arrow[from=1-3, to=3-3]
        \arrow[from=3-1, to=3-3]
    \end{tikzcd}
\end{equation}
which truncates the space of objects by cutting off the dimension of the vector bundles at $d$.

\nid To endow a conically smooth stratified space $M$ with a $\Vect^{\inj, \BSO}_{\leq d}$-framing is to endow every $k$-dimensional stratum $M_k$ with a lift $M_k \longra \sB\sS\sO(k)$ 

\begin{equation}
    \begin{tikzcd}[ampersand replacement=\&]
        \&\& {\sB\sS\sO(k)} \\
        M_{k} \&\& {\sB\sO(k)} \, .
        \arrow[from=1-3, to=2-3]
        \arrow[from=2-1, to=1-3]
        \arrow["\sT_{M_{k}}"', from=2-1, to=2-3]
    \end{tikzcd}
\end{equation}

\smallskip

\nid In particular, every connected stratum $M_k$ is endowed with an orientation and no further compatibility data. 

\nid Having stipulated the form of framing we are interested in, we want to consider stratified spaces equipped with such framings and constructible families thereof. To do that, we make use of the technology from Chapter 4. Let $\cB = \Vect^{\inj, \BSO}_{\leq d}$, we define the $(\oo,1)$-category of oriented stratified topological $d$-manifold basics as the following pullback $(\oo,1)$-category

\begin{equation}
    \begin{tikzcd}[ampersand replacement=\&]
        {\Bsc_{d}^{\cS_{d-1}, \mathsf{or}}} \&\& {\Bsc^{\cS_{d-1}}_{d}} \\
        \\
        {\Bun^{\cB,\emb}} \& {\Bun^{\emb}} \& {\Snglr = \Strat^{\emb}} 
        \arrow[from=1-1, to=1-3]
        \arrow[from=1-1, to=3-1]
        \arrow[hook, from=1-3, to=3-3]
        \arrow[from=3-1, to=3-2]
        \arrow["\simeq", from=3-2, to=3-3]
    \end{tikzcd}
\end{equation}
\end{construction}

\begin{observation}[A sketch of the oriented $2$-dimensional dictionary]
    Let \[
        \Bsc_{2}^{\cS_{1},\mathsf{or},\DD} \longra \Bsc_{2}^{\cS_{1}, \mathsf{or}}
    \]
    be a $2$-dimensional oriented discrete defect datum, and let
    \begin{equation}
        \cF \maps \big(\Bsc_{2}^{\cS_{1},\mathsf{or}}\big)^{\op} \longra \set
    \end{equation}
    be its straightening. Then one expects $\cF$ to determine, and to be determined by, defect data of the form
    \begin{equation}
        \DD = \big(D_0,D_1,D_2;f_0,f_1,f_2\big)
    \end{equation}
    in the sense of \cite{CRS19}.
\end{observation}

\medskip

\nid We will not attempt to formulate this as a theorem in a fully specified comparison category. Instead, we record the dictionary one reads off from the presheaf $\cF$.

\medskip

\nid First, the values of $\cF$ on the oriented local models supply the sets of labels:
\begin{enumerate}
    \item The value
    \begin{equation}
        D_2 \, \overset{\sf def}{=} \, \cF(\RR^{2})
    \end{equation}
    is the set of labels for oriented $2$-strata.
    \item The value
    \begin{equation}
        D_1 \, \overset{\sf def}{=} \, \cF\big(\RR\times \sC(S^0)\big)
    \end{equation}
    is the set of labels for oriented $1$-strata.
    \item For each oriented stratified $1$-sphere $\SS^1$, the value
    \begin{equation}
        D_0(\SS^1) \, \overset{\sf def}{=} \, \cF\big(\sC(\SS^1)\big)
    \end{equation}
    supplies the labels for oriented $0$-strata of singularity type $\SS^1$. Taking the coproduct over isomorphism classes gives
    \begin{equation}
        D_0 \, \overset{\sf def}{=} \, \coprod_{[\SS^1]} D_0(\SS^1),
    \end{equation}
    where the coproduct ranges over isomorphism classes of oriented stratified $1$-spheres.
\end{enumerate}

\nid Next, the morphisms in $\Bsc_{2}^{\cS_{1},\mathsf{or}}$ induce the familiar incidence maps. The two embeddings of an oriented codimension-$1$ basic into an oriented surface basic determine source and target maps
\begin{equation}
    s,t \maps D_1 \longra D_2,
\end{equation}
while, for each oriented stratified circle $\SS^1$, the entrance-paths in $\sC(\SS^1)$ induce maps from $D_0(\SS^1)$ to the adjacent copies of $D_1$ and $D_2$. These organize into the same kind of fold-type compatibility data that appears in \cite[Definition~2.4]{CRS19}.

\medskip

\nid Conversely, given $2$-dimensional oriented defect data in the sense of Carqueville--Runkel--Schaumann, one expects to reconstruct a set-valued presheaf on $\Bsc_{2}^{\cS_{1},\mathsf{or}}$ by assigning these sets to the corresponding oriented basics and using the structure maps to define restriction along morphisms. The point is that, once orientations are built into the local models, the ambiguity present in the unoriented codimension-$1$ case disappears: the two adjacent $2$-strata are no longer merely an unordered pair, but carry a source and a target.

\medskip

\nid In this way, the right-fibration viewpoint appears to recover precisely the same local combinatorics as the defect-data formalism of \cite{CRS19}. Making this into a genuine equivalence statement would require specifying an appropriate comparison category of oriented defect data and checking that the two constructions are inverse. We do not attempt that here.

\begin{observation}
We note that this discussion admits an immediate generalization to any other kind of stable tangential structure 
\begin{equation}
    \sB G \longra \BO.
\end{equation}
Likewise, the oriented $2$-dimensional dictionary sketched above should extend, at least formally, to higher dimensions: one expects the same right-fibration viewpoint to produce the relevant defect data by reading off labels from higher-dimensional oriented basics and the corresponding incidence maps from their entrance-path combinatorics. The obstruction is not so much conceptual as combinatorial. Already in dimension $2$, the bookkeeping becomes nontrivial; in higher dimensions, the resulting local configurations rapidly become harder to write down, visualize, and parse.

\medskip

\nid If we replace $\Vect^{\inj,\BSO}_{\le d}$ with $\Vect^{\inj,\sB G}_{\le d}$, then we obtain an $(\oo,1)$-category $\SMfd_{d}^{\sB G}$. If $G=*$, then $\SMfd_d^{\sB *}$ is the $(\oo,1)$-category of stratified topological $d$-manifolds whose $k$-strata admit a $k$-framing. 
\end{observation}

\begin{definition}
    Let $\sB G \longra \BO$ be a stable tangential structure. We define the $(\oo,1)$-category of $d$-dimensional $\sB G$-framings to be the pullback
    \begin{equation}
        \begin{tikzcd}[ampersand replacement=\&]
            {\Vect^{\inj,\sB G}_{\leq d}} \&\& {\Vect^{\inj, \sB G}} \\
            \\
            {\Vect^{\inj}_{\leq d}} \&\& {\Vect^{\inj}}.
            \arrow[from=1-1, to=1-3]
            \arrow[from=1-1, to=3-1]
            \arrow[from=1-3, to=3-3]
            \arrow[from=3-1, to=3-3]
        \end{tikzcd}
    \end{equation}
    Denoting this pullback $(\oo,1)$-category by $\cB$, we construct the $(\oo,1)$-category of stratified $d$-dimensional $\sB G$-basics to be the pullback $(\oo,1)$-category 
    \begin{equation}
        \begin{tikzcd}[ampersand replacement=\&]
            {\Bsc_{d}^{\cS_{d-1}, \sB G}} \&\& {\Bsc^{\cS_{d-1}}_{d}} \\
            \\
            {\Bun^{\cB,\emb}} \& {\Bun^{\emb}} \& {\Snglr = \Strat^{\emb}} \\
            {.}
            \arrow[from=1-1, to=1-3]
            \arrow[from=1-1, to=3-1]
            \arrow[hook, from=1-3, to=3-3]
            \arrow[from=3-1, to=3-2]
            \arrow["\simeq", from=3-2, to=3-3]
        \end{tikzcd}
    \end{equation}
\end{definition}

\section{Stratified topological quantum field theories}

\begin{definition}
    Let $\cC$ be a symmetric monoidal $(\oo,1)$-category. A \textbf{$d$-dimensional unoriented stratified topological quantum field theory} with values in $\cC$ is a symmetric monoidal functor 
    \begin{equation}
        \sZ \maps \bordstr_{d} \longra \cC.
    \end{equation}
\end{definition}

\begin{observation}
    There is an evident map of interval Segal spaces 
    \begin{equation}
        \bord_d \xhookrightarrow{\,\,\,\,\,\,\,\,} \bordstr_{d}
    \end{equation}
    which, for $I \in \Int$, sends a smooth submanifold 
    \begin{equation}
        \big( M\subset \AA^{\bullet}\times \core(I) \times \RR^{\oo}\big) 
    \end{equation}
    to the same manifold equipped with the trivial stratification. In that case, the constructible tangent bundle of $\big(M \ra *\big)$ is simply the ordinary tangent bundle of $M$, so the stratified tangential datum reduces to the usual rank-$d$ tangent bundle. In other words, one recovers no more data than that of a smooth $(d+\bullet)$-dimensional manifold with boundary.
\end{observation}

\begin{proposition}
    The above map induces a symmetric monoidal faithful functor
    \begin{equation}
        \ho(\bord_d) \xhookrightarrow{\,\,\,\,\,\,\,\,} \ho(\bordstr_{d}),
    \end{equation}
    which for $d>0$ is not full.
\end{proposition}
\begin{proof}
    Symmetric monoidality is immediate, since passing to the trivial stratification commutes with disjoint union.

    To prove faithfulness, we compare with the ordinary bordism categories. By \Cref{ch3:: thm: homotopy category bordism} and \Cref{ch5:: thm: homotopy category stratified bord}, the homotopy categories $\ho(\bord_d)$ and $\ho(\bordstr_d)$ are identified with $\mathsf{Cob}_d$ and $\cS\mathsf{Cob}_d$, respectively. Under these identifications, the induced functor sends the isomorphism class of a smooth bordism to the same bordism endowed with the trivial stratification. If two smooth bordisms become equal in $\cS\mathsf{Cob}_d$, then there is an isomorphism of stratified bordisms between their trivially stratified representatives; but such an isomorphism is exactly a diffeomorphism of smooth bordisms. Hence the original smooth bordism classes were already equal, so the functor is faithful.

    For $d>0$ the functor is not full. Indeed, the cylinder $S^{d-1}\times[0,1]$ admits nontrivial conically smooth stratifications with unstratified incoming and outgoing boundary, for example by declaring the interior hypersurface $S^{d-1}\times\{\frac12\}$ to be a lower-dimensional stratum. This determines a morphism in $\cS\mathsf{Cob}_d$ between unstratified closed $(d-1)$-manifolds which does not lie in the image of any smooth bordism class.
\end{proof}

\begin{observation}
    Every $d$-dimensional unoriented stratified topological quantum field theory 
    \begin{equation}
    \sZ \maps \bordstr_d \longra \cC
    \end{equation} with values in $\cC$ determines, by restriction, a $d$-dimensional unoriented topological quantum field theory with values in $\cC$ 
    \begin{equation}
        \sZ_{\sf sm} \maps \bord_d \xhookrightarrow{\,\,\,\,\,\,\,\,} \bordstr_d \xra{\,\,\,\, \sZ\,\,\,\,} \cC.
    \end{equation}
\end{observation}

More generally, one should expect many other full sub-$(\oo,1)$-categories obtained by imposing local singularity conditions encoded by a right fibration $\cB \ra \Bsc$, provided the allowed singularities retain spherical links of the appropriate codimension.

From this perspective, the present notion of bordism for conically smooth stratified spaces is flexible enough to accommodate examples beyond spaces with an underlying topological $n$-manifold. One might, for instance, try to incorporate finite directed graphs.

%% file: figures/Stcob/SCO3.tikz
\begin{tikzpicture}
	\begin{pgfonlayer}{nodelayer}
		\node [style=none] (4) at (-3, 0) {};
		\node [style=none] (5) at (3, 0) {};
	\end{pgfonlayer}
	\begin{pgfonlayer}{edgelayer}
		\draw [style=geometric] (4.center) to (5.center);
	\end{pgfonlayer}
\end{tikzpicture}

%% file: figures/Stcob/SCO4.tikz
\begin{tikzpicture}
	\begin{pgfonlayer}{nodelayer}
		\node [style=0 strat] (3) at (0, 0) {};
		\node [style=none] (4) at (-3, 0) {};
		\node [style=none] (5) at (3, 0) {};
	\end{pgfonlayer}
	\begin{pgfonlayer}{edgelayer}
		\draw [style=geometric] (4.center) to (5.center);
	\end{pgfonlayer}
\end{tikzpicture}

%% file: Chapters/FO.tex
\cleardoublepage
\chapter{Future outlook}
\phantomsection

We end by offering an overview on how the author expects the work started in this thesis to be continued. In the first section we will explore how the interval category $\Int$ can be generalized beyond dimension $n=1$. The goal is articulating a category $\cube_n$, which when localized along a certain class of maps produces a category equivalent to $\btheta_n$. Using this conjectured localization we define $(\oo,n)$-categories using $\cube_n$, by introducing a generalization of the interval Segal condition \Cref{ch2:: def: interval Segal condition}. Lastly, we put forth a set of submersivity conditions intended to generalize \Cref{ch3:: def: bord_d} and \Cref{ch5:: def: bord_d strat}, and we conjecture that, together with a concordance-resolution step, it produces a model for the $(\oo,n)$-category of fully-extended (stratified) bordisms.

In the last section, we conjecture the existence of a mechanism to produce stratified $d$-dimensional (non-extended) topological quantum field theories using \emph{stratified factorization homology} — a variant of ordinary factorization homology, first introduced in \cite{LurieHA} as topological chiral homology, and later developed in \cite{AFT1}.\footnote{As of this writing, many authors refer to this also as $\alpha$-factorization homology.} The main goal is to bring forth the first touchstone in the attempt to give a stratified generalization of Scheimbauer's thesis \cite{Sch14}, where the author builds a fully-extended topological quantum field theory using $\alpha$-factorization homology. 

\section{Extending intervals to cubes}

The main idea driving this section is to extend the interval Segal model of \((\infty,1)\)-categories to higher dimensions.  For \((\infty,1)\)-categories we used the category \(\Int\), of finite strings of closed intervals as a Segal-type indexing category playing the role of \(\bDelta\).  In higher dimensions, a natural globular replacement for \(\bDelta\) is Joyal–Berger’s category \(\btheta_n\); Rezk’s \(\Theta_n\)-spaces give a well–known model for
\((\infty,n)\)-categories.

Following an idea learned from unpublished work of Manuel Araújo and Pedro Boavida de Brito, we introduce an \(n\)-dimensional cubical indexing category \(\cube_n\) which plays for \(\btheta_n\) the same role that \(\Int\) plays for \(\bDelta\). Objects of \(\cube_n\) are \(n\)-cubes equipped with systems of time coordinates, given by decompositions of each edge into subintervals; morphisms are generated by face maps, degeneracies and refinements of these decompositions. Crucially, \(\cube_n\) is not simply the product \(\Int \times \cdots \times \Int\), since we allow different faces of a cube to carry independent decompositions.

\subsection{\texorpdfstring{$\bTheta_n$}{Theta_n}-spaces}

Before we define $\cube_n$, we start by reviewing $\bTheta_n$, and the notion of Rezk--Segal $\btheta$-space as a model for $(\oo,n)$-categories.

\subsubsection{The category \texorpdfstring{$\bTheta_n$}{Theta_n}}

The literature is inundated with differing models for $(\oo,n)$-categories. This is made especially salient in the wonderful exposition of Chang--Lauda \cite{CL04} and the comparison works of Bergner--Rezk \cite{BR13,BR14}, using iterated wreath products. The wreath construction can be defined in varying degrees of generality, we will follow the variant first described in the work of Berger \cite{BER05} and Bergner--Rezk.

First, we re-introduce Segal's category of finite sets.

\begin{warning}
We will now \emph{re-define} Segal's category and we will add a remark afterwards about the equivalence with the previously defined one.
\end{warning}

\begin{definition}
    Let $\Gamma$ be the category with the following structure:\
    \begin{itemize}
        \item \textbf{Objects}: Finite sets $X,Y,\ldots$.
        \item\textbf{Morphisms}: A morphism from $X$ to $Y$ is given by a function 
        \begin{equation}
            \theta \maps X \longra \cP(Y)
        \end{equation} 
        where $\cP(Y)$ is the set of all subsets of $Y$, with the additional property that $\theta(x_1) \cap \theta(x_2) = \emptyset$ are disjoint for $x_1\neq x_2$. 
        \item \textbf{Composition}: Given a pair of composable morphism $X\xra{\theta} Y \xra{\mu} Z$, we define the composite $\mu \circ \theta$ by sending
        \begin{equation}
            x \longmapsto  \bigcup_{y \in \theta(x)} \mu(y).
        \end{equation} 
    \end{itemize}
\end{definition}
\begin{proposition}
There is a functor, which we will henceforth denote as the \textbf{Segal functor},
\begin{equation}
\cT \maps \bDelta \longra \Gamma    
\end{equation}
which maps $[n]$ to the underlying finite set $\{0,1\ldots,n\}$, and maps a morphism $\alpha\maps [n] \longra [m]$ to the following map in $\Gamma$:
\begin{equation}
    i \longmapsto \{j \in [m]\mkern1mu : \, \alpha(i-1) < j \leq \alpha(i)\}.
\end{equation}
\end{proposition}
\begin{proof}
    The functor is well-defined since $\cT(f)(i)$ must necessarily be disjoint from $\cT(f)(j)$ for $i\neq j$. Moreover, functoriality is direct from the fact that the composite of order-preserving maps is still order-preserving.
\end{proof}
Having defined Segal's category of finite sets, we can now turn to the definition of \emph{wreath product}. 

\begin{definition}\label{ch6:: def: wreath_product}
    Let $\cD$ be an arbitrary category. We define $\Gamma \wr \cD$, which we will call the \textbf{wreath product}, to be the following category:
    \begin{itemize}
        \item\textbf{Objects}: Objects are \emph{symbols} $X(d_{x})_{x \in X}$, where $X$ is a finite set and $(d_{x})_{x\in X}$ is an element of $\mathsf{Obj}(\cD)^{X}$, i.e., an $X$-indexed sequence of objects of $\cD$.
        \item\textbf{Morphisms}: A morphism 
        \begin{equation}
            \theta \maps X(d_{x})_{x \in X} \longra Y(e_{y})_{y \in Y}  
        \end{equation}
        consists of tuples $(\theta_{\Gamma},\{\theta_{xy}\})$; these are a morphism $\theta_\Gamma \maps X \longra Y$ in $\Gamma$, together with a collection of morphisms 
        \begin{equation}
            \theta_{xy} \maps d_x \longra e_y
        \end{equation}
        in $\mathsf{Hom}(\cF)_{\cD}(d_x,e_y)$ whenever $y\in \theta_{\Gamma}(x)$.
        \item \textbf{Composition:} Given a pair of morphisms $\theta \maps X(d_x)_{x \in X} \longra Y(e_y)_{y\in Y}, \phi \maps Y(e_y)_{y\in Y} \longra Z(c_z)_{z\in Z}$; we define its composite morphism to be:
        \begin{enumerate}
            \item Its $\Gamma$ component is defined as 
            \begin{equation}
                (\phi \circ \theta)_{\Gamma} = \phi_\Gamma \circ \theta_{\Gamma}.
            \end{equation}
            \item The collection of morphisms 
            \begin{equation}
                (\phi \circ \theta)_{xz} \maps d_x \longra c_{z}
            \end{equation}
            is defined as containing all composites $\phi_{yz}\circ \theta_{xy}$, subject to $y \in \theta_\Gamma (x)$ and $z \in \phi_\Gamma(y)$.
        \end{enumerate}
    \end{itemize} 
\end{definition}

There is manifestly a forgetful functor $\Gamma \wr \cD \longra \Gamma$, which sends the symbol $X(d_x)_{x\in X}$ to the finite set $X$, and sends a morphism $\theta$ to the underlying map $\theta_\Gamma\colon X \longra Y$ in $\Gamma$. 

\begin{definition}\label{ch6:: def: wreath_prod_general}
    Let $\cF \longra \Gamma$ be an object in $\Cat/\Gamma$. Let the functorial construction:
    \begin{align}
    (-)\wr(-) \maps (\Cat/\Gamma) \times \Cat &\longra \Cat \\
    (\cF \longra \Gamma, \cD) &\longmapsto \cF \times_{\Gamma} (\Gamma \wr \cD) 
    \end{align} 
    which we will denote as $\cF \wr \cD$, and call it the \textbf{wreath product} of $\cF$ with $\cD$.
\end{definition}

Using the Segal functor $\cT \maps \bDelta \longra \Gamma$, we are now allowed to define the wreath product $\bDelta \wr \cF$ with any category $\cF$. Despite the redundancy in exposition, let us carefully unpack its structure, as it is the keystone to building $\bTheta_n$.

\begin{definition}
    Let $\cF$ be an ordinary 1-category. Let $\bTheta \wr \cF$ be the category with the following structure:
    \begin{itemize}
        \item \textbf{Objects}: Tuples $([m],c_1,\ldots,c_m)$, where $[m]\in \bDelta$ and $c_i \in \cF$.
        \item \textbf{Morphisms}: The set of morphisms $[m](c_1,\ldots,c_m) \longra [n](d_1,\ldots,d_n)$ consists of tuples $(\alpha,\{f_{ij}\})$ such that:
        \begin{enumerate}
            \item $\alpha \maps [m] \longra [n]$ is a morphism in $\bDelta$.
            \item For each pair $(i,j)\in [m]\times [n]$ where $\alpha(i-1) < j \leq \alpha(i)$, a morphism 
            \begin{equation}
                f_{ij} \maps c_i \longra d_j
            \end{equation}
            in $\cF$. 
        \end{enumerate}
        \item \textbf{Composition}: Given a pair of composable morphisms:
        \begin{equation}
            ([m],c_1,\ldots,c_m) \xra{(\alpha,\{f_{ij}\})} ([n],d_1,\ldots,d_n) \xra{(\beta,\{g_{jk}\})}
            ([p],e_1,\ldots,e_p)
        \end{equation}
        we define its composite $(\beta\circ \alpha, \{g \circ f\}_{ik})$ as:
        \begin{enumerate}
            \item $(\beta \circ \alpha) \maps [m] \longra [p]$.
            \item For each triple $(i,j,k) \in [m] \times [n] \times [p]$ where:
            \begin{enumerate}
                \item $\alpha(i-1) < j \leq \alpha(i)$, and
                \item $\beta(j-1) < k \leq \beta(j)$,
            \end{enumerate}
            a morphism $(g \circ f)_{ij}$ defined as the composite $g_{jk}\circ f_{ij}$. 
        \end{enumerate}
    \end{itemize} 
\end{definition}
\begin{observation}
There is a useful diagrammatic way of thinking about objects in $\bDelta \wr \cF$, where we depict the tuples as labelled finite directed graphs, and the morphisms as shown below: 
\begin{equation}
    \begin{tikzcd}[ampersand replacement=\&]
        \& 0 \&\& 1 \&\&\&\& 0 \&\& 1 \&\& 2 \&\& 3 \\
        \\
        0 \&\& 1 \&\& 2 \&\& 0 \&\& 1 \&\& 2 \&\& 3
        \arrow[""{name=0, anchor=center, inner sep=0}, "{c_1}"{description}, from=1-2, to=1-4]
        \arrow[between={0.1}{0.8}, maps to, from=1-2, to=3-1]
        \arrow[between={0.1}{0.8}, maps to, from=1-4, to=3-5]
        \arrow[""{name=1, anchor=center, inner sep=0}, "{c_1}"{description}, from=1-8, to=1-10]
        \arrow[between={0.1}{0.9}, maps to, from=1-8, to=3-7]
        \arrow["{c_2}"{description}, from=1-10, to=1-12]
        \arrow[between={0.1}{0.8}, maps to, from=1-10, to=3-11]
        \arrow[""{name=2, anchor=center, inner sep=0}, "{c_3}"{description}, from=1-12, to=1-14]
        \arrow[between={0.1}{0.8}, maps to, from=1-12, to=3-11]
        \arrow[between={0.1}{0.8}, maps to, from=1-14, to=3-13]
        \arrow[""{name=3, anchor=center, inner sep=0}, "{d_1}"{description}, from=3-1, to=3-3]
        \arrow[""{name=4, anchor=center, inner sep=0}, "{d_2}"{description}, from=3-3, to=3-5]
        \arrow[""{name=5, anchor=center, inner sep=0}, "{d_1}"{description}, from=3-7, to=3-9]
        \arrow[""{name=6, anchor=center, inner sep=0}, "{d_2}"{description}, from=3-9, to=3-11]
        \arrow[""{name=7, anchor=center, inner sep=0}, "{d_3}"{description}, from=3-11, to=3-13]
        \arrow["{f_{11}}"{description}, between={0.2}{0.8}, squiggly, from=0, to=3]
        \arrow["{f_{12}}"{description}, between={0.2}{0.8}, squiggly, from=0, to=4]
        \arrow["{f_{11}}"{description}, between={0.2}{0.8}, squiggly, from=1, to=5]
        \arrow["{f_{12}}"{description}, between={0.2}{0.8}, squiggly, from=1, to=6]
        \arrow["{f_{33}}"{description}, between={0.2}{0.8}, squiggly, from=2, to=7]
    \end{tikzcd}
\end{equation}

The example on the left corresponds to the tuple $([1] \xra{\langle 02 \rangle} {[2]},\{f_{11},f_{12}\})$ representing a morphism $[1](c_1) \longra [2](d_1,d_2)$ and the one on the right to the tuple $([3] \xra{\langle 0 2 2 3 \rangle} [3],\{f_{11},f_{12},f_{33}\})$ representing a morphism $[3](c_1,c_2,c_3) \longra [3](d_1,d_2,d_3)$.
\end{observation}

\begin{definition}\label{ch3:: def: theta_n}
    Let $n\geq 0$. We define $\bTheta_{1} = \bDelta$ the terminal category, and let 
    \begin{equation}
        \bTheta_{n} = \bDelta \wr \bTheta_{n-1}.
    \end{equation} 
    The objects of $\Theta_{n}$ are given by $[n](c_1,\ldots,c_n)$ where $[n]\in \bDelta$ and $(c_1,\ldots,c_n)\in \bTheta_{n-1}^{\times n}$. 
\end{definition}
Unpacking the definition of $\bTheta_2$ we observe that its objects admit a familiar interpretation as pasting diagrams. As examples, we showcase the following list:
\begin{equation}
\begin{tikzcd}[ampersand replacement=\&]
	{[0]} \& {=} \& \bullet \\
	\\
	{[1]([0])} \& {=} \& \bullet \&\& \bullet \\
	\\
	{[1]([1])} \& {=} \& \bullet \&\& \bullet \\
	\\
	{[1]([2])} \& {=} \& \bullet \&\& \bullet \\
	\\
	{[3]([1],[0],[1])} \& {=} \& \bullet \&\& \bullet \&\& \bullet \&\& \bullet
	\arrow["0"', shift right=5, draw=none, from=1-3, to=1-3, loop, in=55, out=125, distance=10mm]
	\arrow["0"', shift right=5, draw=none, from=3-3, to=3-3, loop, in=55, out=125, distance=10mm]
	\arrow[from=3-3, to=3-5]
	\arrow["1"', shift right=5, draw=none, from=3-5, to=3-5, loop, in=55, out=125, distance=10mm]
	\arrow["0"', shift right=5, draw=none, from=5-3, to=5-3, loop, in=55, out=125, distance=10mm]
	\arrow[""{name=0, anchor=center, inner sep=0}, curve={height=18pt}, from=5-3, to=5-5]
	\arrow[""{name=1, anchor=center, inner sep=0}, curve={height=-18pt}, from=5-3, to=5-5]
	\arrow["1"', shift right=5, draw=none, from=5-5, to=5-5, loop, in=55, out=125, distance=10mm]
	\arrow["0"', shift right=5, draw=none, from=7-3, to=7-3, loop, in=55, out=125, distance=10mm]
	\arrow[""{name=2, anchor=center, inner sep=0}, from=7-3, to=7-5]
	\arrow[""{name=3, anchor=center, inner sep=0}, curve={height=-18pt}, from=7-3, to=7-5]
	\arrow[""{name=4, anchor=center, inner sep=0}, curve={height=18pt}, from=7-3, to=7-5]
	\arrow["1"', shift right=5, draw=none, from=7-5, to=7-5, loop, in=55, out=125, distance=10mm]
	\arrow["0"', shift right=5, draw=none, from=9-3, to=9-3, loop, in=55, out=125, distance=10mm]
	\arrow[""{name=5, anchor=center, inner sep=0}, curve={height=18pt}, from=9-3, to=9-5]
	\arrow[""{name=6, anchor=center, inner sep=0}, curve={height=-18pt}, from=9-3, to=9-5]
	\arrow["1"', shift right=5, draw=none, from=9-5, to=9-5, loop, in=55, out=125, distance=10mm]
	\arrow[from=9-5, to=9-7]
	\arrow["2"', shift right=5, draw=none, from=9-7, to=9-7, loop, in=55, out=125, distance=10mm]
	\arrow[""{name=7, anchor=center, inner sep=0}, curve={height=18pt}, from=9-7, to=9-9]
	\arrow[""{name=8, anchor=center, inner sep=0}, curve={height=-18pt}, from=9-7, to=9-9]
	\arrow["3"', shift right=5, draw=none, from=9-9, to=9-9, loop, in=55, out=125, distance=10mm]
	\arrow[between={0.2}{0.8}, Rightarrow, from=1, to=0]
	\arrow[between={0.2}{0.8}, Rightarrow, from=3, to=2]
	\arrow[between={0.2}{0.8}, Rightarrow, from=2, to=4]
	\arrow[between={0.2}{0.8}, Rightarrow, from=6, to=5]
	\arrow[between={0.2}{0.8}, Rightarrow, from=8, to=7]
\end{tikzcd}
\end{equation}

\medskip 

\nid It is also instructive to make sense of what morphisms in this category do. Since a lot of the intuition for how $\bTheta_n$ works can be drawn from $n=2$, let us focus on this particular case: the object $[0]$ is terminal, so every other object naturally admits a unique morphism to it. Note that these diagrams come endowed with an implicit choice of framing — namely, the 2-cells are directed top to bottom and, of course, the 1-cells are directed left to right. This corresponds to a choice, and does not reflect any inherent geometric properties of the category.

\begin{observation}
    Note that as we iterate the wreath product, its functoriality provides a generalization of the Segal functor for every $n$, given as
    \begin{equation}
        \cT_n = \alpha \circ (\cT \wr \cT_{n-1})\maps \bDelta \wr \bTheta_{n-1} \longra \Gamma \wr \Gamma \longra \Gamma.
    \end{equation}
    This functor, also called an \emph{assembly functor}, sends a pasting diagram to its set of $n$-cells.
\end{observation}

There is an important factorization system on $\bTheta_n$.

\begin{definition}
    Let $\alpha\colon [n] \longra [m]$ in $\bDelta$ be \textbf{inert} if it is the inclusion of a subinterval in $[m]$. More concretely, for every $i$
    \begin{equation}
        \alpha(i+1) = \alpha(0) + i.
    \end{equation}
    Let $(\alpha,{f}_{ij})$ in $\bTheta_n$ be \textbf{inert} if $\alpha$ is inert in $\bDelta$ and every $f_{ij}$ is inert in $\bTheta_{n-1}$. 
\end{definition}

\begin{definition}
    Let $\alpha \colon [n] \longra [m]$ in $\bDelta$ be \textbf{active} if it preserves the endpoints $\alpha(0) = 0$ and $\alpha(n) = m$. As we did for the inert morphisms above, we define active morphisms in $\bTheta_n$ inductively. Let $(\alpha, f_{ij})$ be \textbf{active} if $\alpha$ is active in $\bDelta$ and every $f_{ij}$ is active in $\bTheta_{n-1}$.
\end{definition}

\subsubsection{Segal \texorpdfstring{$\bTheta_n$}{Theta_n}-spaces}

The Segal conditions for $\bTheta_n$-objects say that, similarly to the $n=1$ case, a functor $\cF\maps \bTheta_n \longra \sset$ is Segal if $\cF(T)$ can be computed as the limit of $\cF$ on the basic $i$-cells.

\begin{definition}
    Let us denote by $c_i$, for $i=0,\ldots,n$, the basic $i$-cells in $\bTheta_n$. We define these objects inductively: for $n=1$, let $c_1 = [1]$ and $c_0 = [0]$; for general $n$, let $c_i = [1](c_{i-1})$.
\end{definition}

\begin{definition}
    Let $\GG_n$ called the \textbf{$n$-globular category} denote the full subcategory of $\bTheta_{n,\mathsf{inr}}$ which contains the $i$-cells $c_0,\ldots,c_n$. This subcategory might be depicted as 
    \begin{equation}
        \begin{tikzcd}[ampersand replacement=\&]
            {c_0} \&\& {c_1} \& {} \& \cdots \&\& {c_n}
            \arrow[shift right=2, from=1-1, to=1-3]
            \arrow[shift left=2, from=1-1, to=1-3]
            \arrow[shift right=2, from=1-3, to=1-5]
            \arrow[shift left=2, from=1-3, to=1-5]
            \arrow[shift right=2, from=1-5, to=1-7]
            \arrow[shift left=2, from=1-5, to=1-7]
        \end{tikzcd}
    \end{equation}  
    Let $i_n \maps \GG_n \hookra \bTheta_{n,\mathsf{inr}}$ denote the inclusion of the full $n$-globular subcategory. Given $T\in \bTheta_n$, we define
    \begin{equation}
        \GG_n /T = \GG_n \underset{\bTheta_{n,\mathsf{inr}}}{\times} \bTheta_{n,\mathsf{inr}}/ T.
    \end{equation} 
    The objects of this category correspond bijectively to the set of $k$-cells of $T$ for every $0\leq k \leq n$. 
\end{definition}

\begin{definition}
    A presheaf $\cF\maps \bTheta_{n}^{\sf op} \longra \sset$ is called a \textbf{Segal $\bTheta_n$-object} if the restriction of $\cF$ to its wide subcategory of \textit{inert} morphisms agrees with $\rkan_{i_{n}}(\cF_{|\GG_n})$. More explicitly, the structure map
    \begin{equation}
        \cF(T) \longra \hoLimit[(c_{i} \xra{\mathsf{inrt}} T)\in \GG_{n}^{\sf op}/T] \cF(c_{i})
    \end{equation}
    is an equivalence. 
\end{definition}

Explicitly, this means that if we evaluate $\cF$ on the following diagram 
\begin{equation}
\begin{tikzcd}[ampersand replacement=\&]
	\bullet \&\& \bullet \&\& \bullet
	\arrow["0"', shift right=5, draw=none, from=1-1, to=1-1, loop, in=55, out=125, distance=10mm]
	\arrow[from=1-1, to=1-3]
	\arrow["1"', shift right=5, draw=none, from=1-3, to=1-3, loop, in=55, out=125, distance=10mm]
	\arrow[""{name=0, anchor=center, inner sep=0}, curve={height=-18pt}, from=1-3, to=1-5]
	\arrow[""{name=1, anchor=center, inner sep=0}, curve={height=18pt}, from=1-3, to=1-5]
	\arrow[""{name=2, anchor=center, inner sep=0}, from=1-3, to=1-5]
	\arrow["2"', shift right=5, draw=none, from=1-5, to=1-5, loop, in=55, out=125, distance=10mm]
	\arrow[between={0.2}{0.8}, Rightarrow, from=0, to=2]
	\arrow[between={0.2}{0.8}, Rightarrow, from=2, to=1]
\end{tikzcd}    
\end{equation} 
we can first evaluate it on 
\begin{equation}
    \begin{tikzcd}[ampersand replacement=\&]
        \bullet \&\& \bullet
        \arrow["0"', shift right=5, draw=none, from=1-1, to=1-1, loop, in=55, out=125, distance=10mm]
        \arrow[from=1-1, to=1-3]
        \arrow["1"', shift right=5, draw=none, from=1-3, to=1-3, loop, in=55, out=125, distance=10mm]
    \end{tikzcd} \qquad \text{ and } \qquad
    \begin{tikzcd}[ampersand replacement=\&]
	\bullet \&\& \bullet
	\arrow["1"', shift right=5, draw=none, from=1-1, to=1-1, loop, in=55, out=125, distance=10mm]
	\arrow[""{name=0, anchor=center, inner sep=0}, curve={height=-18pt}, from=1-1, to=1-3]
	\arrow[""{name=1, anchor=center, inner sep=0}, curve={height=18pt}, from=1-1, to=1-3]
	\arrow[""{name=2, anchor=center, inner sep=0}, from=1-1, to=1-3]
	\arrow["2"', shift right=5, draw=none, from=1-3, to=1-3, loop, in=55, out=125, distance=10mm]
	\arrow[between={0.2}{0.8}, Rightarrow, from=0, to=2]
	\arrow[between={0.2}{0.8}, Rightarrow, from=2, to=1]
    \end{tikzcd} 
\end{equation} separately and taking the fiber product. The rightmost piece of the diagram can now be computed by  evaluating 
\begin{equation}
    \begin{tikzcd}[ampersand replacement=\&]
	\bullet \&\& \bullet
	\arrow["1"', shift right=5, draw=none, from=1-1, to=1-1, loop, in=55, out=125, distance=10mm]
	\arrow[""{name=0, anchor=center, inner sep=0}, from=1-1, to=1-3]
	\arrow[""{name=1, anchor=center, inner sep=0}, curve={height=-18pt}, from=1-1, to=1-3]
	\arrow["2"', shift right=5, draw=none, from=1-3, to=1-3, loop, in=55, out=125, distance=10mm]
	\arrow[between={0.2}{0.8}, Rightarrow, from=1, to=0]
\end{tikzcd} \qquad \text{ and } \qquad
\begin{tikzcd}[ampersand replacement=\&]
	\bullet \&\& \bullet
	\arrow["1"', shift right=5, draw=none, from=1-1, to=1-1, loop, in=55, out=125, distance=10mm]
	\arrow[""{name=0, anchor=center, inner sep=0}, from=1-1, to=1-3]
	\arrow[""{name=1, anchor=center, inner sep=0}, curve={height=18pt}, from=1-1, to=1-3]
	\arrow["2"', shift right=5, draw=none, from=1-3, to=1-3, loop, in=55, out=125, distance=10mm]
	\arrow[between={0.2}{0.8}, Rightarrow, from=0, to=1]
\end{tikzcd}
\end{equation} separately and, again, further taking the fiber product.

\subsection{Category of higher cubes}
Before delving into the definition of $\cube_n$, the main character of this section, let us recall that the interval category $\Int$ admits a canonical functor 
\begin{equation}
    \gamma_1 \colon \Int^{\op} \longra \Gamma
\end{equation}
which is defined as the composition of $\gamma^{\op} \colon \Int^{\op} \longra \bDelta$ and the Segal functor $\cT \colon \bDelta \longra \Gamma$. It sends an interval $I$ to the set of path components of its complement $[0,1]\setminus I$. Predicated on the existence of this functor, we introduce the notion of the $n$-cubes category, a higher dimensional adaptation of $\Int$. 

The definition below, including the iterated wreath-product presentation and the intended localization to $\bTheta_n$, is based on unpublished work of Manuel Araújo and Pedro Boavida de Brito.

\begin{definition}\label{def:: ch6: cube_n category wreath}
    Let $\cube_1^{\op} = \Int^{\op}$. We define the \define{$n$-cubes category} to be the iterated wreath product of $\Int$ with itself $n$ times: 
    \begin{equation}
        \cube_n^{\op} = \cube_1^{\op} \wr \cube_{n-1}^{\op}. 
    \end{equation}
\end{definition}

Alternatively, we can describe the objects of $\cube_n$ as follows: let $[0,1]^{n}$ denote the $n$-dimensional unit cube. The objects of $\cube_n$ can be described as open subsets $R \subseteq [0,1]^{n}$ consisting of:
    \begin{enumerate}
        \item An open collar \[C \simeq [0,\epsilon_l)\times [0,1]^{n-1} \coprod (\epsilon_r,1]\times [0,1]^{n-1},\] of $\partial [0,1] \times [0,1]^{n-1}$, where $0<\epsilon_l< \epsilon_r < 1$.
        \item In the complement of the collar, $R$ consists of a disjoint union of open cubes of the following form:
        \begin{equation}
            \coprod_{\alpha \in A} I_1^{\alpha}\times J_{\alpha}
        \end{equation}
        where $A$ is the set of inner connected components of $I_1 \in \Int_1$ and such that 
        \begin{equation}
            C\simeq I_1^{\sf out}\times [0,1]^{n-1} .
        \end{equation}  
        and $J_{\alpha} \in \cube_{n-1}$.
    \end{enumerate}
Morphisms are given by inclusions of open subsets.

\begin{remark}\label{ch6:: rmk: unpack cube_n}
    We now unpack the definition. For $n=0$, we have that $\cube_0 = \ast$ by definition. For $n=1$, objects are merely finite disjoint unions of open intervals of $[0,1]$ which always must include a collar of the boundary $\{0,1\}\subset [0,1]$. In other words, $\cube_1 = \Int$, i.e., the interval category. For $n=2$, things get immediately more interesting. As a way of illustrating the geometric conditions forced by condition (2), we depict below three examples of finite disjoint unions of open cubes in $[0,1]^2$: 
    \begin{equation}
        \tikzfig{figures/FO/REC1}{\input{figures/FO/REC1.tikz}}
    \end{equation}
    only the left one actually corresponds to an object of $\cube_2$, the one on the right fails to satisfy the collar condition and the one in the middle fails to satisfy both the collar and the product-decomposition condition. To indicate a morphism we might make use of the opacity as a way of distinguishing objects. The next figure illustrates this in a convenient manner, by drawing the overlapping subsets on top of the morphism arrow:
    \begin{equation}
        \tikzfig{figures/FO/REC4}{\input{figures/FO/REC4.tikz}} \,\,\,\,\,\, \xra{\scalebox{0.6}{\tikzfig{figures/FO/REC2}{\input{figures/FO/REC2.tikz}}}} \,\,\,\,\,\,\, \tikzfig{figures/FO/REC3}{\input{figures/FO/REC3.tikz}}
    \end{equation}
\end{remark}

\noindent We now initiate the comparison between $\cube_n$ and $\bTheta_n$, analogously to the way we defined $\gamma\colon \Int \longra \bDelta^{\op}$ in Chapter 2. In fact, our definition will be just as iterative as our definition of $\cube_n$.

One of the advantages of defining $\cube_n$ using an iterated wreath product construction, is that it automatically comes with a functor to $\bTheta_{n}^{\op}$ — defined inductively as follows:
    \begin{equation}
        \gamma_n = \gamma_{1}\wr \gamma_{n-1} \maps \cube_n \longra \bTheta_n^{\op}.
    \end{equation}
This maps 
    \begin{equation}
        \coprod_{\alpha \in A} I_1^{\alpha}\times J_{\alpha} \,\,\longmapsto \,\,\gamma_1(I_1)\big(\gamma_{n-1}(J_1),\ldots, \gamma_{n-1}(J_k)\big).
    \end{equation}

\begin{remark}
    The functor $\gamma_n$, as is defined, might seem opaque and devoid of immediate geometric content, yet the underlying idea underpinning it is wholly geometric. To sharpen our discussion and in order to undress the definition of its inductive nature, let us think about the case of $\gamma_2$. Let us consider the object $R$ in $\cube_2$ which we depicted above:
    \begin{equation}
    \tikzfig{figures/FO/REC4}{\input{figures/FO/REC4.tikz}} 
    \end{equation}
    This object can be written as 
    \begin{equation}
       R = (I_1^{1} \times J_1)\, \coprod\, (I_{1}^{2} \times J_2)
    \end{equation}
    such that
    \begin{align}
        I_1 & = \tikzfig{figures/FO/REC5}{\input{figures/FO/REC5.tikz}}\,, \\ 
        J_1 & = \tikzfig{figures/FO/REC6}{\input{figures/FO/REC6.tikz}}\,, \\ 
        J_2 &= \tikzfig{figures/FO/REC7}{\input{figures/FO/REC7.tikz}}\,. 
    \end{align}
    with the arrows at the bottom highlighting the directions along which $\gamma_2$ assigns $[2]([2],[3])$ to $R$, where the left outer $[2]$ computes $\pi_0\big([0,1]\setminus\pr_1(R)\big)$ with the linear order inherited from the standard orientation on $[0,1]$. The left inner $[2]$ corresponds to $\gamma_{1}(J_1)$ and the right inner $[3]$ corresponds to $\gamma_1(J_2)$. Interpreting this object as a pasting diagram, we note that it corresponds to the following one:
    \begin{equation}
    \begin{tikzcd}
        \bullet && \bullet && \bullet
        \arrow[""{name=0, anchor=center, inner sep=0}, from=1-1, to=1-3]
        \arrow[""{name=1, anchor=center, inner sep=0}, curve={height=-18pt}, from=1-1, to=1-3]
        \arrow[""{name=2, anchor=center, inner sep=0}, curve={height=18pt}, from=1-1, to=1-3]
        \arrow[""{name=3, anchor=center, inner sep=0}, curve={height=-30pt}, from=1-3, to=1-5]
        \arrow[""{name=4, anchor=center, inner sep=0}, curve={height=30pt}, from=1-3, to=1-5]
        \arrow[""{name=5, anchor=center, inner sep=0}, curve={height=12pt}, from=1-3, to=1-5]
        \arrow[""{name=6, anchor=center, inner sep=0}, curve={height=-12pt}, from=1-3, to=1-5]
        \arrow[between={0.2}{0.8}, Rightarrow, from=1, to=0]
        \arrow[between={0.2}{0.8}, Rightarrow, from=0, to=2]
        \arrow[between={0.2}{0.8}, Rightarrow, from=3, to=6]
        \arrow[between={0.2}{0.8}, Rightarrow, from=6, to=5]
        \arrow[between={0.2}{0.8}, Rightarrow, from=5, to=4]
    \end{tikzcd}
    \end{equation}
    This interpretation of the objects of $\cube_n$ is an important one. 

    To give this example some additional flesh, we study also the action of $\gamma_2$ under the morphism depicted in \Cref{ch6:: rmk: unpack cube_n}. The morphism is induced by an inclusion along the first direction which induces a $\pi_0$-isomorphism corresponding to a slight thickening of the first interval component, and two inclusions along the second direction: the first merges the two cubes over the first component whilst the second is just an equality. As seen from the image of $\gamma_2$, this corresponds to the morphism 
    \begin{equation}
        \begin{tikzcd}
            \bullet && \bullet && \bullet
            \arrow[""{name=0, anchor=center, inner sep=0}, curve={height=-18pt}, from=1-1, to=1-3]
            \arrow[""{name=1, anchor=center, inner sep=0}, curve={height=18pt}, from=1-1, to=1-3]
            \arrow[""{name=2, anchor=center, inner sep=0}, curve={height=-30pt}, from=1-3, to=1-5]
            \arrow[""{name=3, anchor=center, inner sep=0}, curve={height=30pt}, from=1-3, to=1-5]
            \arrow[""{name=4, anchor=center, inner sep=0}, curve={height=12pt}, from=1-3, to=1-5]
            \arrow[""{name=5, anchor=center, inner sep=0}, curve={height=-12pt}, from=1-3, to=1-5]
            \arrow[between={0.2}{0.8}, Rightarrow, from=0, to=1]
            \arrow[between={0.2}{0.8}, Rightarrow, from=2, to=5]
            \arrow[between={0.2}{0.8}, Rightarrow, from=5, to=4]
            \arrow[between={0.2}{0.8}, Rightarrow, from=4, to=3]
        \end{tikzcd} \longrightarrow 
        \begin{tikzcd}
            \bullet && \bullet && \bullet
            \arrow[""{name=0, anchor=center, inner sep=0}, from=1-1, to=1-3]
            \arrow[""{name=1, anchor=center, inner sep=0}, curve={height=-18pt}, from=1-1, to=1-3]
            \arrow[""{name=2, anchor=center, inner sep=0}, curve={height=18pt}, from=1-1, to=1-3]
            \arrow[""{name=3, anchor=center, inner sep=0}, curve={height=-30pt}, from=1-3, to=1-5]
            \arrow[""{name=4, anchor=center, inner sep=0}, curve={height=30pt}, from=1-3, to=1-5]
            \arrow[""{name=5, anchor=center, inner sep=0}, curve={height=12pt}, from=1-3, to=1-5]
            \arrow[""{name=6, anchor=center, inner sep=0}, curve={height=-12pt}, from=1-3, to=1-5]
            \arrow[between={0.2}{0.8}, Rightarrow, from=1, to=0]
            \arrow[between={0.2}{0.8}, Rightarrow, from=0, to=2]
            \arrow[between={0.2}{0.8}, Rightarrow, from=3, to=6]
            \arrow[between={0.2}{0.8}, Rightarrow, from=6, to=5]
            \arrow[between={0.2}{0.8}, Rightarrow, from=5, to=4]
        \end{tikzcd}
    \end{equation}  
    which vertically composes the two cells on the left, also encoded as
    \begin{equation}
        [2]([1],[3]) \xra{\id (\langle 0,2 \rangle, \id)} [2]([2],[3])
    \end{equation}
\end{remark}

\begin{remark}
    Although we will mostly use the geometric model, it is helpful to keep in mind the inductive form of morphisms. A morphism
    \begin{equation}
        f\maps I_{1}(J_1,\ldots, J_n) \longra I'_{1}(J'_{1},\ldots,J'_{k})
    \end{equation}
    in $\cube_n$ may be viewed as an inclusion in the first coordinate together with, for each inner connected component of $I_1$, a morphism in $\cube_{n-1}$ from the corresponding $J_i$ to the lower-dimensional piece of the target lying over the image component. In the geometric model this is simply the inclusion of one rectangular decomposition into another.
\end{remark}

\subsection{Segal \texorpdfstring{$n$}{n}-cube-spaces}
The guiding idea is the same as in Chapter 2. Rather than work with all morphisms in $\Int = \cube_1$, we singled out those which ought to become invertible once the model is topologized, and localized along them. For $\cube_n$ we follow the same philosophy: we isolate the morphisms which should play the role of weak equivalences from the geometric point of view.

\begin{remark}
For $n=1$, these are precisely the inclusions $I \longra J$ of open intervals inducing isomorphisms on $\pi_0$. 
\end{remark}

We generalize this set as follows.

\begin{definition}
    Let $\cW_{\pi_{0}}(n)$ denote the subcategory of $\cube_n$ spanned by the morphisms 
    \begin{equation}
    f\colon R \longra S, \quad \text{ such that } \gamma_n(f) = \id.
    \end{equation} 
\end{definition}

This is the higher-dimensional analogue of the class $\cW_{\pi_{0}}$ introduced in Chapter 2. Geometrically, one should think of the morphisms in $\cW_{\pi_{0}}(n)$ as those inclusions obtained by enlarging the cubes of $R$ without creating new overlaps or forcing previously distinct rectangles to merge. This suggests the following conjectural comparison with the usual $\bTheta_n$-model.   

\begin{conjecture} \label{fo:: conj: local.}
    There should be a Quillen equivalence
\begin{equation}
    \begin{tikzcd}[ampersand replacement=\&]
	{\fun(\cube_{n},\sset)_{\mathsf{inj}}} \&\& {\fun(\bTheta_{n}^{\op},\sset)_{\mathsf{inj}}.}
	\arrow[""{name=0, anchor=center, inner sep=0}, "{\gamma_{n}^*}", bend left=10, from=1-1, to=1-3]
	\arrow[""{name=1, anchor=center, inner sep=0}, "{\rkan_{\gamma_{n}}}", bend left=10, from=1-3, to=1-1]
	\arrow["\dashv"{anchor=center, rotate=-90}, draw=none, from=1, to=0]
\end{tikzcd}
\end{equation}
induced by the functor $\gamma_{n} \maps \cube_{n} \longra \bTheta_{n}^{\op}$.
\end{conjecture}

We only indicate why one might expect this to hold.

\begin{proof}[Sketch of the expected argument]
    The blueprint of this proof is identical to the one for $n=1$ (see the proof of \Cref{ch2:: prop: weak contractibility of pi_0}), with an important caveat: that proof relies on the contractibility of ordered configuration spaces on the line. For higher $n$, this ceases to be true even for $n=2$. Evidently, this means we need to adapt our previous argument. We use instead a subspace $\conf^{\sf special}_T \subset ([0,1]^n)^{\pi_{0}(R_{T})}$ consisting of configurations adapted to the tree-like structure of $T$. This is inspired by Ayala--Hepworth's description of configuration-space cells indexed by $n$-orderings, or equivalently by healthy planar level trees \cite{ayalaConfigurationSpaces$Theta_n$2014}; in the present rectangular setting, one expects $\conf^{\sf special}_T$ to be a convex cell.
    \begin{align}
        |\cJ_T| \,\, &\simeq \,\, \hoColimit[R_T \in \cJ_T]\ast \\
        &\overset{(1)}{\simeq} \,\,\hoColimit[R_T \in \cJ_T] \conf(R_T) \cap \conf^{\sf special}_T \\
        &\overset{(2)}{\simeq} \,\,\conf^{\sf special}_T \\
        &\overset{(3)}{\simeq} \,\,\ast 
    \end{align}

    Some remarks about how this proof works: both $\conf(R_T)$ and $\conf^{\sf special}_T$ are contractible spaces; the former being homeomorphic to a product of Euclidean spaces and the latter, by design, to a convex cell in $([0,1]^{n})^{\pi_0(R_T)}$, thereby rendering their intersection contractible, too; this proves $(1)$. To obtain $(2)$, one argues as in the interval case: given a special configuration $\phi \in \conf^{\sf special}_T$, the subposet 
    \begin{equation}
        \big\{R_T \in \cJ_T\mkern1mu :\, \phi \in \conf(R_T)\cap \conf^{\sf special}_T\big\} \subset \cJ_{T}
    \end{equation}
    is filtered, therefore it has contractible geometric realization. Lurie's $\oo$-categorical version of Seifert-Van Kampen then identifies the homotopy colimit of the cover by the contractible pieces $\conf(R_T)\cap \conf^{\sf special}_T$ with $\conf^{\sf special}_T$. Finally, $(3)$ follows from the convexity of $\conf^{\sf special}_T$.
\end{proof}

Assuming \Cref{fo:: conj: local.}, we will call a functor $\cF \maps \cube_n \longra \sset$ a $\cube_n$-space if it sends morphisms in $\cW_{\pi_{0}}(n)$ to weak equivalences. We will call it a \emph{Segal} $\cube_n$-\emph{space} if, in addition, for every $R \in \cube_n$, the value of $\cF(R)$ can be recovered by successively cutting $R$ into single cubes, that is, into the pieces representing $k$-cells, one coordinate direction at a time and in increasing order. A typical such decomposition is illustrated by
\begin{equation}
    \tikzfig{figures/FO/REC14}{\input{figures/FO/REC14.tikz}}
\end{equation}
\medskip

\nid Since $\cF$ sends axis-aligned enlargements to equivalences, many different cube configurations may encode the same $k$-cell of $\bTheta_n$. Herein lies the advantage of working with $\cube_n$ rather than directly with $\bTheta_n^{\op}$ in geometric contexts: one keeps the rigid globular bookkeeping, while retaining the flexibility of rectangular decompositions and the freedom to isotope configurations.   

\subsubsection{\texorpdfstring{$(\oo,d)$}{(infty,d)}-categories of stratified cobordisms}

One of the goals of defining $(\oo,d)$-categories using Segal $\cube_d$-spaces is so that we can construct $(\oo,d)$-categories of fully extended stratified bordisms by slightly adapting the methods used in Chapter 3 and Chapter 5.

\medskip

\nid Before proceeding with this construction, let us first illustrate why naively copying the structure of $\tilbordstr_d$ will fail if we do not begin by addressing the definition of $\cube_n$. The issue is that a general object of $\cube_2$ allows different components in the first direction to carry unrelated decompositions in the second direction. After cutting once, there is then no canonical lower-dimensional cube recording where the next cuts should occur, so the inductive description of composition loses control of the slicing data. Assume $(\oo,d=2)$, and that we are working with 2-dimensional 2-bordisms.

\begin{equation}
    \tikzfig{figures/Stcob/SCO11}{\input{figures/Stcob/SCO11.tikz}}
\end{equation}
\medskip

The point of the picture is precisely that the two pieces want different interval data in the transverse direction, so they do not naturally refine a common lower-dimensional cutting pattern.

In order to accommodate for the construction of the $(\oo,d)$-category of stratified bordisms, we must restrict $\cube_n$ to a full subcategory $\cR\cube_{n}$ spanned by the objects $X \in \cube_n$ such that

\begin{equation}
    X = I_1 (J_1, \ldots, J_n), \quad \gamma_{1}(I_1) = [n]
\end{equation}
and there exists $J \in \cR\cube_{n-1}$ together with refinements $J_{i} \ra J$ for each $1 \leq i \leq n$.

\begin{construction}\label{fo:: def: higher strat bord}
     Fix $V$ to be a finite-dimensional real vector space. We would like to define a functor 
    \begin{equation}
    \tilbordstr_{[0,\ldots,d];V}\maps \man^{\op}\times \cR\cube_{d} \longrightarrow \set    
    \end{equation}
    sending a pair $(M,R)$, where $M  \in \man^{\op}$ and $R\in \cR\cube_d$ such that 
    \begin{equation}
        R = I(J_1, \ldots, J_n),
    \end{equation} 
    to the set of bounded closed subspaces 
    \begin{equation}
   \big( W\subseteq M\times \core(R) \times V \big)
    \end{equation} 
    such that the following conditions are satisfied:
    \begin{enumerate}
    \item \label{id: cor 1} $W$ is endowed with the structure of a conically smooth stratification and a $\SMfd_{d}^{\lag d \rag}$-framing. In particular, it is a conically smooth stratified space with a $\lag d \rag$-corner structure
    \begin{equation}
         (W \ra P),\quad (W \xra{c} \cP(\{1,\cdots,d\}));
    \end{equation} 
    and the embedding must preserve the underlying stratifications i.e., send the boundary of $W$ to the boundary of $M\times \core(R) \times V$.
    \item \label{id: cor 2} The restricted projection 
    \begin{equation}
        \pi\maps \big(W \cap (M\times I^{\mathsf{c}}\times [0,1]^{d-1} \times V) \big) \hookra M\times I^{\mathsf{c}}\times [0,1]^{d-1} \times V \fib  M\times I^{\mathsf{c}}\times [0,1]^{d-1}
    \end{equation} 
    is a conically smooth submersion.
    \item \label{id: cor 3} When restricting $W$ to $J_{i} \times I^{\sf inn}$, for every $t\in I^{\sf inn}$ one expects the resulting slice over $t$ to determine an element of $\tilbordstr_{[0,\ldots,d-1];V}(M,J_{i})$.
    \item \label{id: cor 4} The restricted projection 
    \begin{equation}
    \pi_W\maps W \hookra M\times \core(R) \times V \fib M 
    \end{equation}
    is a conically smooth submersion.
    \end{enumerate} 

    \noindent Now we define how $\tilbordstr_{d,V}$ acts on morphisms:
    \begin{itemize}
    \item[$(\man)$] The morphism $(f,\mathsf{id}_{R})\maps (M,R) \longra (N,R)$ maps the closed submanifold 
    \begin{equation}
    W\subseteq M\times \core(R) \times V
    \end{equation}
    to the pullback along $(f\times \id)$:
    \begin{equation}
        \begin{tikzcd}[ampersand replacement=\&]
            {(f\times \id)^\ast W} \& {W} \\
            \\
            N\times \core(R) \times V \& M\times \core(R) \times V
            \arrow[from=1-1, to=1-2]
            \arrow[hook, "(f\times \id)^\ast \pi_{|W}"',from=1-1, to=3-1]
            \arrow["\lrcorner"{anchor=center, pos=0.125}, draw=none, from=1-1, to=3-2]
            \arrow[hook, from=1-2, to=3-2]
            \arrow["f\times \id"', from=3-1, to=3-2]
        \end{tikzcd}
    \end{equation}
    which is a conically smooth stratified $\big(d + \mathsf{dim}(N)\big)$-dimensional submanifold $f^\ast W \subseteq N \times \core(R) \times V$.
    \item[$(\cR\cube)$] A morphism $(\id_M, R \longra S)$ should act by \emph{forgetting} connected components and restricting the stratified $d$-manifold with corners $(W\ra P)$.
    \end{itemize} 
\end{construction}

\begin{warning}
\Cref{fo:: def: higher strat bord} should be read as a first-pass candidate rather than as a completed definition. The recursive slicing condition captures the geometry one wants, but in its present form it does not yet control where the induced lower-dimensional cuts land after iterating the construction. In particular, this does not yet produce a true functor on the nose.
\end{warning}

The point of \Cref{fo:: def: higher strat bord} is that, if a correct definition of this flavour exists, one expects the $\cube_n$-Segal condition to reduce to \emph{cuts} in a single direction. This is the reason for keeping the construction inductive, and it is precisely this one-directional cutting picture that already appears in the proof of the $(\oo,1)$-case of \Cref{ch5:: thm: stratbord interval segal}.

As a way of wrapping up this section, we will try to motivate our inductive definition by way of examples in $d=2$. Let us start with $R=I() \in \cube_2$ — which corepresents objects in the $\cube_d$ model for $(\oo,d)$-categories. An element of $\tilbordstr_{[0,1,2];V}(\AA^{0},R)$ is a stratified 2-dimensional manifold with corners whose projection to $\core(R) \simeq [0,1]^2$ is a conically smooth submersion as follows:
\medskip
\begin{equation}
    \tikzfig{figures/FO/REC15}{\input{figures/FO/REC15.tikz}}
\end{equation}
\medskip

\nid and it corresponds to the object which is a disjoint union of a pair of unoriented points. If we consider the case of 1-morphisms which are co-represented by objects of the type $R = I(J_1)$ where $J_1 \in \cube_1 = \Int$ such that $\gamma(J_1) = [0]$. An element of $\tilbordstr_{[0,1,2];V}(\AA^0,R)$ is a stratified 2-dimensional manifold with corners whose projection would look as follows:
\begin{equation}
    \tikzfig{figures/FO/REC16}{\input{figures/FO/REC16.tikz}}
\end{equation}
\medskip

\nid Finally, to represent a 2-morphism, one should allow a cube $R = I(J_1)\in \cube_2$ where $J_1 $ has $\gamma(J_1)=[1]$ and whose two first-direction components carry interval data in the transverse direction. An element of $\tilbordstr_{[0,1,2];V}(\AA^0,R)$ should then be thought of as a stratified surface interpolating between 1-morphisms of the previous kind, as in the following picture:

\begin{equation}
    \tikzfig{figures/FO/REC17}{\input{figures/FO/REC17.tikz}}
\end{equation}
\bigskip

\section{Stratified TQFTs via factorization homology}

In the unstratified setting, factorization homology is one of the cleanest mechanisms we have for passing from local algebraic data to topological field theories. The reason is simple: every smooth $d$-manifold is locally modelled on $\RR^d$. One may therefore encode local observables by a single algebraic object --- an $\cE_d$-algebra, or equivalently a $\Disk_d$-algebra --- and then integrate that datum over the manifold. In the fully extended setting, this perspective is already present in Lurie's sketch and is later developed in Scheimbauer's thesis \cite{Lur09,LurieHA,AF1,Sch14}.

For stratified spaces, however, there is no longer a unique local model. Near a point in a codimension-$k$ stratum one sees a basic of the form $\RR^{d-k}\times \sC(L)$, and the link $L$ remembers how the deeper strata sit nearby. So the coefficient system must now keep track not only of the top-dimensional disk, but of all the allowed singularity types and the embeddings among them. This is the conceptual reason for replacing ordinary disk algebras by $\cB$-disk algebras for a suitable category of basics $\cB$.

There is also a second language lurking in the background, which is often better adapted to field-theoretic intuition. Rather than presenting the coefficient system directly as a disk-algebra, one can recast the same local-to-global structure in the language of factorization algebras on opens. In the stratified setting the relevant objects should be \emph{constructible factorization algebras}: roughly, factorization algebras that are locally constant along each stratum and whose variation across strata is controlled by the stratification. We will not try to develop that formalism here, but it is worth keeping it in mind, since it is the natural language for pushing observables along the projection of a stratified bordism to an interval. With that motivation in place, let us recall the operadic definition from \cite{AFT3}.

\begin{definition}\cite[Definition~2.2 and Notation~2.9]{AFT3}
    Let $\cB \longra \Bsc$ be an $(\oo,1)$-category of basics and let $\cV^{\otimes}$ be a symmetric monoidal $(\oo,1)$-category. The symmetric monoidal $(\oo,1)$-category $\Disk(\cB)\subset \Mfd(\cB)$ is the smallest full symmetric monoidal subcategory generated by the basics in $\cB$. A \define{$\cB$-disk algebra in} $\cV$ is a symmetric monoidal functor
    \begin{equation}
        A \maps \Disk(\cB) \longra \cV.
    \end{equation}
\end{definition}

\begin{definition}\cite[Definition~2.14 and Theorem~2.15]{AFT3}
    Let $A$ be a $\cB$-disk algebra in a symmetric monoidal $(\oo,1)$-category $\cV^{\otimes}$ for which factorization homology exists, and let $X\in \Mfd(\cB)$. The \define{factorization homology of $X$ with coefficients in $A$} is the left Kan extension of $A$ along $\Disk(\cB)\hookra \Mfd(\cB)$, and is denoted
    \begin{equation}
        \int_{X} A.
    \end{equation}
    When $\cV^{\otimes}$ is $\otimes$-sifted cocomplete, it is computed by the colimit formula
    \begin{equation}
        \int_{X} A \,\simeq\, \colim_{(U\hookra X)\in \Disk(\cB)_{/X}} A(U).
    \end{equation}
\end{definition}

The picture to keep in mind is the following. Once a defect datum
\begin{equation}
    \Bsc_{d}^{\cS_{d-1},\DD} \longra \Bsc_{d}^{\cS_{d-1}}
\end{equation}
has been fixed, the relevant local coefficients should be $\Bsc_{d}^{\cS_{d-1},\DD}$-disk algebras. In other words, one does not attach algebraic data only to the smooth $d$-disk, but to every allowed local singularity type, together with the compatibilities imposed by embeddings of basics. This is precisely the point where the defect labels introduced in Chapter~5 begin to interact with factorization homology in a meaningful way.

Now let $W$ be a morphism in our interval Segal space $\bordstr_d$, together with its projection $\pi\maps W \longra \core(I)$. If $A$ is such a stratified disk algebra, one expects the assignment
\begin{equation}
    U \subset \core(I) \longmapsto \int_{\pi^{-1}(U)} A
\end{equation}
to define a constructible factorization algebra on $\core(I)$ relative to the finite stratification induced by $I$. If $I\in \Int$ encodes a finite family of cuts, let $S_{I}\subset \core(I)$ denote the corresponding finite set of separating points. Then the expected pushforward should be locally constant on each connected component of $\core(I)\setminus S_{I}$ and only allowed to jump across the strata in the induced stratification
\begin{equation}
    \emptyset \subset S_{I} \subset \core(I).
\end{equation}
This is precisely the kind of one-dimensional constructibility one wants in order to remember the algebraic data assigned to composition. In the unstratified case, Scheimbauer uses factorization algebras on $(0,1)$ which are locally constant with respect to such finite stratifications to model bimodules and, more generally, the levels of the Morita target \cite[Corollary~3.2.25 and Definition~3.2.27]{Sch14}. The same pattern suggests that in the stratified setting the first target should not yet be $\Vect$, but rather an interval Segal space of constructible factorization algebras on stratified intervals. Moreover, if one cuts $W$ along a regular value of $\pi$, the excision theorem of \cite{AFT3} identifies the factorization homology of the whole bordism with the relative tensor product of the factorization homologies of the two pieces. Thus the gluing law required for composition is not added by hand; it is already built into the formalism.

\begin{conjecture}
    Let $\cV^{\otimes}$ be a symmetric monoidal $(\oo,1)$-category with the colimits needed to define stratified factorization homology, and let $A$ be a $\Bsc_{d}^{\cS_{d-1},\DD}$-disk algebra in $\cV$. Then factorization homology with coefficients in $A$ should determine a symmetric monoidal functor out of the stratified bordism category constructed in Chapter~5, landing naturally in an interval Segal space of constructible factorization algebras on stratified intervals. After comparison with the corresponding Morita-type target, this should produce a non-extended stratified topological quantum field theory.
\end{conjecture}

%% file: Chapters/Appendix.tex
\chapter{Auxiliary background}\label{sec:: ap: categorical_background}

This appendix collects a small amount of background used repeatedly in Chapters~2--5. The goal is not to give a self-contained treatment of these topics, but to gather the constructions and lemmas that are invoked later while keeping the main chapters focused on the new material of the thesis.

\section{Categories of elements and finality}

\begin{definition}\label{ap:: def: cat of els}
    Let $\cC$ be a small category and let $\cF \maps \cC \longra \set$ be a functor. The \define{category of elements} of $\cF$, denoted
    \begin{equation}
        \int \cF,
    \end{equation}
    is the category whose:
    \begin{enumerate}
        \item objects are pairs $(c,x)$ with $c\in \cC$ and $x\in \cF(c)$;
        \item morphisms $(c,x)\to (d,y)$ are morphisms $f\maps c\to d$ in $\cC$ satisfying $\cF(f)(x)=y$.
    \end{enumerate}
    There is a tautological forgetful functor
    \begin{equation}
        \int \cF \longra \cC.
    \end{equation}
\end{definition}

\begin{remark}
    If $K\maps \bDelta^{\op}\longra \set$ is a simplicial set, then its category of simplices is precisely the category of elements of $K$. This is the construction used in Chapter~2.
\end{remark}

\begin{definition}\label{ap:: def: initial functor}
    Let $K\maps \cF \longra \cD$ be a functor.
    \begin{enumerate}
        \item We say that $K$ is \define{initial} if for every complete category $\cE$ and every functor $F\maps \cD \longra \cE$, the canonical map
        \begin{equation}
            \Limit[\cD]{F} \longra \Limit[\cF]{F\circ K}
        \end{equation}
        is an isomorphism whenever either side exists.
        \item We say that $K$ is \define{final} if for every cocomplete category $\cE$ and every functor $F\maps \cD \longra \cE$, the canonical map
        \begin{equation}
            \Colimit[\cF]{F\circ K} \longra \Colimit[\cD]{F}
        \end{equation}
        is an isomorphism whenever either side exists.
    \end{enumerate}
\end{definition}

\begin{proposition}
    A functor $K\maps \cF \longra \cD$ is initial if, for every object $d\in \cD$, the comma category $d/K$ is non-empty and connected. Dually, $K$ is final if, for every object $d\in \cD$, the comma category $K/d$ is non-empty and connected.
\end{proposition}

\begin{remark}
    The homotopical version of this criterion replaces connectedness by weak contractibility. This is the guiding principle behind several contractibility arguments in Chapter~2.
\end{remark}

\section{Kan extensions}\label{ap:: subsec: Kan extensions}

\begin{definition}
    Let $F\maps \cF \longra \cE$ and $K\maps \cF \longra \cD$ be functors. A \define{right Kan extension} of $F$ along $K$ is a functor
    \begin{equation}
        \rkan_{K}(F)\maps \cD \longra \cE
    \end{equation}
    together with a natural transformation
    \begin{equation}
        \epsilon \maps \rkan_{K}(F)\circ K \Rightarrow F
    \end{equation}
    which is terminal among such pairs. Dually, a \define{left Kan extension} of $F$ along $K$ is initial among functors out of $\cD$ equipped with a natural transformation from $F$.
\end{definition}

\begin{theorem}\label{ap:: thm: right Kan extension pointwise}
    Assume that $\cE$ is complete. Then the right Kan extension of $F$ along $K$ exists and is computed pointwise by
    \begin{equation}
        \rkan_{K}(F)(d)\cong \Limit[d/K]{\big(d/K \longra \cF \xra{\,\,F\,\,} \cE\big)}.
    \end{equation}
    Dually, if $\cE$ is cocomplete, then left Kan extensions exist and are computed by colimits over the corresponding comma categories.
\end{theorem}

\section{Discrete fibrations and presheaves}

Fibrations enter the thesis in two different ways: first as the ordinary categorical shadow of presheaves, and later as the $(\oo,1)$-categorical language in which tangential and defect data are most naturally expressed.

\begin{definition}
    A functor $p\maps \cE \longra \cB$ is a \define{discrete right fibration} if for every morphism $f\maps x \to y$ in $\cB$ and every object $\widetilde{y}\in \cE$ with $p(\widetilde{y})=y$, there exists a unique morphism $\widetilde{f}\maps \widetilde{x}\to \widetilde{y}$ in $\cE$ such that $p(\widetilde{f})=f$.

    \medskip

    \nid Dually, $p$ is a \define{discrete left fibration} if $p^{\op}\maps \cE^{\op}\longra \cB^{\op}$ is a discrete right fibration.
\end{definition}

\begin{construction}\label{ap:: cons: discrete fibration}
    Let $p\maps \cE \longra \cB$ be a discrete right fibration. We define a presheaf
    \begin{equation}
        \mathsf{Un}(p)\maps \cB^{\op} \longra \set
    \end{equation}
    by sending an object $b\in \cB$ to the fibre
    \begin{equation}
        \mathsf{Un}(p)(b)=\{\widetilde{b}\in \cE : p(\widetilde{b})=b\}.
    \end{equation}
    Given a morphism $f\maps b\to b'$ in $\cB$, the map
    \begin{equation}
        \mathsf{Un}(p)(f)\maps \mathsf{Un}(p)(b') \longra \mathsf{Un}(p)(b)
    \end{equation}
    sends $\widetilde{b'}$ to the source of the unique lift of $f$ with target $\widetilde{b'}$.

    \medskip

    \nid Conversely, if $\cF\maps \cB^{\op}\longra \set$ is a presheaf, then the forgetful functor
    \begin{equation}
        \int \cF \longra \cB
    \end{equation}
    is a discrete right fibration. These two constructions are inverse to one another up to canonical isomorphism.
\end{construction}

\begin{proposition}
    The construction of \Cref{ap:: cons: discrete fibration} defines an equivalence of categories
    \begin{equation}
        \mathsf{dRFib}_{\cB} \,\, \simeq \,\, \Fun(\cB^{\op},\set).
    \end{equation}
\end{proposition}

\begin{definition}
    Let $p\maps \cE \longra \cC$ be a map of quasicategories. We say that $p$ is a \define{right fibration} if it has the right lifting property with respect to every inner horn inclusion
    \begin{equation}
        \Lambda^{n}_{i}\hookra \Delta^{n},
    \end{equation}
    for $0<i\leq n$.
\end{definition}

\begin{theorem}\label{ap:: thm: straightening unstraightening}
    For every $(\oo,1)$-category $\cC$, straightening and unstraightening induce an equivalence of $(\oo,1)$-categories
    \begin{equation}
        \mathsf{RFib}_{\cC} \,\, \simeq \,\, \fun(\cC^{\op},\mathsf{Kan}).
    \end{equation}
\end{theorem}

\begin{remark}
    Here $\mathsf{Kan}$ is a convenient model for the $(\oo,1)$-category of spaces. Restricting to right fibrations with discrete fibres recovers set-valued presheaves. This is the perspective used in Chapter~5 to encode discrete defect data.
\end{remark}

\section{A homotopy covering criterion}

\begin{theorem}\label{ap:: thm: lurie Seifert}
    Let $X$ be a topological space, let $K$ be a small category, and let
    \begin{equation}
        F\maps K \longra \mathcal{U}(X)
    \end{equation}
    be a functor to the poset of open subsets of $X$. Assume that the opens $\{F(k)\}_{k\in K}$ cover $X$ and that, for every point $x\in X$, the full subcategory
    \begin{equation}
        K_{x}=\{k\in K:\, x\in F(k)\}
    \end{equation}
    has weakly contractible nerve. Then the canonical map
    \begin{equation}
        \hoColimit[K]{F} \longra X
    \end{equation}
    is a weak homotopy equivalence.
\end{theorem}

\begin{remark}
    In Chapter~2 we use this criterion to identify certain nerves with contractible spaces by exhibiting them as homotopy colimits of contractible open sets in a contractible ambient space.
\end{remark}

\section{Differential topology}

The bordism constructions of Chapters~3 and~4 repeatedly use the same small package of differential-topological facts: proper submersions are locally trivial, bundles over contractible bases are trivial, and the same remains true in charts adapted to a boundary. We gather those statements here.

\begin{definition}\label{ap:: def: smooth fibration}
    A smooth map $\pi \maps E \longra B$ is a \define{smooth fibre bundle} if every point $b\in B$ admits an open neighbourhood $U\subseteq B$, a smooth manifold $F$, and a diffeomorphism
    \begin{equation}
        \varphi_{U}\maps \pi^{-1}(U) \xra{\,\,\cong\,\,} U\times F
    \end{equation}
    commuting with the projections to $U$.
\end{definition}

\begin{proposition}\label{ap:: prop: triviality contractibility}
    Let $\pi\maps E \longra B$ be a smooth fibre bundle with fibre $F$. If $B$ is contractible, then $\pi$ is isomorphic, as a smooth fibre bundle, to the trivial bundle
    \begin{equation}
        B\times F \longra B.
    \end{equation}
\end{proposition}
\begin{proof}
    This is standard. One way to prove it is to use path lifting to transport the fibre over a chosen basepoint along a contraction of $B$. See \cite[Corollary~10.3]{Hu94}.
\end{proof}

\begin{definition}\label{ap:: def: proper map}
    A map $f\maps M \longra N$ is \define{proper} if the inverse image of every compact subset of $N$ is compact in $M$.
\end{definition}

\begin{proposition}[Ehresmann fibration theorem]\label{ap:: prop: Ehresmann fibration no boundary}
    Let $f\maps M \longra N$ be a smooth proper submersion between smooth manifolds without boundary. Then $f$ is a smooth fibre bundle.
\end{proposition}

\begin{cor}
    Under the hypotheses of \Cref{ap:: prop: Ehresmann fibration no boundary}, all fibres of $f$ are diffeomorphic.
\end{cor}

\begin{proposition}\label{ap:: prop: ehresmann boundary}
    Let $M$ and $N$ be smooth manifolds with boundary, and let $f\maps M \longra N$ be a smooth map such that:
    \begin{enumerate}
        \item $f$ is a proper submersion;
        \item $f(\partial M)\subseteq \partial N$.
    \end{enumerate}
    Then $f$ is a smooth fibre bundle in charts adapted to the boundary. In particular, its fibres are smooth manifolds without boundary and are locally constant up to diffeomorphism.
\end{proposition}
\begin{proof}
    Away from $\partial N$ this is the classical Ehresmann theorem. Near a boundary point of $N$, choose a chart modelled on $\RR_{\geq 0}\times \RR^{n-1}$. The inward-pointing vector field in the normal direction can be lifted locally along the submersion $f$. Because $f(\partial M)\subseteq \partial N$, the lift can be chosen inward-pointing along $\partial M$, so its flow remains inside $M$. Integrating these local lifts yields local trivializations compatible with the boundary, exactly as in the standard proof of Ehresmann's theorem.
\end{proof}

\begin{remark}
    This is the version used in Chapter~3 when one studies projections to intervals in the presence of collars. Chapter~4 proves the corresponding conically smooth stratified analogue.
\end{remark}